%% file: asymptotics_en.tex
\documentclass{amsbook}

\usepackage{amssymb}
\usepackage{amsfonts}
\usepackage{amsmath}
\usepackage[legalpaper,bookmarks=true,colorlinks=true,linkcolor=blue,citecolor=blue]{hyperref}
\usepackage{graphicx}%
\usepackage{fancyhdr}
\usepackage{color}
\usepackage[mathlines]{lineno}

\newtheorem{theorem}{Theorem}
\theoremstyle{plain}

\newtheorem{corollary}{Corollary}

\newtheorem{definition}{Definition}

\newtheorem{lemma}{Lemma}

\newtheorem{proposition}{Proposition}

\numberwithin{equation}{section}

\begin{document}
\Large
\pagenumbering{roman}
\begin{center}

\huge \textbf{Gane Samb LO}\\
\bigskip
\huge \textbf{Modou NGOM and Tchilabalo A. KPANZOU}\\
\vskip 5cm
\Huge \textbf{Weak Convergence (IA)} \\
\bigskip
-\\
\bigskip \textbf{Sequences of Random Vectors}
\vskip 5cm

\huge \textit{\textbf{Statistics and Probability African Society (SPAS) Books Series}.\\
 \textbf{Calgary, Alberta. 2016}}.\\

\bigskip \Large  \textbf{DOI} : http://dx.doi.org/10.16929/sbs/2016.0001\\
\bigskip \textbf{ISBN} 978-2-9559183-1-9
\end{center}

\newpage
\begin{center}
\huge \textbf{SPAS TEXTBOOKS SERIES}
\end{center}

\bigskip \bigskip

\Large

 \begin{center}
 \textbf{GENERAL EDITOR of SPAS EDITIONS}
 \end{center}

\bigskip
\noindent \textbf{Prof Gane Samb LO}\\
gane-samb.lo@ugb.edu.sn, gslo@ugb.edu.ng\\
Gaston Berger University (UGB), Saint-Louis, SENEGAL.\\
African University of Sciences and Technology, AUST, Abuja, Nigeria.\\

\bigskip

\begin{center}
\Large \textbf{ASSOCIATED EDITORS}
\end{center}

\bigskip
\noindent \textbf{KEHINDE DAHUD SHANGODOYIN}\\
shangodoyink@mopipi.ub.bw\\
UNIVERSITY Of BOTSWANA\\

\noindent \textbf{Blaise SOME}\\
some@univ-ouaga.bf\\
Chairman of LANIBIO, UFR/SEA\\
Ouaga I Pr Joseph Ki-Zerbo University.\\

\bigskip
\begin{center}
\Large \textbf{ADVISORS}
\end{center}

\bigskip

\noindent \textbf{Ahmadou Bamba SOW}\\
ahmadou-bamba.sow@ugb.edu.sn\\
Gaston Berger University, Senegal.\\

\noindent \textbf{Tchilabalo Abozou KPANZOU}\\
kpanzout@yahoo.fr\\
Kara University, Togo.\\

\newpage

\noindent \textbf{Library of Congress Cataloging-in-Publication Data}\\

\noindent Gane Samb LO, 1958-, Modou NGOM , Tchilabalo A. KPANZOU \\

\noindent Weak Convergence (IA). Sequences of Random Vectors.\\

\noindent SPAS Books Series, 2016.\\

\noindent Statistics and Probability African Society (SPAS).\\

\noindent \textit{DOI} : 10.16929/sbs/2016.0001\\

\noindent \textit{ISBN}  978-2-9559183-1-9

\newpage

\noindent \textbf{Author : Gane Samb LO}\\
\bigskip

\bigskip
\noindent \textbf{Emails}:\\
\noindent gane-samb.lo@ugb.edu.sn, ganesamblo@ganesamblo.net.\\

\bigskip
\noindent \textbf{Url's}:\\
\noindent www.ganesamblo@ganesamblo.net\\
\noindent www.statpas.net/cva.php?email.ganesamblo@yahoo.com.\\

\bigskip \noindent \textbf{Affiliations}.\\
Main affiliation : Gaston Berger University, UGB, SENEGAL.\\
African University of Sciences and Technology, AUST, ABuja, Nigeria.\\
Affiliated as a researcher to : LSTA, Pierre et Marie Curie University, Paris VI, France.\\

\noindent \textbf{Teaches or has taught} at the graduate level in the following universities:\\
Saint-Louis, Senegal (UGB)\\
Banjul, Gambia (TUG)\\
Bamako, Mali (USTTB)\\
Ouagadougou - Burkina Faso (UJK)\\
African Institute of Mathematical Sciences, Mbour, SENEGAL, AIMS.\\
Franceville, Gabon\\

\bigskip \noindent \textbf{Dedicatory}.\\

\noindent \textbf{To my wife Mbaye Ndaw Fall who is accompanying for decades with love and patience}

\bigskip \noindent \textbf{Acknowledgment of Funding}.\\

\noindent The author acknowledges continuous support of the World Bank Excellence Center in Mathematics, Computer Sciences and Intelligence Technology, CEA-MITIC. His research projects in 2014, 2015 and 2016 are funded by the University of Gaston Berger in different forms and by CEA-MITIC.

\newpage
\noindent \textbf{Author : Modou NGOM}\\
\bigskip

\bigskip

\noindent \textbf{Email}:\\
\noindent ngomodoungom@gmail.com,ngomodoungom@yahoo.fr\\

\bigskip \noindent \textbf{Affiliations}.\\
Main affiliation : Gaston Berger University, UGB, SENEGAL.\\

\noindent \textbf{Modou Ngom} participated in the project of the book while preparing his PhD thesis :
\textbf{On a stochastic process doubly indexed, with margins estimating the extreme value index, and its Gaussian and non Gaussian Asymptotic}. He collaborated with Prof Gane Samb LO update this part of his PhD course in Weak convergence to produce the book.\\

\noindent Modou Ngom has also been a lecturer many years in Mathematics, Probability and Measure, at the University of Gaston Berger of Saint-Louis, Senegal,Saint-Louis, Senegal (UGB)\\

\bigskip \noindent \textbf{Dedicatory of Modou NGOM}.\\

\noindent \textbf{To my family, whose support,  encouragement and faithful prayers arouse my envy to go forward}.\\

\newpage

\noindent \textbf{Author : Tchilabalo A. KPANZOU}\\
\bigskip

\noindent \textbf{Emails}:\\
\noindent kpanzout@gmail.com, kpanzout@yahoo.fr\\

\noindent \textbf{Url's}:\\
\noindent https://sites.google.com/a/aims.ac.za/tchilabalo\\
\noindent http://univi.net/spas/cvf.php?email=kpanzout@yahoo.fr\\
 
\bigskip \noindent \textbf{Affiliations}.\\
Main affiliation: University of Kara, Kara, TOGO.\\ 
\noindent \textbf{Teaches or has taught} at the graduate level in the following universities:\\
University of Kara (UK), TOGO\\
University of Lom\'e (UL), TOGO\\
Ecole Normale Sup\'erieure (ENS), TOGO\\
University of Abomey-Calavi (UAC), BENIN\\

\title{Weak Convergence (IA). Sequences of Random Vectors}

\begin{abstract} \large (\textbf{English}) This monograph aims at presenting the core weak convergence theory for sequences of random vectors with values in $\mathbb{R}^k$. In some places, a more general formulation in metric spaces is provided. It lays out the necessary foundation that paves the way to applications in particular sub-fields of the theory. In particular, the needs of Asymptotic Statistics are addressed. A whole chapter is devoted to weak convergence in $\mathbb{R}$ where specific tools, for example for handling weak convergence of sequences using independent and identically distributed random variables such that the Renyi's representations by means of standard uniform or exponential random variables, are stated. The functional empirical process is presented as a powerful tool for solving a  considerable number of asymptotic problems in Statistics. The text is written in a self-contained approach with the proofs of all used results at the exception of the general Skorohod-Wichura Theorem.\\

\noindent (\textbf{Fran\c{c}ais}) \noindent Cet ouvrage a l'ambition de pr\'esenter le noyau dur de la th\'eorie de la convergence vague de suite de vecteurs al\'eatoires dans $\mathbb{R}^k$. Autant que possible, dans certaines situations, la th\'eorie g\'en\'erale dans des espaces m\'etriques est donn\'ee. Il pr\'epare la voie \`a une sp\'ecialisation dans certains sous-domaines de la convergence vague. En particulier, les besoins de la statistique asymptotique ont \'et\'e satisfaits. Un chapitre de l'ouvrage concerne la convergence vague dans $\mathbb{R}$ avec des outils sp\'ecifiques, par exemple, pour \'etudier les suites de variables al\'eatoires ind\'ependantes et identiquement distribu\'ees tels que la repr\'esentation de Renyi au moyen de variables al\'eatoires uniformes ou exponentielles standard. Le processus empirique fonctionnel est introduit comme un outil puissant pour \'etudier des probl\`emes asymptotiques en Statistiques. 
L'ouvrage est r\'edig\'e dans une approche auto-citante avec toutes les preuves des r\'esultats utilis\'es, \`a l'exception du Th\'eor\`eme de Skorohod-Wichura.\\


\noindent \textbf{AMS 2010 Classification Subjects :} 60XXX; 62G30
\end{abstract}

\maketitle

\frontmatter
\tableofcontents
\mainmatter
\Large
\include{asymptotics_cv_preface_gen_en}

\include{asymptotics_cv_preface_WC_gen_en}

\include{asymptotics_cv_preface_ed1_en}

\include{asymptotics_cv_review_rk_en} 
\include{asymptotics_cv_01_en} 
\include{asymptotics_cv_02_en} 
\include{asymptotics_cv_03_en} 
\include{asymptotics_cv_04_en} 
\include{asymptotics_analysis_01_en} 
\include{asymptotics_analysis_02_en} 

\include{asymptotics_biblio_en}
\end{document}

%% file: asymptotics_cv_preface_gen_en.tex
\chapter*{General Preface}

\noindent \textbf{This textbook} is the first of series whose ambition is to cover broad part of Probability Theory and Statistics.  These textbooks are intended to help learners and readers, of all levels, to train themselves.\\

\noindent As well, they may constitute helpful documents for professors and teachers for both courses and exercises.  For more ambitious  people, they are only starting points towards more advanced and personalized books. So, these textbooks are kindly put at the disposal of professors and learners.

\bigskip \noindent \textbf{Our textbooks are classified into categories}.\\

\noindent \textbf{A series of introductory  books for beginners}. Books of this series are usually accessible to student of first year in 
universities. They do not require advanced mathematics.  Books on elementary probability theory and descriptive statistics are to be put in that category. Books of that kind are usually introductions to more advanced and mathematical versions of the same theory. The first prepare the applications of the second.\\

\noindent \textbf{A series of books oriented to applications}. Students or researchers in very related disciplines  such as Health studies, Hydrology, Finance, Economics, etc.  may be in need of Probability Theory or Statistics. They are not interested by these disciplines  by themselves.  Rather, the need to apply their findings as tools to solve their specific problems. So adapted books on Probability Theory and Statistics may be composed to on the applications of such fields. A perfect example concerns the need of mathematical statistics for economists who do not necessarily have a good background in Measure Theory.\\

\noindent \textbf{A series of specialized books on Probability theory and Statistics of high level}. This series begin with a book on Measure Theory, its counterpart of probability theory, and an introductory book on topology. On that basis, we will have, as much as possible,  a coherent presentation of branches of Probability theory and Statistics. We will try  to have a self-contained, as much as possible, so that anything we need will be in the series.\\

\noindent Finally, \textbf{research monographs} close this architecture. The architecture should be so large and deep that the readers of monographs booklets will find all needed theories and inputs in it.\\

\bigskip \noindent We conclude by saying that, with  only an undergraduate level, the reader will  open the door of anything in Probability theory and statistics with \textbf{Measure Theory and integration}. Once this course validated, eventually combined with two solid courses on topology and functional analysis, he will have all the means to get specialized in any branch in these disciplines.\\

\bigskip \noindent Our collaborators and former students are invited to make live this trend and to develop it so  that  the center of Saint-Louis becomes or continues to be a reknown mathematical school, especially in Probability Theory and Statistics.

%% file: asymptotics_cv_preface_WC_gen_en.tex
\chapter*{General Preface of Our Series of Weak Convergence}

\noindent \textbf{The series Weak convergence} is an open project with three categories.\\

\noindent \textbf{The special series  Weak convergence I} consists of texts devoted to the core theory of weak convergence, each of them concentrated on the handling of one specific class of objects. The texts will have labels $A$, $B$, etc. Here are some examples.\\

\noindent (1) Weak convergence of Random Vectors (IA).\\

\noindent (2) Weak convergence of stochastic processes and empirical processes (IB).\\

\noindent (3) Weak convergence of random measures (IC).\\

\noindent (4) Weak convergence of fuzzy random measures (IC).\\

\bigskip \noindent \textbf{The special series Weak convergence II} consists of textbooks related to the theory of weak convergence, each of them concentrated on one specialized field using weak convergence. Usually, these subfields are treated apart in the literature. Here, we want to put them in our general frame as continuations of the Weak Convergence Series I. Some examples are the following.\\

\noindent (1) Weak laws of sums of independent random variables.\\

\noindent (2) Weak laws of sums of associated random variables.\\

\noindent (3) Univariate Extreme value Theory.\\

\noindent (4) Multivariate Extreme value Theory.\\

\noindent (5) Etc.\\

\bigskip \noindent \textbf{The special series Weak convergence III} consists of textbooks focusing on statistical applications of Parts of the Weak Convergence Series I and Weak Convergence Series II. Examples :\\

\noindent (1) A handbook of Gaussian Asymptotic Distribution Using the Functional Empirical Process.\\

\noindent (2) A handbook of Statistical Estimation of the Extreme Value index.\\

\noindent (1) etc.\\

%% file: asymptotics_cv_preface_ed1_en.tex
\chapter*{Preface of The Series Weak Convergence : Sequence of Random Vectors}

\noindent \textbf{The series Weak convergence (IA)} concerns the theory of weak convergence of sequences of random vectors. Due to the theorem of Kolmgorov, stating broadly that the probability law of any random element is characterized by its finite distribution under the appropriate state spaces, the place of the distributions of random vectors is surely central to Probability Theory.\\

\noindent This motivated us to begin this series by the weak convergence of random vectors as the foundation of all the structure.\\

\noindent Another reason is that the needs of Asymptotic Statistics, which is one of the main motivations of the development of Weak Convergence Theory, generally does not need more than that. This booklet then gives to some readers exactly what they specifically.\\

\noindent This textbook focuses on the study of random elements in $\mathbb{R}^k$, $k\geq 1$. So the properties and the topology of $\mathbb{R}^k$ are used.\\

\noindent But when only the general properties of the metric are used, we prefer to give the results in the general case where the studied sequences have their values in a metric space with a metric $d$.\\

\noindent The concept of tightness is essential in weak convergence theory. In this text, the Helly-Bray method is exclusively used.\\

\noindent This textbook is concluded by a chapter of the functional empirical process. Here, only the weak limits of its finite distributions are treated. We show how to use it for deriving asymptotic results in many research problems. With such tools, even at this somewhat elementary level of weak convergence, it is possible for readers to provide contributions in many research fields in Statistics and in applied related fields.

\noindent I wish you a pleasant reading and hope receiving your feedback.\\

\noindent To my wife Mbaye Ndaw Fall who is accompanying me since decades.\\

\noindent Saint-Louis, Calgary, Abuja, Bamako, Ouagadougou, 2016. 
  
\newpage
\noindent \huge \textbf{Preliminary Remarks and Notations}.\\

\Large
 
\bigskip
\noindent \textbf{WARNINGS}\\

\noindent \textbf{(1)} In all this book, any unspecified limit in presence with the subscripts $n$ are meant as $n\rightarrow +\infty$.\\

\noindent \textbf{(2)} This textbook deals with general distribution functions $F$ on $\mathbb{R}^k$, $k\geq 1$. The Lebesgue-Stieljes measure induced by a general distribution function is not necessarily a probability measure. If this induced Lebesgue-Stieljes is a probability measure, we precise this distribution function as a \textbf{probability distribution function}. As well, a distribution function of a random vector $X$ of Lebesgue-Stieljes is implicitly a probability distribution function although we do not say : the probability distribution function of $X$.\\






%% file: asymptotics_cv_review_rk_en.tex
\chapter{Review of Usual Weak Convergence Results in $\mathbb{R}^k$} \label{ChapRevCvRk}

\section{Introduction} \label{cv.review.sec1}

In this chapter, we will see that most of the readers, actually know a considerable number of weak convergence results, even if they
did not use this concept. What has to be done, on top of this review, is to present these individual results in the frame of a unified theory in the most general setting. This is the target of this book which will be given in the subsequent chapters.\\

\noindent Here, we are going to recall classical convergence results that any student should have encountered from the first courses in probability theory or in Statistics.\\

\noindent We begin to set the general frame of weak convergence in $\mathbb{R}^{k}$, $k\geq 1$. We will admit the statements in the following section. We will be able to establish their validity in Chapter \ref{cv}, in particular in Theorem \ref{cv.theo.portmanteau.rk} of that chapter.\\

\section{Weak Convergence in  $\mathbb{R}^{k}$} \label{cv.review.sec2}

\noindent Let us remind that the probability law of any vector random variable $X:(\Omega ,\mathbb{A},\mathbb{P})\mapsto \mathbb{R}^{k}$ is 
characterized by \\ 

\bigskip \noindent \textbf{(a)} its distribution function: \newline
\begin{equation*}
\mathbb{R}^{k}\ni x\hookrightarrow F_{X}(x)=\mathbb{P}(X\leq x), 
\end{equation*}

\bigskip \noindent \textbf{(b)} its characteristic function (Here, $i$ is the complex number such that $i^{2}=-1$ with positive sinus, and  $<.,.>$ stands for the classical product space on  $\mathbb{R}^{k}$) 
\begin{equation*}
\mathbb{R}^{k}\ni u\hookrightarrow \Phi (u)=\mathbb{E}(\exp(i<u,X>)),
\end{equation*}

\bigskip \noindent \textbf{(c)} its moment generating function (if it exists in a neighborhood of the null vector) 
\begin{equation*}
\mathbb{R}^{k}\ni u\hookrightarrow \Psi _{X}(u)=\mathbb{E}(\exp(<u,X>)).
\end{equation*}

\bigskip \noindent and

\bigskip \noindent \textbf{(d)} its Radon-Nikodym derivative, or probability density function (\textit{pdf}), (if it exists), with respect to (\textit{w.r.t})  a measure $\nu$ on $\mathbb{R}^{k}$ :  
\begin{equation*}
d\mathbb{P}/d\nu =f_{X}.
\end{equation*}

\bigskip \noindent It is interesting that these characteristics also play the main roles in weak convergence through Theorem 
\ref{cv.theo.portmanteau.rk} we will prove in Chapter \ref{cv}.\\

\noindent We have :

\begin{theorem} (\textbf{THEOREM-DEFINITION-LEMMA}) \label{cv.review.rk} Let $X_{n} :(\Omega_{n},\mathcal{A}_{n},\mathbb{P}_{n})\mapsto
(\mathbb{R}^{k},\mathbb{B}(\mathbb{R}^{k}))$ be a sequence of random vectors, $X : (\Omega_{\infty} ,\mathcal{A}_{\infty},\mathbb{P}_{\infty})\mapsto (\mathbb{R}^{k},\mathbb{B}(\mathbb{R}^{k}))$ a random vector. Then the assertions (a) and (b) below are equivalent.\newline

\bigskip \noindent \textbf{(a)} For any $u\in \mathbb{R}^{k}$,\newline

\begin{equation*}
\Phi_{X_{n}}(u)\rightarrow \Phi _{X}(u)\text{ as }n\rightarrow +\infty. 
\end{equation*}

\bigskip \noindent \textbf{(b)} For any continuity point $u\in \mathbb{R}^{k}$ of $F_X$, 
\begin{equation*}
F_{X_{n}}(x)\rightarrow F_{X}(x) \text{ as } n\rightarrow +\infty .
\end{equation*}

\bigskip \noindent If one the assertions (a) or (b) holds, we say that the sequence $X_{n}$ weakly converges to $X$, or $X_{n}$ converges in distributions to $X$ or $X_{n}$ converges in law to $X$, as $n\rightarrow +\infty$, and we denote this by
 
\begin{equation*}
X_{n}\rightsquigarrow X\text{ or }X_{n}\longrightarrow _{d}X\text{ or }X_{n}%
\overset{\mathcal{L}}{\longrightarrow }X\text{ or }X_{n}\overset{w}{%
\longrightarrow }X\text{ or }X_{n}\longrightarrow _{w}X
\end{equation*}

\bigskip The weak limit is unique in distribution, meaning that if $X_{n}$ weakly converges to $X$ and to $Y$, then $X$ and $Y$ have the same distribution, that is $F_X=F_Y$ in the context of $\mathbb{R}^{k}$.\\
 
\bigskip \noindent We also have the following sufficiency weak convergence conditions.\\

\bigskip \noindent \textbf{(c)} If the moment generating functions (\textit{mgf}) $\Psi_{X_{n}}$ exist on  $B_{n}$, $n\geq 1$ and $\Psi _{X}$ exists on $B$, where the $B_{n}$ and $B$ are neighborhoods of $0$ such that $B \subset \cap_{n\geq 0} B_n$, and if for any $x\in B$, 
\begin{equation*}
\Psi _{X_{n}}(x)\rightarrow \Psi _{X}(x)\text{ as }n\rightarrow +\infty, 
\end{equation*}

\noindent then $X_{n}$ weakly converges to $X$.\newline

\bigskip \noindent \textbf{(d)} Finally, suppose that the probability distribution $\mathbb{P}_{n}(\circ)=\mathbb{P}_{n}(X_n \in \circ)$, $n\geq 1$, and $\mathbb{P}_{X}(\circ)=\mathbb{P}_{\infty}(X \in \circ)$ have Radon-Nikodym derivatives with respect to a measure $\nu$ on $\mathbb{R}^k$, denoted by
$$
d\mathbb{P}^{n}/d\nu =f_{X_{n}} \ \ n\geq 1, \ \ d\mathbb{P}_{X}/d\nu =f_{X}.
$$

\noindent If for any  $x\in D_{X}=\{x,f_{X}(x)>0\}$, 
\begin{equation*}
f_{X_{n}}(x)\rightarrow f_{X}(x)\text{ as }n\rightarrow +\infty ,
\end{equation*}

\noindent then  $X_{n}\rightsquigarrow X$.\newline

\noindent  We have the following last point.\\

\noindent \textbf{(e)} Assume that the sequence $\{X_n, \ \ n\geq 1\} \subset \mathbb{R}^k$ weakly converges to $X \in \mathbb{R}^k$, as $n \rightarrow +\infty$ and let $A$ be a real $(m,k)$-matrix with $m\geq 1$. Then $\{AX_n, \ \ n\geq 1\} \subset \mathbb{R}^m$ weakly converges to $AX \in \mathbb{R}^m$.
\end{theorem}

\bigskip \noindent \textbf{Remark}. Point (e) of Theorem \ref{cv.review.rk} above is a consequence of the continuous mapping Theorem \ref{cv.mappingTh} in Chapter \ref{cv}.\\

\bigskip \noindent \textbf{In summary, the weak convergence in $\mathbb{R}^{k}$ holds when the distribution functions, the characteristic functions, the moment generating functions (if they exist) or the probability density functions (if they exist) with respect to the same measure $\nu$, point-wisely converge  to the distribution function, or to the characteristic function or to moment generating function (if it exists), or to the probability density unction (if it exists) with respect to $\nu$ of a probability measure in $\mathbb{R}^{k}$. In the case of point-wise convergence of the distribution functions, only matters the convergence for continuity points of the limiting distribution functions}.\\

\bigskip \noindent  All this is awesome and gives us pretty well tools to deal with weak convergence. The examples given below form the core set of examples you cannot ignore.\\

\noindent But before we proceed to this review, we need a handsome criterion derived from the convergence of characteristic functions.\\

\begin{proposition} \label{cv.wold} \textbf{(Wold Criterion)}. The sequence $\{X_n, \ \ n\geq 1\} \subset \mathbb{R}^k$ weakly converges to $X \in \mathbb{R}^k$, as $n \rightarrow +\infty$ if and only if for any $a \in \mathbb{R}^k$, the sequence $\{<a,X_n>, \ \ n\geq 1\} \subset \mathbb{R}$ weakly converges to $<A,X> \in \mathbb{R}$ as $n \rightarrow +\infty$.
\end{proposition}

\noindent \textbf{Proof}. The proof is quick and uses the notation above. Suppose that $X_n$ weakly converges to $X$ in $\mathbb{R}^k$ as $n \rightarrow +\infty$. By using the convergence of characteristic functions, we have for any $u\in \mathbb{R}^k$
$$
\mathbb{E}\left(\exp(i<X_n,u>) \right) \rightarrow \mathbb{E}\left(\exp(i<X,u>)\right)\ \ as \ \ n \rightarrow +\infty.
$$  

\bigskip \noindent It follows that for any $a\in \mathbb{R}^k$ and for any $t \in \mathbb{R}$, we have

\begin{equation}
\mathbb{E}\left(\exp(it<X_n,a>) \right)\rightarrow \mathbb{E}\left(\exp(it<X,a>) \right) \ \ as \ \ n \rightarrow +\infty, \label{cv.proj}
\end{equation}

\bigskip \noindent that is, by taking $u=ta$ in the formula above, and by denoting $Z_n=<X_n,a>$ and $Z=<X,a>$, we have

$$
\mathbb{E}\left(\exp(itZ_n) \right) \rightarrow \mathbb{E}\left(exp(itZ)\right) \ \ as \ \ n \rightarrow +\infty.
$$

\noindent This means that $Z_n \rightsquigarrow Z$, that is $<a,X_n>$ weakly converges to $<a,X>$.\\

\noindent Conversely, suppose that for any $a \in \mathbb{R}^k$, the sequence $\{<a,X_n>, \ \ n\geq 1\} \subset \mathbb{R}$ weakly converges to 
$<A,X> \in \mathbb{R}$ as $n \rightarrow +\infty$. Then by taking $t=1$ in (\ref{cv.proj}) we get for any $a=u \in \mathbb{R}^k$,

$$
\mathbb{E}\left(\exp(i<X,u>)\right) \rightarrow \mathbb{E}\left(\exp(i<X,u>) \right)\ \ as \ \ n \rightarrow +\infty.
$$

\noindent which means that $X_n \rightsquigarrow X$ as $n\rightarrow +\infty$.\\

\newpage

\section{Examples of Weak Convergence in $\mathbb{R}$} \label{cv.review.sec3}

\subsection{Weak Convergence of a sequence of Hyper-geometric random variables to a Binomial random variable} \label{cv.review.subsec.HypBin} $ $\\

\noindent Let $X_{N}$ be a random variable following a Hyper-geometric law $\mathcal{H}(N,M,n)$
with $M/N\rightarrow p$, $N\rightarrow \infty$, $n$ being fixed. Then $X_{N}$ weakly converges to a Binomial random variable $X$, that 
is $X \sim \mathcal{B}(n,p)$.\\

\noindent \textbf{Proof}. Let us use the probability density functions with respect to the counting measure $\nu$ on $\mathbb{N}$. We have

\begin{equation*}
f_{X_{n}}(k)=\frac{\left( 
\begin{tabular}{l}
$M$ \\ 
$k$%
\end{tabular}%
\right) \left( 
\begin{tabular}{l}
$N-M$ \\ 
$n-k$%
\end{tabular}%
\right)}{\left( 
\begin{tabular}{l}
$N$ \\ 
$n$%
\end{tabular}%
\right) },0\leq k\leq \min (n,M).
\end{equation*}

\noindent Suppose that $M/N\rightarrow p$, $N\rightarrow \infty$. We have 
\begin{eqnarray*}
f_{X_{n}}(k) &=&\frac{M!}{k!(M-k)!}\frac{(N-M)!}{(n-k)!(N-M-(n-k))!}\frac{%
n!(N-n)!}{N!} \\
&=&\frac{n!}{k!(n-k)!}\times \frac{M!}{(M-k)!}\times \frac{(N-M)!}{(N-M-(n-k))!}\times \frac{(M-n)!}{N!} \\
&=&\left( 
\begin{tabular}{l}
$n$ \\ 
$k$%
\end{tabular}%
\right) \times \left\{ \frac{M!}{(M-k)!}\right\} \left\{ \frac{(N-M)!}{%
(N-M-(n-k))!}\right\} \left\{ \frac{(M-n)!}{N!}\right\} .
\end{eqnarray*}

\noindent But
\begin{eqnarray*}
\left\{ \frac{M!}{(M-k)!}\right\}  &=&(M-k+1)(M-k+2)...(M-1)M \\
&=&M^{k}(1-\frac{k-1}{M})(1-\frac{k-2}{M}) \times ...\times (1-\frac{1}{M}) \\
&=&M^{k}(1+o(1))
\end{eqnarray*}

\noindent since $M\rightarrow \infty $ and $k$ is fixed. Next, 

\begin{eqnarray*}
\left\{ \frac{(N-M)!}{(N-M-(n-k))!}\right\}  &=&(N-M-(n-k)+1)\times
...\times (N-M-1)(N-M) \\
&=&(N-M)^{n-k}(1+\frac{n-k-1}{N-M})(1+\frac{n-k-2}{N-M})\\
&\times& ...\times ((1+\frac{1}{N-M}) \\
&=&(N-M)^{n-k}(1+o(1)),
\end{eqnarray*}

\noindent since, also, $N-M=N(1-M/N)\sim N(1-p)\rightarrow \infty $ and $n-k$ is fixed. Finally 
\begin{eqnarray*}
\left\{ \frac{(M-n)!}{N!}\right\}  &=&\frac{1}{(N-n+1)(N-n+2)...(N-1)N} \\
&=&\frac{1}{N^{n}(1-\frac{n-1}{N})(1-\frac{n-2}{N})...(1-\frac{1}{N})} \\
&=&\frac{1}{N^{n}(1+o(1))}.
\end{eqnarray*}

\noindent for similar reasons. In total for any $0\leq k\leq n$%

\begin{equation*}
f_{X_{n}}(k)=\left( 
\begin{tabular}{l}
$n$ \\ 
$k$%
\end{tabular}%
\right) \left( \frac{M}{N}\right) ^{k}\left( \frac{N-M}{N}\right)
^{n-k}(1+o(1))\rightarrow \left( 
\begin{tabular}{l}
$n$ \\ 
$k$%
\end{tabular}%
\right) p^{k}(1-p)^{n-k}.
\end{equation*}

\noindent Hence, for any $k$ in the support set of the \textit{pdf} of a $\mathcal{B}(n,p)$ random variable \textit{w.r.t} to the counting measure $\nu$, denoted 
\begin{equation*}
f_{X}(k)=\left( 
\begin{tabular}{l}
$n$ \\ 
$k$%
\end{tabular}%
\right) p^{k}(1-p)^{n-k},
\end{equation*}

\noindent we have
\begin{equation*}
\forall (1\leq k\leq n),f_{X_{n}}(k)\rightarrow f_{X}(k).
\end{equation*}

\bigskip \noindent The proof is finished.\\

\noindent \textbf{Useful remark in sampling technique theory}. This result allows to treat drawing without replacement (which generates a hyper-geometric law) may be approximated as a drawing with replacement (which gives a Binomial law) when the size of the global population is large. The idea behind this is the following : if we randomly draw a small number of individuals from a large set, it is almost improbable that we draw one individual more than one time.\\

\subsection{Weak Convergence of a sequence of Binomial random variables to a Poisson random variable} \label{cv.review.subsec.BinPois} $ $\\

\noindent
Let $X_{n}$ be a sequence of $\mathcal{B}(n,p)$ random variable with  $p=p_{n}\rightarrow 0$
and  $np_{n}\rightarrow \lambda$, $0<\lambda$, as  $n\rightarrow \infty$. Then $X_{n}$ weakly converges to a Poisson random variable  $X$ with parameter $\lambda$, that is $X \sim \mathcal{\lambda}$.\\

\noindent \textbf{Proof}. Let us use the moment generating functions. Let $X_{n}$ be a sequence of $\mathcal{B}(n,p_n)$-random variable and $X$ be a $\mathcal{P}(\lambda)$ random variable. We have 
\begin{equation*}
\Psi _{X_{n}}(t)=((1-p_{n})+p_{n}e^{t})^{n}, \text{ }n\geq 1;\Psi
_{X}(t)=\exp (\lambda (e^{t}-1)),t\in \mathbb{R}.
\end{equation*}

\noindent Put $\lambda _{n}=np_{n}\rightarrow \lambda$. For any fixed $t$, we have
\begin{equation*}
\Psi _{X_{n}}(t)=\left(\frac{\lambda _{n}}{n}+\left(1-\frac{\lambda _{n}}{n}\right)e^{t}\right)^{n}=\left( 1-\frac{\lambda _{n}(e^{t}-1)}{n}\right) ^{n}\rightarrow
\exp (\lambda (e^{t}-1))=\Psi _{X}(t)
\end{equation*}

\noindent by the following classical results of Calculus courses :  
\begin{equation*}
\left(1+\frac{x_{n}}{n}\right)^{n}\rightarrow e^{x}\text{ as }n\rightarrow +\infty 
\text{ whenever  }x_{n}\rightarrow x\in \mathbb{R}\text{ as }%
n\rightarrow +\infty .
\end{equation*}

\subsection{Weak Convergence of a sequence of Poisson random variable to a Gaussian random variable} \label{cv.review.subsec.PoisGauss} $ $\\

\noindent Let $Z_{\lambda }$ be a Poisson random variable with parameter $\lambda>0$, that is $Z_{\lambda }\sim \mathcal{P}(\lambda)$. Then the random variable 

\begin{equation*}
\frac{Z_{\lambda }-\lambda }{\sqrt{\lambda }}
\end{equation*}

\noindent weakly converges to standard Gaussian random variable $X$, that is $X \sim \mathcal{N}(0,1)$, as $\lambda \rightarrow +\infty$.\\
  
\bigskip \noindent \textbf{Proof}. Let us use the moment generating functions. The moment generating function of $Z_{\lambda }\sim \mathcal{P}(\lambda)$ is 
\begin{equation*}
\Psi _{Z_{\lambda }}(t)=\exp (\lambda (e^{t}-1)).
\end{equation*}

\noindent Set
\begin{equation*}
Y(\lambda )=\frac{Z_{\lambda }}{\sqrt{\lambda }}=\frac{Z_{\lambda}-\mathbb{E}(X)}{\sigma
_{Z_{\lambda}}}.
\end{equation*}

\noindent We have 
\begin{equation*}
\Psi _{Y(\lambda )}(u)=e^{-\sqrt{\lambda }}\times \varphi _{Z}(u/\sqrt{%
\lambda })=e^{-\sqrt{\lambda }}\times \exp (\lambda (e^{u/\sqrt{\lambda }%
}-1)).
\end{equation*}

\noindent As $\lambda \rightarrow \infty$, we may use the following expansion 
\begin{equation*}
\lambda \left(e^{u/\sqrt{\lambda }}-1\right)=\lambda (1+\frac{u}{\sqrt{\lambda }}+\frac{%
u^{2}}{2\lambda }+O(\lambda ^{-3/2})-1
\end{equation*}%
\begin{equation*}
=u\sqrt{\lambda }+\frac{u^{2}}{2}+O(\lambda ^{-1/2}).
\end{equation*}

\noindent Hence
\begin{equation*}
\Psi _{Y(\lambda )}(u)=\exp (\frac{u^{2}}{2}+O(\lambda ^{-1/2}))\rightarrow
\exp (u^{2}/2).
\end{equation*}

\noindent We conclude that 
\begin{equation*}
\frac{Z_{\lambda}}{\sqrt{\lambda }}\rightarrow \mathcal{N}(0,1)
\end{equation*}

\noindent as $\lambda \rightarrow \infty .$

\subsection{Convergence of a sequence of Binomial random variables to a standard Gaussian random variable} \label{cv.review.subsec.Bin} $ $\\

\noindent Let $X_{n}$ be a $\mathcal{B}(n,p)$ random variable with $p\in ]0,1[$ which is fixed and $n\geq 1$. Then, as $n\rightarrow \infty$, 
\begin{equation}
Z_{n}=\frac{X_{n}-np}{\sqrt{npq}}\rightsquigarrow \mathcal{N}(0,1). \label{cv.cltBin}
\end{equation}

\noindent \textbf{Proof}. Let us use the moment generating functions. Let $X\sim \mathcal{B}(n,p)$. We have

\begin{equation*}
\Psi _{X_{n}}(u)=(q+pe^{u})^{n}.
\end{equation*}

\noindent where $q=1-p$. Then 
\begin{equation}
\Psi _{(X_{n}-np)/\sqrt{npq}}(u)=e^{-\sqrt{np/q}}\times \Psi _{X_{n}}(u/%
\sqrt{npq}),  \label{binom00}
\end{equation}

\noindent with 
\begin{equation*}
\Psi _{X}(u/\sqrt{npq})=(q+pe^{u/\sqrt{npq}})^{n}.
\end{equation*}

\noindent The idea behind the coming computations is to use a second order expansion of $e^{u/\sqrt{npq}}$ in the neighborhood of $0$ as 
$n\rightarrow \infty$ and $u$ fixed. We get an expression of the form $1+v_{n}$, where $v_{n}$ tends to zero. Finally an expansion of the logarithm function $\log(1+v_{n})$ of order 2 is operated.\\

\noindent Hence, as $n\rightarrow \infty $ and $u$ is fixed, we have, 
\begin{equation*}
e^{u/\sqrt{npq}}=1+\frac{u}{\sqrt{npq}}+\frac{u^{2}}{2npq}+O(n^{-3/2}).
\end{equation*}

\noindent Next,
\begin{equation*}
(q+pe^{u/\sqrt{npq}})=1+u\sqrt{p/nq}+\frac{u^{2}}{2nq}+O(n^{-3/2})=1+v_{n}
\end{equation*}

\noindent with
\begin{equation*}
v_{n}=u\sqrt{p/nq}+\frac{u^{2}}{2nq}+O(n^{-3/2})\rightarrow 0.
\end{equation*}

\noindent Thus, 
\begin{eqnarray*}
\log \left(1+u\sqrt{p/nq}+\frac{u^{2}}{2nq}+O(n^{-3/2})\right)&=&\log (1+v_{n})\\
&=&v_{n}-\frac{1}{2}v_{n}^{2}+O(v_{n}^{3})\\
&=&u\sqrt{p/nq}+\frac{u^{2}}{2nq}-\frac{pu^{2}}{2nq}+O(n^{-3/2}).
\end{eqnarray*}

\noindent Hence 
\begin{eqnarray*}
\Psi _{X_{n}}(u/\sqrt{npq})&=&(q+pe^{u/\sqrt{npq}})^{n}=\exp (n\log (q+pe^{u/\sqrt{npq}}))\\
&=&\exp \left(n\left(u\sqrt{p/nq}+\frac{u^{2}}{2nq}-\frac{pu^{2}}{2nq}+O(n^{-3/2})\right)\right)\\
&=&\exp \left(u\sqrt{np/q}+\frac{u^{2}}{2q}-\frac{pu^{2}}{2q}+O(n^{-1/2})\right)\\
&=&e^{u\sqrt{np/q}}e^{u^{2}/2+O(n^{-1/2})}.
\end{eqnarray*}

\noindent By going back to (\ref{binom00}), we arrive at  
\begin{equation*}
\Psi _{(X_{n}-np)/\sqrt{npq}}(u)\rightarrow \exp (u^{2}/2).
\end{equation*}

\noindent This is 
\begin{equation*}
(\beta (n,p)-np)/\sqrt{npq}\rightarrow_{w} \mathcal{N}(0,1) \ \ as \ \ n \rightarrow +\infty.
\end{equation*}

\noindent QED.\\

\bigskip \noindent \textbf{Remark}. We will come back for a direct proof of this result using the
central limit theorem stated just below.\\

\subsection{Convergence of a sequence of Negative Binomial random variables to a standard Gaussian random variable} \label{cv.review.subsec.BinNeg} $ $\\

\noindent Let $Y_{k}$ be a sequence of $\mathcal{NB}(k,p)$ random variable with $p\in ]0,1[$ which is fixed and $k\geq 1$. Then, as $k\rightarrow \infty$, 
\begin{equation}
Z_k=\frac{p(Y_k-\frac{k}{p})}{\sqrt{qk}}\rightsquigarrow \mathcal{N}(0,1). \label{cv.cltBinNeg}
\end{equation}

\bigskip \textbf{Proof}. \noindent Let $t \in \mathbb{R}$ fixed. We write

$$
Z_k=\frac{p}{\sqrt{kq}}Y_k-\sqrt{\frac{k}{q}}
$$

\bigskip \noindent and remind that 

$$
\varphi_{Y_k}(t)=\left( \frac{pe^t}{1-qe^t}\right)^k.
$$

\begin{eqnarray*}
\varphi_{Z_k}(t) &=& \exp(-t\sqrt{\frac{k}{q}}) \varphi_{Y_k}\left(\frac{pt}{\sqrt{kq}}\right)\\
&=&  \exp\left(-t\sqrt{\frac{k}{q}}\right) \left(\frac{p\exp(\frac{pt}{\sqrt{kq}})}{1-q\exp(\frac{pt}{\sqrt{kq}})}\right)^k\\
&=&  \exp\left(-t\sqrt{\frac{k}{q}}\right) \exp\left( k\log\left(\frac{p\exp(\frac{pt}{\sqrt{kq}})}{1-q\exp(\frac{pt}{\sqrt{kq}})}\right)\right)\\
&=:& \exp\left(-t\sqrt{\frac{k}{q}}\right)\exp(k\log (B_k))
\end{eqnarray*}

\bigskip \noindent with

\begin{eqnarray*}
B_k&=:&\frac{p\exp(\frac{pt}{\sqrt{kq}})}{1-q\exp(\frac{pt}{\sqrt{kq}})}\\
&=:&\frac{B_{k,1}}{B_{k,2}}.
\end{eqnarray*}

\bigskip \noindent Now, since $pt/\sqrt{kq}\rightarrow 0$ as $k \rightarrow +\infty$, we get the second order expansions of 
$B_{k,1}$ and $B_{k,1}$ as follows :

$$
B_{k,1}=p\left( 1+\frac{pt}{\sqrt{kq}}+\frac{p^2t^2}{2kq}+ O(k^{\frac{-3}{2}})\right)
$$

\bigskip \noindent and

$$
B_{k,2}=1-q\left( \frac{pt}{\sqrt{kq}}+\frac{p^2t^2}{2kq}+ O(k^{\frac{-3}{2}})\right)=p\left( 1-t\sqrt{\frac{q}{k}}-\frac{pt^2}{2k}+ O(k^{\frac{-3}{2}})\right)
$$

\bigskip \noindent It comes that

\begin{eqnarray*}
B_k&=&\frac{B_{k,1}}{B_{k,2}}\\
&=&\frac{p\left( 1+\frac{pt}{\sqrt{kq}}+\frac{p^2t^2}{2kq}+ O(k^{\frac{-3}{2}})\right)}{p\left( 1-t\sqrt{\frac{q}{k}}-\frac{pt^2}{2k}+ O(k^{\frac{-3}{2}})\right)}\\
&=:&\frac{1+a_k}{1-b_k},
\end{eqnarray*}

\noindent with

$$
\log(1+a_k)=a_k-\frac{1}{2}a_k^2+O(a_k^3)
$$

\bigskip \noindent and

$$
\log(1-b_k)=-b_k-\frac{1}{2}b_k^2+O(b_k^3).
$$

\bigskip \noindent Hence, we get 

\begin{eqnarray*}
\log(B_k)&=&\frac{pt}{\sqrt{kq}}+t\sqrt{\frac{k}{q}}+\frac{t^2}{2k}+O(k^{\frac{-3}{2}}).
\end{eqnarray*}

\bigskip \noindent Next, we have

\begin{eqnarray*}
\exp\left(k\log(B_k)\right)&=&\exp\left(pt\sqrt{\frac{k}{q}}+t\sqrt{kq}+\frac{t^2}{2}+O(k^{\frac{-1}{2}})\right)\\
&=&\exp\left((1-q)t\sqrt{\frac{k}{q}}+t\sqrt{kq}+\frac{t^2}{2}+O(k^{\frac{-1}{2}})\right)\\
&=&\exp\left(t\sqrt{\frac{k}{q}}+\frac{t^2}{2}+O(k^{\frac{-1}{2}})\right).
\end{eqnarray*}

\bigskip \noindent Finally, we get 

\begin{eqnarray*}
\varphi_{Z_k}(t)&=&\exp\left(-t\sqrt{\frac{k}{q}}\right) \exp\left( t\sqrt{\frac{k}{q}}+\frac{t^2}{2}+O(k^{\frac{-1}{2}})\right)\\
&=&\exp\left( \frac{t^2}{2}+O\left(k^{-1/2}\right) \right)\\
&\rightarrow&\varphi_{\mathcal{N}(0,1)}(t) \ as \ k\rightarrow +\infty. \ \square
\end{eqnarray*}

\bigskip \subsection{Simple Central Limit Theorem in $\mathbb{R}$.} \label{cv.review.subsec.cltR} $ $\\

\noindent The two last cases are special cases of a more general weak convergence theorem, called the central limit theorem (\textit{CLT}) of Probability Theory. We say that a sequence of real random variables $(X_{n})_{n\geq 1}$, for which each $X_n$ has a positive finite second moment, satisfies the \textit{CTL} property if and only if  
\begin{equation*}
\frac{X_{n}-E(X_{n})}{\sigma _{X_{n}}}
\end{equation*}

\noindent weakly converges to Gaussian standard random variable. This, of course, is not always true. Here, we will see a simple case. Later, we will give a global solution of this problem in $\mathbb{R}$.\\

\noindent Let $X_{1},X_{2},...$ be a sequence of real valued random variables which are independent and identically distributed (\textit{iid}) random variables with common distribution function $F$ with 
\begin{equation*}
E(X_{i})=\mu =\int xdF(x)=0,\sigma _{X_{i}}^{2}=\sigma ^{2}=\int (x-\mu)^{2}dF(x)=1.
\end{equation*}

\noindent Put, for $n\geq 1$, 
\begin{equation*}
S_{n}=X_{1}+...+X_{n}.
\end{equation*}

\noindent We have, as $n\rightarrow \infty$,

\begin{equation*}
\frac{S_{n}}{\sqrt{n}}\rightarrow \mathcal{N}(0,1).
\end{equation*}

\noindent  \textbf{Proof}. Consider the common characteristic function 
\begin{equation*}
\mathbb{R}\ni u\hookrightarrow \Phi _{X_{i}}(u)=E(e^{iuX_{i}})=\Psi (u).
\end{equation*}

\noindent Since the second moment exists, we have the following expansion at order 2, 
\begin{equation*}
\Phi (u)=1+iu\Phi ^{\prime }(0)+\frac{1}{2}u^{2}\Phi ^{\prime \prime
}(0)+O(u^{3})
\end{equation*}%
\begin{equation*}
=1-\frac{1}{2}u^{2}+O(u^{2})
\end{equation*}%

\noindent since
\begin{equation*}
\Phi ^{\prime }(0)=i\text{ }\mathbb{E}(X)=0,\text{ }\Phi ^{\prime \prime
}(0)=-\mathbb{E}(X^{2})=-1.
\end{equation*}

\noindent Thus  
\begin{equation*}
\Phi _{S_{n}/\sqrt{n}}(u)=(\Phi (u/\sqrt{n}))^{n}.
\end{equation*}

\noindent For $u$ fixed, as $n\rightarrow \infty$,

\begin{eqnarray*}
\Phi _{S_{n}/\sqrt{n}}(u) &=&(\Phi (u/\sqrt{n}))^{n}=\exp \left(n\log (1-\frac{%
u^{2}}{2n}+O(n^{-3/2}))\right) \\
&=&\exp (n(-\frac{u^{2}}{n}+O(n^{-3/2})) \\
&=&\exp (-u^{2}/2+O(n^{-1/2})) \\
&\rightarrow &\exp (-u^{2}/2).
\end{eqnarray*}

\noindent We just established  
\begin{equation*}
\frac{S_{n}}{\sqrt{n}}\rightarrow \mathcal{N}(0,1) \ \ as\ n\rightarrow +\infty.
\end{equation*}

\noindent In a more general case of an \textit{iid} sequence of random variables $X_1$, $X_2$, ... with
\begin{equation*}
\mathbb{E}(X_{i})=\mu =\int xdF(x)=\mu ,\sigma _{X_{i}}^{2}=\sigma ^{2}=\int (x-\mu)^{2}dF(x)=\sigma^2,
\end{equation*}

\noindent we apply the former result to the sequence $(X_i-\mu)/\sigma$, $i=1,2,...$ to get

\begin{equation*}
\frac{1}{\sigma \sqrt{n}} (S_{n}-n\mu) \rightarrow \mathcal{N}(0,1).
\end{equation*}

\bigskip \noindent Let us give two examples of applications of the simple central limit theorem on the binomial trials.\\

\noindent \textbf{Example 1 : Weak convergence of the binomial random variable}.\\

\noindent We are going to prove the result (\ref{cv.cltBin}) of Subsetion \ref{cv.review.subsec.Bin} concerning the weak law of a sequence of binomial random variables as the number of trials, $n$, increases  while the probability of success, $p \in ]0,1[$, is fixed. So we keep the notation of that subsection.\\

\noindent We know from the earlier courses on elementary Probability Theory we may find in a considerable number of books, especially in \cite{ept-en2016}, with the current Probability Theory and Statistics Series, in Chapter 2, Lemma 1, that if $X_n \sim \mathcal{B}(n,p)$, then $X_{n}$ is the sum of $n$ independent Bernoulli $\mathcal{B}(p)$ random variables $Y_{1},...,Y_{n}$ such that

\begin{equation*}
X_{n}=Y_{1}+...+Y_{n}.
\end{equation*}

\noindent For each of the $Y_{i}$'s random variables, we have

\begin{equation*}
\mathbb{E}(Y_{i})=p\text{ and }\sigma ^{2}=\mathbb{V}ar(Y_{i})=pq\text{ where }q=1-p.
\end{equation*}

\noindent Then, the random variable $Z_n$ in Formula (\ref{cv.cltBin}) becomes

\begin{equation*}
Z_{n}=\frac{X_{n}-np}{\sqrt{npq}}=\frac{1}{\sigma \sqrt{n}}
\sum_{i=1}^{n}(Y_{i}-\mathbb{E}(Y_{i})).
\end{equation*}

\noindent Hence, the weak convergence of  $Z_{n}$ to $\mathcal{N}(0,1)$  as $n \rightarrow +\infty$, is a consequence of the simple central limit standard on $\mathbb{R}$.\\

\bigskip \noindent \textbf{Remark}. This proof is quick and beautiful. The first proof is still useful. Because, we may be in a position to teach this result at a level where the central limit theorem is not available. Besides, this proof is part of History. In the same spirit, the oldest proof of this result goes back to 1732 by \textit{de} Moivre and to 1801 by Laplace (see Lo\`{e}ve \cite{loeve}, page 23). These historical methods can also be found in \cite{ept-en2016} and in \cite{ept-fr2016} with a writing which is appropriate to beginners of first year of university.\\

\bigskip \noindent \textbf{Example 2 : Negative Binomial Law}.\\

\noindent For a fixed integer $k\geq 1$, a Negative Binomial random variable $X_{k}$ is defined relatively to Bernoulli trails of probability of success $p \in ]0,1[$. The number of repetitions of a Bernoulli experiment of parameter $p$ which is necessary to obtain $k$ successes is
said to follow a Negative Binomial random variable with parameters $k$ and $p$, denoted by  $X_{k} \sim \mathcal{NB}(k,p)$. For $k=1$, it is said that $X_{1}$ follows a geometric law with parameter $p$, denoted $X_1 \sim \mathcal{G}(p)$.\\

\noindent Similarly to the sequence of binomial random variable, we may apply the central limit theorem to the sequence of negative binomial random variables $X_k$, $k\geq 1$ to get the following result

\begin{equation}
Z_{n}=\frac{p(X_{k}-\frac{k}{p})}{\sqrt{nq}} \rightsquigarrow \mathcal{N}(0,1) \text{ as } k \rightarrow +\infty. \label{cv.cltNegBin01}
\end{equation} 

\noindent To this purpose, the reader may find more details in classical elementary books in probability theory, for instance in \cite{ept-en2016} or in \cite{ept-fr2016}, Chapters 2 and 3. In Chapter 2 of these monographs, Lemma 2, ensures that a $\mathcal{NB}(k,p)$ random variables $X_{k}$ is the sum of $k$ independent and geometric $\mathcal{G}(p)$ random variables 
$Y_{1},...,Y_{k}$ such that

\begin{equation*}
X_{k}=Y_{1}+...+Y_{n},
\end{equation*}

\noindent and for each of these random variables $Z_i$'s, we have

\begin{equation*}
\mathbb{E}(Y_{i})=\frac{1}{p}\text{ and }\sigma ^{2}=\mathbb{V}ar(Y_{i})=\frac{q}{p^{2}}\text{ where }q=1-p.
\end{equation*}

\noindent Thus, by the central limit theorem 

\begin{equation*}
Z_{n}=\frac{p(X_{k}-\frac{k}{p})}{\sqrt{nq}}=\frac{1}{\sigma \sqrt{n}}\sum_{i=1}^{n}(Y_{i}-\mathbb{E}(Y_{i}))\rightsquigarrow \mathcal{N}(0,1) \text{ as } k \rightarrow +\infty. \label{cv.cltNegBin02}
\end{equation*}

\noindent which proves (\ref{cv.cltNegBin01}).\\

\subsection{Limit laws in Extreme value Theory} \label{cv.review.subsec.evt} $ $\\

\noindent \bigskip \noindent Consider $X_{1},X_{2}$, ...  a sequence of \textit{iid} random variables with common distribution function $F$. Put for each $n\geq 1$,

\begin{equation*}
M_{n}=\max (X_{1},...,X_{n}).
\end{equation*}

\bigskip \noindent Recall that for any $x \in \mathbb{R}$ 
\begin{equation*}
P(M_{n}\leq x)=F(x)^{n},x\in \mathbb{R}.
\end{equation*}

\noindent The basic problem of extreme value theory is finding sequences $\left(a_{n}>0\right) _{n\geq 1}$ and $\left( b_{n}\right) _{n\geq 1}$ such that  
\begin{equation*}
\frac{M_{n}-b_{n}}{a_{n}}
\end{equation*}

\noindent weakly converges to some random variable $Z$, 
\begin{equation*}
\frac{M_{n}-b_{n}}{a_{n}}\rightsquigarrow Z.
\end{equation*}

\noindent If this holds, we write $F \in D(F_Z)$.\\

\noindent We are going to give three examples corresponding to the three nontrivial cases.\\

\bigskip \noindent \textbf{(a)} Let $\Lambda$ be a Gumbel random variable of distribution function
$$
\Lambda(x)=exp(-e^{-x}), \ \ x \in \mathbb{R}.
$$

\noindent Let the $X_i$'s are standard exponential random variables, $X_i \sim \mathcal{E}(1)$, with 
\begin{equation*}
F(x)=(1-\exp (-x))1_{(x\geq 0)}, \text{ } x\in \mathbb{R}.
\end{equation*}

\noindent We have, as $n\rightarrow +\infty$,

\begin{equation}
M_{n}-\log n \rightsquigarrow \Lambda. \label{portal.example1a}
\end{equation}

\bigskip  \noindent \textbf{Proof}. By using the distribution functions, we want to prove that for any $x\in \mathbb{R}$,

\begin{equation}
\mathbb{P}(M_{n}-\log n \leq x) \rightarrow \Lambda(x). \label{portal.example1b}
\end{equation}

\bigskip \noindent  \textbf{Proof}. We are going to show, by using the distribution functions, that

\noindent Indeed, we have

\begin{equation*}
\mathbb{P}(M_{n}-\log n\leq x)=P(M_{n}\leq x+\log n)=F(x+\log n)^{n}.
\end{equation*}

\bigskip \noindent But for any $x\in \mathbb{R},$ $x+\log n\geq 0$ for $n\geq
\exp (-x).$ Then for large values of $n$, $P(M_{n}\leq x+\log n)=(1-\exp
(-x-\log n))^n$ and next for any $x\in \mathbb{R}$ and for $n$ large enough,

\begin{equation*}
\mathbb{P}(M_{n}-\log n\leq x)=\left(1-\frac{e^{-x}}{n}\right)\rightarrow e^{-e^{-x}}=\Lambda(x).
\end{equation*}

\noindent So (\ref{portal.example1b}) holds and so does (\ref{portal.example1a}), that is : $X \in D(\Lambda)$.\\

\bigskip \noindent \textbf{(b)} Let $FR(\alpha)$ a Fr\'echet random variable with parameter $\alpha>0$, with distribution function

$$
\phi_{\alpha}(x)=exp(-x^{-\alpha})1_{(x\geq 0)},
$$

\noindent where $1_A$ is the indicator function of the set $A$ that assigns the value one to elements of $A$ and zero to elements of the complementary of $A$.\\

\noindent Let the $X_i$'s be Pareto random variables with parameter $\alpha>0$, $X \sim \mathcal{P}ar(\alpha)$, with common distribution function

\begin{equation*}
F(x)=(1-x^{-\alpha }) \ 1_{(x\geq 1)},  \text{ } x\in \mathbb{R}
\end{equation*}

\noindent Then, as $n \rightarrow+\infty$, we have

\begin{equation}
n^{-1/\alpha }M_{n} \rightsquigarrow FR(\alpha). (\label{portal.example2a})
\end{equation}

\bigskip \noindent \textbf{Proof}. We want to prove that for any $x\in \mathbb{R}$, we have as $n\rightarrow +\infty$,

\begin{equation}
\mathbb{P}( n^{-1/\alpha }M_{n} \leq x)  \rightarrow \phi_{\alpha}(x). \label{portal.example2b}
\end{equation}

\bigskip \noindent The observations $X_i$'s are non-negative since the support of a $\mathcal{P}ar(\alpha)$ law is $\mathbb{R}_{+}$. So the maxima $M_{n}$ are non-negative for any $n\geq 1$. We may discuss two cases.\newline

\bigskip \noindent Case $x\leq 0$. In this case, we have
\begin{equation*}
\mathbb{P}(n^{-1/\alpha }M_{n}\leq 0)=0=\phi_{\alpha }(x),
\end{equation*}

\bigskip \noindent and then (\ref{portal.example2b}) holds.\\

\bigskip \noindent Case $x>0$. In this case 
\begin{equation*}
P(n^{-1/\alpha }M_{n}\leq x)=P(M_{n}\leq n^{1/\alpha }x).
\end{equation*}

\bigskip \noindent For large values of $n$, we have $n^{1/\alpha }x>1$ (take $n\geq x^{-\alpha}$, to ensure that) and for these values, 
\begin{eqnarray*}
\mathbb{P}(n^{-1/\alpha }M_{n} &\leq &x)=F(n^{1/\alpha }x)^{n}=(1-(n^{1/\alpha
}x)^{-\alpha })^{n} \\
&=&\left(1-\frac{x^{-\alpha }}{n}\right)^{n}\rightarrow \exp (-x^{-\alpha })=\phi _{\alpha }(x).
\end{eqnarray*}

\bigskip \noindent So  (\ref{portal.example2b}) holds for $x>0$. But putting together the two cases, we have $F \in D(FR(\alpha))$.\\

\bigskip \noindent \textbf{(c)} Let $W(\beta)$ be a Weibull random variable with parameter $\beta>0$, with distribution function

$$
\psi_{\alpha}(x)=exp(-(-x)^{\beta})1_{x\leq 0} + 1_{(x>0)}.
$$

\noindent Let the $X_i$'s be uniformly distributed on $(0,1)$ with probability distribution function : 
\begin{equation*}
F(x)=x1_{(0\leq x\leq 1)}+1_{(x\geq 1)}, \text{ } x\in \mathbb{R}.
\end{equation*}

\noindent We have

\begin{equation}
n(M_{n}-1)\overset{d}{\rightarrow }W(1) \text{ as } n\rightarrow +\infty. \label{portal.example3a}
\end{equation}

\bigskip \noindent \textbf{Proof}. We have to prove that for any $x\in \mathbb{R}$, as $n\rightarrow +\infty$,

\begin{equation}
\mathbb{P}(n(M_{n}-1)\leq x)=F\left(1+\frac{x}{n}\right)^{n} \rightarrow \psi_{1}(x). \label{portal.example3b}
\end{equation}

\noindent We have two cases.\newline

\bigskip \noindent Case $x\geq 0$. We see that $1+x/n$ is non-negative for $n\geq 1$ and 
\begin{equation*}
P(n(M_{n}-1)\leq x)=F(1+\frac{x}{n})^{n}=1=\psi _{1}(x)
\end{equation*}

\bigskip \noindent and we see that (\ref{portal.example3b}) holds for $x\geq 0$.\\

\noindent Case $x<0$. For large values of $n$, we have $0\leq 1+x/n\leq 1$ (take $0\geq -x \geq n$, to get it) and for these values of $n$, 
\begin{eqnarray*}
P(n(M_{n}-1) &\leq &x)=F\left(1+\frac{x}{n}\right)^{n} \\
&=&\left(1+\frac{x}{n}\right)^{n}\rightarrow e^{x}=\psi _{1}(x).
\end{eqnarray*}

\bigskip \noindent Then (\ref{portal.example3b}) also holds for $x<0$ and then (\ref{portal.example3b}) holds for any $x\in \mathbb{R}$,
 
\begin{equation*}
\mathbb{P}(n(M_{n}-1)\leq x)\longrightarrow \psi _{1}(x).
\end{equation*}

\noindent Conclusion : $F \in D(W(1))$.\\

\bigskip \noindent \textbf{Summary}. In Uni-variate Extreme Value Theory (UEVT), it is proved that the three non-degenerated possible limits are the three we gave above. You will have the opportunity to go deep in that theory in the book of this series \cite{wc-ib-sp-ang}.

\section{Examples of Convergence in $\mathbb{R}^{k}$} \label{cv.review.sec4}

\subsection{Simple Central Limit in $\mathbb{R}^{k}$} \label{cv.review.subsec.clt} $ $\\

\noindent We now move to the Central Limit Theorem in $\mathbb{R}^{k}$ in the \textit{iid} case. Let $X_{1},X_{2},....$ be centered \textit{iid} $\mathbb{R}^k$-random variables with common finite variance-covariance matrix $\Sigma =(\sigma _{ij})_{1\leq i\leq k,1\leq j\leq k}$, that is

\begin{equation*}
\sigma _{ij}=Cov(X_{i},X_{j}) \in \mathbb{R}, \ \ 1\leq i,j\leq k.
\end{equation*}

\bigskip \noindent Set the partial sums 
\begin{equation*}
S_{n}=X_{1}+X_{2}+...+X_{n}, \ \ n\geq 1.
\end{equation*}

\bigskip \noindent We have the following central limit theorem on $\mathbb{R}^{k}$, 
\begin{equation*}
S_{n}/\sqrt{n}\rightsquigarrow \mathcal{N}(0,\Sigma) \ \ as \ \ n\rightarrow +\infty.
\end{equation*}

\bigskip \noindent \textbf{Proof}. The matrix $\Sigma$ is symmetrical and non-negative in the sense that for any
 $u\in \mathbb{R}^{k}$
\begin{equation*}
\text{ }^{t}u\Sigma u=\text{ }^{t}u\mathbb{E}(XX^{\prime })u=\mathbb{E}%
((\text{ }^{t}Xu)(\text{ }^{t}Xu))=\mathbb{E}((\text{ }^{t}Xu)^{2})\geq 0.
\end{equation*}

\noindent  By the matrices theory, $\Sigma$ has $k$ non-negative eigenvalues $\lambda_{1}$, $\lambda_{2}$,...,$\lambda_{k}$ and there exists
a orthogonal $(k,k)$-matrix $T$ such that
\begin{equation*}
\text{ }^{t}T\Sigma T=diag(\lambda _{1},\lambda _{2},...,\lambda _{n})=\Lambda .
\end{equation*}

\noindent Set  
\begin{equation*}
Y_{i}=\text{ }^{t}TX_{i}.
\end{equation*}

\noindent The random variables $Y_{i}$ are centered, \textit{iid} and have common variance-covariance matrix equal to
\begin{equation*}
\Sigma_{Y}=\text{ }^{t}T\Sigma T=\Lambda .
\end{equation*}

\noindent This means that the components of each $Y_{i}$ are uncorrelated and have variances equal to  $\lambda _{1},\lambda
_{2},...,\lambda _{n}$. Set  
\begin{equation}
M_{n}=\frac{1}{\sqrt{n}}(Y_{1}+Y_{2}+...+Y_{n})=\text{ }^{t}T\left(\frac{S_{n}}{\sqrt{n}}\right).  \label{rk00}
\end{equation}

\bigskip \noindent For any $A=\text{ }^{t}(a_{1},a_{2},...,a_{k})\in \mathbb{R}^{k}$,
\begin{equation*}
<A,M_{n}>=\frac{1}{\sqrt{n}}\sum_{i=1}^{n}<A,Y_{i}>.
\end{equation*}

\noindent The variables $<A,Y_{i}>$ then are \textit{iid} and have common variance 

\begin{equation*}
\mathbb{E}<A,Y_{i}>^{2}=\sum_{i=1}^{n}a_{i}^{2}\lambda _{i}=\text{ }%
^{t}A\Lambda A,
\end{equation*}

\noindent because of the mack of correlation between the components of each  $Y_{i}$. We may apply the central limit theorem in $\mathbb{R}$ to get 
 
\begin{equation*}
<A,M_{n}>\rightarrow \mathcal{N}(0,\sum_{i=1}^{i=n}a_{i}^{2}\lambda _{i})=\mathcal{N}(0,\text{ }^{t}A\Lambda A)
\end{equation*}

\noindent But $\mathcal{N}(0,\text{ }^{t}A\Lambda A)$ is the law of a Gaussian random variable that is the linear transform $\text{ }^{t}AZ=<A,Z>$ of $Z$, where 
$Z$ follows the $\mathcal{N}(0,\Lambda )$ law. Then
 
\begin{equation*}
\forall A\in \mathbb{R}^{k},<A,M_{n}>\rightsquigarrow <A,Z>.
\end{equation*}

\noindent In terms of characteristic functions, we have for any  $t\in \mathbb{R}$ and for any $A\in \mathbb{R}^{k}$,

\begin{equation*}
\mathbb{E}\exp (it<A,M_{n}>)\rightarrow \mathbb{E}\exp (it<A,Z>).
\end{equation*}

\noindent For $t=1$, we have for any  $A\in \mathbb{R}^{k}$%
\begin{equation*}
\Phi _{M_{n}}(A)=\mathbb{E}\exp (i<A,M_{n}>)\rightarrow \Phi _{Z}(A)=\mathbb{E}\exp (i<A,Z>).
\end{equation*}

\noindent This means that  
\begin{equation*}
M_{n}\rightsquigarrow Z.
\end{equation*}

\noindent This, Point (e) of Theorem \ref{cv.review.rk} and (\ref{rk00}) together implies that

\begin{equation*}
S_{n}/\sqrt{n}=TM_{n}\rightsquigarrow TZ
\end{equation*}

\noindent and then

\begin{equation*}
\text{ }^{t}TZ\sim \mathcal{N}(0,T\Lambda \text{ }^{t}T)=\mathcal{N}(0,\Sigma ).
\end{equation*}

\noindent Hence, finally, as $n\rightarrow +\infty$, 
\begin{equation*}
S_{n}/\sqrt{n}\rightsquigarrow \mathcal{N}(0,\Sigma).
\end{equation*}

\subsection{Weak Convergence of the Multinomial Law} \label{cv.review.subsec.multinomial} $ $\\

\noindent A $k$-tuple  $X_n=(X_{1,n},...,X_{k,n})$ follows a multinomial law with parameters
$n\geq 1$ and  $p=(p_{1},p_{2},...,p_{k})$, with 
\begin{equation*}
\forall (1\leq i\leq k),p_{i}>0\text{ \ }and\text{ \ }\sum_{1\leq i\leq k}p_{i}=1,
\end{equation*}

\noindent denoted $X \sim \mathcal{M}_{k}(n,p)$, if and only if its probability law is given by  
\begin{equation*}
\mathbb{P}(X_{1,n}=n_{1},...,X_{k,n}=n_{k})=\frac{n!}{n_{1}!\times ...\times
n_{k}!}p_{1}^{n_{1}}\times p_{2}^{n_{2}}\times ...\times p_{k}^{n_{k}},
\end{equation*}

\noindent where $(n_{1},...,n_{k})$ satisfies 
\begin{equation*}
\forall (1\leq i\leq k),\text{ }n_{i}\geq 0\text{ \ et \ }\sum_{1\leq i\leq
k}n_{i}=n.
\end{equation*}

\noindent A random variable following the $\mathcal{M}_{k}(n,p)$ law may be generated as follows :\\

\noindent Consider a random experiment with $k$ possible outcomes $E_{i},1$\ $\leq i\leq k$, each of them occurring with $p_{i}>0$. After $n$ repetitions, the number of occurrences $X_{i,n}$ of each $E_i$ is observed for $i=1,...,k$. The resulting random vector follows the $\mathcal{M}_{k}(n,p)$ law. Each individual coordinate $X_{i,n}$ follows the Binomial law $\mathcal{B}(n,p_{i})$.\\

\noindent We have the following weak convergence result.\\

\noindent Put 
\begin{equation}
Z_{n}=\text{ }^{t}\left(\frac{X_{1,n}-np_{1}}{\sqrt{np_{1}}},...,\frac{X_{k,n}-np_{k}}{\sqrt{np_{k}}}\right) \rightsquigarrow \mathcal{N}_{k}(0,\Sigma) \ \ as \ \ n\rightarrow +\infty, \label{cv.tclMultiNom01}
\end{equation}

\noindent where $\Sigma$ is a $(k,k)$-matrix with elements $\Sigma_{i,i}=1-p_{i}$ and $\Sigma_{i,j}=-\sqrt{p_{i}p_{j}}$, $1\leq i,j\leq k$.\\ 

\noindent \textbf{Important remark}. This result has a significant number of applications. We may cite two of them. It is used to have the finite-distribution function of the empirical process. It also serves as the foundations of Chi-square statistical tests that will be studied latter in one the books of this series.\\

\noindent \textbf{Proof}. We have at least two ways of proving the result. The first is based on the use of the moment generating function on logarithm expansions. The second exploits the central limit theorem in $\mathbb{R}^k$ we just proved.\\

\noindent \textbf{First proof}. We already know from \cite{bmtp} that its moment generating function is
\begin{equation*}
\phi_{X_n}(u)=\left(\sum_{1\leq i\leq k}p_{i}e^{u_{i}}\right)^{n}.
\end{equation*}

\noindent We have
\begin{equation*}
Z_n=AX+B,
\end{equation*}

\noindent where $A$ is the diagonal matrix

\begin{equation*}
A=diag\left(\frac{1}{\sqrt{np_{1}}}, \frac{1}{\sqrt{np_{2}}}, ...,\frac{1}{\sqrt{np_{k}}}\right)
\end{equation*}

\noindent and 
\begin{equation*}
B=\left( 
\begin{array}{c}
-\sqrt{np_{1}} \\ 
-\sqrt{np_{2}} \\ 
... \\ 
-\sqrt{np_{k}}%
\end{array}%
\right) .
\end{equation*}

\noindent Thus 
\begin{eqnarray*}
\phi _{Z_{n}}(u)&=&\exp (<B,u>)\times \phi _{X}(\text{ }^{t}Au)\\
&=& \left( \exp \left(\sum_{1\leq i\leq k}-\sqrt{np_{i}}u_{i}\right) \right) \times \left(\sum_{1\leq i\leq
k}p_{i}e^{u_{i}/\sqrt{np_{i}}}\right)^{n}.
\end{eqnarray*}

\noindent Let $u$ be fixed. For each fixed $i$, $1\leq i \leq k$, $u_{i}/\sqrt{np_{i}}
\rightarrow +\infty$\ as $n\rightarrow \infty $\ since $p_{i}>0$. We have the expansion 
\begin{equation*}
e^{u_{i}/\sqrt{np_{i}}}=1+u_{i}/\sqrt{np_{i}}+\frac{1}{2}\frac{u_{i}^{2}}{%
np_{i}}+O(n^{-3/2}).
\end{equation*}

\noindent Next
\begin{eqnarray*}
A&=&\left(\sum_{1\leq i\leq k}p_{i}e^{u_{i}/\sqrt{np_{i}}}\right)^{n}=\exp \left(n\log
\left(\sum_{1\leq i\leq k}p_{i}e^{u_{i}/\sqrt{np_{i}}}\right)\right).\\
&=&\exp `\left(n\log \left(1+\sum_{1\leq i\leq k}u_{i}\sqrt{p_{i/}n_{i}}+\sum_{1\leq
i\leq k}\frac{1}{2}\frac{u_{i}^{2}}{n_{i}}+O(n^{-3/2})\right)\right).
\end{eqnarray*}

\noindent Set 
\begin{equation*}
a=\sum_{1\leq i\leq k}u_{i}\sqrt{p_{i/}n}+\sum_{1\leq i\leq k}\frac{1}{2}%
\frac{u_{i}^{2}}{n}\rightarrow 0\text{ as }n\rightarrow \infty .
\end{equation*}

\noindent We have
\begin{equation*}
A=\exp (n\log (1+a)).
\end{equation*}

\noindent Let us expand $log(1+a)$ at the second order $2$. We obtain

\begin{eqnarray*}
A&=&\exp (n(a-\frac{1}{2}a^{2}+O(a^{3})).\\
&=&\exp \left(n\left(\sum_{1\leq i\leq k}u_{i}\sqrt{p_{i/}n}+\sum_{1\leq i\leq k}\frac{1%
}{2}\frac{u_{i}^{2}}{n}-\frac{1}{2}\left(\sum_{1\leq i\leq k}u_{i}\sqrt{p_{i/}n}\right)^{2}+O(n^{-3/2})\right)\right)\\
&=&\exp \left(\sum_{1\leq i\leq k}u_{i}\sqrt{np_{i}}+\sum_{1\leq i\leq k}\frac{1}{2}%
u_{i}^{2}-\frac{1}{2}\left(\sum_{1\leq i\leq k}u_{i}\sqrt{p_{i}}\right)^{2}+O(n^{-1/2})\right)\\
&=&\exp \left(\sum_{1\leq i\leq k}u_{i}\sqrt{np_{i}}\right)\times \exp \left(\sum_{1\leq i\leq
k}\frac{1}{2}u_{i}^{2}-\frac{1}{2}\left(\sum_{1\leq i\leq k}u_{i}\sqrt{p_{i}}\right)^{2}+O(n^{-1/2})\right).
\end{eqnarray*}

\noindent Putting all this together, we get 
\begin{equation*}
\phi _{Z_{n}}(u)=\exp \left(\sum_{1\leq i\leq k}\frac{1}{2}u_{i}^{2}-\frac{1}{2}
\left(\sum_{1\leq i\leq k}u_{i}\sqrt{p_{i}}\right)^{2}+O(n^{-1/2})\right)
\end{equation*}
\begin{equation*}
\rightarrow \phi _{Z}(u)=\exp \left(\sum_{1\leq i\leq k}\frac{1}{2}u_{i}^{2}-%
\frac{1}{2}\left(\sum_{1\leq i\leq k}u_{i}\sqrt{p_{i}}\right)^{2}\right).
\end{equation*}

\noindent and
\begin{equation}
\phi _{Z}(u)=\exp \left\{ \sum_{1\leq i\leq k}\frac{1}{2}%
(1-p_{i})u_{i}^{2}-\sum_{1\leq i,j\leq k}u_{i}u_{j}\sqrt{p_{i}p_{j}}\right\}, 
\label{cv7001}
\end{equation}

\noindent which is the moment generating function of a $k$-dimensional centered Gaussian vector $Z$ whose 
variance-covariance matrix is $\Sigma$. The first proof finishes here.\\

\bigskip \noindent \textbf{Second proof}. At the $i$-th repetition of the experiment, $i\in \{1,...,n\}$, we have a random vector 
\begin{equation*}
Z^{(i)}=\left( 
\begin{tabular}{l}
$Z_{1}^{(i)}$ \\ 
... \\ 
$Z_{k}^{(i)}$%
\end{tabular}%
\right) 
\end{equation*}

\noindent defined as follows : for each $1\leq r \leq k$ 

\begin{equation*}
Z_{r}^{(i)}=\left\{ 
\begin{tabular}{lll}
$1$ & if  & the outcome $E_{r}$ occurs at the $i^{th}$ experiment and any
other did not \\ 
$0$ & if & a different outcome occurs at the $i^{th}$ experiment\  
\end{tabular}%
\right. 
\end{equation*}

\noindent It is clear that each $Z^{(i)}$ is distributed as a multivariate $\mathcal{M}_{k}(1,k)$ random variable, and that the $Z^{(i)}$'s are independent.\\

\noindent Further, for a fixed $i\in \{1,...,n\}$, each $Z_{r}^{(i)}$, $1\leq r \leq k$, follows a Bernoulli law of parameter 
$p$ and only one of the $Z_{r}^{(i)}$'s ($1\leq r\leq k$) takes the value one, the others being null. This implying that 
\begin{equation*}
Z_{r}^{(i)}Z_{s}^{(i)}=0\text{ for }1\leq r\neq s\leq k,1\leq i\leq n.
\end{equation*}

\noindent We also have 
\begin{equation*}
Z_{1}^{(i)}+...+Z_{n}^{(i)}=1.
\end{equation*}

\noindent Then for each $i\in \{1,...,n\},$

\begin{equation*}
\mathbb{E}(Z_{r}^{(i)})=p_{i}\text{ and }\mathbb{V}ar(Z_{r}^{(i)})=p_{i}(1-p_{i}),\text{ }%
1\leq r\leq k
\end{equation*}

\noindent and for $1\leq r\neq s\leq k$
\begin{equation*}
cov(Z_{r}^{(i)},Z_{s}^{(i)})=\mathbb{E}(Z_{r}^{(i)}Z_{s}^{(i)})-\mathbb{E}(Z_{r}^{(i)})\mathbb{E}(Z_{s}^{(i)})=-p_{r}p_{s},
\end{equation*}

\noindent since $Z_{r}^{(i)}Z_{s}^{(i)}=0.$ So, each $Z^{(i)}$ has the
variance-covariance matrix

\begin{equation*}
\Sigma ^{0}=\left( 
\begin{tabular}{lllll}
$p_{1}(1-p_{1})$ & $-p_{1}p_{2}$ & ... & $-p_{1}p_{k-1}$ & $-p_{1}p_{k}$ \\ 
$-p_{2}p_{1}$ & $p_{2}(1-p_{2})$ & ... & $-p_{2}p_{k-1}$ & $-p_{2}p_{k1}$ \\ 
... & ... & ... & ... & ... \\ 
$-p_{k-1}p_{1}$ & $-p_{k-1}p_{2}$ & ... & $p_{k-1}(1-p_{k-1})$ & $%
-p_{k-1}p_{k}$ \\ 
$-p_{k}p_{1}$ & $-p_{k}p_{2}$ & ... & $-p_{k}p_{k-1}$ & $-p_{k}(1-p_{k})$%
\end{tabular}%
\right) 
\end{equation*}

\noindent or, in a different notation,

\begin{equation*}
\Sigma _{0}=(\sigma _{ij}^{0})_{1\leq i,j\leq k}\text{ with }\sigma
_{ij}^{0}=\left\{ 
\begin{tabular}{lll}
$p_{i}(1-p_{i})$ & if & $i=j$ \\ 
$-p_{i}p_{j}$ & if & $i\neq j$%
\end{tabular}
\right. .
\end{equation*}

\noindent After $n$ repetitions of the experiment, the random variables $%
Z^{(1)},...,Z^{(n)}$ which are independent and $\mathcal{M}_{k}(1,k)$ random
vectors add up to $X_{n},$ which means that 
\begin{equation*}
X_{n}=Z^{(1)}+...+Z^{(n)}.
\end{equation*}

\noindent By the multivariate standard central limit theorem, we have

\begin{equation*}
S_{n}=\frac{1}{\sqrt{n}}\sum_{i=1}^{n}\left( Z^{(i)}-\mathbb{E}(Z^{(i)}\right)
\rightsquigarrow Z_{0}\sim \mathcal{N}_{k}(0,\Sigma _{0}).
\end{equation*}

\noindent We easily check that

\begin{equation*}
S_{n}=\frac{1}{\sqrt{n}}\sum_{i=1}^{n}\left( Z^{(i)}-\mathbb{E}(Z^{(i)}\right)
=\text{ }^{t}\left( \frac{X_{1,n}-np_{1}}{\sqrt{n}},\frac{X_{2,n}-np_{2}}{\sqrt{n}}%
,...,\frac{X_{k,n}-np_{k}}{\sqrt{n}}\right) .
\end{equation*}

\noindent And then, we have the matrix relation 

\begin{equation*}
DS_{n}=Z_{n},
\end{equation*}

\noindent where $D$ is the diagonal matrix

\begin{equation*}
D=diag(1/\sqrt{p_{1}},...,1/\sqrt{p_{k}}).
\end{equation*}

\noindent By the continuous mapping theorem (Point (e) of Theorem \ref{cv.review.rk}),

\begin{equation*}
Z_{n}=DS_{n}\rightsquigarrow DZ_{0}\sim \mathcal{N}_{k}(0,D\Sigma _{0}D),
\end{equation*}

\noindent since $D$ is a symmetrical matrix. It remains to compute

\begin{equation*}
\Sigma =D\Sigma _{0}D=(\sigma _{ij})_{1\leq i,j\leq k}.
\end{equation*}

\noindent For $1\leq h,j\leq k,$ $(\Sigma _{0}D)_{hj}$ is the matrix product of the $%
h^{th}$ row of $\Sigma _{0}$ by the $j^{th}$ column of $D.$ By using the fact
that $D$ is diagonal, we get for $1\leq h,j\leq k,$

\begin{equation*}
(\Sigma _{0}D)_{hj}=\sigma _{hj}^{0}/\sqrt{p_{j}}.
\end{equation*}

\noindent Next  $\sigma _{ij}=(D\Sigma _{0}D)_{ij}$ is the product of $i^{th}$
row of $D$ by the $j^{th}$ column of $(\Sigma _{0}D)^{(j)}=\text{ }^{t}((\Sigma
_{0}D)_{1j},(\Sigma _{0}D)_{2j},...,(\Sigma _{0}D)_{kj})$ and then, by using
the diagonal property of $D$, we have
 
\begin{equation*}
(D\Sigma _{0}D)_{ij}=\frac{1}{\sqrt{p_{i}}}(\Sigma _{0}D)_{ij},
\end{equation*}

\noindent and then
\begin{equation*}
\sigma _{ij}=(D\Sigma _{0}D)_{ij}=\frac{1}{\sqrt{p_{i}p_{j}}}\sigma
_{ij}^{0}=\left\{ 
\begin{tabular}{lll}
$\sigma _{ii}^{0}/p_{i}=1-p_{i}$ & if & $i=j$ \\ 
-$\sqrt{p_{i}p_{j}}$ & if & $i\neq j$%
\end{tabular}%
\right. ..
\end{equation*}

\noindent We get again that
\begin{equation*}
Z_{n}\rightsquigarrow \mathcal{N}_{k}(0,\Sigma ),
\end{equation*}

\noindent where $\Sigma$ is defined in the line following Formula (\ref{cv.tclMultiNom01}) in head part of this subsection. This ends the second proof.\\

\bigskip \noindent We may conclude in a form of a proposition.\\

\begin{proposition} \label{cv.multinomial}
Le $X(n)=\text{ }^{t}(X_{1}(n),...,X_{k}(n))$, $n\geq 1$, be a sequence of $k$-dimensional random vectors such that each $X(n)$ follows a multinomial law with parameters  $n\geq 1$\ and $p=(p_{1},p_{2},...p_{k})$ with 
\begin{equation*}
\forall (1\leq i\leq k),p_{i}>0\text{ \ and \ }\sum_{1\leq i\leq k}p_{i}=1.
\end{equation*}

\noindent Then, as $n\rightarrow +\infty$, 
\begin{equation*}
Z_{n}=\text{ }^{t}(\frac{X_{1}-np_{1}}{\sqrt{np_{1}}},...,\frac{X_{1}-np_{k}}{\sqrt{%
np_{k}}})
\end{equation*}

\noindent weakly converges to a $k$-dimensional Gaussian vector of variance-covariance matrix $\Sigma$ whose elements are 
\begin{equation}
\Sigma _{ii}=(1-p_{i}) \label{cv7002}
\end{equation}

\noindent and
\begin{equation}
\Sigma _{ij}=-\sqrt{p_{i}p_{j}}. \label{cv7003}
\end{equation}
\end{proposition}

\subsection{Finite dimensional weak limits of the uniform empirical process} \label{cv.review.subsec.ep} $ $\\

\noindent Let $U_{1},$\ $U_{2},...$ be a sequence of independent and standard uniformly distributed random variables on 
$(0,1)$ defined on the same probability space $(\Omega,\mathcal{A},\mathbb{P})$, with common distribution function $F(s)=s1_{(0\leq s\leq
1)}+1_{(s\geq 1)}$. For each $n\geq 1$, we may define the empirical distribution function associated with $U_{1}$, $U_{2}$,...,$U_{n}$ : 
\begin{equation*}
\begin{array}{ccc}
\mathbb{R}\ni x & \mapsto  & U_{n}(s)=\frac{1}{n}Card\{i,\text{ }1\leq i\leq
n\text{, \ }U_{i}\leq s\}%
\end{array}%
\end{equation*}

\noindent The empirical process associated with $U_{1}$, $U_{2}$,...,$U_{n}$ is defined as follows
\begin{equation*}
\alpha _{n}(s)=\sqrt{n}(U_{n}(s)-s),0\leq s\leq 1.
\end{equation*}

\bigskip \noindent Consider $0=t_{0}<t_{1}<...<t_{k}<t_{k+1}=1$ a partition of $(0,1)$ and set
\begin{equation*}
Y_{n}=\text{ }^{t}(\alpha _{n}(t_{1}),...,\alpha _{n}(t_{k+1})).
\end{equation*}

\noindent We have :

\begin{proposition} \label{cv.ep.df} 
Any finite distribution of the uniform empirical process of the form
$$
\text{ }^{t}(\alpha _{n}(t_{1}),...,\alpha _{n}(t_{k+1}))
$$ 

\noindent with
$$
0=t_{0}<t_{1}<...<t_{k}<t_{k+1}=1,
$$

\noindent weakly converges to a $k$-dimensional Gaussian random variable with variance-covariance matrix
$$
\left(\min (t_{i},t_{j})-t_{i}t_{j}\right) _{1\leq i,j\leq k},
$$

\noindent that is,
  
\begin{equation*}
\text{ }^{t}(\alpha _{n}(t_{1}),...,\alpha _{n}(t_{k+1}))\rightarrow \mathcal{N}%
_{k}(0,\left(\min (t_{i},t_{j})-t_{i}t_{j}\right) _{1\leq i,j\leq k})
\end{equation*}
\end{proposition}

\bigskip \noindent \textbf{Proof}. Set 
\begin{equation}
Z_{n}=\text{ }^{t}\left(\frac{\alpha _{n}(t_{1})}{\sqrt{t_{1}}},\frac{\alpha
_{n}(t_{2})-\alpha _{n}(t_{1})}{\sqrt{t_{2}-t_{1}}}...,\frac{\alpha
_{n}(t_{k+1})-\alpha _{n}(t_{k})}{\sqrt{t_{k+1}-t_{k}}}\right).  \label{cv730}
\end{equation}

\noindent Let us remark that
\begin{equation*}
N_{n}=(nF_{n}(t_{1}),nF_{n}(t_{2})-nF_{n}(t_{1}),...,nF_{n}(t_{k+1})-nF_{n}(t_{k}))
\end{equation*}

\noindent follows a multinomial law with outcomes probabilities $t_{1}$, $t_{2}-t_{1}$,...,$t_{k+1}-t_{k}$. Indeed, we have that
\begin{equation*}
nF_{n}(t_{j})-nF_{n}(t_{j-1})
\end{equation*}

\noindent is the number of observations falling in $]t_{j-1}, t_{j}]$ and for each $j$, the probability that one observation falls in 
$]t_{j-1}, t_{j}]$ is $p_{j}=t_{j}-t_{j-1}$.\\

\noindent We may apply the weak convergence of the multinomial law we established in Subsection \ref{cv.review.subsec.multinomial}.\\

\noindent Let us define $Z_{n}$ by centering each $j$th component of $N_{n}$ at $n(t_{j}-t_{j-1})$ and normalizing it by $\sqrt{n(t_{j}-t_{j-1})}$.\\

\noindent Remind that $Y_{n}=\text{ }^{t}(\alpha _{n}(t_{1}),...,\alpha _{n}(t_{k+1}))$. We have the matrix relation

$$
Z_n=AY_n \Leftrightarrow Y_n=BZ_n
$$

\noindent where the relation $y=Bz$ is the following correspondence  

$$
y_i=\sqrt{t_1}x_{1} + \sqrt{t_{2}-t_{1}}x_{2} + ... + \sqrt{t_{2}-t_{1}}x_{i}, \ \ 1\leq i \leq k+1.
$$

\noindent By the weak convergence of the  multinomial law, $Z_n$ weakly converges to a centered Gaussian vector $Z=(Z_{1},Z_{2},...,Z_{k+1})$ such that 

\begin{equation*}
\mathbb{E}(Z_{j}^{2})=1-(t_{j}-t_{j-1})
\end{equation*}

\noindent and 
\begin{equation*}
\mathbb{E}(Z_{i}Z_{j})=-\sqrt{(t_{i}-t_{i-1})(t_{j}-t_{j-1})}.
\end{equation*}

\noindent By the continuous mapping theorem (Point (e) of Theorem \ref{cv.review.rk}), $Y_n=BZ_n$ weakly converges to $Y=BZ$, where

$$
Y_i=\sqrt{t_1}Z_{1} + \sqrt{t_{2}-t_{1}}Z_{2} + ... + \sqrt{t_{2}-t_{1}}Z_{i}, \ \ 1\leq i \leq k+1.
$$

\bigskip \noindent Let $T=(T_{1},...,T_{k+1})$ be defined by $(t_{j}-t_{j-1})Z_{j}=T_{j}$, $1\leq i \leq k$, that is 

\begin{equation*}
Z=\text{ }^{t}\left(\frac{T_{1}}{\sqrt{t_{1}}},\frac{T_{2}}{\sqrt{(t_{2}-t_{1})}},...,\frac{T_{j}}{\sqrt{(t_{j}-t_{j-1})}},...,\frac{T_{k+1}}{\sqrt{(t_{k+1}-t_{k})}}\right)
\end{equation*}

\bigskip \noindent We have
\begin{equation*}
\mathbb{E}(T_{j}^{2})=\mathbb{E}\left(\left(Z_{j}\sqrt{(t_{j}-t_{j-1})}%
)^{2}=(t_{j}-t_{j-1})(1-(t_{j}-t_{j-1}\right)\right).
\end{equation*}

\bigskip \noindent and
\begin{equation*}
\mathbb{E}(T_{i}T_{j})=\sqrt{(t_{j}-t_{j-1})(t_{i}-t_{i-1})}\mathbb{E}%
(Z_{i}Z_{j})=-(t_{j}-t_{j-1})(t_{i}-t_{i-1}).
\end{equation*}

\noindent Before we compute the covariance of $Y_{i}$ and $Y_{j}$, we check that for $t_i \leq t_j$,%
\begin{eqnarray*}
t_i t_j&=&\left(\sum_{h=1}^{i} (t_h-t_{h-1})\right)\left(\sum_{r=1}^{j} (t_r-t_{r-1})\right)\\
&=& \left(\sum_{h=1}^{i} (t_h-t_{h-1})\right) \left(\sum_{r=1}^{i} (t_r-t_{r-1})+\sum_{r=i+1}^{j} (t_r-t_{r-1})\right)\\
&=&\left(\sum_{h=1}^{i} (t_h-t_{h-1})\right)^2+ \sum_{h=1}^{i} \sum_{r=i+1}^{j} (t_h-t_{h-1})(t_r-t_{r-1})\\
&=&\sum_{h=1}^{i} (t_h-t_{h-1})^2+\sum_{1\leq h\neq r \leq i} (t_h-t_{h-1})(t_r-t_{r-1})\\
&-& \sum_{h=1}^{h=i} \sum_{r=i+1}^{r} (t_h-t_{h-1})(t_r-t_{r-1})
\end{eqnarray*}

\noindent By putting together all these points, we are going to compute the variance-covariance matrix of $Y$. For $1\leq i\leq j\leq 1$, we have%
\begin{eqnarray*}
Y_{i}Y_{j}&=&\left(\sum_{h=1}^{i}T_{h}\right)^{2}+\sum_{h=1}^{h=i}\sum_{r=i+1}^{j}T_{h}T_{k}\\
&=& \sum_{h=1}^{i}T_{h}^2+\sum_{1\leq h\neq r \leq i}T_{h}T_{r} +\sum_{h=1}^{i}\sum_{r=i+1}^{j}T_{h}T_{r}.
\end{eqnarray*}

\noindent Finally, we get

\begin{eqnarray*}
\mathbb{E}(Y_{i}Y_{j})&=& \sum_{h=1}^{i}(1-(t_h-t_{h-1})) -\sum_{1\leq h\neq r \leq i} (t_h-t_{h-1})(t_r-t_{r-1})\\
&-&\sum_{h=1}^{h=i}\sum_{r=i+1}^{r=j} (t_h-t_{h-1})(t_r-t_{r-1})\\
&=& \sum_{h=1}^{i} (t_h-t_{h-1})-\sum_{h=1}^{i}(t_h-t_{h-1})^2\\
&-& \sum_{1\leq h\neq r \leq i} (t_h-t_{h-1})(t_r-t_{r-1}) -\sum_{h=1}^{i}\sum_{r=i+1}^{j} (t_h-t_{h-1})(t_r-t_{r-1})\\
&=& t_i -\sum_{h=1}^{i}(t_h-t_{h-1})^2 -\sum_{1\leq h\neq r \leq i} (t_h-t_{h-1})(t_r-t_{r-1})\\
&-& \sum_{h=1}^{i}\sum_{r=i+1}^{j} (t_h-t_{h-1})(t_r-t_{r-1})\\
&=& t_i - t_i t_j=\min(t_i,t_j)-t_i t_j.
\end{eqnarray*}

\noindent This completes the proof.\\

\section{Invariance principle} \label{cv.review.sec5}

Let $X_{1}$, $X$,.... be a sequence of \textit{iid} centered random variables with finite variances, that is $E\left\vert X_{i}\right\vert ^{2}<\infty$. For each $n\geq 1$, set

\begin{equation*}
S_{n}=X_{1}+\cdots+X_{n}
\end{equation*}

\noindent For $0\leq t\leq 1$ and $n\geq 1$, put  
\begin{equation*}
S_{n}(t)=\frac{S_{\left[ nt\right] }}{\sqrt{n}}
\end{equation*}

\noindent where, for any real $u$, $[u]$ stands for the integer part of $u$, which is the greatest integer less or equal to $u$.\\

\noindent we are going to explore the weak convergence of the finite distributions if the stochastic process \{$S_{n}(t),$\ $0\leq
t\leq 1\}$.\\

\noindent For this purpose, let $0=t_{0}<t_{2}<...<t_{k}=1,$\ $\ k\geq 1$. We have :

\begin{proposition} \label{cv.review.ip} The sequence of  finite distributions
\begin{equation*}
\left( \frac{S_{\left[ nt_{j}\right] }}{\sqrt{n}},1\leq j\leq k\right), \ \ n\geq 1,
\end{equation*}

\noindent weakly converges to $k$-dimensional centered Gaussian vector with variance-covariance matrix
\begin{equation*}
\left( \min (t_{i},t_{j})\right) _{1\leq i,j\leq k}.
\end{equation*}
\end{proposition}

\bigskip \noindent \textbf{Proof}.   We have
\begin{equation*}
\left\{ 
\begin{array}{c}
Y_{n}(t_{1})=X_{n}(t_{1})-X_{n}(t_{0})=\frac{1}{\sqrt{n}}\sum_{[nt_{0}]<j%
\leq \lbrack nt_{1}\}}X_{j} \\ 
\vdots \\ 
Y_{n}(t_{i})=X_{n}(t_{i})-X_{n}(t_{i-1})=\frac{1}{\sqrt{n}}%
\sum_{[nt_{i-1}]<j\leq \lbrack nt_{i}\}}X_{j} \\ 
\vdots\\ 
Y_{n}(t_{k})=X_{n}(t_{k})-X_{n}(t_{k-1})=\frac{1}{\sqrt{n}}%
\sum_{[nt_{k-1}]<j\leq \lbrack nt_{k}\}}X_{j}%
\end{array}%
\right. .
\end{equation*}

\bigskip \noindent We easily see that the random variables $Y_{n}(t_{i})$ are independent and that for each $1\leq i\leq k$, we apply the central limit theorem in $\mathbb{R}$ to get,  
\begin{equation*}
Y_{n}(t_{i})=\frac{1}{\sqrt{n}}\sum_{[nt_{i-1}]<j\leq \lbrack
nt_{i}\}}X_{j}\rightarrow \mathcal{N}(0,t_{i}-t_{i-1})
\end{equation*}

\noindent Hence, for any $u=(u_{1},...,u_{k})\in \mathbb{R}^{k}$, 
\begin{equation*}
\mathbb{E}\left(\exp \left(\sum_{1\leq i\leq 1}Y_{n}(t_{i})u_{i}\right) \right)=\prod_{1\leq i\leq 1}%
\mathbb{E}(\exp (Y_{n}(t_{i})u_{i})\rightarrow \prod_{1\leq i\leq 1}e^{\frac{1}{2}u_{i}^{2}/(t_{i}-t_{i-1})}.
\end{equation*}

\bigskip \noindent Thus, the vector $Y_{n}=text{ }^{t}(Y_{n}(t_{i}),  1\leq i\leq k)$ weakly converges to a Gaussian random vector $Z$, which has independent components and for each $1\leq i \leq k$, the $i$-th component has the variance $t_{i}-t_{i-1}$.\\

\noindent  The vector $X_{n}=text{ }^{t}(X_{n}(t_{i}),$\ $1\leq i\leq k)$ is the linear transform of $Y_{n}$ of the form  
\begin{equation*}
X_{n}=AY_{n}=\left( 
\begin{array}{cccc}
1 & 0 & ... & 0 \\ 
1 & 1 & ... & 0 \\ 
1 & ... & 1 & 0 \\ 
1 & ... & ... & 1%
\end{array}%
\right) Y_{n}
\end{equation*}

\noindent with

\begin{equation*}
A_{ij}=1_{(i\leq j)}.
\end{equation*}

\bigskip \noindent Then $X_{n}$ weakly converges to  $V=AZ$, whose components satisfy  
\begin{equation*}
V_{i}=Z_{1}+...+Z_{i}
\end{equation*}

\noindent and
\begin{equation*}
Z_{i}=V_{i}-V_{i-1}.
\end{equation*}

\noindent Then for any $1\leq i\leq k$,
\begin{equation*}
\mathbb{E}(V_{i}^{2})=\sum_{1\leq j\leq i}\mathbb{E}(Z_{j}^{2})=\sum_{1\leq
j\leq i}(t_{j}-t_{j-1})=t_{i}.
\end{equation*}

\noindent And for any 1$\leq i\leq j\leq k$,
\begin{eqnarray*}
\mathbb{E}(V_{i}V_{j})&=&\mathbb{E}(V_{i}(V_{i}+(V_{j}-V_{i})\\
&=&\mathbb{E}(V_{i}^{2})+\mathbb{E}(V_{i}(V_{j}-V_{i})).
\end{eqnarray*}

\noindent Since
\begin{equation*}
V_{i}=Z_{1}+...+Z_{i}
\end{equation*}

\noindent and since the random variables
  
\begin{equation*}
Z_{i}=V_{i}-V_{i-1}
\end{equation*}

\noindent are independent and centered, we get
\begin{equation*}
\mathbb{E}(V_{i}V_{j})=\mathbb{E}(V_{i}^2)=t_{i}=t_{i}\wedge t_{j}.
\end{equation*}

\noindent This suffices to conclude.\\

\bigskip \noindent \textit{Terminology}. The result presented in Proposition \ref{cv.review.ip} is the first step of what is called \textit{invariance principle} in Probability Theory.
 

%% file: asymptotics_cv_01_en.tex
\chapter{Weak Convergence Theory} \label{cv}

\section{Introduction}

In this chapter, we treat a unified theory of weak convergence by its functional characterization. We want to have complete theory of limits of sequences of probability measures on $(\mathbb{R}^{k},\mathcal{B}(\mathbb{R}^{k}))$, where $\mathcal{B}(\mathbb{R}^{k})$ is the Borel $\sigma$-algebra of $\mathbb{R}^{k}$.\\

\noindent However, the handling of the fundamental results only uses the metric structure of $\mathbb{R}^{k}$. This is why, whenever possible, we deal with sequences of probability measures on a metric spaces $(S,d)$, endowed with its Borel $\sigma $-algebra $\mathcal{B}(S)$.\\

\noindent But when dealing with limits of sub-sequences of sequences if random variables or probability measures, we essentially remain in $\mathbb{R}^{k}$ by making profit of the Helly-Bray theorem.\\

\noindent As in any theory on limits, we will have to deal with the uniqueness of limits, and convergence criteria, and relative compactness. Here, we will speak of weak compactness or simply tightness or uniform tightness.\\

\section{Definition, Uniqueness and Portmanteau Theorem}

\bigskip

\begin{definition} \label{cv.cvDEF} The sequence of measurable applications $X_{n}:(\Omega _{n},\mathcal{A}_{n},\mathbb{P}_{n})\mapsto (S,B(S))$ weakly converges to the measurable application  $X : (\Omega_{\infty},\mathcal{A}_{\infty},\mathbb{P}_{\infty})\mapsto (S,\mathcal{B}(S))$ if and only for any continuous and bounded function $f:S\mapsto \mathbb{R}$, (denoted  $f\in \mathcal{C}_{b}(S)$), we have 
\begin{equation}
\mathbb{E}f(X_{n})\rightarrow \mathbb{E}f(X) \text{ as } n \rightarrow +\infty. \label{cv21}
\end{equation}
\end{definition}

\bigskip \noindent We notice that the spaces on which the applications  $X_{n}$ are defined have no importance here. Only matter their probability laws on $(S,d)$. Indeed, denote  $L=\mathbb{P}_{X}=\mathbb{P}_{\infty}\circ X^{-1}$, the probability law of $X$ defined by  
\begin{equation*}
\forall \text{ }B\in \mathcal{B}(S),\text{ }L(B)=\mathbb{P}_{\infty}(X^{-1}(B))=\mathbb{P}_{\infty}(X\in B).
\end{equation*}

\noindent and for each $n\geq 1$, $\mathbb{P}^{(n)}$ the probability law of $X_n$ defined by  
\begin{equation*}
\forall \text{ }B\in \mathcal{B}(S),\text{ } \mathbb{P}^{(n)}(B)=\mathbb{P}_{n}(X_{n}^{-1}(B))=\mathbb{P}_{n}(X_n \in B).
\end{equation*}

\bigskip \noindent The definition says that $X_n$ weakly converges to $X$ if and only if for any $f\in \mathcal{C}_{b}(S)$,
\begin{equation*}
\int_{S}\text{ }f(x)\text{ }d\mathbb{P}^{(n)}(x) \rightarrow \int_{S}\text{ }f(x)\text{ }dL(x)  \text{ as } n \rightarrow +\infty.
\end{equation*}

\noindent We might also replace (\ref{cv21}) by 
\begin{equation}
\mathbb{E}f(X_{n })\rightarrow \int_{S}\text{ }f\text{ }dL \text{ as } n \rightarrow +\infty,  \label{cv22}
\end{equation}

\noindent and only say that $(X_{n})_{n\geq 1\ }$ weakly converges to the probability measure $L$. In the sequel, we will use both terminologies. 

\bigskip \textbf{Warning}. It is also important to see that the expectation symbols in (\ref{cv21}) depend of the probability measures that they use, and consequently, they should be labeled accordingly as  

$$
\mathbb{E}_{\infty}(f(X))=\int f(X) d\mathbb{P}_{\infty}, \text{  }\mathbb{E}_{n}(f(X_n))=\int f(X_n) d\mathbb{P}_{n}, \text{  }n\geq 1.
$$

\noindent But, for sake of simplicity, we choose not to put the subscripts $n$ and $\infty$ to keep the writing simple and to use them only when necessary.\\

\bigskip \noindent \textbf{Notation}. When $(X_{n})_{n\geq 1}$ weakly converges $X$ as $n\rightarrow +\infty$, we mainly use the notation%
\begin{equation*}
X_{n}\rightsquigarrow X \text{ as } \rightarrow +\infty,
\end{equation*}

\noindent but we may also use $X_{n}\rightarrow _{w}X$ ($w$ standing for \textit{weakly}) or $X_{n}\rightarrow _{d}X$ ($d$ standing for : \textit{in distribution}).\\

\noindent We are going to show that the limit we have defined is unique, \textbf{but in distribution}, in the following sense.\\

\bigskip 

\begin{proposition} \label{tm_cv_prop1} Let  $X_{n}:(\Omega _{n},\mathcal{A}_{n},\mathbb{P}_{n})\mapsto (S,B(S))$ be a sequence of measurable applications and, $\mathbb{P}_{1}$ and $\mathbb{P}_{2}$ two probability measures on $(S,\mathcal{B}(S))$. Suppose that  $X_{n}$
weakly converges to $\mathbb{Q}_{1}$ and to $\mathbb{Q}_{2}$. Then, we  necessarily have  
\begin{equation*}
\mathbb{Q}_{1}=\mathbb{Q}_{2}.
\end{equation*}
\end{proposition}

\noindent This means that if $X_{n}$ weakly converges to $X$ and to $Y$, then they have the same probability measure, meaning that they are equal in distribution.\\

\noindent \textbf{Proof}. Suppose that $X_{n}$ weakly converges to $\mathbb{P}_{1}$ and to $\mathbb{P}_{2}$. We want to show that 
$\mathbb{Q}_{1}=\mathbb{Q}_{2}$. But it suffices to show that the two probability measures coincide on the class $\Theta$ of open sets of $(S,d)$. Indeed, the class $\Theta $ is a $\pi$-system (that is : a class which is closed under finite intersection), which generates $\mathcal{B}(S)$. Then, by the $\lambda- \pi$ lemma, two probability measures on  $(S,\mathcal{B}(S))$ that coincide on $\Theta$ are equal on $\mathcal{B}(S).$\\

\noindent Now let $G$ be an open set of $S$. For any integer number $m\geq 1$, set the function $f_{m}(x)=\min
(m$\ $d(x,G^{c}),1)$, $x\in S$. We may see that for any $m\geq 1$, $f_{m}$ has values in $[0,1]$, and is bounded. Since $G^{c}$ is closed, we have
\begin{equation*}
d(x,G^{c})=\left\{ 
\begin{array}{c}
>0\text{ if }x\in G \\ 
0\text{ if }x\in G^{c}\text{ }%
\end{array}%
\right. .
\end{equation*}

\bigskip \noindent It is clear that $f_m=0$ on the border $\partial G$of $G$. We will not use this fact in what follows. But, we surely do use it later in this chpater.\\
 
\noindent Let us show that $f_{m}$ is a Lipschitz function. Let us handle $\left\vert f_{m}(x)-f_{m}(y)\right\vert $ through
three cases.\\

\noindent Case 1. $(x,y)\in (G^{c})^{2}$. Then 
\begin{equation*}
\left\vert f_{m}(x)-f_{m}(y)\right\vert =0\leq m\text{ }d(x,y).
\end{equation*}

\bigskip \noindent Case 2. $x\in G$\ and  $y\in G^{c}$\ (including also the case where the roles of $x$ and $y$ are switched). We have 
\begin{equation*}
\left\vert f_{m}(x)-f_{m}(y)\right\vert =\left\vert \min
(md(x,G^{c}),1)\right\vert \leq m\text{ }d(x,G^{c})\leq m\text{ }d(x,y),
\end{equation*}

\bigskip \noindent by the very definition of $d(x,G^{c})=\inf \{d(x,z),$\ z$\in G^{c}\}.$\newline

\noindent Case 3. $(x,y) \in G^{2}$. We use Property (\ref{annexe2}) in the Annexe Section (\ref{cv.annexe}) below and get, 
\begin{equation*}
\left\vert f_{m}(x)-f_{m}(y)\right\vert =\left\vert \min
(md(x,G^{c}),1)-\min (md(y,G^{c}),1)\right\vert \leq \left\vert
md(x,G^{c})-md(y,G^{c})\right\vert ,
\end{equation*}%
\begin{equation*}
\leq m\text{ }d(x,y)
\end{equation*}

\bigskip \noindent by the second triangle inequality. Then  $f_{m}$ is a Lipschitz function with coefficient $m$. Now, let us show that
\begin{equation*}
f_{m}\uparrow 1_{G}\text{ as m}\uparrow \infty .
\end{equation*}

\bigskip \noindent Indeed, if $x\in G^{c}$, we obviously have $f_{m}(x)=0\uparrow
0=1_{G}(x)$. If $x\in G$, that $d(x,G^{c})>0$ and $md(x,G^{c})\uparrow \infty$ as $m\uparrow \infty$. Then for $m$ large enough, 
\begin{equation}
f_{m}(x)=1\uparrow 1_{G}(x)=1 \text{ as} m\uparrow \infty.  \label{limfm}
\end{equation}

\bigskip \noindent In summary, each function $f_{m}$ is a non-negative and bounded Lipschitz function, that implies that 
$f_m \in \mathcal{C}_{b}(S)$, $m\geq 1$.\\

\noindent Now let us apply the definition of the weak convergence. The assumption implies that for any $f \in \mathcal{C}_{b}(S)$, we have as  $n\rightarrow +\infty$,
\begin{equation}
\mathbb{E}f(X_{n})\rightarrow \int f\text{ }d\mathbb{Q}_{1}\text{ \ and \ }%
\mathbb{E}f(X_{n})\rightarrow \int f\text{ }d\mathbb{Q}_{2}.
\end{equation}

\noindent By the uniqueness of real limits in $\mathbb{R}$, we get
\begin{equation*}
\forall (f\in C_{b}(S)),\int f\text{ }d\mathbb{Q}_{1}=\int f\text{ }d\mathbb{Q}_{2.}
\end{equation*}

\noindent Now, we apply this to the $f_m$, $m\geq 1$ to say 
\begin{equation*}
\forall (m\geq 1),\text{ }\int f_{m}\text{ }d\mathbb{Q}_{1}=\int f_{m}\text{ 
}d\mathbb{Q}_{2.}
\end{equation*}

\noindent Next, as $m$ increases to $+\infty$, we use (\ref{limfm}) and apply the Monotone Convergence Theorem to conclude that
\begin{equation*}
\int 1_{G}\text{ }d\mathbb{Q}_{1}=\int 1_{G}\text{ }d\mathbb{Q}_{2},
\end{equation*}

\noindent that is 
\begin{equation*}
\mathbb{Q}_{1}(G)=\mathbb{Q}_{2}(G).
\end{equation*}

\noindent Since $G$ is arbitrary fixed, this equality holds for all open sets of $S$. We conclude that $\mathbb{Q}_1=\mathbb{Q}_{2}$. $\square$\\

\noindent From that proof, we use specific functions $f_m$ to justify that ${Q}_{1}=\mathbb{Q}_2$. By using only the properties of the $f_m$'s, we have the following general laws.

\begin{proposition} \label{tm_cv_prop1b} Let us consider two probability measures $\mathbb{Q}_1$ and $\mathbb{Q}_{2}$ on $(S,\mathcal{B}(S))$. Let use define the assertions
$$
\mathbb{Q}_1(A)=\mathbb{Q}_{2}(B) \ \ (E1)
$$

\noindent and 

$$
\int f \ d\mathbb{Q}_1=\int f \ d\mathbb{Q}_{2},
$$

\noindent depending respectively on a measurable subset of $S$ and a real-valued measurable map defined on $S$. Then we have the equivalence between the following assertions :

\noindent (a) $\mathbb{Q}_1(A)=\mathbb{Q}_{2}$.\\

\noindent (b) Formula (E1) holds for any open set $A$.\\

\noindent (c) Formula (E1) holds for any closed set $A$.\\

\noindent (d) Formula (E1) holds for any continuous and bounded mapping $f$.\\

\noindent (e) Formula (E1) holds for any continuous and bounded mapping $f$ vanishing outside an open set.\\

\noindent (f) Formula (E1) holds for any Lipschitz and bounded mapping $f$.\\

\noindent (g) Formula (E1) holds for any Lipschitz and bounded mapping $f$ vanishing outside an open set.\\

\noindent If $E=\mathbb{R}^k$, $k\geq 1$, \\

\noindent (h) we may replace the phrase \textit{vanishing outside an open set} by \textit{vanishing outside an bounded open set} in (e) and (g).\\

\noindent Furthermore, the assertion :\\

\noindent (i) Formula (E1) holds for function $f$ of the form

$$
f(x)=\prod_{j=1}^{k} f_{j}(x_j), \ x=(x_1,\cdots,x_k) \in \mathbb{R}^k,
$$ 

\noindent, where each $f_j$ is a real-valued mapping defined on $\mathbb{R}$, which is Lipschitz, bounded and vanishing outside a bounded open set of $\mathbb{R}$,\\

\noindent is equivalent to any of the assertion (a) - (g). 
\end{proposition}

\noindent \textbf{Remark}. Only the last line of the proposition has to be justified, since the others are merely easy deduction from the proof of the preceding proposition. But on $S=\mathbb{R}^k$, the class $\mathcal{O}_b$ of bounded open intervals of the form

$$
]a,b[=\prod_{j=1}^{k} ]a_j,b_j[
$$

\noindent is a $\pi$-system and $\mathbb{R}^k$ is an increasing limit of elements of $\mathcal{O}_b$. By a slightly modified form of the $\pi-\lambda$ rule, $\mathcal{O}_b$ is a determining class of probability measures. For each $]a,b[$, for each $1\leq j \leq k$, we construct the sequence (in $m$) of functions $f_{j,m}$ based on  $]a_j,b_j[$, that is $f_{j,m}(x)=\min(md(x,]a_j,b_j[^c),1)$ to get that

$$
f_m(x)=\prod_{j=1}^{k} f_{j,m}(x_j) \uparrow 1_{]a,b[} \ as \ m\uparrow +\infty, \ x=(x_1,\cdots,x_k).
$$ 

\noindent Since each $f_m$ is Lipschitz, bounded (by one), vanishing outside a bounded open interval, we apply the same reasoning of the proof of the preceding proposition justify the two statements (h) and (i). $\blacksquare$\\

\bigskip \noindent \textbf{Notation}. When $(X_{n})_{n\geq 1}$ weakly converges to $X$, we use the following main notation
\begin{equation*}
X_{n}\rightsquigarrow X \text{ as } n\rightarrow +\infty.
\end{equation*}

\noindent But we will also use other notations like : $X_{n}\rightarrow _{\mathcal{L}}X$ (for convergence in law) or $X_{n}\rightarrow _{d}X$ (or convergence in distribution) or $X_{n}\rightarrow _{w}X$ (weak convergence).\\ 

\bigskip \noindent Next, we need to characterize the weak convergence using several criteria. This will furnish a rich set of tools for establishing weak convergence results.

\begin{theorem} \label{cv.theo.portmanteau} The sequence of measurable mappings $X_{n }:(\Omega _{n },\mathcal{A}_{n },P_{n })\mapsto (S,B(S))$ weakly converges to the probability measure $L$ if and only if one of these assertions holds.

\bigskip \noindent (ii) For any open set $G$ of  $S$ , 
\begin{equation*}
\liminf_{n\rightarrow +\infty} \mathbb{P}_{n}(X_{n }\in G)\geq L(G).
\end{equation*}

\bigskip \noindent (iii) For any closed set $G$ of  $S$, we have
\begin{equation*}
\limsup_{n\rightarrow +\infty} \mathbb{P}_{n}(X_{n }\in F)\leq L(F).
\end{equation*}

\bigskip \noindent (iv) For any lower semi-continuous and bounded below function $f$, we have
\begin{equation*}
\liminf_{n\rightarrow +\infty} \mathbb{E}f(X)\geq \int f \text{ }dL.
\end{equation*}

\bigskip \noindent (v) For any upper semi-continuous and bounded above function $f$, we have
\begin{equation*}
\limsup_{n\rightarrow +\infty} \mathbb{E}f(X_n)\leq \int f dL.
\end{equation*}

\bigskip \noindent (vi) For any Borel set $B$ of $S$ that is $L$-continuous, that is $L(\partial B)=0$, we have
\begin{equation*}
\lim_{n\rightarrow +\infty} \mathbb{P}(X_{n }\in B)=\lim_{n\rightarrow +\infty} \mathbb{P}_{n}(X_{n }\in B)=L(B).
\end{equation*}

\bigskip \noindent (vii) For any non-negative and bounded Lipschitz function $f$, we have : 
\begin{equation*}
\liminf_{n\rightarrow +\infty} \mathbb{E}f(X_n)\geq \int f \\
dL.
\end{equation*}
\end{theorem}

\bigskip

\bigskip \noindent Before we begin the proof, we recall that $\partial B$ is the boundary of the set $B$. If  $L(\partial B)=0$, it is said that
$B$ est $L$-continuous. As to the semi-continuous functions, we will give a reminder in the Annexe below.\\

\bigskip \noindent \textbf{Proof}. To unify the notation, we denote Formula (\ref{cv21}) as by Point (i) of the definition of weak convergence. From now, we break the proof into points.\\

\noindent (1) $(ii) \Leftrightarrow (iii)$. This is achieved by complementation.\\

\noindent (2) $(iv) \Leftrightarrow (v)$. This is achieved by moving from $f$ to  $-f$ and by remarking that opposite of upper semi-continuous functions are lower semi-continuous and vice-versa.\\
 
\noindent (3) $(i)\Rightarrow (vii)$. This is obvious since a Lipschitz function is continuous.\\

\noindent  $(vii)\Rightarrow (ii)$. Let $G$ be an open subset of $S$. For any $m\geq 1$, set $f_{m}(x)=\min (m$\ $d(x,G^{c}),1)$. We already knew from the proof of Proposition \ref{tm_cv_prop1} that for each $m\geq1$, $f_{m}$ is a non-negative and bounded Lipschitz function such that 
\begin{equation*}
f_{m}\uparrow 1_{G}\text{ as } m\uparrow \infty .
\end{equation*}

\bigskip \noindent We have for any $n\geq 1$ and for any  $m\geq 1,$ 
\begin{equation*}
\mathbb{E}(1_{G}(X_{n}))\geq \mathbb{E}f_{m}(X_{n}).
\end{equation*}

\bigskip \noindent Let us apply $(vii)$ to get  
\begin{equation}
\liminf_{n\rightarrow +\infty }\text{ }\mathbb{E}(1_{G}(X_{n}))\geq \lim
\inf_{n\rightarrow +\infty }\mathbb{E}f_{m}(X_{n})\geq \int f_{m}\text{ }dL.
\label{cv23}
\end{equation}

\noindent But for any measurable set $B$ and for any probability measure $\mathbb{Q}$,
\begin{equation*}
\mathbb{E}_{\mathbb{Q}}(1_{B})=\mathbb{Q}(B)
\end{equation*}

\bigskip \noindent For  $B=1_{X_{n}^{-1}(G)}=1_{(X_{n}\in G)}$, we let $m$ increase to $+\infty$ and use the Monotone Convergence Theorem to (\ref{cv23}), and get  
\begin{equation*}
\liminf_{n\rightarrow +\infty }\mathbb{P}(X_{n}\in G)\geq \int 1_{G}\text{ 
}dL=L(G).
\end{equation*}

\noindent Thus $(ii)$ holds true.\\

\bigskip \noindent (4) $(ii)\Rightarrow (iv)$. Assume $(ii)$ is true. Let $f$ be an lower semi-continuous function bounded below, say by $M$. In a first step, we are going to prove $(iv)$ for $f-M=g$, which is nonnative and lower semi-continuous. Then the sets  $(g\leq c)$ are closed by Proposition \ref{cv.annexe.SC} in the Annexe Section \ref{cv.annexe}. Set for $m\geq 1$ fixed,
 
\begin{equation*}
G_{i}=\{g>i/m\},\text{ i}\geq 1
\end{equation*}

\noindent and 
\begin{equation*}
g_{m}=\frac{1}{m}\sum_{i=1}^{m^{2}}1_{G_{i}}
\end{equation*}

\noindent The sets  $G_{i}$ are open since $g$ is lower semi-continuous. Let us remark that 
\begin{equation}
g_{m}(x)=\frac{i}{m}\text{ for  }\text{ }\frac{i}{m}<g(x)\leq \frac{i+1}{m},%
\text{ for } i=1,...,m^{2}-1  \label{cv24}
\end{equation}

\noindent  and 
\begin{equation*}
g_{m}(x)=m\text{ for } g(x)>m.
\end{equation*}

\noindent Then  
\begin{equation*}
g_{m}\leq g.
\end{equation*}

\noindent Further, by (\ref{cv24}) 
\begin{equation*}
\left| g_{m}(x)-g(m)\right| \leq 1/m\text{ }for\text{ }g(x)\leq m.
\end{equation*}

\noindent This implies  
\begin{equation*}
g(X_{n })\geq g_{m}(X_{n })=\frac{1}{m}\sum_{i=1}^{m^{2}}1_{G_{i}}(X_{n })=%
\frac{1}{m}\sum_{i=1}^{m^{2}}1_{(X_{n }\in G_{i})}
\end{equation*}

\noindent and next

\begin{equation}
\mathbb{E}g(X_{n })\geq \mathbb{E}g_{m}(X_{n })=\frac{1}{m}\mathbb{E}%
\sum_{i=1}^{m^{2}}1_{(X_{n }\in G_{i})}.  \label{cv25}
\end{equation}

\bigskip \noindent Then (\ref{cv25}) yields 
\begin{equation*}
\mathbb{E}g(X_{n })\geq \mathbb{E}g_{m}(X_{n })\geq \frac{1}{m}%
\sum_{i=1}^{m^{2}}\mathbb{E}1_{(X_{n }\in G_{i})} = \frac{1}{m}%
\sum_{i=1}^{m^{2}}\mathbb{P}(X_{n }\in G_{i}).
\end{equation*}

\noindent By letting $n$ go to $+\infty$ and by applying $(ii)$, we get  
\begin{equation*}
\liminf_{n \rightarrow +\infty} \mathbb{E}g(X_{n })\geq \liminf_{n \rightarrow +\infty} \text{ }\mathbb{E}g_{m}(X_{n
})\geq \frac{1}{m}\sum_{i=1}^{m^{2}}L(G_{i})=\int g_{m}\text{ }dL\geq
\int_{(g\leq m)}g_{m}\text{ }dL
\end{equation*}

\begin{equation*}
\geq \int_{(g\leq m)}g\text{ }dL+\int_{(g\leq m)}(g_{m}-g)\text{ }dL.
\end{equation*}

\noindent Now, as \noindent $m\rightarrow \infty$, we have 
\begin{equation*}
\int_{(g\leq m)}g\text{ }dL\rightarrow \int g\text{ }dL
\end{equation*}

\noindent and next,
\begin{equation*}
\left\vert \int_{(g\leq m)}(g_{m}-g)\text{ }dL \right\vert \leq L(S)/m\rightarrow 0.
\end{equation*}

\noindent Hence 
\begin{equation*}
\liminf_{n \rightarrow +\infty} \text{ }\mathbb{E}g(X_{n })\geq \int g\text{ }dL.
\end{equation*}

\bigskip \noindent Now, we come back to $f$ and see that by replacing $g$ by $f-M$ in $(vi)$, the formula remains true for $f$ by simplification of the finite number $M$. Hence $(iv)$ holds.

\bigskip \noindent (5) $(ii)\Rightarrow (vi)$. Recall that the boundary $\partial B$ of $B$ is the difference of interior $B$ from it 
adherence (closure), denoted as $\partial B=\overline{B}-int(B)$. Since 
$$
int(B)\subseteq B\subseteq \overline{B},
$$

\noindent we have

\begin{equation}
L(\partial B)=L(int(B))-L(\overline{B})=0 \Rightarrow L(int(B))=L(\overline{B%
})=L(B).  \label{cv26a}
\end{equation}

\bigskip \noindent Since $int(B)$ is open and $\overline{B}$ is closed, we may apply both (ii) and (iii) to get
\begin{equation}
L(int(B)) \leq \liminf_{n\rightarrow +\infty} \mathbb{P}_n(X_{n}\in int(B))\leq \limsup_{n\rightarrow +\infty} \mathbb{P}_n(X_{n }\in B),  \label{cv26b}
\end{equation}

\begin{equation}
\leq \mathbb{P}_n(X_{n}\in \overline{B}) \leq L(int(B)).  \label{cv26X}
\end{equation}

\bigskip \noindent Thus, by (\ref{cv26a}),

\begin{equation*}
L(B)=\liminf_{n\rightarrow +\infty} \mathbb{P}_n(X_{n }\in B)=\lim \mathbb{P}_n(X_{n}\in B),
\end{equation*}

\bigskip \noindent which was the target.\newline

\bigskip \noindent (6) $(vi)\Rightarrow (iii)$. Assume $(vi)$ holds and let $F$ be a closed subset of $S$. Set $F(\epsilon)=\{x,$\ $d(x,F)\leq \epsilon \}$\ for $\epsilon \geq 0$. We have  
\begin{equation*}
F\subseteq F(\epsilon)
\end{equation*}

\noindent and, since $F$ is closed, 
\begin{equation*}
F(\epsilon )\downarrow F \text{ as } \epsilon \downarrow 0
\end{equation*}

\bigskip \noindent Now, $\partial F(\epsilon )\subseteq \{x,\ d(x,F)=\epsilon \}$ and the sets $\{x,\ d(x,F)=\epsilon \}$ are disjoint. So the sets $\partial F(\epsilon)$ are disjoint. So they have null probabilities except eventually for a countable number of values of $\epsilon$, that is 
$$
L(\partial F(\epsilon))=0,
$$

\noindent except eventually for a countable number of values of $\epsilon$. (See Proposition \ref{cv.annexe.FDS} in the Annexe Section \ref{cv.annexe}). Then, we may easily find a sequence  $\epsilon _{n}\downarrow 0$ such that for any $p\geq 1$, 
\begin{equation*}
L(\partial F(\epsilon _{p}))=0.
\end{equation*}

\noindent For $p$ fixed, $F\subseteq F(\epsilon _{p})$ and this implies

\begin{equation*}
\limsup_{n\rightarrow +\infty} \mathbb{P}_n(X_{n }\in F)\leq \limsup_{n\rightarrow +\infty} \mathbb{P}_n(X_{n }\in
F(\epsilon _{p}))
\end{equation*}

\bigskip \noindent Next, by applying $(vi)$ 
\begin{equation*}
\limsup_{n\rightarrow +\infty} \mathbb{P}_n(X_{n}\in F)\leq \limsup_{n\rightarrow +\infty} \mathbb{P}_n(X_{n }\in
F(\epsilon _{p}))\leq L(F(\epsilon _{p})).
\end{equation*}

\bigskip \noindent Finally by letting $p$ go to infinity, we arrive at
\begin{equation*}
\lim \sup \mathbb{P}(X_{n }\in F)\leq L(F),
\end{equation*}

\bigskip \noindent and this is $(iii)$.\newline

\bigskip \noindent (7) $(iv)\Rightarrow (i)$. Assume $(iv)$ is true. Then $(v)$ is also true. Then for any bounded and continuous function $f$, it is lower semi-continuous and bounded below and upper semi-continuous and bounded above. We may apply both $(iv)$ and $(v)$ to have

\begin{equation*}
\int fdL\leq \liminf_{n\rightarrow +\infty} \mathbb{E}f(X_{n })\leq \limsup_{n\rightarrow +\infty} \mathbb{E}f(X_{n })\leq \int f \ dL.
\end{equation*}

\bigskip \noindent Thus 
\begin{equation*}
\int fdL=\liminf_{n\rightarrow +\infty} \mathbb{E}f(X_{n })=\limsup_{n\rightarrow +\infty} \mathbb{E}_{\ast }f(X_{n }).
\end{equation*}

\bigskip \noindent In summary, we have proved the Theorem through the following graph. We may check that each point implies all the others by using the right path in : 
\begin{equation*}
\begin{array}{ccccccc}
(i) & \Rightarrow & (vii) & \Rightarrow & (ii) & \Leftrightarrow & (iii) \\ 
\Uparrow &  &  &  & \Downarrow &  & \Uparrow \\ 
(v) & \Leftrightarrow & (iv) & = & 
\begin{array}{cc}
(iv) & (vi)%
\end{array}
& = & (vi)%
\end{array}%
\end{equation*}

\bigskip \noindent And this shows that the six assertions are equivalent.\\

\noindent  \textbf{Some extensions}. We may and do remark that the implications $(vii) \Rightarrow (ii)$ uses functions $f_m$ vanishing outside open sets $G$. So the weak convergence is also equivalent to the two other assertions, which are implied par Assertion (i).

\begin{corollary} \label{CvExtension} The sequence of measurable mappings $X_{n }:(\Omega _{n },\mathcal{A}_{n },P_{n })\mapsto (S,B(S))$ weakly converges to the probability measure $L$ if and only if \\

\bigskip \noindent (viia) For any bounded Lipschitz function $f$, we have : 
\begin{equation*}
\lim_{n\rightarrow +\infty} \mathbb{E}f(X_n) = \int f dL,
\end{equation*}

\noindent if and only if \\

\bigskip \noindent (viib) For bounded Lipschitz function $f$ vanishing outside an open set , we have : 
\begin{equation*}
\lim_{n\rightarrow +\infty} \mathbb{E}f(X_n) = \int f dL.
\end{equation*}

\noindent if and only if \\

\bigskip \noindent (viic) For bounded continuous function $f$ vanishing outside an open set , we have : 
\begin{equation*}
\lim_{n\rightarrow +\infty} \mathbb{E}f(X_n) = \int f \ dL.
\end{equation*}
\end{corollary}

\section{Continuous Mapping Theorem}

Let $(X_{n})$, $n\geq 1$, be a sequence of measurable applications with values in the metric space $(S, d)$ converging to the measurable application $X$ with values in $(S,d)$. Suppose we have a mapping of $(S,d)$ into another metric space $(E,r)$. The natural question we may ask ourselves is the following : Does the sequence $(g(X_{n}))$, $n\geq 1$, weakly converge to $g(X)$?\\

\noindent The answer is easy if $g$ is continuous. To make the ideas clear, denote $Y_n=g(X_{n})$, $n\geq 1$, and $Y=g(X)$. Then for any 
$f\in C_{b}(E)$, we have $h=(f\circ g) \in C_{b}(S)$ and for any $n\geq 1$, 
\begin{equation*}
\mathbb{E}f(Y_{n}) = \mathbb{E}(f\circ g)(X_n)=\mathbb{E}(h(X_n))
\end{equation*}

\noindent and

\begin{equation*}
\mathbb{E}f(Y) = \mathbb{E}(f\circ g(X))=\mathbb{E}(h(X)).
\end{equation*}

\bigskip\noindent Then $\mathbb{E}f(Y_{n})$ converges to $\mathbb{E}f(Y_n)$, whenever $(X_{n})$ weakly converges to $X$. Thus, we may conclude that  $g(X_{n})$, $n\geq 1$, weakly converges to $g(X)$.\\

\noindent This result is a particular case of a more general answer given below. Define by $D_g$ the set of all discontinuity points of $g$. The continuity of $g$ means that $D_g$ is empty. The generalization of the result given below requires that the function $g$ be $\mathbb{P}_X$-continuous, that is $\mathbb{P}_X(D_g)=0$. But we cannot write $\mathbb{P}_X(D_g)$ unless we are sure that $D_g$ is measurable. Fortunately, by Lemma \ref{cv.annexe.Discont} in the Appendix Section \ref{cv.annexe} below, it is a surprising fact that $D_g$ is measurable whatever be $g$. We have the following more general result.\\

\begin{proposition} \label{cv.mappingTh} Let $X_{n }:(\Omega _{n }, \mathcal{A}_{n },\mathbb{P}_{n})\mapsto (S,B(S))$ be a sequence of measurable applications weakly to converging to a measurable application $X :(\Omega_{\infty} ,\mathcal{A}_{\infty},\mathbb{P}_{\infty})\mapsto (S,B(S))$ (or to the probability measure $L$) and let $g$ be a mapping of $(S,d)$ into the metric space $(E,r)$ such that $g$ is $\mathbb{P}_X$-continuous (or $L$-continuous), then the sequences $g(X_{n})$ weakly converges to $g(X)$ (or to $L\circ g^{-1}$).
\end{proposition}

\bigskip \noindent \textbf{Proof}. Suppose that $X_{n }\rightarrow _{w}L$ with $L(discont(g))=0$. Let $F$ be a closed subset of $E$. Let us show that the Point (iii) of Portmanteau Theorem \ref{cv.theo.portmanteau} holds. Let us first show that,
  
\begin{equation}
\overline{g^{-1}(F))}\text{ }\subseteq g^{-1}(F)\cup discont(g).  \label{cv28}
\end{equation}

\noindent where $\overline{g^{-1}(F)}$ is the closure of $g^{-1}(F)$.  Indeed, let $x\in \overline{g^{-1}(F)}$. Then there exists a sequence $(y_{n})_{n\geq 1}\in $\ $g^{-1}(F)$ such that $y_{n}\rightarrow x$ and for each $n \geq 1,$\ $g(y_{n})\in F$. From here, we have two cases.\\

\noindent  Either $x\in discont(g)$ and then $x\in g^{-1}(F)\cup discont(g)$.\\

\noindent Or $x\notin  discont(g)$, that is $g$ is continuous at $x$. Then, since, $y_{n}\rightarrow x$, we have $g(y_{n}) \rightarrow g(x)$. Since the sequence $g(y_{n})$ is in $F$, which is closed, then $g(x)\in F$. This is equivalent to $x\in g^{-1}(F)$ and finally : $x\in g^{-1}(F)\cup discont(g)$.\\

\noindent We conclude that (\ref{cv28}) is true by combining both cases.\\

\noindent Now, let us use (\ref{cv28}) in the following way. We have 
\begin{eqnarray*}
\limsup_{n\rightarrow +\infty} \mathbb{P}_n(g(X_{n })\in F)&=&\limsup_{n\rightarrow +\infty} \mathbb{P}_{n}(X_{n }\in g^{-1}(F))\\
&\leq& \limsup_{n\rightarrow +\infty} \mathbb{P}_n(X_{n }\in \overline{g^{-1}(F))}),
\end{eqnarray*}

\noindent and subsequently,
 
\begin{equation*}
\limsup_{n\rightarrow +\infty} \mathbb{P}_n(X_{n }\in \overline{g^{-1}(F)}\text{ })\leq L(%
\overline{g^{-1}(F))}\text{ )}\leq L(g^{-1}(F))+L(discont(g))
\end{equation*}

\noindent This concludes the proof by
 
\begin{equation*}
\limsup_{n\rightarrow +\infty} \mathbb{P}_n(g(X_{n })\in F\text{ })\leq L\circ g^{-1}(F).
\end{equation*}

\bigskip \noindent In that proof, we used general properties of the metric. But when we have specific knowledge of the metric, we may go deeper and get particular criteria. Here, we are going to exploit the metrics of $\mathbb{R}^k$, $k\geq 1$. The combination of the Portmanteau theorem with the characterization results of probability measures in $\mathbb{R}^k$ leads to stunning and fine results.\\

\section{Space $\mathbb{R}^{k}$}

In this section we focus on the particular metric space $S=\mathbb{R}^{k}$. Before we begin, let us make some reminder on the characterization of the distributions in $\mathbb{R}$.\\

\bigskip \noindent Let $X=\left[ 
\begin{array}{c}
X_{1} \\ 
\cdot \cdot \cdot \\ 
X_{k}%
\end{array}
\right] ,$\ X$_{n}=\left[ 
\begin{array}{c}
X^{(n)}_{1} \\ 
\cdot \cdot \cdot \\ 
X^{(n)}_{k}%
\end{array}%
\right]$,\\

\noindent $n\geq 1$, be random vectors of dimension $k\geq 1$.\\

\noindent \textbf{Terminology}. By random vectors in $\mathbb{R}^{k}$, we mean measurable applications defined on some measurable space with values in $\mathbb{R}^{k}$.\\

\noindent Before we proceed further, we need some adaptations of the General Portmanteau Theorem \ref{cv.theo.portmanteau} to prepare more precise rules of weak convergence on $\mathbb{R}^d$. On $E=\mathbb{R}^{k}$, Point (ii) of that theorem may be restricted to bounded intervals of the form

$$
G=\prod_{1\leq j \leq k} ]a_j,b_j[, \ \forall j \in \{1,...,k\}, \ -\infty<a_j<b_j<+\infty.
$$

\noindent If so, the extension corollary \ref{CvExtension} (See page \pageref{CvExtension}) may also be adapted and based on functions vanishing outsides compacts sets. We we the following extension $E=\mathbb{R}^{k}$.

\begin{corollary} \label{CvExtensionRK} The sequence of measurable mappings $X_{n }:(\Omega _{n },\mathcal{A}_{n },P_{n })\mapsto (\mathbb{R}^k,\mathcal{B}(\mathbb{R}^k))$ weakly converges to a probability measure $L$ if and only if  \\

\bigskip \noindent (iiA) For any bounded open set $G$, we have  

\begin{equation*}
\lim_{n\rightarrow +\infty} \mathbb{E}f(X_n) = \int f \ dL,
\end{equation*}

\noindent if and only if\\

\bigskip \noindent (viiA) For ant Lipschitz and bounded function $f$, we have 
\begin{equation*}
\lim_{n\rightarrow +\infty} \mathbb{E}f(X_n) = \int f \ dL,
\end{equation*}

\noindent if and only if\\

\bigskip \noindent (viiB) For any $f$ Lipschitz and bounded function $f$, vanishing outside a compact set, we have 
\begin{equation*}
\lim_{n\rightarrow +\infty} \mathbb{E}f(X_n) = \int f \ dL,
\end{equation*}

\noindent if and only if\\

\bigskip \noindent (viiC) For any continuous and bounded function $f$, vanishing outside a compact set, we have 
\begin{equation*}
\lim_{n\rightarrow +\infty} \mathbb{E}f(X_n) = \int f \ dL
\end{equation*}
\end{corollary}

\noindent \textbf{Proof}. We only prove Point (iiA) by showing its equivalence with Assertion (ii) of the Portmanteau Theorem. Obviously (ii) implies (iiA). Now, suppose that (iiA) is true. Any open 
$G$ set is a countable union of bounded intervals $G_j=]a^{(j)}, b^{(j)}[$. So for each $p\geq 1$,

$$
\mathbb{P}_n(X_n \in G) \geq \mathbb{P}_n\left(X_n \in \bigcup_{1}^{p} G_j\right). 
$$ 

\noindent Since $\cup_{1\leq j \leq p} G_j$ is a bounded open set, we may apply (iiA) in the latter formula to have 

$$
\liminf_{n\rightarrow +\infty} \mathbb{P}_n(X_n \in G) \geq L \left(X_n \in \bigcup_{1}^{p} G_j\right),
$$ 

\noindent and by letting $p \uparrow +\infty$ and by using the continuity of $L$, we conclude by

$$
\liminf_{n\rightarrow +\infty} \mathbb{P}_n(X_n \in G) \geq L(G).
$$

\bigskip \noindent We recall that the probability law of a random vector $X$ of $\mathbb{R}^{k}$ is characterized by its \textit{distribution function}, defined by  
\begin{equation*}
(t_{1},t_{2},...,t_{k})^{T}\mapsto F_{X}(t_{1},t_{2},...,t_{k})=\mathbb{P}(X_{1}\leq t_{1},X_{2}\leq t_{2},...,X_{k}\leq t_{k})
\end{equation*}

\bigskip \noindent or by its characteristic function 
\begin{equation*}
(u_{1},u_{2},...,u_{k})^{T}\mapsto \Phi _{X}(u_{1},u_{2},...,u_{k})=\mathbb{E}\left( \exp \left( \sum_{j}^{k}i\text{ u}_{j}X_{j}\right) \right)
\end{equation*}

\noindent or by its moment generating function (whenever its exists) defined by  
\begin{equation*}
(u_{1},u_{2},...,u_{k})^{T}\mapsto \Psi _{X}(u_{1},u_{2},...,u_{k})=\mathbb{E}\left(\exp \left(\sum_{j}^{k}\text{ u}_{j}X_{j}\right) \right)
\end{equation*}

\noindent or by its probability density function whenever it exists. And it exists with respect to the Lebesgue measure for instance if and only if \begin{equation*}
(t_{1},t_{2},...,t_{k})^{T}\mapsto f_{X}(t_{1},t_{2},...,t_{k})=\frac{%
\partial ^{(k)}F_{X}(t_{1},t_{2},...,t_{k})}{\partial t_{1}\partial
t_{2}\cdot \cdot \cdot \partial t_{k}},
\end{equation*}

\noindent \textit{a.e.} with respect to the Lebesgue measures.\\

\noindent It remarkable that these characteristics also play important roles in the theory of weak convergence of random vectors.\\

\noindent We have the following characterizations and criteria.\\

\begin{proposition} \label{cv.FRDIR} Let $X_{n}:\ (\Omega _{n }, \mathcal{A}_{n },P_{n})\mapsto (\mathbb{R}^{k},\mathcal{B}(\mathbb{R}^{k}))$, $n\geq 1$ be a sequence of random vectors and $X : (\Omega_{\infty} ,\mathcal{A}_{\infty},\mathbb{P}_{\infty})\mapsto (\mathbb{R}^{k},\mathcal{B}(\mathbb{R}^{k}))$.\\

\noindent If $X_{n}$ weakly converges to $X$, then for any continuity point $t=(t_{1},t_{2},...,t_{k})$ of $F_{X}$, we have, as $n \rightarrow +\infty$,

\begin{equation}
\mathbb{P}_n\left(X_{n}\in \prod_{i=1}^{k}\left] -\infty ,t_{i}\right]
\right)\rightarrow F_{X}(t_{1},t_{2},...,t_{k}).  \label{cv32}
\end{equation}
 \end{proposition}

\bigskip \noindent \textbf{Proof}. Consider the distribution function of $X$ 
\begin{equation*}
F_{X}(t_{1},t_{2},...,t_{k})=\mathbb{P}(X_{1}\leq t_{1},X_{2}\leq
t_{2},...,X_{k}\leq t_{k})
\end{equation*}%
\begin{equation*}
=\mathbb{P}(X\in \prod_{i=1}^{k}\left] -\infty ,t_{i}\right] )
\end{equation*}

\noindent Denote $t=(t_{1},...,t_{k})$ and $t(n)=(t_{1}(n),t_{2}(n),...,t_{k}(n))$, $n\geq 1$. We have : $t(n)\uparrow t$\ (resp $t(n)\downarrow t)$ as $n \rightarrow +\infty$ if and only if
\begin{equation*}
\forall (1\leq i\leq k),\text{ }t_{i}(n)\uparrow t_{i}\text{ }(resp.\text{ }%
t_{i}(n)\downarrow t_{i}) \text{ as } n \rightarrow +\infty.
\end{equation*}

\noindent Set $A(t)= \prod_{i=1}^{k}\left] -\infty ,t_{i}\right] $. We have as $n\uparrow £\infty$,
\begin{equation*}
A(t(n))\downarrow A(t),
\end{equation*}

\noindent and by using the Monotone Convergence Theorem,
\begin{equation*}
F_{X}(t(n))=\mathbb{P}(X\in A(t(n))\downarrow \mathbb{P}(X\in A(t))=F_{X}(t),
\end{equation*}

\noindent as $n \rightarrow +\infty$. Then $F_{X}$ is right continuous at each point $t$. But 
\begin{equation*}
A(t(n))\uparrow A^{+}(t)=\prod_{i=1}^{k}\left] -\infty ,t_{i}\right[,
\end{equation*}

\noindent as $n \rightarrow +\infty$, and next, still by the Monotone Convergence Theorem,  
\begin{equation*}
F_{X}(t(n))=\mathbb{P}(X\in A(t(n))\uparrow \mathbb{P}(X\in A^{+}(t)),
\end{equation*}

\noindent as $n \rightarrow +\infty$. But we have
 
\begin{eqnarray}
D(t)&=&A(t)\text{ }\backslash \text{ }A^{+}(t)\\
&=&\{x=(x_{1},...,x_{k})\in A(t),\exists 1\leq i\leq k,\text{ }x_{i}=t_{i}\}.
\end{eqnarray}

\bigskip \noindent To better understand this formula, let us have a look at it for $k=1$ : 
\begin{equation*}
]-\infty ,\text{ }a]\text{ }\backslash \text{ }]-\infty ,\text{ }a[=\{a\}
\end{equation*}

\noindent and for $k=2$ (a diagram would help) : 
\begin{equation*}
]-\infty ,\text{ }a]\text{ }\times \text{ }]-\infty ,\text{ }b]\text{ }%
\backslash \text{ }]-\infty ,\text{ }a[\text{ }\times \text{ }]-\infty ,%
\text{ }b[
\end{equation*}%
\begin{equation*}
=\{(x,y)\in ]-\infty ,\text{ }a]\text{ }\times \text{ }]-\infty ,\text{ }b],%
\text{ }x=a\text{ or }y=b\}
\end{equation*}

\noindent Hence, if 
\begin{equation}
\mathbb{P}_{\infty}(X\in D(t))=L(D(t))=0,  \label{cv30}
\end{equation}

\noindent we get, as $n\rightarrow \infty$, 

\begin{eqnarray*}
F_{X}(t(n))=\mathbb{P}(X\in A(t(n)))\uparrow \mathbb{P}(X\in A^{+}(t))&=&\mathbb{P}(X\in A(t))-\mathbb{P}(X\in D(t))\\
&=&F_{X}(t).
\end{eqnarray*}

\bigskip \noindent We conclude that (\ref{cv30}) is the condition for $t$ to be a continuity point of $F_{X}$. Further, $D(t)$ is the boundary of 
$A(t)$, that is  
\begin{equation}
\partial A(t)=D(t)
\end{equation}

\bigskip \noindent To see this, just check that $A(t)$ is closed and that the interior of $A(t)$ is $A^{+}(t)$. By Point $(vi)$ of Portmanteau
Theorem \ref{cv.theo.portmanteau}, we get that for any continuity point $t$ of $F_X$,

\begin{equation*}
F_{X_n}(t)=\mathbb{P}(X_n \in A(t)) \rightarrow F_{X}(t)=\mathbb{P}(X \in A(t)) \text{ as } n\rightarrow +\infty.
\end{equation*}

\bigskip \noindent This ends the proof. Conversely, we will have :\\

\begin{proposition} \label{cv.FRINV}
Let $X_{n}:\ (\Omega _{n }, \mathcal{A}_{n },P_{n})\mapsto (\mathbb{R}^{k},\mathcal{B}(\mathbb{R}^{k}))$, $n\geq 1$ be a sequence of random vectors and $X : (\Omega_{\infty},\mathcal{A}_{\infty},P_{\infty})\mapsto (\mathbb{R}^{k},\mathcal{B}(\mathbb{R}^{k}))$. Suppose that for any continuity point $t$ of $F_X$, we have

\begin{equation}
F_{X_n}(t)=\mathbb{P}(X_n \in A(t)) \rightarrow F_{X}(t)=\mathbb{P}(X \in A(t)) \text{ as } n\rightarrow +\infty. \label{cv.conFr}
\end{equation}
  
\noindent Then $X_{n}$ weakly converges to $X$.
\end{proposition}

\noindent \textbf{Warning}. The proof of this proposition below is lengthy and very technical. It is stated only for people who are training to be a researcher in fundamental mathematics, probability or Statistics. If you are not among these people, you may skip it.\\

\noindent \textbf{Proof}. Suppose that for any $t=(t_{1},t_{2},...,t_{k})$ continuity point of $F_{X}$ and $F_{X_{n}}(t)\rightarrow F_{X}(t)$, as $n\rightarrow +\infty$.\\

\noindent To show that $X_{n}$ weakly converges to $X$, we are going to use Point $(ii)$ of Portmanteau Theorem \ref{cv.theo.portmanteau}, that is, for any open set $G$ on  $\mathbb{R}^{k}$, we have 
\begin{equation*}
\liminf_{n\rightarrow +\infty} \mathbb{P}_n(X_{n}\in G)\geq \mathbb{P}(X\in G).
\end{equation*}

\noindent Let $G$ be an arbitrary open set in $\mathbb{R}^{k}$. By using Proposition \label{cv.GFcontinuous}, $G$ is a countable union of $F_X$-continuous intervals in the form 
\begin{equation*}
G=\bigcup_{j\geq 1}]a^{j},b^{j}],
\end{equation*}

\noindent where for all $j\geq 1$, all the points $c$ defined by 
\begin{equation*}
c_{i}=a_{i}^{(j)}\text{ or }c_{i}=b_{i}^{(j)},
\end{equation*}

\noindent are continuity points of $F_{X}$. Following the notation in Formula (\ref{sec_append_not}) of the Appendix Section \ref{cv.annexe}, these points may be parametrized as

$$
c=b+\varepsilon* (a-b).
$$

\noindent In the sequel, $\mathcal{U}$ denotes the set of all bounded $F_X$-intervals.\\

\noindent Now, by the continuity of the probability measure $\mathbb{P}_{X}$, we can find for any 
$\eta >0$, an integer $m$ such that  
\begin{equation}
\mathbb{P}_{X}(G)-\eta \leq \mathbb{P}_{X}\left(\bigcup_{j=1}^{m}]a^{j},b^{j}]\right). \label{cv31b}
\end{equation}

\noindent We set $A_{j}=]a^{j},b^{j}]$ and use the Poincarr\'{e} formula, that is the inclusion-exclusion formula, that gives
\begin{eqnarray*}
\mathbb{P}_{X}\left(\bigcup_{j=1}^{m}A_{j}\right)&=&\sum \mathbb{P}_{X}(A_{j})-\sum \mathbb{P}_{X}(A_{i}A_{j}) \ \ \ (FP1)\\
&&+\sum \mathbb{P}_{X}(A_{i}A_{j}A_{k})+...+(-1)^{n+1}\mathbb{P}_{X}(A_{1}A_{2}...A_{n})
\label{FP1}
\end{eqnarray*}

\noindent and  
\begin{eqnarray*}
\mathbb{P}_{X_{n}}\left(\bigcup_{j=1}^{m}A_{j}\right)&=&\sum \mathbb{P}_{X_{n}}(A_{j})-\sum \mathbb{P}_{X_{n}}(A_{i}A_{j}) \ \ \ (FP2)\\
&&+\sum \mathbb{P}_{X_{n}}(A_{i}A_{j}A_{k})+...+(-1)^{n+1}\mathbb{P}_{X_{n}}(A_{1}A_{2}...A_{n}).  \label{FP2}
\end{eqnarray*}

\noindent We are going to handle each of these terms of the expressions above. Let us take one of th terms
\begin{equation*}
\mathbb{P}_{X}(A_{i_{1}}A_{i_{2}}...A_{i_{k}}).
\end{equation*}

\noindent As showed in Subsection \ref{cv.subsFcontinuous} in the Annexe Section \ref{cv.annexe} below, the class $\mathcal{U}$ of $F_{X}$-continuous intervals is stable under finite intersection. Thus, any set $A_{i_{1}}A_{i_{2}}...A_{i_{k}}$, which is of the type $]a,b]$, is in $\mathcal{U}$. It is a $F_{X}$-continuous interval. The Lebesgue-Stieljes Formula, gives

\begin{equation*}
\mathbb{P}_{X}(A_{i_{1}}A_{i_{2}}...A_{i_{k}})=\Delta_{a,b}F,
\end{equation*}

\noindent with

\begin{equation} \label{sec2.DeltaFormula}
\Delta_{a,b}F =\sum\limits_{\varepsilon \in \{0,1\}^{k}}(-1)^{(\sum_{1\leq i\leq k}\varepsilon _{i})}F_{X}(b+\varepsilon
\ast (a-b)).
\end{equation}

\noindent We similarly get that
\begin{equation*}
\mathbb{P}_{X_{n}}(A_{i_{1}}A_{i_{2}}...A_{i_{k}})=\sum\limits_{\varepsilon
\in \{0,1\}^{k}}(-1)^{(\sum_{1\leq i\leq k}\varepsilon
_{i})}F_{X_{n}}(b+\varepsilon \ast (a-b)).
\end{equation*}

\noindent And we are able to apply the assumption of the convergence of $F_{X_{n}}$ to $F_X$ for continuity points of $F_X$ to have, as $n \rightarrow +\infty$,
$$
\mathbb{P}_{X_{n}}(A_{i_{1}}A_{i_{2}}...A_{i_{k}})\rightarrow \mathbb{P}_{X}(A_{i_{1}}A_{i_{2}}...A_{i_{k}}).
$$

\noindent By operating term by term in Formulas (FP1) and in Formula (FP2), we conclude that, as $n \rightarrow +\infty$,
\begin{equation*}
\mathbb{P}_{X_{n}}\left(\bigcup_{j=1}^{m}A_{j}\right)\rightarrow \mathbb{P}
_{X}(\bigcup_{j=1}^{m}A_{j}).
\end{equation*}

\noindent Then 
\begin{eqnarray*}
\liminf_{n\rightarrow +\infty} \mathbb{P}_n(X_{n}\in G)&=&\liminf_{n\rightarrow +\infty} \mathbb{P}_n\left(X_{n}\in \bigcup_{j\geq 1}]a^{j},b^{j}]\right)\\
&&\geq \lim_{n\rightarrow +\infty} \mathbb{P}_n(X_{n}\in \bigcup_{j=1}^{m}]a^{j},b^{j}])\geq \mathbb{P}_{\infty}(X\in G)-\eta, 
\end{eqnarray*}

\noindent and this for an arbitrary $\eta >0$. Then, by letting $\eta \downarrow 0$, we arrive at  
\begin{equation*}
\lim \inf \mathbb{P}_n(X_{n}\in G)\geq \mathbb{P}_{\infty}(X\in G),
\end{equation*}

\noindent for any open set $G$ in $\mathbb{R}^k$. We finally conclude that  
\begin{equation*}
X_{n}\rightarrow _{w}X  \text{ as } n\rightarrow +\infty.
\end{equation*}

\bigskip \noindent We are moving to characteristic functions. We have the following characterization.\newline

\begin{proposition} \label{cv.FC} Let $X_{n}:\ (\Omega _{n},\mathcal{A}_{n},\mathbb{P}_{n})\mapsto (\mathbb{R}^{k},\mathcal{B}(\mathbb{R}^{k}))$, $n\geq 1$ be a sequence of random
vectors and $X:(\Omega _{\infty },\mathcal{A}_{\infty },\mathbb{P}_{\infty})\mapsto (\mathbb{R}^{k},\mathcal{B}(\mathbb{R}^{k}))$ another random vector. 
Then \noindent $X_{n}$ weakly converges to $X$ as $n\rightarrow +\infty ,$ if and only if for any point $(u_{1},u_{2},...,u_{k})^{T}\in 
\mathbb{R}^{k}$, 
\begin{equation*}
\Phi _{X_{n}}(u_{1},u_{2},...,u_{k})\mapsto \Phi _{X}(u_{1},u_{2},...,u_{k})%
\text{ as }n\rightarrow +\infty .
\end{equation*}
\end{proposition}

\bigskip \noindent \textbf{Remark} The proof we are proposing here is based Corollary \ref{CvExtensionSC} above and on a of a version of the  Stone-Weierstrass Theorem which is an important theorem in spaces of continuous functions defined on a compact set. Version of that theorem are recalled  in Section \ref{cv.annexe} (See page \pageref{cv.annexe}). Another proof, that is more beautiful to us, is provided in Theorem \ref{cv.tension.ConvFC} in Chapter \ref{cv.tensRk}. This latter is based on the concept of tightness and the Levy continuity theorem. But since the proof of Corollary \ref{CvExtensionSC} is based on the tightness, which by the way, is the key to both methods.

\bigskip \noindent \textbf{Proof}. Recall the definition of the characteristic
function : 
\begin{equation*}
(u_{1},u_{2},...,u_{k})^{T}\mapsto \Phi_{X}(u_{1},u_{2},...,u_{k})
=\mathbb{E}\left(\exp \left(i \sum_{j}^{k} u_{j} X_{j}\right)\right),
\end{equation*}

\bigskip \noindent which can written as follows. 
\begin{equation*}
(u_{1},u_{2},...,u_{k})^{T}\mapsto \exp \left(\sum_{j}^{k}i\text{ u}%
_{j}X_{j}\right)=\cos \left(\sum_{j}^{k} u_{j}X_{j}\right)+i \sin \left(\sum_{j}^{k} u_{j}X_{j}\right).
\end{equation*}

\bigskip \noindent This is a complex function whose components are bounded and
continuous functions of $X$ and by definition, we have 
\begin{equation*}
\mathbb{E}\exp \left(\sum_{j}^{k}i\text{ u}_{j}X_{j}\right)=\mathbb{E}\cos
\left(i\sum_{j}^{k} u_{j}X_{j}\right)+i \mathbb{E}\sin \left(\sum_{j}^{k}\text{u}_{j}X_{j}\right).
\end{equation*}

\noindent Hence, by the very definition of weak convergence, for any point $%
^{t}(u_{1},u_{2},...,u_{k})\in \mathbb{R}^{k},$%
\begin{equation}
\Phi _{X_{n }}(u_{1},u_{2},...,u_{k})\mapsto \Phi
_{X}(u_{1},u_{2},...,u_{k}).  \label{cv33}
\end{equation}

\bigskip \noindent This proves the direct implication of our proposition. To prove the indirect one, we appeal to an extension of the Stone-Weierstrass Theorem (See Corollary \ref{sec_EF_cor_05} in the appendix, page \pageref{sec_EF_cor_05}, in Subsection \ref{cv.subsec.annexe.SW1} in Section \ref{cv.annexe} below) and to Corollary \ref{CvExtensionSC}. According to that corollary, we need to prove that 
 $f \in \mathcal{L}_b$,  

\begin{equation}
\mathbb{E}(f(X_{n}))\rightarrow \mathbb{E}(f(X))\text{ as }n\rightarrow +\infty. \ \ \label{targ}
\end{equation}

\noindent * Let $f \in \mathcal{L}_b$, vanishing outside $[-r,r]^k$, $r>0$. Consider any $a>r$. Now let us consider the class $\mathcal{H}$ of finite linear combinations of functions of the form

\begin{equation}
\prod_{j=1}^{d} \exp\biggr(i n_j\pi x_j/r\biggr), \label{EX}
\end{equation}

\noindent where $n_j\in \mathbb{Z}$ is a constant and $i$ is the normed complex of angle $\pi/2$ and let $\mathcal{H}_{a}$ be the class of the restrictions $h_{a}$ of elements $h \in \mathcal{H}$ on 
$K_r=[-a,a]^k$.\\

\noindent It is clear that $\mathcal{H}_{a}$ is a sub-algebra of $C_b(K_a)$, $K_a=[-a,a]^k$, with the following properties.\\

\noindent (a) $f \equiv 0$ on $\partial K_a$.\\

\noindent (b) for each $h \in \mathcal{H}$, the uniform norm of $h$ on $\mathbb{R}^k$ is equal to the uniform norm of $h$ on $K_a$, that is

$$
\|f\|_{\infty}=\sup_{x \in \mathbb{R}^k} |h(x)|= \sup_{x \in K_a} |h(x)|=\|f\|_{K_a}. 
$$ 

\noindent This comes from that remark that $h$ is a finite linear combination of functions of the form in Formula \ref{EX} above and each factor $\exp\left(i n_j\pi x_j/r\right)$ is a 
$2r$-periodic function.\\

\noindent (c) $\mathcal{H}_{a}$ separates the points of $K_a\setminus \partial K_a$ and separates points of $K_r\setminus \partial K_a$ from points of $\partial K_a$. Indeed, if $x$ and $y$ are two points in $K_a$, at the exception where both of them are edge points of $K_a$ of the form

$$
(x,y) \in \{(s_1,...,s_d) \in K_a, \ \forall j\in\{1,...,d\}, \ s_j=a- \ or \ s_j=a\}^2,
$$

\noindent \noindent there exists $j_0 \in \{1,\cdots,d\}$ such that $0<|x_{j_0}-y_{j_0}|<2r$ that is $|(x_{j_0}-y_{j_0})/a|<2$ and the function 
$$
h_{r}(x)=\exp(i\pi x_{j_0}/a)
$$ 

\noindent separates $x$ and $y$ since $h_{a}(x)=h_{a}(y)$ would imply $\exp( i\pi (x_{j_0}-x_{j_0})/a)=1$, which  would imply $x_{j_0}-x_{j_0}=2\ell a$, $\ell \in \mathbb{Z}$. The only possible value of $\ell$ would be zero and this is impossible since $x_{j_0}-y_{j_0} \neq 0$.\\

\noindent (d) $\mathcal{H}_{r}$ contains all the constant functions.\\

\noindent We may then apply Corollary 2 in \cite{loSW2018} (Corollary \ref{sec_EF_cor_05} in the appendix, page \pageref{sec_EF_cor_05}) to get that : there exists $h_{a} \in \mathcal{H}_{a}$ such that

\begin{equation*}
\|f-h_{a}\|_{K_a} \leq \varepsilon/3. \label{approx19}
\end{equation*}

\noindent By the remark in , we have $\|h\|_{+\infty}=\|h_{a}\|_{K_a}\leq \|f\|_{K_a} +\|f-h_{a}\|_{K_a}$, so that

\begin{equation}
\|h\|_{+\infty} \leq \|f\|_{+\infty} +\varepsilon/3 \leq \|f\|_{+\infty}+1. \label{approx19}
\end{equation}

\noindent By (\ref{cv33}), we have 
\begin{equation*}
\mathbb{E}(h(X_{n}))\rightarrow \mathbb{E}(h(X))\text{ as }n\rightarrow +\infty .
\end{equation*}

\noindent Let $n_{0}\geq $ such that, for any $n\geq n_{0},$ 
\begin{equation}
\normalsize
\left\vert \mathbb{E}(h(X_{n}))\rightarrow \mathbb{E}(h(X))\right\vert
=\left\vert \int h\text{ }d\mathbb{P}_{n}\circ X_{n}^{-1}-\int h\text{ }d%
\mathbb{P}\circ X^{-1}\right\vert \leq \varepsilon /3.  \label{cv34g}
\end{equation}

\bigskip \noindent We have 
\begin{eqnarray*}
\mathbb{E}(f(X_{n}))-\mathbb{E}(f(X)) &=&\left(\int f\text{ }d\mathbb{P}_{n}\circ
X_{n}^{-1}-\int h\text{ }d\mathbb{P}_{n}\circ X_{n}^{-1}\right)  \notag \\
&&+\left(\int h\text{ }d\mathbb{P}_{n}\circ X_{n}^{-1}-\int h\text{ }d%
\mathbb{P}\circ X^{-1}\right) \notag \\
&&+\left(\int h\text{ }d\mathbb{P}\circ X^{-1}-\int f\text{ }d\mathbb{P}%
\circ X^{-1}\right). \notag
\end{eqnarray*}

\bigskip \noindent The first term satisfies 

\begin{eqnarray}
\mathbb{E}\left\vert \int f\text{ }d\mathbb{P}_{n}\circ X_{n}^{-1}-\int h\text{ }d\mathbb{P}_{n}\circ X_{n}^{-1}\right\vert  &\leq &\int_{K_a}\left\vert f-h\right\vert \text{ }d\mathbb{P}_{n}\circ X_{n}^{-1} \notag \\
&&+\int_{K_a^c}\left\vert f-h\right\vert \text{ }d\mathbb{P}_{n}\circ X_{n}^{-1} \notag \\
&\leq &\varepsilon /3+(\left\Vert f\right\Vert +\left\Vert h\right\Vert)
\mathbb{P}_{n}(X_{n}\in K_a^c). \label{cv34d}
\end{eqnarray}
\noindent 

\bigskip \noindent By treating the third term in the same manner, we also get \\ 
\begin{equation}
\normalsize
\mathbb{E}\left\vert \int f\text{ }d\mathbb{P}_{\infty }\circ X^{-1}-\int h\text{ }d\mathbb{P}_{\infty }\circ X^{-1}\right\vert \leq \varepsilon
/3+(\left\Vert f\right\Vert +\left\Vert h\right\Vert )\text{ }\mathbb{P}
_{\infty }(X\in K_a^c)  \label{cv34e}
\end{equation}

\bigskip \noindent By putting together Formulas (\ref{cv34g}), (\ref{cv34d}) and (\ref{cv34e}), and by using Formula \ref{approx19}, we get for each fixed $n\geq n_{0}$,
 
\begin{equation*}
\left\vert \mathbb{E}(f(X_{n}))-\mathbb{E}(f(X))\right\vert \leq \varepsilon +(2\left\Vert f\right\Vert +1)(\mathbb{P}_{n}(X_{n}\in K_a^c)+\mathbb{P}_{\infty }(X\in K_a^c)).
\end{equation*}

\bigskip \noindent For each fixed $n\geq n_{0},$ by letting $a\uparrow +\infty$,  $\mathbb{P}_{n}(X_{n}\in K_a^c)+\mathbb{P}_{\infty }(X\in K_a^c)\downarrow 0$. Then for each $n\geq
n_{0}$, we have

\begin{equation*}
\left\vert \mathbb{E}(f(X_{n}))-\mathbb{E}(f(X))\right\vert \leq \varepsilon. 
\end{equation*}

\noindent Taking the superior limit as $n \rightarrow +\infty$, and next letting $\varepsilon \downarrow 0$ make us reach the target. $\blacksquare$\\

\bigskip \noindent By putting together (\ref{cv.FRDIR}), (\ref{cv.FRINV}) and (\ref{cv.FC}), we have the full Portmanteau Theorem in $\mathbb{R}^k$.\\

\begin{theorem} \label{cv.theo.portmanteau.rk} Let $k$ be a positive integer. The sequence of random vectors $X_{n}:(\Omega _{n },\mathcal{A}_{n },
\mathbb{P}_{n })\mapsto (\mathbb{R}^k,\mathcal{B}(\mathbb{R}^k))$, $\geq 1$, weakly converges to the random vector $X:(\Omega _{\infty},\mathcal{A}_{\infty},\mathbb{P}_{\infty})\mapsto (\mathbb{R}^k,\mathcal{B}(\mathbb{R}^k))$ if and only if one of these assertions holds.\\

\bigskip \noindent (i) For any real-valued continuous and bounded function $f$ defined on $\mathbb{R}^k$,
\begin{equation*}
\lim_{n\rightarrow +\infty} \mathbb{E}f(X_n) =\mathbb{E}f(X).
\end{equation*}

\bigskip \noindent (ii) For any open set $G$ in $\mathbb{R}^k$, 
\begin{equation*}
\liminf_{n\rightarrow +\infty} \mathbb{P}_n(X_{n}\in G)\geq \mathbb{P}_{\infty}(X \in G).
\end{equation*}

\bigskip \noindent (iii) For any closed set $F$ of  $S$, we have
\begin{equation*}
\limsup_{n\rightarrow +\infty} \mathbb{P}_{n}(X_{n}\in F)\leq \mathbb{P_{\infty}}(X \in F).
\end{equation*}

\bigskip \noindent (iv) For any lower semi-continuous and bounded below function $f$, we have
\begin{equation*}
\liminf_{n\rightarrow +\infty} \mathbb{E}f(X_n)\geq \mathbb{E}f(X).
\end{equation*}

\bigskip \noindent (v) For any upper semi-continuous and bounded above function $f$, we have
\begin{equation*}
\limsup_{n\rightarrow +\infty} \mathbb{E}f(X_n)\leq \mathbb{E}f(X).
\end{equation*}

\bigskip \noindent (vi) For any Borel set $B$ of $S$ that is $\mathbb{P}_{X}$-continuous, that is $\mathbb{P}_{\infty}(X \in \partial B)=0$, we have
\begin{equation*}
\lim_{n\rightarrow +\infty} \mathbb{P}_n(X_{n }\in B)=\mathbb{P}_{X}(B)=\mathbb{P}_{\infty}(X \in B).
\end{equation*}

\bigskip \noindent (vii) For any non-negative and bounded Lipschitz function $f$, we have 
\begin{equation*}
\liminf_{n\rightarrow +\infty} \mathbb{E}f(X_n)\geq \mathbb{E}f(X).\\
\end{equation*}

\bigskip \noindent (viii) For any continuity point $t=(t_{1},t_{2},...,t_{k})$ of $F_{X}$, we have,

\begin{equation*}
F_{X_n}(t) \rightarrow F_{X}(t) \text{ as } n\rightarrow +\infty.
\end{equation*}

\noindent where for each $n\geq 1$, $F_{X_n}$ is the distribution function of $X_n$ and $F_{X}$ that of $X$.\\

\noindent (ix) For any point $u=(u_{1},u_{2},...,u_{k})\in \mathbb{R}^{k}$, 
\begin{equation*}
\Phi _{X_{n}}(u)\mapsto \Phi _{X}(u) \text{ as }n\rightarrow +\infty,
\end{equation*}

\noindent where for each $n\geq 1$, $\Phi_{X_n}$ is the characteristic function of $X_n$ and $\Phi_{X}$ is that of $X$
\end{theorem}

\bigskip \noindent The characteristic function as a tool of weak convergence is also used through the following criteria.\\

\noindent \textbf{Wold Criterion}. The sequence $\{X_n, \ \ n\geq 1\} \subset \mathbb{R}^k$ weakly converges to $X \in \mathbb{R}^k$, as $n \rightarrow +\infty$ if and only if for any $a \in \mathbb{R}^k$, the sequence $\{<a,X_n>, \ \ n\geq 1\} \subset \mathbb{R}$ weakly converges to $X \in \mathbb{R}$ as $n \rightarrow +\infty$.

\noindent \textbf{Proof}. The proof is quick and uses the notation above. Suppose that $X_n$ weakly converges to $X$ in $\mathbb{R}^k$ as $n \rightarrow +\infty$. By using the convergence of characteristic functions, we have for any $u\in \mathbb{R}^k$
$$
\mathbb{E}\left(\exp(i<X_n,u>)\right) \rightarrow \mathbb{E}\left(\exp(i<X,u>)\right) \ \ as \ \ n \rightarrow +\infty.
$$  

It follows for any $a\in \mathbb{R}^k$ and for any $t \in \mathbb{R}$, we have

\begin{equation}
\mathbb{E}\left(\exp(it<X_n,a>)\right) \rightarrow \mathbb{E}\left(\exp(it<X,a>)\right) \ \ as \ \ n \rightarrow +\infty. \label{cv.proj}
\end{equation}

\noindent that is, by taking $u=ta$ in the formula above, and by denoting $Z_n=<X_n,a>$ and $Z=<X,a>$

$$
\mathbb{E}\left(\exp(itZ_n)\right) \rightarrow \mathbb{E}\left(\exp(itZ)\right) \ \ as \ \ n \rightarrow +\infty.
$$

\noindent This means that $Z_n \rightsquigarrow Z$, that is $<a,X_n>$ weakly converges t0 $<a,X>$.\\

\noindent Conversely, suppose that for any $a \in \mathbb{R}^k$, the sequence $\{<a,X_n>, \ \ n\geq 1\} \subset \mathbb{R}$ weakly converges to $X \in \mathbb{R}$ as $n \rightarrow +\infty$. Then by taking $t=1$ in (\ref{cv.proj}) we get for any $a=u \in \mathbb{R}^k$,

$$
\mathbb{E}\left(exp(i<X,u>)\right) \rightarrow \mathbb{E}\left(exp(i<X,u>)\right) \ \ as \ \ n \rightarrow +\infty.
$$

\noindent which means that $X_n \rightsquigarrow +\infty$ as $n\rightarrow +\infty$.\\

\section{Theorem of Scheff\'{e}}

In the previous section, we linked the weak convergence to some characteristics of random vectors distributions, in particular the distribution functions and the characteristic functions. Now, what happens for the probability density functions? The theorem of Sheff\'e goes beyond the particular case of $\mathbb{R}^k$ and gives a very general answer as follows.\\

\begin{theorem} \label{cv.scheffe}. Let $\lambda$ be a measure on some measurable space $(E,B)$. Let $p$, $(p_{n})_{n\geq 1}$ be probability densities with respect to $\lambda $, that are real-valued, non-negative and measurable functions defined on $E$ such that
\begin{equation}
\forall n\geq 1,\int p_{n}\text{ }d\lambda =\int p\text{ }d\lambda =1. \label{scheffe1}
\end{equation}

\bigskip \noindent Suppose that
\begin{equation*}
p_{n}\rightarrow p,\text{ }\lambda -a.e.
\end{equation*}

\bigskip \noindent Then
\begin{equation}
\sup_{B\in \mathcal{B}}\left\vert \int_{B}p_{n}\text{ }d\lambda -\int_{B}p%
\text{ }d\lambda \right\vert =\frac{1}{2}\int \left\vert p_{n}-p\right\vert 
\text{ }d\lambda \rightarrow 0.  \label{scheffe2}
\end{equation}
\end{theorem}

\bigskip \noindent \textbf{Proof}. Suppose $p_{n}\rightarrow p,$\ $\lambda -a.e.$  Set $\Delta
_{n}=p-p_{n}.$ Then  (\ref{scheffe1}) implies 
\begin{equation*}
\int \Delta _{n}\text{ }d\lambda =0.
\end{equation*}

\noindent Then, for  $B\in \mathcal{B},$%
\begin{equation*}
\int_{B^{c}}\Delta _{n}\text{ }d\lambda =\int \Delta _{n}\text{ }d\lambda
-\int_{B}\Delta _{n}\text{ }d\lambda =-\int_{B}\Delta _{n}\text{ }d\lambda .
\end{equation*}

\noindent \noindent Thus, 
\begin{eqnarray}
2\left| \int_{B}\Delta _{n}\text{ }d\lambda \right| &=&\left| \int_{B}\Delta_{n}\text{ }d\lambda \right| +\left| \int_{B^{c}}\Delta _{n}\text{ }d\lambda \right| \notag\\
&\leq& \int_{B}\left| \Delta _{n}\right| \text{ }d\lambda
+\int_{B^{c}}\left| \Delta _{n}\right| \text{ }d\lambda \leq \int \left|
\Delta _{n}\right| \text{ }d\lambda ,  \label{scheffe3}
\end{eqnarray}

\noindent meaning that 
\begin{equation}
\left| \int_{B}\Delta _{n}\text{ }d\lambda \right| \leq \frac{1}{2}\int
\left| \Delta _{n}\right| \text{ }d\lambda .  \label{scheffe4}
\end{equation}

\bigskip \noindent By taking $B=(\Delta _{n}\geq 0)$ in (\ref{scheffe3}%
), we get
\begin{equation*}
2\left| \int_{B}\Delta _{n}\text{ }d\lambda \right| =\left| \int_{B}\Delta
_{n}^{+}\text{ }d\lambda \right| +\left| \int_{B^{c}}-\Delta _{n}^{-}\text{ }%
d\lambda \right| =\int \Delta _{n}^{+}d\lambda +\int \Delta _{n}^{-}d\lambda
=\int \left| \Delta _{n}\right| d\lambda .
\end{equation*}

\bigskip \noindent By putting together the two last formulas, we have
\begin{equation}
\sup_{B\in \mathcal{B}}\left| \int_{B}p_{n}\text{ }d\lambda -\int_{B}p\text{ 
}d\lambda \right| =\frac{1}{2}\int \left| p_{n}-p\right| \text{ }d\lambda .
\label{scheffe6}
\end{equation}

\bigskip \noindent Now we get, 
\begin{equation*}
0\leq \Delta _{n}^{+}=\max (0,p-p_{n})\leq p.
\end{equation*}

\bigskip \noindent Besides, we have 
\begin{eqnarray*}
\int \Delta _{n}^{+}\text{ }d\lambda &=&\int_{(\Delta _{n}\geq 0)}\Delta _{n}\text{ }d\lambda\\
&=&\int \Delta _{n}\text{ }d\lambda -\int_{(\Delta _{n}\leq 0)}\Delta _{n}\text{ }d\lambda\\
&=&\int_{(\Delta _{n}\leq 0)}-\Delta _{n} \ \ d\lambda =\int \Delta _{n}^{-}\text{ }d\lambda ,
\end{eqnarray*}

\bigskip \noindent so that
\begin{equation}
\int \left| \Delta _{n}\right| \text{ }d\lambda =2\int \Delta _{n}^{+}\text{ 
}d\lambda  \label{scheffe7}
\end{equation}

\bigskip \noindent Here, we apply the Fatou-Lebesgue Dominated Theorem to  
$$
0\leq \Delta _{n}^{+}\leq \left| \Delta _{n}\right|
\rightarrow 0 \text{ }\lambda -a.e. \text{ as } n\rightarrow +\infty, \text{ and } 0\leq \Delta _{n}^{+}\leq p
$$. 

\noindent We get  
\begin{equation*}
\int \Delta _{n}^{+}\text{ }d\lambda \rightarrow 0,
\end{equation*}

\bigskip \noindent in virtue of (\ref{scheffe6}), 
\begin{equation*}
\sup_{B\in \mathcal{B}}\left| \int_{B}p_{n}\text{ }d\lambda -\int_{B}p\text{ 
}d\lambda \right| =\frac{1}{2}\int \left| p_{n}-p\right| \text{ }d\lambda
=\int \Delta _{n}^{+}\text{ }d\lambda \rightarrow 0.
\end{equation*}

\bigskip \noindent which puts and end to the proof.\\

\bigskip \noindent The Theorem of Scheff\'{e} may be applied to probability densities in $\mathbb{R}^{k}$ with respect to the Lebesgue measure or to a counting measure.\\

\begin{proposition} \label{cv.scheffeExt} These two assertions hold.\\

\noindent (A) Let $X_{n}: (\Omega _{n },\mathcal{A}_{n},\mathbb{P}_{n}) \mapsto (\mathbb{R}^{k},\mathcal{B}(\mathbb{R}^{k}))$ be random vectors
and $X : (\Omega_{\infty},\mathcal{A}_{\infty},\mathbb{P}_{\infty})\mapsto (\mathbb{R}^{k},\mathcal{B}(\mathbb{R}^{k}))$ another random vector, all of them absolutely continuous with respect to the Lebesgues measure denoted as $\lambda_k$. Denote $f_{X_{n}}$ the probability density function of $X_n$,  
$n\geq 1$ and by $f_{X}$ the probability density function of $X$. Suppose that we have 

\begin{equation*}
f_{X_{n}}  \rightarrow f_{X}, \text{ }\lambda_k-a.e., \text{ as } n\rightarrow +\infty.
\end{equation*}

\noindent Then $X_{n}$ weakly converges to $X$ as $n\rightarrow +\infty$.\\

\noindent (B) Let $X_{n }: (\Omega _{n },\mathcal{A}_{n},\mathbb{P}_{n })\mapsto (\mathbb{R}^{k},\mathcal{B}(\mathbb{R}^{k}))$ be discrete random vectors and $X$ : $(\Omega_{\infty},\mathcal{A}_{\infty},\mathbb{P}_{\infty})\mapsto (\mathbb{R}^{k},\mathcal{B}(\mathbb{R}^{k}))$ another discrete random vector. For each $n$, define $D_n$ the countable support of $X_n$, that 
$$
\mathbb{P}_n(X_n \in D_n)=1 \text{ and for each } x \in D_n, \text{ } \mathbb{P}(X_n=x)\neq 0,
$$

\noindent and $D_{\infty}$ the countable support of $X$. Set  $D=D_{\infty} \cup (\cup_{n \geq 1} D_n)$ and denote by $\nu$ as the counting measure on $D$. Then the probability densities of the $X_n$ and of $X$ with respect to $\nu$ are defined on $D$ by
$$
f_{X_n}(x)=\mathbb{P}_{n}(X_n=x), \text{ } n\geq 1, \text{ } f_{X}(x)=\mathbb{P}_{\infty}(X=x), \text{ }x\in D.
$$

\noindent If 
\begin{equation*}
(\forall x\in D), f_{X_{n }}(x)  \rightarrow f_{X}(x),
\end{equation*}

\noindent then $X_{n}$ weakly converges to $X$.
\end{proposition}

\bigskip \noindent We will finish by giving a very refined complement the Portmanteau Theorem on $\mathbb{R}^k$. Denote by $\mathcal{L}_b$ the class of functions of the form

$$
f(x)=\prod_{j=1}^{k} f_j(x_j), \ x=(x_1,...,x_k)^t, \ (PF)
$$

\noindent where each $f_j$ is continuous bounded and vanishing outside a compact set. denote also by $\mathcal{L}{bl}$ the class of functions of the form (PF) where each $f_j$ is bounded, Lipschitz and 
vanishing outside a compact set.\\

\noindent We have the following extension.

\begin{corollary} \label{CvExtensionSC} The sequence of measurable mappings $X_{n }:(\Omega _{n },\mathcal{A}_{n },P_{n })\mapsto (\mathbb{R}^k,\mathcal{B}(\mathbb{R}^k))$ weakly converges to a probability measure $L$ if and only if for any $f \mathcal{L}_b$,

$$
\mathbb{E}f(X_n) \rightarrow \int f \ dL,
$$

\noindent if and only if for any $f \mathcal{L}_{bl}$,

$$
\mathbb{E}f(X_n) \rightarrow \int f \ dL,
$$
\end{corollary}

\noindent \textbf{Proof}. The proof of that important extension uses the tightness which will be addressed in Chapter \ref{cv.tensRk}, Subsection \ref{proofCvExtensionSC}, page \pageref{proofCvExtensionSC}. $\square$

\section[Weak convergence of sequences on the same space]{Weak Convergence and Convergence in Probability on one Probability Space} \label{cv.CvCp}

\noindent In this section, we place the weak convergence limit in the
general frame of the convergence of random variables defined on the same
probability space ($\Omega ,\mathcal{A},\mathbb{P})$ with values in a metric
space $(S,d).$ We already saw that the weak convergence of random variables
does not require from them and from the weak limit random variable that they
are defined on a common probability space. In the particular case where this
happens, and only in this case, we are able to have interesting relations
with other types of convergences.\\

\noindent Conversely, the powerful theorem of theorem of Skorohod-Wichura-Dudley
allows to transform any weak convergence, under specific conditions on the
space $S$, to an almost-sure convergence of versions of the sequences and on
the limit. In this text, this theorem is only proved when $S$ is the
real line $\mathbb{R}$ in Chapter \ref{cv.R}. The proof is expected
in a more general book on weak convergence.\\

\noindent Let us begin with the definitions.\\

\subsection{Definitions} $ $\\

\noindent In all this section, except in the Subsection \ref{cv.skorohod}, the random variables  $(X_{n})_{n\geq 0},$ $(Y_{n})_{n\geq 0}$, etc., and
the random variables $X,$ $Y,$ $etc.$ are defined on the same probability
space ($\Omega ,\mathcal{A},\mathbb{P})$ \ and have their values in the
metric space  $(S,d)$. We will also have to use constants $c$ in  $S$. So you will not find
probability measures $\mathbb{P}_{\infty}$, and $\mathbb{P}_n$n, $n\geq 1$, here.\\

\bigskip \noindent \textbf{(a) Almost-sure convergence}.\\

\noindent The sequence $(X_{n})_{n\geq 0}$ converges almost-surely to $X$ as 
$n\rightarrow \infty $, denoted as $X_{n}\longrightarrow X,$ $a.s.$ as $%
n\rightarrow \infty ,$ if and only if the subspace of $\Omega $ on which $%
(X_{n})_{n\geq 0}$\ fails to converge to $X$\ is a $\mathbb{P}$-null set,
that is%
\begin{equation*}
\mathbb{P}(\{\omega \in \Omega ,X_{n}\nrightarrow X\})=\mathbb{P}(\{\omega
\in \Omega ,d(X_{n},X)\nrightarrow 0\})=0.
\end{equation*}

\noindent This may be expressed as 
\begin{equation*}
(X_{n}\nrightarrow X)=\bigcup\limits_{k\geq 1}\bigcap\limits_{n\geq
0}\bigcup\limits_{p\geq n}(d(X_{p},X)>k^{-1}).
\end{equation*}

\noindent and this is surely measurable because of the continuity of the
metric $d$. This leads to the new definition :  $(X_{n})_{n\geq 0}$
almost-surely converges to $X$ as $n\rightarrow \infty ,$ if and only if : 
\begin{equation}
\forall k\geq 1,\text{ }\mathbb{P}\left( \bigcap\limits_{n\geq
0}\bigcup\limits_{p\geq n}(d(X_{p},X)>k^{-1})\right) =0.  \label{cps}
\end{equation}

\bigskip \noindent \textbf{(b) Convergence in probability}.\\

\noindent The sequence $(X_{n})_{n\geq 0}$ converges in probability to $X$, as $n\rightarrow +\infty$, denoted as  $X_{n}\longrightarrow _{\mathbb{%
\mathbb{P}}}X$, if and only if  
\begin{equation*}
\forall \varepsilon >0,\text{ }\lim_{n\longrightarrow +\infty }\mathbb{P}%
\left( d(X_{n},Y)>\varepsilon \right) =0.
\end{equation*}

\noindent Now, we are going to make a brief comparison between these two types of convergence. The following proposition is already known to
the reader in the case where $S$ is $\mathbb{R}$.

\bigskip

\begin{proposition} \label{cv.CvCp.01}
If $X_{n}\longrightarrow X$ $a.s.$ as $n\rightarrow +\infty ,$ then $%
X_{n}\longrightarrow _{\mathbb{P}}X$ as $n\rightarrow +\infty$.
\end{proposition}

\noindent \textbf{Proof}. The proof is the same as in $\mathbb{R}$. Suppose
that $X_{n}\longrightarrow X$ $a.s.$ as $n\rightarrow +\infty .$\  We
have to prove (\ref{cps}). We have for $k\geq 1,$%
\begin{equation*}
(d(X_{n},X)>k^{-1})\subset \bigcup\limits_{p\geq
n}(d(X_{p},X)>k^{-1})=:B_{n,k}.
\end{equation*}

\noindent But the sequence $B_{n,k}$ is non-decreasing in $n$ to 
\begin{equation*}
\bigcap\limits_{n\geq 0}\bigcup\limits_{p\geq n}(d(X_{p},X)>k^{-1})=:B_{k}
\end{equation*}

\noindent and for any $n\geq 0$ and  $k\geq 1$, 
\begin{equation}
\mathbb{P}(d(X_{n},X)>k^{-1})\leq \mathbb{P}\left( B_{n,k}\right). \label{cv.AS1}
\end{equation}

\noindent By the continuity of the probability, 

\begin{equation*}
\lim \sup_{n\rightarrow \infty }\mathbb{P}(d(X_{n},X)>k^{-1})\leq
\lim_{n\rightarrow \infty }\mathbb{P}\left( B_{n,k}\right) =\mathbb{P}\left(
B_{k}\right) =0,
\end{equation*}

\noindent where we applied (\ref{cps}) to the left member of (\ref{cv.AS1}).\\

\bigskip \noindent Now we are going to give a number of relations between the convergence in
probability and the weak convergence

\subsection{Weak Convergence and Convergence in Probability} $ $\\

\noindent Before we step in the comparison results, we have to enrich the  Portmanteau
Theorem \ref{cv.theo.portmanteau} by this supplementary point.\\

\begin{lemma}
\label{cvcp.lem01} The sequence $(X_{n})_{n\geq 0}$ weakly converges to $X$ as $n \rightarrow +\infty$ if and only if \\

\noindent (viia) For any bounded Lipschitz function  $f:S\longrightarrow 
\mathbb{R}$,
\begin{equation*}
\mathbb{E}f(X_{n})\rightarrow f(X),
\end{equation*}

\noindent as $n \rightarrow +\infty$.
\end{lemma}

\noindent \textbf{Proof}. Let us place ourselves in the proof of Portmanteau
Theorem \ref{cv.theo.portmanteau}.  Now $(vii)$ is a sub-case of $(viia),$
and then $(viia)\Longrightarrow (vii).$ Now if  $(vii)$ holds, we may take
the infimum $A$ and the supremum $B$ of a bounded Lipschitz function $f.$ By
applying Point $(vii)$ to  $f-A$ and to  $-f+B,$ we get $(viia).$ Then we
have $(vii)\Longleftrightarrow (viia)$.\\ 

\noindent We are going to state a number of properties.\\

\noindent In the sequence, all limits in presence of subscripts $n$ are meant as $%
n\rightarrow +\infty $ unless the contrary is specified.

\bigskip \noindent \textbf{(a) The convergence in probability implies the weak convergence}

\begin{proposition} \label{cvcp.prop2} If\ $X_{n}\longrightarrow _{\mathbb{P}}X$ as $n \rightarrow +\infty$, then 
$X_{n}\rightsquigarrow X$ as $n \rightarrow +\infty$.
\end{proposition}

\noindent \textbf{Proof}. Suppose that  $X_{n}\longrightarrow _{\mathbb{P}}X.
$ Let us show that $X_{n}\rightsquigarrow X$ by using Point $(viia)$ of
Lemma \ref{cvcp.lem01} above. Let $f$ be a Lipschitz bounded function of
coefficient $\ell >0$ and of bound $M.$ We have for any  $n\geq 0,$%
\begin{equation*}
\left\vert f(X_{n})-f(X)\right\vert \leq \ell d(X_{n},X).
\end{equation*}

\noindent We have for any $n\geq 0$ and for any $\varepsilon>0$,

\begin{eqnarray*}
\left\vert \mathbb{E}f(X_{n})-\mathbb{E}f(X)\right\vert  &\leq &\mathbb{E}%
\left\vert f(X_{n})-f(X)\right\vert  \\
&\leq &\int_{(d(X_{n},X)\leq \varepsilon )}\left\vert f(X_{n})-f(X)\right\vert d\mathbb{P}\\
&+&\int_{(d(X_{n},X)>\varepsilon )}\left\vert f(X_{n})-f(X)\right\vert d\mathbb{P}.
\end{eqnarray*}

\noindent But for any $n\geq 0$ and for any  $\varepsilon >0,$%
\begin{equation*}
\int_{(d(X_{n},X)\leq \varepsilon )}\ell d(X_{n},X)d\mathbb{P}\leq \ell
\varepsilon .
\end{equation*}

\noindent Furthermore, for any $n\geq 0$ and for $\varepsilon >0,$%
\begin{equation*}
\int_{(d(X_{n},X)>\varepsilon )}\left\vert f(X_{n})-f(X)\right\vert d\mathbb{%
P}\leq \int_{(d(X_{n},X)>\varepsilon )}2Md\mathbb{P}\leq 2M\text{ }\mathbb{P}%
(d(X_{n},X)>\varepsilon ).
\end{equation*}

\noindent Then for any $n\geq 0$ and for any  $\varepsilon >0,$%
\begin{equation*}
\left\vert \mathbb{E}f(X_{n})-\mathbb{E}f(X)\right\vert \leq \ell
\varepsilon +2M\text{ }\mathbb{P}(d(X_{n},X)>\varepsilon ).
\end{equation*}

\noindent Then for any $\varepsilon >0,$%
\begin{equation*}
\lim \sup_{n\rightarrow \infty }\left\vert \mathbb{E}f(X_{n})-\mathbb{E}%
f(X)\right\vert \leq \ell \varepsilon .
\end{equation*}

\noindent By letting $\varepsilon \downarrow 0$, we get 
\begin{equation*}
\mathbb{E}f(X_{n})\longrightarrow \mathbb{E}f(X),
\end{equation*}

\noindent which finishes the proof.\newline

\bigskip \noindent \textbf{(b) Weak convergence and convergence in probability to a constant are equivalent}.\\

\begin{proposition} \label{cv.CvECp} We have the following equivalence : $X_{n}\longrightarrow _{\mathbb{P}}c$ as $n \rightarrow +\infty$
if and only if $X_{n}\rightsquigarrow c$ as $n \rightarrow +\infty$.
\end{proposition}

\noindent \textbf{Proof}. The implication ($X_{n}\rightarrow _{\mathbb{P}%
}c)\Rightarrow \left( X_{n}\rightsquigarrow c\right) $ comes from
Proposition \ref{cvcp.prop2}. Let us prove that $\left(
X_{n}\rightsquigarrow c\newline
\right) \Rightarrow (X_{n}\rightarrow _{\mathbb{P}}c).$ Suppose that $\left(
X_{n}\rightsquigarrow c\right) .$ Let $\varepsilon >0$. Point ($ii)$ of
Portmanteau Theorem \ref{cv.theo.portmanteau} gives%
\begin{eqnarray*}
\lim \inf_{n\rightarrow +\infty }\mathbb{P}(d(X_{n},c) &<&\varepsilon )=\lim
\inf_{n\rightarrow +\infty }\mathbb{P}(X_{n}\in B(c,\varepsilon ))\leq 
\mathbb{P}(c\in B(c,\varepsilon )) \\
&=&\mathbb{P}(d(c,c)<\varepsilon ) \\
&=&\mathbb{P}(\Omega )=1.
\end{eqnarray*}

\noindent Then 
\begin{equation*}
\lim \sup_{n\rightarrow +\infty }\mathbb{P}(d(X_{n},c)\geq \varepsilon
)=1-\lim \inf_{n\rightarrow +\infty }\mathbb{P}(d(X_{n},c)\leq \varepsilon
)\leq 1-1=0.
\end{equation*}

\noindent Then for any $\varepsilon >0,$%
\begin{equation*}
\lim \sup_{n\rightarrow +\infty }\mathbb{P}(d(X_{n},c)>\varepsilon )\leq
\lim \sup_{n\rightarrow +\infty }\mathbb{P}(d(X_{n},c)\geq \varepsilon )=0.
\end{equation*}%
Hence $X_{n}\rightarrow _{\mathbb{P}}c$.\\

\bigskip \noindent \textbf{(c) Two equivalent sequences  in probability weakly converge to
the same limit if one of them does.}\\

\begin{proposition} \label{cvcp.prop4} If $X_{n}\rightsquigarrow X$ \ and $%
d(X_{n},Y_{n})\longrightarrow _{\mathbb{P}}0$ as $n \rightarrow +\infty$
, then $Y_{n}\rightsquigarrow X$.
\end{proposition}

\noindent \textbf{Proof}. Suppose $X_{n}\rightsquigarrow X$ \ and $%
d(X_{n},Y_{n})\longrightarrow _{\mathbb{P}}0.$ Let us prove that $%
Y_{n}\rightsquigarrow X$ by using Point $(viia)$ of Proposition \ref{cvcp.prop2}
above. Let  $f$ be a bounded Lipschitz function with coefficient $\ell >0$ and bound  
$M.$ We have for any $n\geq 0$ and $\varepsilon >0,$ \newline
\ 
\begin{eqnarray*}
\left\vert \mathbb{E}f(Y_{n})-\mathbb{E}f(X)\right\vert  &\leq &\mathbb{E}%
\left\vert f(Y_{n})-f(X)\right\vert  \\
&\leq &\mathbb{E}\left\vert f(X_{n})-f(X)\right\vert +\mathbb{E}\left\vert
f(Y_{n})-f(X_{n})\right\vert .
\end{eqnarray*}

\noindent By applying Point $(viia)$ of Lemma \ref{cvcp.lem01} above and by
using the weak limit $X_{n}\rightsquigarrow X$, we get
\begin{equation*}
\lim \sup_{n\rightarrow +\infty }\left\vert \mathbb{E}f(Y_{n})-\mathbb{E}%
f(X)\right\vert \leq \lim \sup_{n\rightarrow +\infty }\mathbb{E}\left\vert
f(Y_{n})-f(X_{n})\right\vert .
\end{equation*}

\noindent Now we use the same method used in the proof of Proposition \ref{cvcp.prop2} to have 
\begin{eqnarray*}
\mathbb{E}\left\vert f(Y_{n})-f(X_{n})\right\vert  &\leq&\int_{(d(Y_{n},X_{n})\leq \varepsilon )}\left\vert f(Y_{n})-f(X_{n})\right\vert d\mathbb{P}\\
&+&\int_{(d(Y_{n},X_{n})>\varepsilon )}\left\vert f(Y_{n})-f(X_{n})\right\vert d\mathbb{P} \\
&\leq &\ell \varepsilon +2M\text{ }d(Y_{n},X_{n}),
\end{eqnarray*}

\noindent which tends to zero as $n\rightarrow +\infty $ and next  $
\varepsilon \downarrow 0.$ We conclude that
\begin{equation*}
\lim \sup_{n\rightarrow +\infty }\left\vert \mathbb{E}f(Y_{n})-\mathbb{E}%
f(X)\right\vert =0.
\end{equation*}

\bigskip \noindent \textbf{(d) Slutsky's Theorem}.\\

\noindent We have the following important and yet simple tool in weak convergence.

\begin{proposition} \label{cv.slutsky}
\bigskip\ If \ $X_{n}\rightsquigarrow X$ \ and $Y_{n}\rightsquigarrow c$,
then $(X_{n},Y_{n})\rightsquigarrow (X,c)$
\end{proposition}

\noindent \textbf{Proof}. Let  $X_{n}\rightsquigarrow X$ \ and $%
Y_{n}\longrightarrow _{\mathbb{P}}c.$ We want to show that $%
(X_{n},Y_{n})\rightsquigarrow (X,c).$ We first remark that \ $%
Y_{n}\longrightarrow _{\mathbb{P}}c$ since $Y_{n}\rightsquigarrow c.$\ Next,
on  $S^{2}$ endowed with the euclidean metric,

\begin{equation*}
d_{e}((x^{\prime },y^{\prime }),(x^{\prime \prime },y^{\prime \prime }))=%
\sqrt{d(x^{\prime },x^{\prime \prime })^{2}+d(y^{\prime },y^{\prime \prime
})^{2}},
\end{equation*}

\noindent we have 
\begin{equation*}
d_{e}((X_{n},Y_{n}),(X_{n},c))=d(Y_{n},c).
\end{equation*}

\noindent It comes that for any $\varepsilon >0$,
\begin{equation*}
\lim \sup_{n\rightarrow +\infty }\mathbb{P}(d_{e}((X_{n},Y_{n}),(X_{n},c))>%
\varepsilon )=\lim \sup_{n\rightarrow +\infty }\mathbb{P}(d(Y_{n},c)>%
\varepsilon )=0,
\end{equation*}

\noindent since $Y_{n}\longrightarrow _{\mathbb{P}}c$. Then $%
d_{e}((X_{n},Y_{n}),(X_{n},c))\longrightarrow _{\mathbb{P}}0.$ By
Proposition \ref{cvcp.prop4}, it is enough to have the weak limit of $(X_{n},c)$
which will be that of $(X_{n},Y_{n})$.\\

\noindent To show the weak convergence of $(X_{n},c)$ to $(X,c)$, \ we consider a real
bounded and continuous function $g(\cdot ,\cdot )$ defined on $S^{2}$ and
try to show that $\mathbb{E}g(X_{n},c)\rightarrow \mathbb{E}g(X,c).$ But it
comes from that  $c$ is fixed \ and the function  $f(x)=g(x,c)$ is bounded
and continuous and then $\mathbb{E}f(X_{n})\rightarrow \mathbb{E}f(X)$ since 
$X_{n}\rightsquigarrow X$. But  $\mathbb{E}f(X_{n})\rightarrow \mathbb{E}f(X)
$ is  $\mathbb{E}g(X_{n},c)\rightarrow \mathbb{E}g(X,c).$ This finishes the
proof.\newline

\bigskip \noindent \textbf{(e) Coordinate-wise convergence in probability}.\\

\begin{proposition} \label{cv.CPcoordinates} $X_{n}\longrightarrow _{\mathbb{P}}X$ \ and $%
Y_{n}\longrightarrow _{\mathbb{P}}Y$ if and only if  $(X_{n},Y_{n})%
\longrightarrow _{\mathbb{P}}(X,Y)$.
\end{proposition}

\bigskip \noindent \textbf{Proof}. Suppose that $X_{n}\longrightarrow
_{\mathbb{P}}X$ \ and $Y_{n}\longrightarrow _{\mathbb{P}}Y.$ Let us use the
Manhattan distance on $S^{2}:$

\begin{equation*}
d_{m}((x^{\prime },y^{\prime }),(x^{\prime \prime},y^{\prime \prime }))=d(x^{\prime },x^{\prime \prime })+d(y^{\prime
},y^{\prime \prime }).
\end{equation*}

\noindent For any  $\varepsilon >0,$ $\lim \sup_{n\rightarrow +\infty }%
\mathbb{P}(d_{m}((X_{n},Y_{n}),(X,Y))>\varepsilon )$ is  
\begin{eqnarray*}
&=&\limsup_{n\rightarrow +\infty }\mathbb{P}(d(X_{n},X)+d(Y_n,Y))>%
\varepsilon ) \\
&\leq &\limsup_{n\rightarrow +\infty} \biggr(\mathbb{P}(d(X_{n},X)>\varepsilon
/2)+\mathbb{P}(d(Y_{n},Y)>\varepsilon /2) \biggr) \\
&\leq &\lim \sup_{n\rightarrow +\infty }\mathbb{P}(d(X_{n},X)>\varepsilon
/2)+\lim \sup_{n\rightarrow +\infty }\mathbb{P}(d(Y_{n},Y)>\varepsilon /2) \\
&=&0.
\end{eqnarray*}

\bigskip \noindent Conversely, suppose that $(X_{n},Y_{n})\longrightarrow _{\mathbb{P}%
}(X,Y).$ Then for any $n\geq 1,$

\begin{equation*}
d(X_{n},X)\leq d(X_{n},X)+d(Y_{n},Y)=d_{m}((X_{n},Y_n),(X,Y))\rightarrow _{%
\mathbb{P}}0\text{ as }n\rightarrow +\infty .
\end{equation*}

\noindent Then $d(X_{n},X)\rightarrow _{\mathbb{P}}0$ as $n\rightarrow +\infty $ and,
in the same manner, $d(Y_{n},Y)\rightarrow _{\mathbb{P}}0.$

\subsection{Skorohod-Wichura Theorem} \label{cv.subsec.skorohod} $ $\\

\noindent \noindent We only state this result in a complete and separable metric space.\\

\begin{theorem} \label{cv.skorohodWichura} Let $(X_{n})_{n\geq 0}$, and $X$ be of measurable applications with values in $(S,d)$, a complete and separable space, not necessarily defined on the same probability space.\\

\noindent If $X_{n}\rightsquigarrow X$, then there exists a probability space ($\Omega ,\mathcal{A},\mathbb{P})$ holding measurable
applications $\ (Y_{n})_{n\geq 0}$ and $Y$  such that  
\begin{equation*}
\mathbb{P}_{X}=\mathbb{P}_{Y}\text{ and }\left( \forall n\geq 0,\mathbb{P}%
_{X_{n}}=\mathbb{P}_{Y_{n}}\right)
\end{equation*}

\noindent and
\begin{equation*}
Y_{n}\rightarrow Y,\text{ }a.e.
\end{equation*}
\end{theorem}

\bigskip

\noindent This theorem is powerful and may reveal itself very usefull in a great number of situations. You will find a proof of it
for $S=\mathbb{R}$ in Chapter \ref{cv.R}, Theorem \ref{cv.skorohod}.\\

\bigskip

\section{Appendix} \label{cv.annexe} 

\subsection{$F$-continuous intervals, where $F$ is a distribution function} \label{cv.subsFcontinuous} $ $\\

\noindent Let $\mathbb{P}$ be a probability measure $\mathbb{P}$ on $(\mathbb{R}^{k},%
\mathcal{B}(\mathbb{R}^{k}))$.  Consider its distribution function 
\begin{equation*}
(x_{1},...,x_{k})\hookrightarrow F(x_{1},...,x_{k})=P\left(
\prod\limits_{i=1}^{k}]-\infty ,x_{i}]\right) .
\end{equation*}

\subsubsection{$F$-continuous intervals}

\bigskip \noindent Let  
\begin{equation*}
]a,b]=\prod\limits_{i=1}^{k}]a_{i},b_{i}]
\end{equation*}

\noindent be and interval of $R^{k}$. Define 
\begin{equation*}
E(a,b)=\{c=(c_{1},...,c_{k})\in \mathbb{R}^{k},\text{ }\forall 1\leq i\leq
k,(c_{i}=a_{i}\text{ ou }c_{i}=b_{i})\}.
\end{equation*}

\noindent We may use extra-notations to get compact forms of $E(a,b)$. Define the product of $k$-tuples term by term
$$
(x_1,...,x_k)*(y_1,...,y_k)=(x_{1}y_{1},...,x_{k}y_{k}).
$$

\noindent We also have 

\begin{equation} \label{sec_append_not}
E(a,b)=\{b+\varepsilon*(a-b), \varepsilon=(\varepsilon_1,...,\varepsilon_k) \in \{0,1\}^k\}.
\end{equation}

\bigskip \noindent  We say that the interval $(a,b)$ is $F$-continuous if and only if $(a,b)$ is bounded and each element of $E(a,b)$ is a continuity point of $F$, that is
\begin{equation*}
\forall c\in E(a,b),\mathbb{P}(\partial ]-\infty ,c])=0.
\end{equation*}

\bigskip \noindent Let $\mathcal{U}$ be the class of all $F$-continuous intervals. By convention, we say that the empty set is an $F$-continuous interval. Here are some properties of $\mathcal{U}$.

\subsubsection{$\mathcal{U}$ is stable by finite intersection} \label{subsubsecUstable}

\bigskip \noindent  Take $]a,b]=\prod\limits_{i=1}^{k}]a_{i},b_{i}]\in \mathcal{U}$ and $]c,d]=\prod\limits_{i=1}^{k}]c_{i},d_{i}] \in 
\mathcal{U}$. We have
\begin{equation*}
]a,b]\cap ]c,d]=\prod\limits_{i=1}^{k}]a_{i}\vee c_{i},b_{i}\wedge
d_{i}]=]\alpha ,\beta ],
\end{equation*}

\bigskip \noindent where $x\vee y$ and $x\wedge y$ respectively stand for the maximum and the minimum of $x$ and $y$, and \ \ $\alpha =(a_{1}\vee c_{1},..,a_{k}\vee c_{k})$ and $\beta =(b_{1}\wedge d_{1},...,b_{k}\wedge d_{k})$. If $]a,b]\cap ]c,d]$ is empty, it is in $\mathcal{U}$. Otherwise, none of the factor $]a_{i}\vee c_{i},b_{i}\wedge d_{i}]$ is empty. We are going to show that : 
\begin{equation}
\forall e\in E(\alpha ,\beta ),\partial ]-\infty ,e]\subset
\bigcup\limits_{z\in E(a,b)\cup E(c,d)}\partial ]-\infty ,z].
\label{unionEab}
\end{equation}

\bigskip \noindent Indeed, take $e\in E(\alpha ,\beta )$. We have
\begin{equation*}
e_{i}=a_{i}\vee c_{i}\text{ or }b_{i}\wedge d_{i},\text{ }1\leq i\leq k.
\end{equation*}

\noindent Take $t\in \partial ]-\infty ,e]$. This means that
\begin{equation*}
(t_{i}\leq e_{i},\text{ }1\leq i\leq k)\text{ and }(\exists i_{0}\text{,}%
t_{i_{0}}=e_{i_{0}})
\end{equation*}

\bigskip \noindent Since $]\alpha ,\beta ]$ is included in $]a,b]$ and in $]c,d]$,  $t$ satisfies
\begin{equation*}
t_{i}\leq b_{i}\text{ and }t_{i}\leq d_{i},1\leq i\leq k.
\end{equation*}

\bigskip \noindent Now, let us consider $i_{0}$ such that $t_{i_{0}}=c_{i_{0}}$. We have four cases
\begin{equation*}
\left\{ 
\begin{tabular}{lll}
$t_{i_{0}}=e_{i_{0}}=a_{i_{0}}\vee c_{i_{0}}=a_{i_{0}}$ & $\Longrightarrow $
& $t_{i_{0}}=a_{i_{0}}$ and $t_{i}\leq b_{i},1\leq i\leq k$ \\ 
$t_{i_{0}}=e_{i_{0}}=a_{i_{0}}\vee c_{i_{0}}=c_{i_{0}}$ & $\Longrightarrow $
& $t_{i_{0}}=c_{i_{0}}$ and $t_{i}\leq d_{i},1\leq i\leq k$ \\ 
$t_{i_{0}}=e_{i_{0}}=b_{i_{0}}\wedge d_{i_{0}}=b_{i_{0}}$ & $\Longrightarrow 
$ & $t_{i_{0}}=b_{i_{0}}$ and $t_{i}\leq b_{i},1\leq i\leq k$ \\ 
$t_{i_{0}}=e_{i_{0}}=b_{i_{0}}\wedge d_{i_{0}}=d_{i_{0}}$ & $\Longrightarrow 
$ & $t_{i_{0}}=d_{i_{0}}$ and $t_{i}\leq d_{i},1\leq i\leq k$%
\end{tabular}%
.\right. 
\end{equation*}

\bigskip \noindent We are going to conclude by considering each line of the formula above.\\

\noindent First line : $t\in
\partial ]-\infty ,z_{1}]$ where $%
z_{1}=(b_{1},...,b_{i_{0}-1},a_{i_{0},}b_{i_{0}+1},b_{k})\in E(a,b)$.\\

\noindent Second line : $t\in \partial ]-\infty ,z_{2}]$ where $z_{2}=(d_{1},...,d_{i_{0}-1},c_{i_{0},}d_{i_{0}+1},d_{k})\in E(c,d)$.\\

\noindent Third line : $t\in \partial ]-\infty ,b]$ and of course $b\in E(a,b)$.\\ 

\noindent Fourth line : $t\in \partial ]-\infty ,d]$ and of course $d\in E(c,d).$ So $t$ is one of
the  $\partial ]-\infty ,z]$ with $z\in E(a,b)\cup E(c,d)$.\\

\noindent  So \ref{unionEab} holds, and since the union is a finite union of null sets, we have

\begin{equation*}
\forall e\in E(\alpha ,\beta ),P(\partial ]-\infty ,e])=0.
\end{equation*}

\noindent Therefore, $\mathcal{U}$ is stable by finite intersection.\\

\begin{lemma} \label{cv.annexe.Inc}
Every neighborhood of an arbitrary point $x$ includes a $F$-continuous interval $]a,b]$ containing $x$.
\end{lemma}

\bigskip \noindent Let $V$ be a neighborhood of $x$. There exists an interval such that $]a,b[$
\begin{equation*}
x\in \prod\limits_{i=1}^{k}]a_{i},b_{i}[.
\end{equation*}

\bigskip \noindent Set
\begin{equation*}
\varepsilon _{0}=\min (x_{i}-a_{i},1\leq i\leq k)\wedge \min
(b_{i}-x_{i},1\leq i\leq k),
\end{equation*}

\bigskip \noindent denote by  $\delta =(1,...,1)$ the vector of $\mathbb{R}^{k}$ whose all components are equal to one. Then for $0<\varepsilon <\varepsilon _{0}$, we have  
\begin{equation*}
]a+\varepsilon \delta ,x+\varepsilon \delta ]\subset ]a,b[.
\end{equation*}

\bigskip \noindent Each point $e$ of $E(a+\varepsilon \delta ,x+\varepsilon \delta )$ is of the form
\begin{equation*}
t(\varepsilon)=(t_{1}+\varepsilon ,t_{2}+\varepsilon ,...,t_{k}+\varepsilon,
)
\end{equation*}

\bigskip \noindent with, of course, $t_{i}=a_{i}$ or $t_{i}=x_{i}$. For any choice of these $t=(t_{1},...,t_{k})$, the sets
$\partial ]-\infty ,t(\varepsilon )]$ are disjoint. Then, by Proposition \ref{cv.annexe.FDS} below we have
\begin{equation*}
\mathbb{P}(\partial ]-\infty ,t(\varepsilon )])>0,
\end{equation*}

\noindent except, eventually, when $\varepsilon$ is out of countable set $D_t \subset ]0,\varepsilon _{0}[$. But $D=\cup _{t}D(t) \subset ]0,\varepsilon _{0}[$ is is countable, since it is at most a union of  $2^{k}$ countable sets. Hence, surely, we may pick a value of $\varepsilon$ out of $]0,\varepsilon _{0}[$, such that for any vector $e$ satisfying
\begin{equation*}
e_{i}=a_{i}+\varepsilon \text{ or }x_{i}+\varepsilon 
\end{equation*}

\bigskip \noindent we have  
\begin{equation*}
P(\partial ]-\infty ,t(\varepsilon )])=0
\end{equation*}

\bigskip \noindent and 
\begin{equation*}
x\in ]a+\varepsilon \delta ,x+\varepsilon \delta ]\subset ]a,b[.
\end{equation*}

\bigskip \noindent We just proved that there exists $]A_{x},B_{x}[=]a+\varepsilon \delta
,x+\varepsilon \delta /2[$ and  $]a_{x},b_{x}]=]a+\varepsilon \delta
,x+\varepsilon \delta ]$ such that
\begin{equation}
x\in ]A_{x},B_{x}[\subset ]a_{x},b_{x}]\subset V.  \label{cleEab}
\end{equation}

\noindent Let us use this to show that \textit{any open set $G$ of  $\mathbb{R}^{k}$ is a countable union of $F$-continuous intervals.}

\bigskip \noindent  Indeed, by (\ref{cleEab}), any open set $G$ may be written as
\begin{equation*}
G=\bigcup\limits_{x\in G}]A_{x},B_{x}[.
\end{equation*}

\bigskip \noindent Since $\mathbb{R}^{k}$ is a separable space, this open cover reduces to a countable cover, that is, there exists a sequence $(x_{j})_{j\geq 0}\subset G$ such that 
\begin{equation*}
G=\bigcup\limits_{j\geq 0}]A_{x_{j}},B_{x_{j}}[.
\end{equation*}

\bigskip \noindent We finally get 
\begin{equation*}
G=\bigcup\limits_{j\geq 0}]a_{x_{j}},b_{x_{j}}],
\end{equation*}

\bigskip \noindent where the $]a_{x_{j}},b_{x_{j}}]$ are $F$-continuous intervals. We have this proposition.

\begin{proposition} \label{cv.GFcontinuous}
Let $F$ be any probability distribution function on $\mathbb{R}^{k}$, $k\geq 1$. Then any open $G$ set in $\mathbb{R}^{k}$ is a countable union of $F$-continuous intervals of the form $]a,b]$ or $]a,b[$, where by definition, an interval $(a,b)$ is $F$-continuous if and only if, for any 
$$
\varepsilon=(\varepsilon_1, \varepsilon_2, ..., \varepsilon_k) \in \{0,1\}^k
$$, 

\noindent the point 
$$
b+\varepsilon*(a-b)=(b_1+\varepsilon_1 (a_1-b_1), b_2+\varepsilon_2 (a_2-b_2), ..., b_k+\varepsilon_k (a_k-b_k))
$$ 

\noindent is a continuity point of $F$.
\end{proposition}

\subsection{Semi-continuous Functions} \label{cv.subsec.annexe.semic} $ $\\

\noindent A function $f : S \mapsto \overline{\mathbb{R}}$, where $S$ is a metric space, is continuous if and only if \\

\noindent (i) For any $x \in \mathbb{R}$, for any $\varepsilon >0$, there exists a neighborhood $V$ of $x$ such that 
\begin{equation*}
y\in V\Rightarrow f(y)\in ]f(x)-\varepsilon ,f(x)+\varepsilon \lbrack .
\end{equation*}

\bigskip \noindent In this formula, we use the whole interval $]f(x)-\varepsilon ,f(x)+\varepsilon \lbrack$\ in the definition. But we might
be interested only by one the half intervals. This gives semi-continuous functions. Precisely, a real-valued function $f$ is upper semi-continuous (u.s.c for short) if and only if \\

\bigskip \noindent (ii) For any $x\in \mathbb{R}$, for any $\varepsilon >0$, there exists a neighborhood $V$ of $x$ such that 
\begin{equation*}
y\in V\Rightarrow f(y)<f(x)+\varepsilon,
\end{equation*}

\noindent and a real-valued function $f$ is lower semi-continuous (l.s.c for short) if and only if \\

\bigskip \noindent (iii) For any $x\in \mathbb{R}$, for any $\varepsilon >0$, there exists a neighborhood $V$ of $x$ such that 
\begin{equation*}
y\in V\Rightarrow f(y)>f(x)-\varepsilon.
\end{equation*}

\bigskip \noindent We have two immediate remarks.\\

\noindent (a) A real function is continuous if and only if it is both \textit{u.s.c} and  \textit{l.s.c}.\\

\noindent (b) A real function $f$ is \textbf{u.s.c.} if and only if its opposite $-f$ is \textit{l.s.c}.\\

\noindent Here is a characterization of real-valued semi-continuous functions.\\

\begin{proposition}\label{cv.annexe.SC} We have the following properties :\\

\noindent (1) $f : S \mapsto \overline{\mathbb{R}}$ is upper semi-continuous function if and only if the set $(f\geq c)$ is closed for any real number $c \in \mathbb{R}$.\\

\noindent (2) $f$ is lower semi-continuous if and only if the set $(f\leq c)$ is closed for any real number $c \in \mathbb{R}$.\\

\noindent (3) If $f$ is \textit{u.s.c} or \textit{l.s.c}, it is measurable.\\
\end{proposition}

\bigskip \noindent \textbf{Proof}. Proof of Point (1). Let us begin by the direct implication. Let $f$ be a \textit{u.s.c} function from $S$ to $\mathbb{R}$. Let us show that the set $(f\geq c)$ is closed by showing that the set $(f < c)$ is open. Let $x \in G^{c}=(f < c)$, that is $f(x)<c$. Let us take $\varepsilon=c-f(x) >0$. Since $f$ is \textit{u.s.c}, there exists a neighborhood $V$ of $x$ such that 

\begin{equation*}
y\in V\Rightarrow f(y) < f(x)+\epsilon=c,
\end{equation*}

\noindent which may be written as
\begin{equation*}
y\in V \Rightarrow f(y) < c,
\end{equation*}

\noindent which means that $V\subseteq G^c$. We proved that $G^c$ contains each of its elements with one of their neighborhood. Then $G^c$ is open. This proves the direct sens.\\

\noindent Now suppose $(f\geq c)$ is closed for any real number $c$. Fix $x$ in $S$. Then for any $\epsilon >0$, the set $G=(f<f(x)+\epsilon)$ is open and $x \in G$. Then, there is a neighborhood of $x$ such that $x \in V \subset G$. We conclude that : for any $x\in S$, for any $\epsilon>0$, there exists a neighborhood of $x$ such that 
\begin{equation*}
y\in V\Rightarrow f(y)\leq f(x)+\epsilon.
\end{equation*}

\noindent So $f$ is \textbf{u.s.c.}. This completes the proof of Point (1).\\

\noindent Point (2) is proved by applying Point (1) to $-f$.\\

\noindent Point (3) is a consequence of Points (1) and (2) and classical measurability criteria for real-valued functions.

\subsection{Probabilistic property of a non-countable family of disjoint events} \label{cv.subsec.annexe.Disjoint} $ $\\

\noindent 
\begin{proposition} \label{cv.annexe.FDS} Let $(B_{\lambda })_{\lambda \in \Gamma }$ be a family of disjoint measurable sets in a probability space $(\Omega,\mathcal{A},\mathbb{P})$. Then at most, a countable number of them are not null sets, or equivalently, the cardinality of the elements $\lambda$ of $\Gamma$ for which $\mathbb{P}(B_{\lambda})>0$, is at most countable.
\end{proposition}

\bigskip \noindent \textbf{Proof}. Define
\begin{equation*}
D=\{\lambda \in \Gamma ,\text{ }L(B_{\lambda })>0\}.
\end{equation*}

\noindent and for any integer $k\geq 1$, 
 
\begin{equation*}
D_{k}=\{\lambda \in \Gamma ,\text{ }L(B_{\lambda })>1/k\}.
\end{equation*}

\noindent It is clear that we have
\begin{equation*}
D=\cup _{k\geq 1}D_{k},
\end{equation*}

\noindent We are going to prove that each $D_{k}$ is finite. Indeed, suppose we can find $r\geq 1$ elements in $D_{k}$
denoted as $\lambda _{1},\lambda _{2},...,\lambda _{r}$. Since the $B_{\lambda}$'s are disjoint, we have
\begin{equation*}
1\geq \mathbb{P}\left(\bigcup_{1}^{r}B_{\lambda _{j}}\right)=\sum_{1}^{r} \mathbb{P}(B_{\lambda _{j}})\geq
r/k.
\end{equation*}

\noindent Then 
\begin{equation*}
r\leq k.
\end{equation*}

\noindent This means that we cannot choose more that $k$ points in $D_{k}$. Hence $D$ is finite, that cardinality of $D_k$ is less than $k$. Thus, $D$ is at most countable as a countable union of finite sets.

\subsection{Measurability of the set discontinuity points in a metric space} \label{cv.subsec.annexe.Discontinuity} $ $\\

\noindent Here is an amazing result, that is the sets of discontinuity points of a function defined from a metric space to another metric space is measurable whatever be the function.\\

\begin{lemma} \label{cv.annexe.Discont}
Let $g$ be a function $g$ from the metric space $(S, d)$ to the metric space $(D, r)$. Denote by $discont(g)$, the set of discontinuity points of $g$. We have

\begin{equation}
discont(g)=\bigcup_{s=1}^{\infty }\bigcap_{t=1}^{\infty }B_{s,t}, \label{annexe4}
\end{equation}%

\noindent where for each couple of positive integers $(s,t)$ 
\begin{equation*}
B_{s,t}=\left\{ x\in S,\exists (y,z)\in S^{2},\text{ }d(x,y)<1/t,d(z,x)<1/t,%
\text{ }r(g(y),g(z))\geq 1/s\right\} .
\end{equation*}

\noindent is an open set.
\end{lemma}

\bigskip \noindent From this lemma, we see that $discont(g)$ is measurable as countable unions and intersections of open sets. But we have to prove the lemma.\\

\noindent \textbf{Proof of the lemma}. Let us show that  
\begin{equation*}
\bigcup_{s=1}^{\infty }\bigcap_{t=1}^{\infty }B_{s,t}\subseteq discont(g).
\end{equation*}

\bigskip \noindent Let  $x\in \bigcup_{s=1}^{\infty }\bigcap_{t=1}^{\infty
}B_{s,t}$. Then there exists an integer $s\geq 1$ fixed  such that for any integer $t\geq 1$, there exist $y_{t}$ and $z_{t}$ such that
\begin{equation*}
d(x,y_{t})<1/t,
\end{equation*}

\bigskip \noindent and
\begin{equation*}
d(x,z_{t})<1/t, 
\end{equation*}

\bigskip \noindent and
\begin{equation}
\forall \text{ }t\geq 1,\text{ }r(g(y_{t}),g(z_{t}))\geq 1/s  \label{cv29}
\end{equation}

\bigskip \noindent Since $g$ is continuous at $x$, we get, as $t \rightarrow +\infty$,  
\begin{equation*}
r(g(y_{t}),g(z_{t}))\leq r(g(y_{t}),g(x))+r(g(x_{t}),g(z_{t}))\rightarrow 0,
\end{equation*}

\bigskip \noindent which is in contradiction with (\ref{cv29}). Then $x$ is a discontinuity point of $g$.\\

\bigskip \noindent Reversely, we have to show that
\begin{equation*}
discont(g)\subseteq \bigcup_{s=1}^{\infty }\bigcap_{t=1}^{\infty }B_{s,t}.
\end{equation*}

\bigskip \noindent Let $x$ be a discontinuity point of $g$. By the negation of the definition of the continuity, we have, 
\begin{equation*}
\exists \text{ }\epsilon >0,\forall \text{ }\eta >0,\exists \text{ }y\in S,%
\text{ }d(x,y)<\eta ,\text{ \ }r(g(y),g(x))\geq \varepsilon .
\end{equation*}

\bigskip \noindent Let $s$ be an integer such that $\varepsilon \geq 1/s$. Then for any $1/t$ where $t$ is a positive integer, we have
\begin{equation*}
\exists \text{ }y\in S,\text{ }d(x,y)<1/t,\text{ \ }r(g(y),g(x))\geq 1/s.
\end{equation*}

\bigskip \noindent Putting $z=x$, leads to  
\begin{equation*}
d(x,z)<1/t,\text{ }d(x,y)<1/t,\text{ \ }r(g(y),g(x))\geq 1/s.
\end{equation*}

\bigskip \noindent Then $x\in \bigcup_{s=1}^{\infty }\bigcap_{t=1}^{\infty
}B_{s,t}$.\\

\noindent By combining the two steps, we get the equality.\\

\bigskip \noindent Let us prove that for each couple of positive integers $(s,t)$, $B_{s,t}$ is an open set. Fix $s\geq 1$, $t\geq 1$. Put
$a=1/s>0$ and $b=1/t>0$. Let $x\in B_{s,t}$. Then  
\begin{equation*}
\exists (y,z)\in S^{2},\text{ }d(x,y)<b,d(z,x)<b,\text{ }r(g(y),g(z))\geq a
\end{equation*}

\bigskip \noindent Set $c=min(b-d(x,y),b-d(z,x))>0$ and take  $x^{\prime }\in
B(x,c)$. Then 
\begin{equation*}
d(x^{\prime },y)\leq d(x^{\prime },x)+d(x,y)<c+d(x,y)\leq b
\end{equation*}

\bigskip \noindent and next, 
\begin{equation*}
d(x^{\prime },z)<d(x^{\prime },x)+d(x,z)\leq c+d(x,z)\leq b
\end{equation*}

\bigskip \noindent and
\begin{equation*}
r(g(y),g(z))\geq a.
\end{equation*}

\bigskip \noindent Thus,  $x^{\prime }\in B_{s,t}$. Hence  
\begin{equation*}
x\in B(x,c)\subseteq B_{s,t}.
\end{equation*}

\bigskip \noindent Therefore each $B_{s,t}$ contains each of its point with an open ball. Hence $B_{s,t}$ is an open set.\\

\noindent We finished the proof of the lemma, which proves the measurability of $g$.\\

\subsection{Stone-Weierstrass Theorem} \label{cv.subsec.annexe.SW1} $ $\\

\noindent \textbf{I - Stone-Weierstrass Theorem}.\\

\noindent Here are two forms of Stone-Weierstrass Theorem. The second is more general and is the one we use in this text.

\begin{proposition} \label{proba02_rv_append_sw_prop1}
Let $(S, d)$ be a compact metric space and $H$ a non-void subclass of the class $\mathcal{C}(S, \mathbb{R})$ of all real-valued continuous functions defined on $S$. Suppose that $H$ satisfies the following conditions.\\

\bigskip \noindent (i) $H$ is \textit{lattice}, that is, for any couple $(f,g)$ of elements of $H$,  $f\wedge g$\ et $f\vee g$ are in $H$
\newline

\bigskip \noindent (ii) For any couple $(x,y)$ of elements of $S$ and for any couple $(a,b)$ of real numbers such that $a=b$ if $x=y$, there exists a couple $(h,k)$ of elements of $H$ such that
\begin{equation*}
h(x)=a\text{ and }k(y)=b.
\end{equation*}

\bigskip \noindent Then $H$ is dense in $\mathcal{C}(S, \mathbb{R})$ endowed with the uniform topology, that is each continuous
function  from $S$ to $\mathbb{R}$ is the uniform limit of a sequence of elements in $H$.
\end{proposition}

\begin{theorem} \label{proba02_rv_append_sw_prop2}
Let $(S, d)$ be a compact metric space and $H$ a non-void subclass of the class $\mathcal{C}(S, \mathbb{C})$ of all real-valued continuous functions defined on $S$. Suppose that $H$ satisfies the following conditions.\\

\noindent (i) $H$ contains all the constant functions.\\

\noindent (ii) For all $(h,k)\in H^{2}$, $h+k\in H,h\times k\in H,\overline{u}\in H$.\\

\noindent (iii) $H$ separates the points of $S$, i.e., for two distinct elements of $S$, $x$ and $y$, that is $x\neq y$, there exists $h\in H$ such that 
\begin{equation*}
h(x)\neq h(y).
\end{equation*}

\noindent Then $H$ is dense in $\mathcal{C}(S, \mathbb{C})$ endowed with the uniform topology, that is each continuous
function  from $S$ to $\mathbb{C}$ is the uniform limit of a sequence of elements in $H$.
\end{theorem}

\bigskip \noindent \textbf{Remark}.\newline

\noindent If we work in $\mathbb{R}$, the condition on the conjugates - $\overline{u}\in H$\ - becomes needless.\\

\noindent But here, these two classical versions do not apply. We use the following extension.

\begin{corollary} \label{sec_EF_cor_05} Let $K$ be a non-singleton compact space and $\mathcal{A}$ be a non-empty sub-algebra of $C(K,\mathbb{C})$. Let $f \in C(K,\mathbb{C})$. Suppose that there exists $K_0 \subset K$ such that $k\setminus K_0$ has at least two elements and $f$ is constant on $K_0$. Suppose that the following assumption hold.

\noindent (1) $\mathcal{A}$ separates the points of $K\setminus K_0$ and separates any point of $K_0$ from any point of $K\setminus K_0$.\\

\noindent (2) $A$ contains all the constant functions.\\

\noindent (3) For all $f \in \mathcal{A}$, its conjugate function $\bar{f}=\mathcal{R}(f) - i \mathcal{Im}(f) \in \mathcal{A}$,\\

\noindent  Then 

$$
f \in \overline{\mathcal{A}}.
$$
\end{corollary}

\noindent A proof if it available in \cite{loSW2018}.\\

\subsection{A useful remark} \label{cv.subsec.annexe.divers} $ $\\

\noindent \noindent \textbf{The min function is Lipschitz}. We have for any real numbers $x$, $y$, $X$, and $Y$, 
\begin{equation}
\left\vert \min (x,y)-\min (X,Y)\right\vert \leq \left\vert x-X\right\vert
+\left\vert y-Y\right\vert.  \label{annexe2}
\end{equation}

\noindent To see that, let us have a look at the four possible cases.\\

\noindent Case 1 : $min(x,y)=x$ and $min(X,Y)=X$. We have
\begin{equation*}
\left\vert \min (x,y)-\min (X,Y)\right\vert \leq \left\vert x-X\right\vert.
\end{equation*}

\noindent Case 2 : $min(x,y)=x$ and $min(X,Y)=Y$. If $x\leq Y$, since $Y\geq X$, we have
\begin{equation*}
0 \leq \min (X,Y) - \min (x,y)=Y-x \leq X-x.
\end{equation*}

\noindent If $x>Y$, since $X\geq Y$, we have
\begin{equation*}
0 \leq \min (x,y)-\min (X,Y)=x-Y\leq y-Y.
\end{equation*}

\noindent Case 3 : $min(x,y)=y$ and $min(X,Y)=Y$. We have
\begin{equation*}
\left\vert \min (x,y)-\min (X,Y)\right\vert \leq \left\vert y-Y\right\vert.
\end{equation*}

\noindent Case 4 : $min(x,y)=y$ and $min(X,Y)=X$. This case id handled as for Case 2 by permuting the roles of $(x,y)$ and $(X,Y)$.\\

\noindent We get (\ref{annexe2}) by putting together the results of the four cases.
 

%% file: asymptotics_cv_02_en.tex
\chapter{Uniform Tightness and Asymptotic Tightness} \label{cv.tensRk}

\section{Introduction}

Any limit theory deals with the notion of compactness through the existence or not for sequences of sub-sequences converging in the sense of
the defined limit. This corresponds to the Bolzano-Weierstrass for real sequences. For the weak convergence, the condition of the existence of such sub-sequences is called \textsl{tightness}. When dealing with weak convergence for general metric spaces, \textit{tightness} leads to the general Prohorov theorem which establishes, under eventually other assumptions, that every \textit{uniformly tight} sequence of measurable applications of a metric space $(S,d)$ has at least a weakly converging sub-sequence.\\

\noindent In this chapter, we focus on weak convergence in $\mathbb{R}^{k}$. And there exists a specific handling of \textit{weak compactness} that is very different from the treatment in the general case. In $\mathbb{R}^{k}$, the major role is played by the theorem of Helly-Bray that directly makes use of the Bolzano-Weierstrass theorem in $\mathbb{R}$.\newline

\noindent Since, we deal with compact sets of $\mathbb{R}^{k}$, just remind two properties which we are going to use. The first is that compact sets of $\mathbb{R}^{k}$ are closed and bounded sets. The second is that $\mathbb{R}^{k}$ is a complete and separable metric space.\newline

\noindent Here, we will be mainly dealing with the $max$-norm defined for $x=(x_{1},...,x_{k})\in \mathbb{R}^{k}$ by
\begin{equation*}
\left\Vert x\right\Vert =\max_{1\leq i\leq k}\left\vert x_{i}\right\vert .
\end{equation*}

\noindent The open balls $B(x,r)$ and the closed balls $B^{f}(x,r)$ with respect to this norm are 
\begin{equation*}
B(x,r)=\{x\in \mathbb{R}^{k},\left\Vert x\right\Vert
<r\}=\prod\limits_{i=1}^{k}]x_{i}-r,x_{i}+r[
\end{equation*}

\noindent for $x=(x_{1},...,x_{k})$ and $r>0$, and  
\begin{equation*}
B^{f}(x,r)=\{x\in \mathbb{R}^{k},\left\Vert x\right\Vert \leq
r\}=\prod\limits_{i=1}^{k}[x_{i}-r,x_{i}+r]
\end{equation*}

\noindent for $r\geq 0$.\

\bigskip \noindent Before we begin, let us make some notation.\newline

\noindent Let $a=(a_{1},...,a_{k})$ and $b=(b_{1},...,b_{k})$. We define the following order relations :

\begin{equation*}
(a\leq b)\Longleftrightarrow (\forall (1\leq i\leq k),\text{ }a_{i}\leq
b_{i}),
\end{equation*}

\noindent and, 
\begin{equation*}
(a < b)\Longleftrightarrow (\forall (1\leq i\leq k),\text{ }a_{i}\leq
b_{i},\text{ }\exists (1\leq i_{0}\leq k),\text{ }a_{i_{0}}<b_{i_{0}})
\end{equation*}

\noindent and finally, 
\begin{equation*}
(a \prec b)\Longleftrightarrow (\forall (1\leq i\leq k),\text{ }a_{i}<b_{i},)
\end{equation*}

\noindent with its symmetrical counterpart,
 
\begin{equation*}
(a \succ b)\Longleftrightarrow (\forall (1\leq i\leq k),\text{ }a_{i}>b_{i},)
\end{equation*}

\bigskip \noindent Also, let us define the following classes of compact sets.\newline

\bigskip \noindent For $A=(A_{1},...,A_{k}) \prec V=(B_{1},...,B_{k})$, denote
\begin{equation*}
K_{A,B}=\prod\limits_{i=1}^{k}[A_{i},B_{i}].
\end{equation*}

\noindent For $A=(A_{1},...,A_{k}) \succ 0$, set
\begin{equation*}
K_{A}=\prod\limits_{i=1}^{k}[-A_{i},A_{i}].
\end{equation*}

\noindent For $M\in \mathbb{R}$, $M>0$, put
\begin{equation*}
K_{c,M}=[-M,M]^{k}.
\end{equation*}

\noindent The sets $K_{A,B},$ $K_{A}$ and $K_{c,M}$, are compact and will be used to characterize the tightness of sequences. The next proposition paves the way for the statements of different and equivalent conditions for tightness.\\

\begin{proposition}
\label{tensTensprop1} Let $\{\mathbb{P}_{n},n\geq 1\}$ be a sequence of probability measures on $(\mathbb{R}^{k},\mathbb{B}(\mathbb{R}^{k}))$. The following propositions are equivalent.\newline

\noindent \textbf{(1a)} For any $\varepsilon >0$, there exists a compact set $K$ in $\mathbb{R}^{k}$ such that  
\begin{equation*}
\inf_{n\geq 1}\mathbb{P}_{n}(K)\geq 1-\varepsilon .
\end{equation*}

\noindent \textbf{(2a)} For any $\varepsilon >0$, there exists a real number $M>0$ such that
\begin{equation*}
\inf_{n\geq 1}\mathbb{P}_{n}(K_{c,M})\geq 1-\varepsilon .
\end{equation*}

\noindent \textbf{(3a)} For any $\varepsilon >0$, there exists a vector $A=(A_{1},...,A_{k}) \succ 0$ of $\mathbb{R}^{k}$ such that
\begin{equation*}
\inf_{n\geq 1}\mathbb{P}_{n}(K_{A})\geq 1-\varepsilon .
\end{equation*}

\noindent \textbf{(4a)} For any $\varepsilon >0$, there exist two vectors $A=(A_{1},...,A_{k}) \prec B=(B_{1},...,B_{k})$ of $\mathbb{R}^{k}$ such that 
\begin{equation*}
\inf_{n\geq 1}\mathbb{P}_{n}(K_{A,B})\geq 1-\varepsilon .
\end{equation*}

\noindent \textbf{(1b)} For any $\varepsilon >0$, there exists a compact set $K$ of $\mathbb{R}^{k} $ such that 
\begin{equation*}
\liminf_{n\rightarrow \infty }\mathbb{P}_{n}(K)\geq 1-\varepsilon .
\end{equation*}

\noindent \textbf{(2b)} For any $\varepsilon >0$, there exists a real number $M>0$ such that 
\begin{equation*}
\liminf_{n\rightarrow \infty }\mathbb{P}_{n}(K_{c,M})\geq 1-\varepsilon .
\end{equation*}

\noindent \textbf{(3b)} Pour tout $\varepsilon >0$, there exists a vector $A=(A_{1},...,A_{k}) \succ 0$ of \ $\mathbb{R}^{k}$ such that
\begin{equation*}
\lim \inf_{n\rightarrow \infty }\mathbb{P}_{n}(K_{A})\geq 1-\varepsilon .
\end{equation*}

\noindent \textbf{(4b)} For any $\varepsilon>0$, there exist two vectors $A=(A_{1},...,A_{k}) \prec B=(B_{1},...,B_{k})$ of $\mathbb{R}^{k}$ such that 
\begin{equation*}
\lim \inf_{n\rightarrow \infty }\mathbb{P}_{n}(K_{A,B})\geq 1-\varepsilon .
\end{equation*}
\end{proposition}

\noindent \textbf{Proof}. We have two groups of formulas : $(1a)-(4a)$ and $(1b)-(4b)$. In fact, we are going to prove that the different points
of each group are equivalent and next, that the two first points of the two groups also are.\newline

\noindent \textbf{Equivalence between the points of the first group  (1a)-(4b)}: Let $\varepsilon >0$ be fixed. Let us show :\\

\noindent $(1a)\Longrightarrow (2a)$. Let $K$ be a compact set such that $\sup_{n\geq 1}\mathbb{P}
_{n}(K)\geq 1-\varepsilon$. Since $K$ is compact, it is bounded. Then, it is included in a set of the form $\{x,\left\Vert x\right\Vert \leq M\}=K_{c,M}$ and then 
\begin{equation*}
\inf_{n\geq 1}\mathbb{P}_{n}(K_{c,M})\geq \inf_{n\geq 1}\mathbb{P}%
_{n}(K)\geq 1-\varepsilon .
\end{equation*}

\noindent $(2a)\Longrightarrow (3a)$. This is obvious since $K_{c,M} $ is  equal to $K_{A}$ with $A=(M,M,...,M).$\newline

\noindent $(3a)\Longrightarrow (4a)$. This is also obvious since a set of the form $K_{A}$, for $A=(A_{1},...,A_{k}) \succ 0$, is exactly  $K_{-A,A}$.\newline

\noindent $(4a)\Longrightarrow (1a)$. This is also obvious since $K_{A,B}$ is a compact set of $\mathbb{R}^k$.\newline

\noindent \textbf{Equivalence between the points of the group (1b)-(4b)}. The proof is exactly the same as for the first group.\newline

\noindent \textbf{Equivalence between the two groups}. It will be enough to prove that : $(1a)\Longleftrightarrow (1b)$.\newline

\noindent \textbf{If (1a) holds}, then for any $\varepsilon>0$, there exists a compact set $K$ of $\mathbb{R}^{k}$ such that, for any 
$n\geq 1$, 
\begin{equation*}
\mathbb{P}_{n}(K)\geq 1-\varepsilon .
\end{equation*}

\noindent Therefore, we have
\begin{equation*}
\lim \inf_{n\rightarrow \infty}\mathbb{P}_{n}(K)\geq 1-\varepsilon.
\end{equation*}

\noindent This leads to $(1b)$.\newline

\noindent \textbf{If (1b) holds}, then for any $\varepsilon >0$, there exists a compact set $K$ such that  
\begin{equation*}
\left\{ \sup_{n\geq 1}\inf_{p\geq n}\mathbb{P}_{p}(K)\right\} \geq
1-\varepsilon /2.
\end{equation*}

\noindent Then, there exists $N\geq 1$, such that  
\begin{equation*}
\inf_{p\geq N+1}\mathbb{P}_{p}(K)\geq 1-\varepsilon,
\end{equation*}

\noindent that is for any $n>N,$%
\begin{equation*}
\mathbb{P}_{n}(K)\geq 1-\varepsilon .
\end{equation*}

\noindent Since $K$ is a compact set, it is in a set of the form $K_{c,M_{\infty}}$, where $M_{\infty}>0$, and then, for any $n>N$, 
\begin{equation*}
\mathbb{P}_{n}(K_{c,M_{\infty }})\geq 1-\varepsilon .
\end{equation*}

\noindent Now, for each fixed $j$, $1\leq j\leq N$, the set $(\left\Vert x\right\Vert \leq M)=K_{c,M}$ increases with $M$ to $\mathbb{R}^{k}$ and then, $\mathbb{P}(\left\Vert X_{j}\right\Vert \leq M)\uparrow 1$. Thus, for any $1\leq j\leq N$, there exists a real number $M_{j}>0$, 
\begin{equation*}
\mathbb{P}_{j}(K_{c,M_{j}})\geq 1-\varepsilon .
\end{equation*}

\noindent By passing, we just demonstrated that each probability measure $\mathbb{P}^{(0)}$
on $(\mathbb{R}^{k},\mathbb{B}(\mathbb{R}^{k}))$ is \textbf{tight}, that is for any $\varepsilon >0$, there exists a compact set  $K^{(0)}=K_{c,M^{(0)}}$ in $\mathbb{R}^{k}$ such that 
\begin{equation}
\mathbb{P}^{0}(K^{(0)})\geq 1-\varepsilon .  \label{tensTindiv}
\end{equation}

\noindent Coming back to our proof, we may take 
\begin{equation*}
M=\max (M_{1},...,M_{N},M_{\infty }),
\end{equation*}

\noindent and see that the sets $K_{c,M_{j}},$ $1\leq j\leq M$ and $K_{c,M_{\infty}}$ are all in $K_{c,M}$ and then for $n\geq 1,$%
\begin{equation*}
\mathbb{P}_{n}(K_{c,M})\geq 1-\varepsilon 
\end{equation*}

\noindent and thus
\begin{equation*}
\inf_{n\geq 1}\mathbb{P}_{n}(K_{c,M})\geq 1-\varepsilon,
\end{equation*}

\noindent which is $(1a)$, since $K_{c,M}$ is a compact set.\newline

\bigskip \noindent In a new step, we provide a link between the formulas  and distribution functions. For a reminder, recall that the probability distribution function associated with a probability measure $\mathbb{P}$ is defined by
\begin{equation*}
F_{\mathbb{P}}(x)=\mathbb{P}(]-\infty ,x])=\mathbb{P}\left(\prod\limits_{i=1}^{k}]-\infty
,x_{i}]\right), \ x=(x_{1},...,x_{k})\in \mathbb{R}^{k}.
\end{equation*}

\noindent This probability distribution function, in turn, determines the probability measure $\mathbb{P}$ as the Lebesgues-Stieljes probability measure defined by : for any  $(a,b)\in \mathbb{R}^{k}\times \mathbb{R}^{k},$ $a\leq b$,

\begin{equation*}
\mathbb{P}(]a,b])=\Delta _{a,b}F=\sum\limits_{\varepsilon \in
\{0,1\}^{k}}(-1)^{s(\varepsilon )}F(b+\varepsilon \ast (a-b))\geq 0,
\end{equation*}

\noindent where  for $\varepsilon =(\varepsilon _{1},...,\varepsilon _{k})\in \{0,1\}^{k},$ $%
s(\varepsilon )=\varepsilon _{1}+...+\varepsilon _{k},$ pour $x=(x_{1},...,x_{k})\in 
\mathbb{R}^{k}$, $y=(y_{1},...,y_{k}),$ $x\ast y=(x_{1}y_{1},...,x_{k}y_{k})$.\\
 
\bigskip  \noindent We are going to use the Lebesque-Stieljes probability measures to deal with uniform tightness. The reader is directed to  \cite{bmtp} or specially to the Chapter 1 of \cite{ips}.\\

\noindent For now, we need this notation. Denote for $M>0$.
\begin{equation*}
L_{M}=\{x,\exists (1\leq i\leq k),x_{i}\leq -c\}
\end{equation*}

\bigskip \noindent We have the following proposition.

\begin{proposition} \label{tensTensprop2} Let $\{\mathbb{P}_{n},n\geq 1\}$ be e sequence of probability measures on $(\mathbb{R}^{k},\mathbb{B}(\mathbb{R}^{k}))$ and consider the sequence of their probability distribution functions $\{F_{n}\geq 1\}$ with $F_{\mathbb{P}_{n}}=F_{n}$ for $n\geq 1$. Then the three following points are equivalent.\newline

\noindent (1c) For any $\varepsilon>0$, there exist a vector $0<C\in \mathbb{R}^{k}$
and a real number $c>0$ such that
\begin{equation*}
\inf_{n\geq 1}F_{n}(C)\geq 1-\varepsilon 
\end{equation*}

\noindent and
\begin{equation*}
\inf_{n\geq 0}\mathbb{P}_{n}(L_{c})\leq \varepsilon .
\end{equation*}

\noindent (2c) For any $\varepsilon >0$, there exists $0<c$ such that for $c^{(k)}=(c,...c)$, there exists $M>0$ such that 
\begin{equation*}
\inf_{n\geq 1}F_{n}(c^{(k)})\geq 1-\varepsilon 
\end{equation*}%

\noindent and
\begin{equation*}
\sup_{n\geq 0}\mathbb{P}_{n}(L_{M})\leq \varepsilon .
\end{equation*}

\noindent (3c) For any $\varepsilon>0$, there exists $M>0$ such that
\begin{equation*}
\inf_{n\geq 1}P_{n}(K_{c,M})\geq 1-\varepsilon .
\end{equation*}

\noindent Since Point $(3c)$ is also Point $(2c)$ of Proposition \ref{tensTensprop1}, then Points (3a) and (3b) are equivalent to all points of that proposition.
\end{proposition}

\noindent \textbf{Proof}. Let us proceed to the proofs of the different equivalence assertions.\\

\noindent \textbf{(a)} $\mathbf{(1c)\Longrightarrow (2c)}$. For any $\varepsilon >0$, there exists $0<C\in \mathbb{R}^{k}$ such that%
\begin{equation*}
\inf_{n\geq 1}F_{n}(C)\geq 1-\varepsilon .
\end{equation*}

\noindent Set $c=max\{C_{i},1\leq i\leq k\}$. We have $]-\infty ,C]\subset ]-\infty ,c^{(k)}]$ and  $F_{n}(c^{(k)})\geq F_{n}(C),$ 
\begin{equation*}
\inf_{n\geq 1}F_{n}(c^{(k)})\geq 1-\varepsilon .
\end{equation*}

\noindent This finishes the proof of this step (a), since the second formula implies all the others.\\

\noindent \textbf{(b)}$\mathbf{(2c)\Longrightarrow (3c)}$. From $(2c)$, we find a vector $d^{(k)}=(d,...,d)$, with $d>0$, such that 
\begin{equation*}
\inf_{n\geq 1}F_{n}(d^{(k)})\geq 1-\varepsilon/2
\end{equation*}

\noindent and real number $e>0$ such that 
\begin{equation*}
\sup_{n\geq 1}\mathbb{P}_{n}(L_{e})\leq \varepsilon/2.
\end{equation*}

\noindent By putting $M=\max (d,e)$, we get 
\begin{equation*}
\inf_{n\geq 1}F_{n}(M^{(k)})\geq 1-\varepsilon /2.
\end{equation*}%

\noindent and next
\begin{equation*}
\sup_{n\geq 1}\mathbb{P}_{n}(L_{M})\leq \varepsilon /2.
\end{equation*}

\noindent Now, let us split $\mathbb{R}^{k}$ as $\mathbb{R}^{k}=L_{M}+L_{M}^{c}$,  with
\begin{equation*}
L_{M}^{c}=\{x,\forall (1\leq i\leq k),x_{i}\geq -M\},
\end{equation*}

\noindent which itself may be decomposed as 

\begin{eqnarray*}
L_{M}^{c} &=&\{x,\forall (1\leq i\leq k),-M\leq x_{i}\leq M\}\\
&+&\{x,\forall (1\leq i\leq k),x_{i}\geq -M\text{ et }\exists (1\leq i\leq k),x_{i}>M\} \\
&=&K_{c,M}+B,
\end{eqnarray*}

\noindent where, obviously, 
\begin{equation*}
B\subset ]-\infty ,M^{(k)}]^{c}.
\end{equation*}

\noindent Therefore, we infer from $\mathbb{R}^{k}=L_{M}+K_{c,M}+B$ that 
\begin{equation}
K_{c,M}^{c}=L_{M}+B.  \label{decomKcM}.
\end{equation}

\noindent Thus, for any $n\geq 1,$%
\begin{equation*}
\mathbb{P}_{n}(K_{c,M})=\mathbb{P}_{n}(L_{M})+\mathbb{P}_{n}(B)\leq
\varepsilon /2+\varepsilon /2=\varepsilon ,
\end{equation*}

\noindent since $B\subset ]-\infty ,M^{(k)}]^{c}$. Hence for any $n\geq 1,$ 
\begin{eqnarray*}
\mathbb{P}_{n}(B) &\leq &\mathbb{P}_{n}(]-\infty ,M^{(k)}]^{c}) \\
&\leq &1-\mathbb{P}_{n}(]-\infty ,M^{(k)}]) \\
&\leq &1-F_{n}(M^{(k)})\leq \varepsilon /2.
\end{eqnarray*}

\noindent This ends the proof of this step (b).\\

\noindent \textbf{(c)}$\mathbf{(3c)\Longrightarrow (1c)}$. Suppose that $(3c)$ holds : 
for any $\varepsilon >0$, there exists $M>0$ such that 
\begin{equation*}
\inf_{n\geq 1}\mathbb{P}_{n}(K_{c,M})\geq 1-\varepsilon .
\end{equation*}

\noindent Then we have  
\begin{equation*}
\inf_{n\geq 1}F_{n}(M^{(k)})=\inf_{n\geq 1}\mathbb{P}_{n}(]-\infty
,M^{(k)}])\geq \inf_{n\geq 1}\mathbb{P}_{n}(K_{c,M})\geq 1-\varepsilon .
\end{equation*}

\noindent Next, because of (\ref{decomKcM}), we get
\begin{equation*}
\mathbb{P}_{n}(L_{M})\leq \mathbb{P}_{n}(K_{c,M}^{c})\leq \varepsilon .
\end{equation*}

\noindent Then $(1c)$ holds. Proposition is entirely proved.\newline

\bigskip  \noindent We move to the study of the tightness concept.\newline

\section{Tightness}

\subsection{Simple tightness}  $ $\\

\noindent In our particular case, each probability measure on $(\mathbb{R}^{k},\mathbb{B}(\mathbb{R}^{k}))$ is tight in the following meaning.\\

\begin{definition} \label{sompTensDef} A probability measure $\mathbb{P}$ on a metric space $(S, d)$ is tight if and only if for any $\varepsilon>0$, there exists a compact set in $S$ such that 
\begin{equation*}
\mathbb{P}(K)\geq 1-\varepsilon .
\end{equation*}
\end{definition}

\noindent We get the following proposition from Formula (\ref{tensTindiv}) above.
 
\begin{proposition} \label{tensTensprop3} A probability measure $\mathbb{P}$ on $(\mathbb{R%
}^{k},\mathbb{B}(\mathbb{R}^{k}))$ is tight.
\end{proposition}

\noindent This result is extensible to complete and separable metric spaces, more generally to totally bounded metric spaces.\\

\subsection{Asymptotic tightness. Uniform tightness}  $ $\\

\noindent Let us begin by the following definitions.

\begin{definition}
\label{tensTensDef1}

\noindent (a) A sequence of probability measures $\{\mathbb{P}_{n},n\geq 1\}$\
on $(\mathbb{R}^{k},\mathcal{B}(\mathbb{R}^{k}))$ is asymptotically tight or is uniformly tight if and only if for $\varepsilon>0$, there exists a compact set $K$ in $\mathbb{R}^{k}$ such that  
\begin{equation}
\inf_{n\geq 1}\mathbb{P}_{n}(K)\geq 1-\varepsilon   \label{tensTensUnif}
\end{equation}

\noindent or, equivalently, 
\begin{equation}
\lim \inf_{n\rightarrow \infty }\mathbb{P}_{n}(K)\geq 1-\varepsilon .
\label{tensTensAsymp}
\end{equation}

\noindent (b) A sequence of random vectors  $X_{n}:(\Omega_{n },\mathcal{A}_{n },\mathbb{P}_{n})\mapsto (\mathbb{R}^k,\mathcal{B}(\mathbb{R}^k))$, $n\geq 1$, is asymptotically tight or is uniformly tight if and only if  the sequence of the probability laws $\{\mathbb{P}_{X_{n}},n\geq 1\}$ is asymptotically tight or is uniformly tight, that is for any $\varepsilon$, there exists a $K$ in $\mathbb{R}^{k}$ such that  
\begin{equation*}
\inf_{n\geq 1}\mathbb{P}_{n}(X_{n}\in K)\geq 1-\varepsilon 
\end{equation*}

\noindent or equivalently, 
\begin{equation*}
\lim \inf_{n\rightarrow \infty }\mathbb{P}_{n}(X_{n}\in K)\geq 1-\varepsilon .
\end{equation*}

\noindent (c) A sequence of probability distribution functions $\{F_{n},n\geq 1\}$ on  $(\mathbb{R}^{k},\mathcal{B}(\mathbb{R}^{k}))$ 
is asymptotically tight or is uniformly tight if and only the sequence of their Lebesgue-Stieljes probability measures $\{\mathbb{P}_{n},n\geq 1\}$ is asymptotically tight or is uniformly tight, or equivalently for any $\varepsilon>0$, there exist $0<C\in \mathbb{R}^{k}$ and a real number $c>0$ such that 
\begin{equation*}
\inf_{n\geq 1}F_{n}(C)\geq 1-\varepsilon 
\end{equation*}

\noindent and
\begin{equation*}
\sup_{n\geq 0}\mathbb{P}_{n}(L_{c})\leq \varepsilon,
\end{equation*}

\noindent that is, if and only if, there exists for any $\varepsilon >0$, a real number $c>0$ such that we have for for $c^{(k)}=(c,...c)$ 
\begin{equation*}
\inf_{n\geq 1}F_{n}(c^{(k)})\geq 1-\varepsilon .
\end{equation*}
\end{definition}

\bigskip \noindent In $\mathbb{R}^{k}$, uniform tightness (\ref{tensTensUnif}) is equivalent to asymptotic tightness because of 
Proposition \ref{tensTensprop1}. Thus from now, we speak only about tightness of sequences of probability measures, or of random vectors, or of probability distribution functions.\newline

\noindent Before, we come to the Helly-Bray theorem, we are going to give three important properties of tightness.\\

\subsection{Tightness and continuous mapping}  $ $\\

\noindent The tightness is preserved by continuous mapping in the following sense.

\begin{proposition} \label{tensTensprop5} Let $X_{n}:(\Omega _{n },\mathcal{A}_{n},\mathbb{P}_{n })\mapsto (\mathbb{R}^k,\mathcal{B}(\mathbb{R}^k))$, $n\geq 1$, be a \textit{tight sequence} of random vectors and let $g:\mathbb{R}^{k}\longmapsto \mathbb{R}^{m},$ $m\geq 1$, be a continuous mapping. Then the sequence $\{g(X_{n}),n\geq 1\}$ is tight.
\end{proposition}

\bigskip
\noindent \textbf{Proof}. Let $\{X_{n},n\geq 1\}$ be tight and $g:\mathbb{R}^{k}\longmapsto \mathbb{R}^{m}$ continuous. For any  $\varepsilon>0$, there exists a compact set $K$ in $\mathbb{R}^{k}$ such that  
\begin{equation}
\inf_{n\geq 1}\mathbb{P}(X_{n}\in K)\geq 1-\varepsilon .
\end{equation}

\noindent But $(X_{n}\in K)\subset (g(X_{n})\in g(K))$ where 
\begin{equation*}
K_{0}=g(K)=\{g(x),x\in K\}
\end{equation*}

\noindent is the direct image of $K$ by $g$, and is a compact set. Indeed, let $\{g(x_{n}),x_{n}\in K,n\geq 1\}$ be a sequence in 
$K_{0}$. Since $K$ is a compact set, the sequence $(x_{n})_{n\geq 1}$, which is in $K$, has a sub-sequence $x_{n(k)}\rightarrow x\in K$ converging, as $k \rightarrow +\infty$, to a point $x$ which is in $K$ since $K$ is closed. Since $g$ is continuous, then  $g(x_{n(k)})_{k\geq 0}$ converges to $g(x) \in K_{0}$ as $k \rightarrow +\infty$. It follows that $K_{0}$ is a compact set in $\mathbb{R}^{m}$ and 

\begin{equation*}
\inf_{n\geq 1}\mathbb{P}(g(X_{n})\in K_{0})\geq \inf_{n\geq 1}\mathbb{P}%
(X_{n}\in K)\geq 1-\varepsilon .
\end{equation*}

\noindent This ends the proof.\newline

\subsection{Characterization of the tightness by that of the components} $ $\\

\noindent In the particular case of $\mathbb{R}^k$, we have 

\begin{proposition}
\label{tensTensprop6} A sequence of random vectors $X_{n}:(\Omega _{n },\mathcal{A}_{n },
\mathbb{P}_{n })\mapsto (\mathbb{R}^k,\mathcal{B}(\mathbb{R}^k))$, $n\geq 1$, is tight if and only if each sequence of components, $\{X_{n}^{(i)},n\geq 1\}$, $1\leq i\leq k$, is tight.
\end{proposition}

\bigskip
\noindent \textbf{Proof}. Let $X_{n}:(\Omega _{n },\mathcal{A}_{n },
\mathbb{P}_{n })\mapsto (\mathbb{R}^k,\mathcal{B}(\mathbb{R}^k))$, $n\geq 1$, be a sequence of random vectors.\newline

\noindent Suppose that this sequence is tight. By Proposition \ref{tensTensprop5}, each sequence of components $\{X_{n}^{(i)},n\geq 1\}=\{\pi
_{i}(X_{n}),n\geq 1\}$, $1\leq i\leq k$, is tight, as continuous transformations of a tight sequence, that is as the $i$-th projection $\pi _{i}$ of a tight sequence.\\

\noindent Suppose that for each $1\leq i\leq k$, $\{X_{n}^{(i)},n\geq 1\}$ is tight. Then for  $1\leq i\leq k$, for any $\varepsilon>0$, there exists a real number $A_{i}>0$ such%
\begin{equation*}
\inf_{n\geq 1}\mathbb{P}(X_{n}^{(i)}\in \lbrack -A_{i},A])\geq 1-\varepsilon
/k.
\end{equation*}

\noindent By setting $A=(A_{1},...,A_{k})$, we have for $A>0$
\begin{equation*}
\bigcap\limits_{i=1}^{k}\left( X_{n}^{(i)}\in \lbrack -A_{i},A]\right)
=\left( X_{n}\in \prod\limits_{i=1}^{k}[-A_{i},A_{i}]\right) ,
\end{equation*}

\noindent It follows that for any $n\geq 1,$%
\begin{eqnarray*}
\mathbb{P}\left( X_{n}\notin \prod\limits_{i=1}^{k}[-A_{i},A_{i}]\right) &=&%
\mathbb{P}\left( \bigcup\limits_{i=1}^{k}\left( X_{n}^{(i)}\notin \lbrack
-A_{i},A]\right) \right) \\
&\leq &\sum\limits_{i=1}^{k}P\left( X_{n}^{(i)}\notin \lbrack
-A_{i},A]\right) \leq \varepsilon,
\end{eqnarray*}

\noindent and then for any $n\geq 1,$%
\begin{equation*}
\mathbb{P}\left( X_{n}\in K_{A}\right) \geq 1-\varepsilon .
\end{equation*}

\noindent Hence, the sequence $\{X_{n},n\geq 1\}$ is tight. The proof is complete.

\subsection{Tightness of a weakly convergent sequence}  $ $\\

\noindent We have the following result. 

\begin{proposition} \label{tensTensprop7} Any sequence $X_{n}:(\Omega _{n },\mathcal{A}_{n }, \mathbb{P}_{n })\mapsto (\mathbb{R}^k,\mathcal{B}(\mathbb{R}^k))$, $n\geq 1$, of random vectors that weakly converges is tight.
\end{proposition}

\noindent \textbf{Proof}. Suppose that $X_{n}$ weakly converges to the probability $\mathbb{P}$. This probability is tight. So for any 
$\varepsilon>0$, there exists a compact $K_{A}=[-A,A]$ of $K$ such that
\begin{equation*}
\mathbb{P}(K)\geq 1-\varepsilon .
\end{equation*}

\noindent Let $0<\delta <1$ and set $A+\delta =(A_{1}+\delta
,...,A_{k}+\delta )$. We have for any $0<\delta <1$,  
\begin{equation*}
\overset{o}{K}_{A+\delta }=\prod\limits_{i=1}^{k}]-A_{i}-\delta ,A_{i}+\delta \lbrack .
\end{equation*}

\noindent Since $X_{n}\rightsquigarrow X$ and $\overset{o}{K}_{A+\delta }$ is open, we use Point $(ii)$ of Portmanteau Theorem \ref{cv.theo.portmanteau} to show that
\begin{equation*}
\liminf_{n\rightarrow \infty }\mathbb{P}(X_{n}\in K_{A+1})\geq \lim
\inf_{n\rightarrow \infty }P(X_{n}\in \overset{o}{K}_{A+\delta })\geq 
\mathbb{P}(\overset{o}{K}_{A+\delta }),
\end{equation*}

\noindent for any $0<\delta <1$ and next
\begin{equation*}
\liminf_{n\rightarrow \infty }\mathbb{P}(X_{n}\in K_{A+1})\geq \mathbb{P}(\overset{o}{K}_{A+\delta }),
\end{equation*}

\noindent for any $0<\delta <1,$. By letting $\delta \downarrow 0$, we have $\overset{o}{K}_{A+\delta }\downarrow 
\overline{K}=K$ since $K$ is a closed set. Applying this in the last formula gives
\begin{equation*}
\liminf_{n\rightarrow \infty }\mathbb{P}(X_{n}\in K_{A+1})\geq \mathbb{P}(K)\geq 1-\varepsilon .
\end{equation*}

\noindent It follows that the sequence $\{X_{n},n\geq 1\}$ is tight.\\

\bigskip \noindent \textbf{Remark}. Actually, we may see that the sequence has inherited the tightness of the weak limit. This result still holds for complete and separable spaces where any probability measure is tight.\newline

\bigskip \noindent In the new section, we are going to deal with the fundamental theorem of tightness.\\

\section{Compactness Theorem for weak convergence in $\mathbb{R}^{k}$}

\noindent This theorem is a kind of inverse of Proposition \ref{tensTensprop7}, concerning the convergence of sub-sequence.\\

\begin{theorem} \label{tensTheo2} (Prohorov - Helly-Bray) Let $X_{n}:(\Omega _{n },\mathcal{A}_{n }, \mathbb{P}_{n })\mapsto (\mathbb{R}^k,\mathcal{B}(\mathbb{R}^k))$, $n\geq 1$, be a tight sequence of random vectors. Then it contains a weakly converging sub-sequence.
\end{theorem}

\bigskip \noindent This theorem may be directly proved, as done in \cite{billingsley} and van
der Vaart and Wellner \cite{vaart}. The proof in Billinsgley is very lengthy. That of van
der Vaart and Wellner is very much simpler and more general. But in this context, we are going to use the Helly-Bray
approach as in van der vaart \cite{vaart_asymp} and Lo\`{e}ve \cite{loeve}. Here, we give a more detailed proof.\newline

\noindent Here, the proof of Theorem \ref{tensTheo2} is based on the following Helly-Bray Theorem in which the hard work is done.\newline

\begin{theorem} \label{tensTheo1} (Helly-Bray) Any sequence $\{F_{n},n\geq 1\}$\ of probability distribution function on $(\mathbb{R}^{k},\mathbb{B}(\mathbb{R}^{k}))$ has a sub-sequence $\{F_{n(k)},k\geq 1\}$ weakly converging to a distribution function $F$, which is not necessarily a probability distribution function.
\end{theorem}

\bigskip \noindent \textbf{Proof}. Let $\{F_{n},n\geq 1\}$ be a sequence of probability distribution functions $(\mathbb{R}^{k},\mathbb{B}(%
\mathbb{R}^{k}))$. Let $\mathbb{Q}^{k}$ be the set of all elements of $\mathbb{R}^{k}$ with rational components. $\mathbb{Q}^{k}$ is everywhere dense in $\mathbb{R}^{k}$. Let us enumerate $\mathbb{Q}^{k}$ as  $\mathbb{Q}^{k}=\{q_{1},q_{2},...\}$ and proceed by steps.\\

\noindent \textbf{Step  1}. We are going to find a sub-sequence $(F_{n(j)})_{j\geq 1}$ of $(F_{n})_{n\geq 1}$ point-wisely converging to some function $G$ on $\mathbb{Q}^{k}$ by using the diagonal sequence method. We have that :  $(F_{n}(q_{1}))_{n\geq 1}$ $\subset \lbrack 0,1]$. Using Bolzano-Weierstrass property on $\mathbb{R}$, we get a sub-sequence $(F_{1,n}(q_{1}))_{n\geq 1}$ of $\{F_{n}-q_1),n\geq 1\}$ converging to $G(q_{1}).$\newline

\noindent Next, we apply the sub-sequence $(F_{1,n})_{n\geq 1}$ to  $q_{2}$ in this way : $(F_{1,n}(q_{2}))_{n\geq 1}$ $\subset \lbrack 0,1]$. 
We find a sub-sequence $(F_{2,n}(q_{2}))_{n\geq 1}$ of $(F_{1,n}(q_{2}))_{n\geq 1}$ that converges to a real number $G(q_{2})$. We proceed so forth and get sub-sequences $(F_{j,n})_{n\geq 1},$ $j=1,2,...$ satisfying : \newline

\noindent (a) For each $j\geq 1,$ $(F_{j+1,n})_{n\geq 1}$ is a sub-sequence of any of the sub-sequences $(F_{i,n})_{n\geq 1}1\leq i\leq j$.\newline

\noindent \bigskip (b) For any $j\geq 1$, for any $1\leq j\leq i,$ $%
F_{j,n}(q_{i})\rightarrow G(q_{i})$.

\noindent Next, we take the diagonal sequence $(F_{j,j})_{j\geq 1}$. We may use a simple graph, as below, to see this : for any fixed $i\geq 1$, the sequence $\{F_{j,j},j\geq i\}$ is a sub-sequence of $(F_{i,n})_{n\geq i}$ and then
\begin{equation*}
F_{j,j}(q_{i})\rightarrow G(q_{i}).
\end{equation*}

\noindent To read this graph, one has to notice that the sequence in one line is a sub-sequence of those in the previous lines. From this, it becomes clear that $\mathbf{F}_{j,j}$ is an element of all the lines from $1$ to $j$.
 
\begin{equation*}
\begin{tabular}{lllllllllll}
$\mathbf{F}_{1,1}$ & $F_{1,2}$ & $F_{1,3}$ & $F_{1,4}$ & $F_{1,5}$ & $%
F_{1,6} $ & $F_{1,7}$ & $F_{1,8}$ & $F_{1,9}$ & $F_{1,10}$ & ... \\ 
& $\mathbf{F}_{2,2}$ & $F_{2,3}$ & $F_{2,4}$ & $F_{2,5}$ & $F_{2,6}$ & $%
F_{1,7}$ & $F_{1,8}$ & $F_{1,9}$ & $F_{1,10}$ & ... \\ 
&  & $\mathbf{F}_{3,3}$ & $F_{3,4}$ & $F_{3,5}$ & $F_{3,6}$ & $F_{3,7}$ & $%
F_{3,8}$ & $F_{3,9}$ & $F_{3,10}$ & ... \\ 
&  &  & $\mathbf{F}_{4,4}$ & $F_{4,5}$ & $F_{4,6}$ & $F_{4,7}$ & $F_{4,8}$ & 
$F_{4,9}$ & $F_{4,10}$ & ... \\ 
&  &  &  & ... & ... & ... & ... & ... & .. & ... \\ 
&  &  &  &  & $\mathbf{F}_{j,j}$ & $F_{j,j+1}$ & $F_{j,j+2}$ & $F_{j,j+3}$ & 
$F_{j,j+5}$ & ...%
\end{tabular}%
\end{equation*}

\noindent We conclude that the diagonal sub-sequence $\left( F_{j,j}\right)
_{j\geq 1}$, written as $(F_{n(j)})_{j\geq 1}$, satisfies 
\begin{equation*}
\forall q\in \mathbb{Q}^{k},\text{ } F_{n(j)}(q)\rightarrow G(q) \text{
as }j\rightarrow +\infty.
\end{equation*}

\bigskip \noindent \textbf{Step 2}. Properties of $G$ on $\mathbb{Q}^{k}$.\\

\noindent (2.1) For any $(a,b)\in \mathbb{Q}^{k}\times \mathbb{Q}^{k}$, as
$j\rightarrow +\infty ,$%
\begin{equation*}
0\leq \Delta _{a,b}F_{n(j)}=\sum\limits_{\epsilon \in
\{0,1\}^{k}}(-1)^{s(\epsilon )}F_{n(j)}(b+\epsilon \ast (a-b))
\end{equation*}
\begin{equation*}
\rightarrow \Delta _{a,b}G=\sum\limits_{\epsilon \in
\{0,1\}^{k}}(-1)^{s(\epsilon )}F(b+\epsilon \ast (a-b))\geq 0,
\end{equation*}

\bigskip
\noindent Since all the points $b+\epsilon \ast (a-b)$ are in $\mathbb{%
Q}^{k}$, it follows that $G$ assigns non-negative volume to cuboids of $\mathbb{Q}^{k}$.%
\newline

\noindent $G$ is non-decreasing on $\mathbb{Q}^{k}$ as inherited from the non-decreasingness of the $F_{n(j)}$, $j\geq 1$, on  $\mathbb{Q}^{k}$,.\newline

\noindent \textbf{step 3}. Define $F$ on  $\mathbb{J}^{k}=(\mathbb{R} \setminus \mathbb{Q})^{k}$ by 
\begin{equation*}
F(x)=\inf \{G(q),q\in \mathbb{Q}^{k},x \prec q\}\in \lbrack 0,1].
\end{equation*}

\noindent for  $x\in \mathbb{J}^{k}$. It is obvious that $F$ is well-defined on $\mathbb{J}^{k}$. It is also sure that $F$ is non-decreasing.\\

\noindent \textbf{(a) Let us show that $F$  is right-continuous}. Let $x\in \mathbb{J}^{k}$ and let $\varepsilon >0$. By definition of the finite infimum, there exists $q\in \mathbb{Q}^{k}$ such  that $x \prec q$ and $G(q)<F(x)+\varepsilon$. For any $y \in \mathbb{J}^{k}$, $x \prec<y<q$, we have $F(y)\leq G(q)$ and $\varepsilon >G(q)-F(x)\geq F(y)-F(x)$. Then 
\begin{equation}
\forall \varepsilon >0,\text{ }\exists q>x, \text{ }x<y<q\Longrightarrow 0\leq
F(y)-F(x)<\varepsilon .  \label{cadG}
\end{equation}

\noindent Then $F$ is right-continuous.\newline

\noindent \textbf{(c) Let us show that $F_{n(j)}(x)\rightarrow F(x)$ for continuity points of $x \in x\in \mathbb{J}^{k}$ of $G$}.\newline

\noindent Let $x$ be a continuity of $F$ on $\mathbb{J}^{k}$. For any $\varepsilon
>0$, we may find $(y^{\prime},y^{\prime \prime })$ in $\left(\mathbb{J}^{k}\right)^2$ such that $y^{\prime}<x<y^{\prime \prime}$ and 
$F(y^{\prime\prime})-F(y^{\prime })<\varepsilon /2$. Let $(q^{\prime },q^{\prime \prime}) \in \mathbb{Q}^{k}$ such that  $y^{\prime }<q^{\prime }<x<q^{\prime \prime }<y^{\prime \prime }$.  Then $ G(q^{\prime \prime })-G(q^{\prime })\leq
F(y^{\prime \prime })-F(y^{\prime })\leq \varepsilon$. Next

\begin{eqnarray*}
F(y^{\prime }) \leq G(q^{\prime })=\lim F_{n(j)}(q^{\prime }) &\leq& \liminf F_{n(j)}(x)\leq \limsup F_{n(j)}(x)\\
&\leq& \limsup F_{n(j)}(q^{\prime \prime })=\lim G(q^{\prime \prime })\leq F(y^{\prime \prime }).
\end{eqnarray*}

\bigskip \noindent Then  $\liminf_{j\rightarrow +\infty} F_{n(j)}(x)$, $\limsup_{j\rightarrow +\infty} F_{n(j)}(x)$ and $F(x)$
are in the interval $[F(y^{\prime }),F(y^{\prime \prime })]$ with length at most equal to $\varepsilon$. This implies that

\begin{equation}
\max (\left\vert F(x)-\lim \inf F_{n(j)}(x)\right\vert ,\left\vert F(x)-\lim
\sup F_{n(j)}(x)\right\vert )\leq \varepsilon ,  \label{cadGG}
\end{equation}

\bigskip \noindent for any $\varepsilon >0$. Therefore, we arrive at 
\begin{equation*}
F_{n(j)}(x)\rightarrow F(x)\text{ as }j\rightarrow +\infty .
\end{equation*}

\bigskip \noindent \textbf{(d) $F$ assigns non-negative volume to cuboids}. For any $(a,b)\in \mathbb{J}^{k}\times \mathbb{J}^{k}$, let $q^{\prime
}\downarrow a$ and $q^{\prime \prime }\downarrow b$ with  $q^{\prime } \succ a $ and $q^{\prime \prime } \succ b$ and $(q^{\prime },q^{\prime \prime })\in \mathbb{Q}^{k}\times \mathbb{Q}^{k}$. By monotone limit, and by the definition of $G$,
\begin{equation}
0\leq \Delta _{q^{\prime },q^{\prime \prime }}G\rightarrow \Delta
_{a,b}F\geq 0.  \label{vpG}
\end{equation}

\bigskip \noindent \textbf{Partial conclusion}. $F$ is a distribution function on $\mathbb{J}^{k}$ and $F_{n(j)}(x)\rightarrow F(x)$ for continuity points $x \in \mathbb{J}^{k}$ of $F$.\\ 

\bigskip 
\noindent \textbf{Step  4}. Now, we may extend $F$ on $\mathbb{J}_{c}^{k}=\mathbb{R}^{k} \setminus \mathbb{J}^{k}$ by 
\begin{equation*}
F(x)=\inf \{F(y), y\in \mathbb{J}^{k}, \ \  x \prec y\} \in \lbrack 0,1],\text{ }x\in \mathbb{J}_{c}^{k}.
\end{equation*}

\noindent First, with a very little effort, we see that $F$ is non-decreasing on $\mathbb{R}^{k}$. Next, we have to prove that $F$ is right-continuous at any point $x\in \mathbb{R}^{k}$. Let us fix $\varepsilon>0$. If $x\in \mathbb{J}^{k}$, by right-continuity of $F$ on $\mathbb{J}^{k}$, we can find $y_1 \succ x$, for any $z\in \mathbb{J}^{k}$ and $x\leq z \prec y_1$, we have : $F(x) \leq F(z) \leq F(y_1)<F(x)+\varepsilon$. This is also true for $z \in \mathbb{J}_{c}^{k}$ since, by construction, $F(z) \leq F(y_1)$.\\

\noindent If $x \in \mathbb{J}_{c}^{k}$, by definition of $F(x)$ we can find $y_1 \succ x$, such that for any $z\in \mathbb{J}^{k}$ and $x\leq z \prec y_1$, we have $G(x) \leq G(z) \leq F(y_1)<F(x)+\varepsilon$. We conclude as in the first case by using in addition the increasingness of $F$. We conclude that for any $\mathbb{R}^{k}$, for any $\varepsilon>0$, we can find $y \succ x$ in $\mathbb{J}^{k}$ such that
\begin{equation*}
(x \leq z \prec y) \text{ and } 0<F(z)-F(x) < F(y)-F(x)<\varepsilon.
\end{equation*}

\noindent Thus, $F$ is right-continuous.\\

\noindent Finally, to show that $F$ assigns non-negative volume to cuboids, we use the right continuity of $F$ and the fact that the property holds on $\mathbb{J}^{k}$.\\

\noindent Finally, we have to prove that $F_{n(j)}(x)\rightarrow F(x)$ for any continuity points of $F$. We may repeat the same
technique that led to (\ref{cadGG}). We are going to give only the beginning. \\

\noindent Let $x$ be a continuity point of $F$ on $\mathbb{R}^{k}$. For any $\varepsilon
>0$, we may find $(y^{\prime},y^{\prime \prime })$ in $\mathbb{R}^{k}$ such that $y^{\prime} \prec x \prec y^{\prime \prime}$ and 
$F(y^{\prime\prime})-F(y^{\prime })<\varepsilon /2$.\\

\noindent From there, we may find $(z^{\prime },z^{\prime \prime}) \in \mathbb{J}^{k}$ such that  $y^{\prime }<z^{\prime }<x<z^{\prime \prime }<y^{\prime \prime}$ and such that $z^{\prime}$ and $z^{\prime \prime}$ are continuity points of $F$ on $\mathbb{J}^{k}$. If we are able to do that, we may re-conduct the same lines that led to (\ref{cadGG}), by replacing $(y^{\prime },y^{\prime \prime })$ by $(z^{\prime },z^{\prime \prime})$.

\bigskip \noindent Now, we may find points of the form $z_{\varepsilon}=y^{\prime}+(\varepsilon)\delta$ in $]y^{\prime}, x[ \cap \mathbb{J}^{k}$, where $\delta=(1,...,1)$ and $0< \varepsilon < \varepsilon_0$. The boundaries of the intervals $]-\infty, z_{\varepsilon}]$ are disjoint. So, for $m_{F,\mathbb{J}}$ being the Lebesgue-Stieljes measure associated with $F$ on $\mathbb{J}^{k}$, we may have $m_{F,\mathbb{J}}(\partial ]-\infty, z_{\varepsilon}])>0$ only for - at most - a countable number of $z_{\varepsilon}$. Then we may easily pick a value of $\varepsilon$ such that $m_{F,\mathbb{J}}(\partial ]-\infty, z_{\varepsilon}])>0$, that is $z^{\prime}=z_{\varepsilon}$ is a continuity point of $F$ on $\mathbb{J}^{k}$. We find $z^{\prime \prime}$ in the same manner.\\

\noindent This completely finishes the proof.\\

\noindent \textbf{Remark} We wanted to give a complete proof with all the necessary details. Our step 4 is needless
if it is possible to prove that $G$ is right-continuous on $\mathbb{Q}^{k}$. If this is the case, one should stop at Step 3 and take $F=G$.\\

\bigskip \noindent Now let us move to the proof of Theorem \ref{tensTheo2}.\newline

\noindent \textbf{Proof of Theorem \ref{tensTheo2} of Prohorov}.
Suppose that the sequence of probability distribution functions $\{F_{n},n\geq 1\}$
is tight, that is the sequence of their Lebesgue-Stieljes measures $\{\mathbb{P}_{n}(]a,b])=\Delta _{a,b}F_{n},n\geq 1\}$ is tight. By 
Proposition \ref{tensTensprop2}, for any $\varepsilon >0$, we may find a vector $C \succ 0$, $C\in \mathbb{R}^{k}$ such that $n\geq 1,$%
\begin{equation*}
F_{n}(C)\geq 1-\varepsilon .
\end{equation*}

\noindent By Theorem \ref{tensTheo1}, there exists a sub-sequence $\left( F_{n(j)}\right) _{j\geq 1}$ of $\left( F_{n}\right)
_{n\geq 1}$ that weakly converges to a distribution function $F$ associated to a measure $L$ defined by  $L(]a,b])=\Delta _{a,b}F$
and bounded by the unity.\\

\noindent  Consider the family $\{C_{h}=C+h^{(k)}, h>0\}$. These points are such that their boundaries are $\partial ]-\infty ,F=C_{h}]$ are disjoint. So, we may choose a sequence $C_{h_{p}}$ such that $L(\partial ]-\infty,C_{h_{p}}])=0$ for any $p\geq 1$ and  $C_{h_{p}}\uparrow (+\infty )^{(k)}$ as $p\uparrow +\infty$. These points are continuity ones of $F$ and are greater than $C$. Then for any fixed $p\geq 1,$%
\begin{equation*}
F_{n(j)}(C_{h_{p}})\geq 1-\varepsilon .
\end{equation*}

\noindent By letting $j\rightarrow \infty$, we get 
\begin{equation*}
F(C_{h_{p}})=L(]-\infty ,C_{h_{p}})\geq 1-\varepsilon .
\end{equation*}

\noindent Next by letting $p\uparrow +\infty$ and next $\varepsilon\downarrow 0$, we get 
\begin{equation}
F((+\infty )^{(k)})=1.  \label{condFRPINF}
\end{equation}

\bigskip \noindent On the other hand, for any $\varepsilon >0$, there exists $M>0$, such that   
\begin{equation*}
\sup \mathbb{P}_{n}(L_{M})\leq \varepsilon .
\end{equation*}

\noindent We have to prove that
\begin{equation}
\lim_{\exists (1\leq i\leq k),x_{i}\rightarrow -\infty }F(x)=0,
\label{condFRMINF}
\end{equation}

\noindent which is equivalent to saying that for any $\varepsilon >0$, there exists $M>0$ such that 

\begin{equation*}
\exists (1\leq i\leq k),x_{i}<-M\Longrightarrow F(x)\leq \varepsilon .
\end{equation*}

\noindent But 
\begin{equation*}
\exists (1\leq i\leq k),\ \ (x_{i}<-M)\ \ \Rightarrow \ \ (]-\infty ,x]\subset L_{M}),
\end{equation*}

\noindent and then for any $n\geq 1$,%
\begin{equation*}
\exists (1\leq i\leq k), \ \ (x_{i}<-M) \ \ \Rightarrow \ \ (F_{n}(x)\leq \varepsilon).
\end{equation*}

\noindent Now, let $x$ be fixed such that :  $\exists (1\leq i\leq
k),x_{i}<-M$. Let $x(h)=x+h^{(k)}$, with $0<h<-(M+x_{i})$. By the now classical method we used just above, we may find a sequence $x(h_{p})$, $p\geq 1,$ of continuity points of $F$ with $h_{p}\downarrow 0$. Then for any fixed $p\geq 1$, for any $j\geq 1,$%
\begin{equation*}
\text{\ }F_{n(j)}(x(h_{p}))\leq \varepsilon.
\end{equation*}

\noindent By letting $j\rightarrow \infty$, we get 
\begin{equation*}
\text{\ }F(x(h_{p}))\leq \varepsilon.
\end{equation*}

\noindent Now, by right continuity, we get, as $p\uparrow+\infty$,
\begin{equation*}
\text{\ }F(x)\leq \varepsilon.
\end{equation*}

\noindent We conclude that for any $\varepsilon>0$,

\begin{equation*}
\exists (1\leq i\leq k),x_{i}<-M\Longrightarrow F(x)\leq \varepsilon .
\end{equation*}

\noindent And this proves (\ref{condFRMINF}). This combined to (\ref{condFRPINF}) shows that $F$ is a probability distribution function on $\mathbb{R}^k$.\\

\section{Applications}

\subsection{Continuity Theorem of L\'{e}vy}  $ $\\

\noindent We have this important property.

\begin{theorem} \label{cv.levy} Let $\psi _{n}$, $n\geq 1$, be a sequence of characteristic functions on $\mathbb{R}$
that converges point-wisely to a function $\psi$ which is continuous at zero. Then $\psi$ is a characteristic function.
\end{theorem}

\bigskip \noindent \textbf{Proof}. We necessarily have $\psi (0)=0$ since  $\psi _{n}(0)=0$ for all $n\geq 1$. We may suppose that the $\psi _{n}$ are the characteristic functions of random variables $X_{n}$, that is for any $,\geq 1$,
\begin{equation*}
\psi _{n}(t)=E(e^{itX_{n}}),t\in \mathbb{R}.
\end{equation*}

\noindent By Fact 1 in Section \ref{funct.facts} in Chapter \ref{funct}, we have that for  $\left\vert \sin a\right\vert \leq a$
for $\left\vert a\right\vert \geq 2$. Then 
\begin{equation*}
1_{(\left\vert \delta x\right\vert >2)}\leq 2\left(1-\frac{\sin \delta x}{\delta x%
}\right)
\end{equation*}

\noindent and by the right equality easily proved, 
\begin{equation*}
1_{(\left\vert \delta x\right\vert >2)}\leq 2\left(1-\frac{\sin \delta x}{\delta x%
}\right)=\frac{1}{\delta }\int_{-\delta }^{\delta }(1-\cos tx)dt.
\end{equation*}

\noindent Let us apply this formula to $X_{n}$ to get 
\begin{equation*}
1_{(\left\vert X_{n}\right\vert >2/\delta )}\leq \frac{1}{\delta }%
\int_{-\delta }^{\delta }(1-\cos tX_{n})dt
\end{equation*}

\noindent and by taking expectations and by applying Fubini Theorem for integrable functions, 
\begin{equation*}
P(\left\vert X_{n}\right\vert >\frac{2}{\delta })\leq \frac{1}{\delta }%
\int_{-\delta }^{\delta }R_{e}(1-Ee^{itX_{n}})dt.
\end{equation*}

\noindent By applying the Monotone Convergence Theorem, to $R_{e}(1-Ee^{itX_{n}})\rightarrow R_{e}(1-\psi (t))$, we obtain
\begin{equation}
\lim \inf_{n\rightarrow \infty }P(\left\vert X_{n}\right\vert >\frac{2}{%
\delta })\leq \frac{1}{\delta }\int_{-\delta }^{\delta }R_{e}(1-\psi (t))dt.
\label{tensionLevy1}
\end{equation}

\noindent The \textit{real part} function $R_{e}(\cdot )$ is continuous and by the assumptions, $R_{e}(1-\psi (t))\rightarrow 0$ as $\delta \rightarrow 0$. Therefore, we may let $\delta \rightarrow 0$ in (\ref{tensionLevy1}) to get 
\begin{equation*}
\lim \inf_{n\rightarrow \infty }P(\left\vert X_{n}\right\vert >\frac{2}{%
\delta })=0.
\end{equation*}

\noindent This implies that the sequence is tight. Then there exists a sub-sequence $X_{n_{k}}$ weakly converging to $X$. By Theorem , we have

\begin{equation*}
\psi _{n_{k}}(t)=E(\exp (itX_{n_{k}}))\rightarrow E(\exp (itX))=\psi _{0}(t).
\end{equation*}

\noindent By the uniqueness of limits in $\mathbb{R}$
\begin{equation*}
\psi =\psi _{0}.
\end{equation*}

\noindent We conclude that $\psi$ is a characteristic function.\\

\bigskip \noindent By applying this, we will have another proof of the characterization of weak convergence by characteristic functions.\\

\subsection{Another proof of the characterization of weak convergence by characteristic functions}  $ $\\

\noindent Here is the beautiful proof we already signaled in the remark after the statement of Proposition \ref{cv.FC} in Chapter \ref{cv}.

\begin{theorem} \label{cv.tension.ConvFC}
A sequence $X_{n}$ of random vectors with values in $\mathbb{R}^{k}$ weakly converges to the random vector  $X\in \mathbb{R}^{k}$ if and only if  $u\in \mathbb{R}^{k},$ $\mathbb{E}(\exp (i<u,X_{n}>))\rightarrow \mathbb{E}(\exp
(i<u,X>))$ \text{ as } $n\rightarrow +\infty$.
\end{theorem}

\noindent \textbf{Proof}. The direct implication comes from the application of the Dominated Convergence Theorem. Let us prove the indirect implication. Suppose that for any $u\in \mathbb{R}^{k}$ 
\begin{equation*}
\psi _{n}(u)=\mathbb{E}(\exp (i<u,X_{n}>))\rightarrow \mathbb{E}(\exp (i<u,X>))=\psi _{n}(u).\end{equation*}

\noindent For any fixed $i$, $1\leq i\leq k$, the sequence of the $i$-th components, $X_{n}^{(i)}$, satisfies, for any $t\in \mathbb{R}$, 
\begin{equation*}
\psi _{n^{(i)}}(t)=\psi _{n}(\underset{i-\text{i-th place}}{\underbrace{%
0,..,t,..0}})=\mathbb{E}(\exp (itX_{n}^{(i)}))\rightarrow \psi (\underset{i-\text{i-th place}}{\underbrace{0,..,t,..0}})=\psi ^{(i)}(t).
\end{equation*}

\noindent The function $\psi_i$ is continuous since it is the characteristic function of $X_i$. Then each $X_{n}^{(i)}$ is tight.  By Proposition \ref{tensTensprop6}, the sequence $X_{n}$ is tight.\\ 

\noindent We may conclude in two steps.\\

\noindent Step 1 : Each sub-sequence of $X_{n}$ contains a sub-sequence weakly converging to a probability $L$. By assumption, $X_{n}$ weakly converges to $\mathbb{P}_{X}$. By the characterization of the probability law by its characteristic function and by the uniqueness of weak limits, we have $L=\mathbb{P}_{X}$. Then, there exists a probability measure, $L_0=\mathbb{P}_{X}$, such that each sub-sequence of $X_{n}$ contains a sub-sequence weakly converging $L_0$.\\

\noindent Step 2 : Let $f : \mathbb{R}^k \mapsto \mathbb{R}$ be a continuous and bounded function. Consider a sub-sequence 
$\mathbb{E}f(X_{n_{j}})$, $j\geq 1$, of the sequence of real numbers $\mathbb{E}f(X_n)$, $n\geq 1$.\\

\noindent This sub-sequence $X_{n_{j}}$, $j\geq 1$, contains a sequence $X_{n_{j_{\ell}}}$, $\ell \geq 1$, weakly converging to $L_0$ as $\ell \rightarrow +\infty$. Then $\mathbb{E}f(X_{n_{j_{\ell}}})$ converges to $\int f dL_0$. So $A=\int f dL_0$ is real number such that each sub-sequence of $\mathbb{E}f(X_n)$, $n\geq 1$ has a sub-sequence converging to $A$.\\

\noindent By Prohorov's Criterion Exercise 4 in Section \ref{funct.sec.1} of Chapter \ref{funct}, $\mathbb{E}f(X_n)$ converges to $\int f dL_0$. Then $X_n$ weakly converges to $L_0=\mathbb{P}_{X}$.

 

%% file: asymptotics_cv_03_en.tex
\chapter{Specific Tools for  Weak Convergence in $\mathbb{R}$} \label{cv.R}

This chapter focuses on tools which are specific to convergence of
sequences of real random variables. For such random variables, we may use
Renyi's representations through uniform or exponential random variables,
especially for sequences of independent and identically distributed random
variables. Such representations use the generalized inverse functions on
which concentrates the first section. Besides, in relation with Section \ref{cv.CvCp} and Theorem \ref{cv.skorohodWichura} in Chapter \ref{cv},  working on weak convergence in the same probability space may become a computation matter.
This chapter gives tools for such an orientation.\\

\section{Generalized inverses of monotone functions} \label{cv.sec2}

\bigskip \noindent This theory is done for non-decreasing and right-continuous functions. It may be done for non-increasing and
left-continuous functions.\\

\noindent Sometimes, left or right continuity is not required (see Point 9 below).\\

\noindent  Let $F$ be a non-decreasing and right-continuous function from $\mathbb{R}$ to $\mathbb{R}$. Let us define the generalized
inverse of $F$ as :

\begin{equation*}
F^{-1}(u)=\inf \{x\in \mathbb{R},F(x)\geq u\}, u\in \mathbb{R}.
\end{equation*}

\bigskip \noindent Because of the importance of this transformation for
univariate extreme value theory, we are going to expose important facts of
generalized inverses. Since we want them to be known by heart, we expose all of
them before we provide their proofs.\newline

\bigskip \noindent \textbf{A - List of most important properties of the generalized inverses}.\\

\noindent \textbf{Point (1).} For any $u\in \mathbb{R}$ and for any 
$t\in \mathbb{R}$ 
\begin{equation}
F(F^{-1}(u))\geq u  \tag{A}
\end{equation}

\noindent and

\begin{equation}
F^{-1}(F(x))\leq x.  \tag{B}
\end{equation}

\bigskip \noindent \textbf{Point (2).} For any $(u,t)\in \mathbb{R}^{2},$

\begin{equation}
(F^{-1}(u)\leq t)\Longleftrightarrow (u\leq F(t))  \tag{A}
\end{equation}

\noindent and

\begin{equation}
(F^{-1}(u)>t)\Longleftrightarrow (u>F(t)).  \tag{B}
\end{equation}

\bigskip
\bigskip \noindent \textbf{Point (3).} $F^{-1}$ is non-decreasing and
left-continuous.\newline

\bigskip \noindent \textbf{Point (4).} The weak convergence for
non-decreasing distribution functions is available by itself and is defined still by
Formula (9) above. Then we have the following implication.

\begin{equation*}
(F_n \rightsquigarrow F) \Rightarrow (F^{-1}_n \rightsquigarrow F^{-1}). 
\end{equation*}

\bigskip \noindent \textbf{Point (5).} Let us suppose that $F_{n}$ and $F$
are \textbf{distribution functions of real random variables} and that $%
F_{n}\rightsquigarrow F$. If $F$ is \textbf{continuous}, that we have the
uniform convergence

\begin{equation*}
\sup_{ x\in \mathbb{R}} |F_n(x)-F(x)| \rightarrow 0 \text{ as } n \rightarrow +\infty 
\end{equation*}

\bigskip
\bigskip \noindent \textbf{Point (6).} A distribution function $F$ on $\mathbb{R}$ has at most a countable points of discontinuity.\newline

\bigskip \noindent \textbf{Point (7).} Let $\mathbb{P}$ be any probability measure
on $\mathbb{R}$ with support $(a,b)$, meaning that 
\begin{equation*}
a=\inf \{x,\text{ }\mathbb{P}(]-\infty ,x])>0\}\text{ and }b=\inf \{x,\text{ 
}\mathbb{P}(]-\infty ,x])=1\}.
\end{equation*}

\noindent which $\mathbb{P}((a,b)^{c})=0$. Then, for $0<\varepsilon <1$, there exists a finite number
partition of $(a,b)$, 
\begin{equation*}
a=t_{0}<t_{1}<t_{2}<...<t_{k}<t_{k+1}=b 
\end{equation*}

\noindent such that for $0<i<k$,

\begin{equation*}
\mathbb{P}(]t_{i},t_{i+1}[)\leq \varepsilon. 
\end{equation*}

\bigskip \noindent We always can extend the bounds to

\begin{equation*}
-\infty \leq t_{0}<t_{1}<t_{2}<...<t_{k}<t_{k+1}\leq +\infty 
\end{equation*}

\bigskip \noindent since $\mathbb{P}(]-\infty ,a[)=0$ and $\mathbb{P}(]b,+\infty
\lbrack )=0$.\\

\bigskip \noindent \textbf{Point (8)} Let $F$ and $G$ be two distribution
function both non-increasing or both non-decreasing. If neither of them is
degenerated, then there exist two continuity points of both $F$ and $G$,  $x_{1}$ and $x_{2}$ such that $x_{1}<x_{2}$ and 
\begin{equation*}
F(x_{1})<F(x_{2})\text{ and }G(x_{1})<G(x_{2}).
\end{equation*}

\bigskip \noindent \textbf{Point (9)} Let $F$ be simply non-decreasing from $%
\mathbb{R}$ to $[a,b]$ without assumption of left or right-continuity. Then for any $y\in ]a,b[$,
\begin{equation*}
F(F^{-1}(y)-0)\leq y\leq F(F^{-1}(y)+0),
\end{equation*}

\noindent where $F(\cdot+)$ and $F(\cdot-)$ respectively stand the right and the left limit at $x$.\\

\noindent If the function $F$ is non-increasing, the generalized inverse is defined by
\begin{equation*}
F^{-1}(y)=\inf \{x\in R,F(x)\leq y\},y\in (a,b).
\end{equation*}

\noindent and we have the formula for $x\in (a,b)$
\begin{equation*}
F(F^{-1}(y)+0)\leq y\leq F(F^{-1}(y)-0)
\end{equation*}

\bigskip \noindent \textbf{B - Proofs of the points}.\\

\noindent \textbf{Proof of Point 1}. Part (A). Set

\begin{equation*}
A_{u}=\left\{ x\in \mathbb{R},F(x)\geq u\right\}, u\in \mathbb{R}.
\end{equation*}

\noindent Since $F^{-1}(u)=\inf A_{n}$, there exists a sequence $%
(x_{n})_{n\geq 0}\in A_{u}$ such that 
\begin{equation*}
\left\{ 
\begin{array}{ccc}
F(x_{n}) & \geq & u \\ 
&  &  \\ 
x_{n} & \downarrow & F^{-1}(u).
\end{array}%
\right.
\end{equation*}

\bigskip \noindent By right-continuity of $F$ we have 
\begin{equation*}
F(F^{-1}(u))\geq u.
\end{equation*}
\bigskip

\noindent This proves Formula (A). As for the Formula (B), consider $%
x\in \mathbb{R}$ and set 
\begin{equation*}
F^{-1}(F(x))=\inf A_{F(x)}.
\end{equation*}

\noindent Let us split $A_{F(x)}$ into 
\begin{equation*}
A_{F(x)}=\left[ -\infty ,x\right[ \cap A_{F(x)}+\left[ x,+\infty \right]
\cap A_{F(x)}
\end{equation*}

\begin{equation*}
=:A_{F(x)}(1)+A_{F(x)}(2).
\end{equation*}

\bigskip
\noindent By \textbf{Fact} 1, stated at the end of this section, we have
 
\begin{equation*}
\inf A_{F(x)}=\min (\inf A_{F(x)}(1),\inf A_{F(x)}(2)).
\end{equation*}

\noindent But

\begin{equation*}
y\in A_{F(x)}(1)\Longrightarrow y\leq x,\text{ then }\inf A_{F(x)}(1)\leq x.
\end{equation*}

\bigskip \noindent Next we obviously have 
\begin{equation*}
\inf A_{F(x)}(2)=\{y\geq x, F(y)\geq F(x) \}=x.
\end{equation*}

\noindent Thus

\begin{equation*}
\inf A_{F(x)}\leq x.
\end{equation*}

\noindent That is :

\begin{equation*}
F^{-1}(F(x))=\inf A_{F(x)}\leq x.
\end{equation*}

\bigskip \noindent This closes the proof of Point 1.\newline

\bigskip \noindent \textbf{Proof of Point 2}. It is obvious that each of
Formulas (A) and (B) is derived from the other by taking complementary. So,
we may only prove one of them, say (B). Suppose $\left( u>F(t)\right)$. By
right-continuity of $F$ at $t$, we can find $\varepsilon$ such that 
\begin{equation*}
u>F(t+\varepsilon ).
\end{equation*}

\noindent Now, for $x\in A_{u}$ we surely have 
\begin{equation*}
x>t+\varepsilon.
\end{equation*}

\noindent Otherwise, we would get, 
\begin{equation*}
x\leq t+\varepsilon \Longrightarrow F(x)\leq F(t+\varepsilon )<u,
\end{equation*}

\noindent and this would lead to the conclusion $x \notin A_{u}$, which is in
contradiction with the assumption. So, $x>t+\varepsilon$ for all $x\in A_{u}$.
This implies that

\begin{equation*}
\inf A_{u}=F^{-1}(u)\geq t+\varepsilon >t.
\end{equation*}

\noindent We proved the direct sens of the first formula. To prove the
indirect sense, consider $F^{-1}(u)>t$. Next, suppose that $u>F(t)$ does not
hold. This implies that $F(t) \geq u$, which is in contradiction with $t\in
A_{u}$ and next, 
\begin{equation*}
\inf A_{u}=F^{-1}(u)\leq t.
\end{equation*}

\noindent This is impossible. Then $u>F(t)$.

\bigskip \noindent \textbf{Proof of Point 3}. We begin to establish that $%
F^{-1}$ is non-decreasing. We have

\begin{equation*}
\forall u\leq u^{\prime },A_{u^{\prime }}\leq A_{u}\Longrightarrow \inf
A_{u^{\prime }}\leq \inf A_{u}.
\end{equation*}

\noindent This implies

\begin{equation*}
F^{-1}(u^{\prime })\geq F^{-1}(u)
\end{equation*}

\noindent Next, we have to prove that $F^{-1}$ is left-continuous. Let $u\in 
\mathbb{R}$. We have for any $h\geq 0$, 
\begin{equation*}
F^{-1}(u-h)\leq F^{-1}(u).
\end{equation*}

\noindent Thus 
\begin{equation*}
\underset{h\downarrow 0}{\lim }F^{-1}(u-h)\leq F^{-1}(u).
\end{equation*}

\noindent Suppose that 
\begin{equation*}
\underset{h\downarrow 0}{\lim }F^{-1}(u-h)=\alpha <F^{-1}(u).
\end{equation*}

\noindent We can find $\varepsilon >0$ such that $\alpha +\varepsilon <F^{-1}(u)$.
Now, for all $h\geq 0$,

\begin{equation*}
F^{-1}(u-h)<\alpha +\varepsilon.
\end{equation*}

\noindent By definition of the infimum, there exists $x$ such that 
\begin{equation*}
F(x)\geq u-h\text{ and }F^{-1}(u-h)<\alpha +\varepsilon.
\end{equation*}

\noindent By Formula (A) of Point 1 and Formula (B) of Point (2), we have

\begin{equation*}
F^{-1}(u-h)<\alpha +\varepsilon \Longrightarrow u-h\leq F(\alpha +\varepsilon ).
\end{equation*}

\noindent Then we get as $h\downarrow 0$%
\begin{equation*}
u\leq F(\alpha +\varepsilon).
\end{equation*}

\noindent Since this is true for any $\varepsilon >0$, we let $\varepsilon
\downarrow 0$ to get

\begin{equation*}
u\leq F(\alpha).
\end{equation*}

\bigskip \noindent But, by Formula (B) of Point (2) and by using the hypothesis, we
arrive at

\begin{equation*}
\left( \alpha <F^{1}(u)\right) \Leftrightarrow \left( F^{-1}(u)>\alpha
\right) \Leftrightarrow \left( u>F(\alpha )\right) .
\end{equation*}

\noindent This is clearly a contradiction. We conclude that 
\begin{equation*}
\underset{h\downarrow 0}{\lim }F^{-1}(u-h)=F^{-1}(u).
\end{equation*}

\noindent And next, $F^{-1}$ is left-continuous.\newline

\bigskip \noindent \textbf{Proof of Point 4}. Suppose that $F_{n}\rightsquigarrow F$. Let $y\in \mathbb{R}$ and
let $\varepsilon >$. Since the number of discontinuity of $F$ is at most
countable, we can find a continuity point $x$ of $F$ in the open interval $%
(F^{-1}(y)-\varepsilon, F^{-1}(y))$. By Point 2, $(F^{-1}(y))$ is equivalent
to $(F(x) < y)$.Since $x\in C(F)$, $F_{n}(x)\rightarrow F(x)$. Then for
values of $n$ large enough, we have $F_{n}(x)<y$ and then $x<F_{n}^{-1}(y)$.
We get 
\begin{equation*}
F^{-1}(y)-\varepsilon \leq x<F_{n}^{-1}(y)
\end{equation*}

\noindent that is for any $\varepsilon >0$, 
\begin{equation*}
F_{n}^{-1}(y)>F^{-1}(y)-\varepsilon .
\end{equation*}

\noindent We let go to $+\infty$ and $\varepsilon$ to
decrease to $0$, and we get for any $y\in \mathbb{R}$,
\begin{equation*}
\liminf_{n\rightarrow \infty }F_{n}^{1}(y)\geq F^{-1}(y).
\end{equation*}

\noindent Now let $y$ be a continuity of $F^{-1}$. For any $y^{\prime}>y$,
we can find a continuity point $x$ of $F$ such that

\begin{equation}
F^{-1}(y^{\prime })<x<F^{-1}(y^{\prime })+\varepsilon.  \label{d2}
\end{equation}

\bigskip
\noindent By Point 1, $x>F^{-1}(y^{\prime})\Longrightarrow F(x)\geq
F(F^{-1}(y))\geq y^{\prime }$. Then 
\begin{equation*}
y<y^{\prime }\leq F(x).
\end{equation*}

\noindent Since $x\in C^{\prime}F)$, we have $F_{n}(x)\rightarrow F(x)$, and then for large
values of $n$, $y<F_{n}(x)$ and by Formula (A) of Point 2, $F_{n}^{-1}(y)\leq x,$. By combining this with Formula  (\ref{d2}), we get 
\begin{equation*}
F^{-1}(y^{\prime })\geq x\geq F_{n}^{-1}(y).
\end{equation*}

\noindent Now let $n\rightarrow +\infty$ to obtain 
\begin{equation*}
\lim \sup_{n\rightarrow \infty }F_{n}^{-1}(y)\leq F^{-1}(y^{\prime }).
\end{equation*}

\noindent Next, let $y^{\prime} \downarrow y$, and get $F^{-1}(y^{\prime
})\downarrow F^{-1}(y)$ by continuity of $F^{-1}$ at $y$. We arrive at 
\begin{equation*}
\lim \sup_{n\rightarrow \infty }F_{n}^{-1}(y)\leq F^{-1}(y).\leq \lim
\inf_{n\rightarrow \infty }F_{n}^{1}(y).
\end{equation*}

\noindent We finally conclude that 
\begin{equation*}
F^{-1}(y)=\lim \sup_{n\rightarrow \infty }F_{n}^{-1}(y)=\lim
\inf_{n\rightarrow \infty }F_{n}^{-1}(t).
\end{equation*}

\bigskip \noindent \textbf{Proof of Point 5}. Since $F$ is non-decreasing, $x
$ is a discontinuity point of $F$ if and only if the discontinuity jump $%
F(x+)-F(x-)$ is positive. Denote by $D$ the set of all discontinuity points
of $F$, and for any $k\geq 1$, denote by $D_{k}$ the set of discontinuity
points such that $F(x+)-F(x-) > 1/k$ and by $D_{k,n}$ the set of
discontinuity points in the interval $[-n,n]$ such that $F(x+)-F(x-) > 1/k$.
We are going to show that $D_{k,n}$ is finite.\newline

\noindent Let us suppose we can find $m$ points $x_1$, ..., $x_m$ in $D_{k,n}
$. Since $F$ is non-decreasing, we may see that the sum of the discontinuity
jumps is less than $F(n)-F(-n)$. You may make a simple drawing for $m=3$ and
project the jumps to the y-axis to see this easily. So

\begin{equation*}
\sum_{1\leq j \leq m} F(x+)-F(x-) \leq F(n)-F(-n). 
\end{equation*}

\bigskip \noindent Since each of these jumps exceeds $1/k$, we have

\begin{equation*}
\sum_{1\leq j \leq m} (1/k) \leq \sum_{1\leq j \leq m} F(x+)-F(x-) \leq
F(n)-F(-n), 
\end{equation*}

\noindent and thus,

\begin{equation*}
m/k \leq F(n)-F(-n),
\end{equation*}

\noindent that is

\begin{equation*}
m \leq k(F(n)-F(-n)). 
\end{equation*}

\noindent We conclude by saying that we cannot have more that $%
[k(F(n)-F(-n))]$ points in $D_{k,n}$, so $D_{k,n}$ is finite. Since

\begin{equation*}
D=\cup_{n\geq 1} \cup_{k\geq 1} D(k,n), 
\end{equation*}

\noindent we see that $D$ is countable. This puts an end to the proof.\\

\bigskip \noindent \textbf{Proof of Point 6}. Let $0<\varepsilon <1$. Let $%
F(t)=\mathbb{P}(]-\infty ,t])$. This is a distribution function such that $%
F(\infty)=0$ and $F(+\infty)=1$. Fix $0<\varepsilon<1$. Set $%
k=[1/\varepsilon]$, where $[t]$ stands for the greatest integer less than or equal
to $t$. We then have 
\begin{equation*}
k\varepsilon \leq 1\leq k\varepsilon +\varepsilon
\end{equation*}

\noindent and denote 
\begin{equation*}
s_{i}=i\varepsilon \text{, for }i=1,...,k \text{ and }s_{k+1}=1.
\end{equation*}

\noindent Put

\begin{equation*}
t_{i}=F^{-1}(s_{i})=\inf \{u,\text{ G(u)}\geq s_{i}\}.
\end{equation*}

\noindent By Point 1, 
\begin{equation}
F(t_{i})\geq s_{i}. \label{tens_a}
\end{equation}

\bigskip
\noindent Next, for any $1\leq i<k,$ 
\begin{equation*}
F(t_{i+1}-)=\lim_{h\downarrow 0}F(t_{i+1}-h).
\end{equation*}

\noindent By definition of $t_{i+1}$, which the supremum of the values $u$
such that $F(u) \geq (i+i)\varepsilon$, we surely have, 
\begin{equation*}
F(t_{i+1}-h)<(i+1)\varepsilon .
\end{equation*}

\bigskip
\noindent By letting $h\downarrow 0$, we get 
\begin{equation}
B(t_{i+1}-)\leq (i+1)\varepsilon.  \label{tens_b}
\end{equation}

\bigskip
\noindent By putting together, (\ref{tens_a}) et (\ref{tens_b}), we have 
\begin{equation*}
\mathbb{P}(]t_{i},t_{i+1}[)=F(t_{i+1}-)-F(t_{i})\leq (i+1)\varepsilon
-i\varepsilon =\varepsilon ,
\end{equation*}

\bigskip
\noindent for $i=1,..,k$. For $i=k$, we have $F(t_{k+1})=1$ and 
\begin{equation*}
\mathbb{P}(]t_{k},t_{k+1}[)=1-F(t_{k})\leq 1-k\varepsilon \leq \varepsilon .
\end{equation*}

\bigskip
\noindent For $i=0$, since $F(t_{0})\geq 0$, we have 
\begin{equation*}
\mathbb{P}(]t_{i},t_{i+1}[)=F(t_{1}-)-F(t_{0})\leq F(t_{1}-)\leq \varepsilon.
\end{equation*}

\bigskip
\noindent We just proved that $0\leq i\leq k,$%
\begin{equation*}
\mathbb{P}(]t_{i},t_{i+1}[)=F(t_{i+1}+)-F(t_{i})\leq
(i+1)\varepsilon-i\varepsilon=\varepsilon .
\end{equation*}

\bigskip \noindent \textbf{Proof of Point 7}. We are going to apply Point 6.
Let us consider the Lebesgue-Stieljes probability measure generated by $F$
and characterized by 
\begin{equation*}
\mathbb{P}(]u,v])=F(v)-F(u), u \leq v. 
\end{equation*}

\bigskip
\noindent In particular, we have $\mathbb{P}(]-\infty,v])=F(v)-F(u)$. Fix $%
\varepsilon >0$ and consider a subdivision 
\bigskip
\begin{equation*}
-\infty=t_{0}<t_{1}<t_{2}<...<t_{k}<t_{k+1}=+\infty 
\end{equation*}

\bigskip
\noindent such that for any $0\leq j \leq k$, 
\begin{equation*}
F(t_{j+1}-)-F(t_{j})=\mathbb{P}_{X}(]t_{i},t_{i+1}[)\leq \varepsilon .
\end{equation*}

\noindent Now we want to prove the uniform convergence. Let $x$ be one of
the $t_{j}$'s. We have 
\begin{equation*}
F_{n}(x)-F(x)\leq \sup_{0\leq j\leq k+1}\left| F_{n}(t_{j})-F(t_{j})\right|.
\end{equation*}

\noindent Any other $x$ is in one of the intervals $]t_{j},t_{j+1}[$. Use
the non-decreasingness of $F$ and $F_{n}$ to have

\begin{eqnarray*}
F_{n}(x)-F(x)&\leq& F_{n}(t_{j+1}-)-F(x)\\
&\leq &F_{n}(t_{j+1}-)-F(t_{j+1}-)+F(t_{j+1}-)-F(x)\\
&\leq & F_{n}(t_{j+1}-)-F(t_{j+1}-)+F(t_{j+1}-)-F(t_{j})\\
&\leq & \sup_{0\leq j\leq k+1}\left| F_{n}(t_{j}-)-F(t_{j}-)\right| +\varepsilon
\end{eqnarray*}

\noindent and

\begin{eqnarray*}
F(x)-F_{n}(x)&\leq & F(x)-F_{n}(t_{j})\\
&\leq & F(x)-F(t_{j})+F(t_{j})-F_{n}(t_{j})\\
&\leq & F(t_{j+1}-)-F(t_{j})+F(t_{j})-F_{n}(t_{j})\\
&\leq & \sup_{0\leq j\leq k+1}\left|F_{n}(t_{j})-F(t_{j})\right| +\varepsilon
\end{eqnarray*}

\noindent At the arrival, we have for any point $x$ different from $t_j$,
 
\begin{equation*}
\left| F(x)-F_{n}(x)\right| \leq \max \left(\sup_{0\leq j\leq k+1}\left|
F_{n}(t_{j}-)-F(t_{j}-)\right| ,\sup_{0\leq j\leq k+1}\left|
F_{n}(t_{j})-F(t_{j})\right| \right)+\varepsilon).
\end{equation*}

\noindent Then

\begin{eqnarray*}
\sup_{x\in \mathbb{R}}\left| F_{n}(x)-F(x)\right| &=&\left\| F_{n}-F\right\|_{\infty }\\
&\leq& \max \left(\sup_{0\leq j\leq k+1}\left|F_{n}(t_{j}-)-F(t_{j}-)\right| ,\sup_{0\leq j\leq k+1}\left|
F_{n}(t_{j})-F(t_{j})\right|\right)\\
&+&\varepsilon.
\end{eqnarray*}

\noindent At this step, we have the more general conclusion. If for all real 
$x$, $F_n(x) \rightarrow F(x)$ and $F_n(x-) \rightarrow F(x-)$, then we may
conclude that

\begin{equation*}
\limsup_{n\rightarrow +\infty} \sup_{ x\in \mathbb{R}} |F_n(x)-F(x)| \leq \varepsilon, 
\end{equation*}

\noindent for an arbitrary $\varepsilon$. Thus we have

\begin{equation*}
\lim_{n\rightarrow +\infty} \sup_{ x\in \mathbb{R}} |F_n(x)-F(x)| =0. 
\end{equation*}

\bigskip \noindent To extend this conclusion to the case $F$ is continuous and $%
F_n(x) \rightarrow F(x)$, we have to prove $F_n(x-) \rightarrow F(x)$ for
any $x$.\newline

\noindent To prove this, fix an arbitrary $x$ and let $0\leq h_p \downarrow 0
$ as $p \uparrow +\infty$. We have for each $n$,

\begin{equation*}
F_n(x-)-F(x) \leq F_n(x)-F(x) \leq |F_n(x-)-F(x)|, 
\end{equation*}

\noindent and

\begin{equation*}
F(x)-F_n(x-) \leq F(x)-F_n(x-h_p) \leq |F(x)-F(x-h_p)|+|F(x-h_p)-F_n(x-h_p)|.
\end{equation*}

\bigskip \noindent By combining these two points, we have 
\begin{equation*}
|F(x)-F_n(x-)| \leq max(|F_n(x-)-F(x)|,
|F(x)-F(x-h_p)|+|F(x-h_p)-F_n(x-h_p)|). 
\end{equation*}

\bigskip
\noindent \noindent Now fix $p$ and let $n\rightarrow +\infty$ to get

\begin{equation*}
\limsup_{n\rightarrow +\infty}|F(x)-F_n(x-)| \leq |F(x)-F(x-h_p)|. 
\end{equation*}

\bigskip
\noindent Finally, let $p\rightarrow +\infty $ to get the conclusion by
continuity of $F$.\newline

\bigskip \noindent \textbf{Proof of Point 8}. $F$ is degenerated if and only
if it has a unique point of increase, say $a$, at which it presents a
discontinuity jump. It is in the form : $F(x)=c_{1}$ for $x<a$ and $F(x)=c_{2}$ for $x\geq a$, with $c_{1}< c_{2}$. So if $F$ is non-degenerated, it has at least  two points of increase. Hence, we can find three continuity points of $F$ : $a_{1}<a_{2}<a_{3}$ such that $F(a_{1})<F(a_{2})<F(a_{3})$. If $G$ is also 
non-degenerated, we also find three continuity points of $G$ : $b_{1}<b_{2}<b_{3}$ such that $F(b_{1})<F(b_{2})<F(b_{3})$. We consider two cases.\newline

\noindent Case 1. The intervals $[a_1, a_3]$ and $[b_1, b_3]$ are disjoint
or have an intersection of one point, which is necessarily a continuity point. Suppose for example that $a_3 \leq b_1$%
. Take $x_1=a_1$ and $x_2=b_3$. We have 
\begin{equation*}
F(x_1)<F(a_3)\leq F(b_1) \leq F(b_3)=F(x_2) 
\end{equation*}

\noindent and 
\begin{equation*}
G(x_1) < G(a_3) \leq G(b_1) < G(b_3)=G(x_2). 
\end{equation*}

\bigskip
\noindent Case 2. The intervals $[a_1, a_3]$ and $[b_1, b_3]$ overlap at least on a non-empty open interval. Take $%
t$ in the intersection. Surely we have $F(a_1)<F(t)$ or $F(t)<F(a_3)$.
Otherwise, we would have $F(a_1)=F(a_3)$, which violates what is above.
Similary $G(b_1)<G(t)$ or $G(t)<G(b_3)$. Now, take $x_1=min(a_1,b_1)$ and $x_x=min(a_3,b_3)$.\newline

\noindent If $F(a_1)<F(t)$, we have

\begin{equation*}
F(x_1)\leq F(a_1) < F(t) \leq F(a_3) \leq F(x_2). 
\end{equation*}

\bigskip

\noindent If $F(t)<F(a_3)$, we have

\begin{equation*}
F(x_1)\leq F(a_1)  < F(t) <  F(a_3) \leq F(x_2) 
\end{equation*}

\bigskip
\noindent We conclude that

\begin{equation*}
F(x_1) < F(x_2). 
\end{equation*}

\noindent We prove similarly that

\begin{equation*}
G(x_1) < G(x_2). 
\end{equation*}

\noindent By Point 5, we know that the discontinuity points of $F$ and $G$ are at
most countable, we may adjust $x_{1}$ and $x_{2}$ to be continuity points of
both $F$ and $G$.

\bigskip 

\noindent \textbf{Proof of Point 9}.\\

\noindent Let us begin by the first case where $F$ est non-decreasing. By definition of the generalized inverse, 
we have for any $h>0$,
\begin{equation*}
F(F^{-1}(y)+h)\geq y
\end{equation*}

\noindent and
\begin{equation*}
F(F^{-1}(y)-h)<y.
\end{equation*}

\noindent By letting $h$ decrease to zero, we get
\begin{equation*}
F(F^{-1}(y)-)\leq y\leq F(F^{-1}(y)+), \ \ y \in (a,b).
\end{equation*}

\noindent Similarly, if $F$ is non-increasing, we have for any $h>0$
\begin{equation*}
F(F^{-1}(y)+h)\leq y
\end{equation*}

\noindent and
\begin{equation*}
F(F^{-1}(y)-h)>y.
\end{equation*}

\noindent By letting $h$ decrease to zero, we get

\begin{equation*}
F(F^{-1}(y)+)\leq y\leq F(F^{-1}(y)-).
\end{equation*}

\bigskip \noindent \textbf{Fact 1}. Let $A$ and $B$ be two disjoint subsets of $R$. We have 
\begin{equation*}
\inf A \cup B = min(\inf A, \inf B). 
\end{equation*}

\noindent Indeed, clearly, $\inf A \cup B$ is less than $\inf A$ and less
than $\inf B$, and then $\inf A \cup B \leq min(\inf A, \inf B)$. Now
suppose that we do not have the equality, that is

\begin{equation*}
\inf A \cup B < min(\inf A, \inf B). 
\end{equation*}

\bigskip \noindent There exists a sequence $(z_n)_{n\leq 0}$ of points of $A \cup B$
decreasing to $\inf (A \cup B)$. Surely for $n$ large enough, $z_n$ will be
less than $\inf A$ and less than $\inf B$. And yet, it is either in $A$ or
in $B$. This is absurd. We conclude that we have the equality.\\

\bigskip  

\section{Applications of Generalized functions}

The first application is the representation of any real random variable by a
standard uniform random variable $U \sim \mathcal{U}(0,1)$ associated with the
distribution function $G(x)=0$ for $x<0$, $G(x)=x$ for $x\in (0,1)$ and $%
G(x)=1$ for $x>0$. We have :

\begin{lemma}
\label{lemmatooldf10} Let F be a distribution function such that $%
F(-\infty)=0$ and $F(+\infty)=1$. Let $U \sim \mathcal{U}(0,1)$, defined on
some probability space $(\Omega,\mathcal{A},\mathbb{P})$. Then $X=F^{-1}(U)$
has the distribution function $F$.
\end{lemma}

\noindent \textit{Proof}. We have by Formula (A) of Point 2 above that, for any $x\mathbb{R}$,
\begin{equation*}
\mathbb{P}(X \leq x)=\mathbb{P}(F^{-1}(U) \leq x)=\mathbb{P}(U
\leq F(x))=F(x). 
\end{equation*}

\noindent $\blacksquare$.\\

\bigskip \noindent A second application is this simple form of Skorohod-Wichura's Theorem.

\begin{theorem} \label{cv.skorohod} Let $F_n \rightsquigarrow F$, where $F_n$ and $F$ are
distribution functions such that $F_n(-\infty)=0$ and $F_n(+\infty)=1$, for $%
n\geq 0$, $F(-\infty)=0$ and $F(+\infty)=1$. Then, there exists a
probability space $(\Omega,\mathcal{A},\mathbb{P})$ holding a sequence
of real random variables $X_n$ and a random variable $X$ such that for any $%
n\geq 0$, $F_n$ is the distribution function of $X_n$, that is $F_n(.)=%
\mathbb{P}(X_n\leq .)$, and $F$ is the distribution function of $X$, that is 
$F(.)=\mathbb{P}(X \leq .)$ and

\begin{equation*}
X_n \rightarrow X \text{a.s as } n\rightarrow +\infty. 
\end{equation*}
\end{theorem}

\noindent \textbf{Proof}. Let us consider $([0,1],\mathcal{B}([0,1]),\lambda)
$ where $\lambda$ is the Lebesgue-measure on $[0,1]$, which is a
probability measure. Consider the identity function  $U$ : $([0,1],\mathcal{%
B}([0,1]),\lambda) \mapsto ([0,1],\mathcal{B}([0,1]),\lambda)$. Then $U$
follows a standard uniform law since for any $x\in (0,1)$,

\begin{equation*}
\lambda(U \leq x)=\lambda(U^{-1}(]-\infty,x])) 
\end{equation*}

\bigskip
\noindent where $U^{-1}$ is the inverse of $U$ and then $%
U^{-1}(]-\infty,x])=]-\infty,x]$. Thus 
\begin{equation*}
\lambda(U \leq x)=\lambda(U^{-1}(]-\infty,x]))=\lambda(]-\infty,x])=x. 
\end{equation*}

\bigskip
\noindent Consider $X_n=F_n^{-1}(U)$, $n\geq 1$, and $X=F^{-1}(U)$. In virtue of Lemma \ref{lemmatooldf10} above, each $F_n$, $n\geq 1$, is
the distribution function of $X_n$ and $F$ is the distribution function of $X
$. Let us show that $X_n$ converges to $X$ almost-surely. By using Point 4
above, we have $F_n^{-1} \rightsquigarrow F^{-1}$. Then

\begin{equation*}
1 \geq \lambda(X_n \rightarrow X) = \lambda(\{u\in [0,1], X_n(u) \rightarrow
X(u) \}) 
\end{equation*}

\begin{equation*}
=\lambda(\{u\in [0,1], F_n^{-1}(u) \rightarrow F^{-1}(u)\}) 
\end{equation*}

\begin{equation*}
\geq \lambda(\{u \in [0,1], \text{u is a continuity point of } F \})=1, 
\end{equation*}

\noindent since the complement of $\{u \in [0,1], \text{u is a continuity
point of } F \}$ is countable and countable sets are null-sets with respect
to the Lebesgues measure.

\section{Representation of Renyi for $iid$ sequences of random variables} \label{cv.R.renyi}

\noindent This section is intended to provide representations of order
statistics $X_{1,n}\leq ...\leq X_{n,n}$ , $n\geq 1,$ of any $n$ independent
random variables $X_{1},...,X_{n}$ with common distribution function $F$ in
that of standard uniform or exponential independent random variables.\\

\noindent We remind again that in this section, all the random variables are defined on the same probability space $\left( \Omega ,\mathcal{A},\mathbb{P}\right).$\\

\noindent We begin by recalling the density probability function of the order statistics from a density probability function $h$.

\subsection{Density of the order statistics} $ $\\

\noindent Let us begin with this lemma.\\

\begin{lemma} \label{cv.R.malmquist01P} Let $Z_{1},Z_{2},...,Z_{n}$ \ be $n$ independent
copies of an absolutely continuous random variable $Z$ of probability
density function $h$ and probability distribution function $H,$ defined on
the same probability space $\left( \Omega ,\mathcal{A},\mathbb{P}\right) .$
Let $1\leq r\leq n$, $1\leq n_{1}<n_{2}<...<n_{r}$. Then the $r$ order
statistics $Z_{n_{1},n}<Z_{n_{2},n}<_{.}..<Z_{n_{r},n}$ have the joint the
joint probability density function in $(z_{1},...,z_{r})$,

\begin{equation}
n!\prod_{j=1}^{r+1}\frac{h(z_{j})(F(z_{j})-F(z_{j-1}))^{n_{j}-n_{j-1}-1}}{%
(n_{j}-n_{j-1}-1)!}1_{(z_{1}<...<z_{r}),}  \label{cv.densityExtreme03}
\end{equation}

\noindent with by convention $n_{0}=0$ and $n_{r}=n+1,z_{0}=-\infty =z_{0}$ and $%
z_{r+1}=+\infty .$
\end{lemma}

\bigskip \noindent \textbf{Proof}. 
Suppose that the assumptions of the proposition holds. Let us find the joint
density probability functions of $r$ order statistics $Z_{n_{1},n}\leq
Z_{n_{2},n}\leq ...\leq Z_{n_{r},n}$, with $1\leq r\leq n$, $1\leq
n_{1}<n_{2}<...<n_{r}$ .\newline

\noindent Since $Z$ is an absolutely continuous random variable, the
observations are distinct almost surely and we have $%
Z_{n_{1},n}<Z_{n_{2},n}<_{.}..<Z_{n_{r},n}$. Then for $dz_{i}$ small
enough and for $z_{1}<z_{2}<...<z_{r}$, the event 
\begin{equation*}
(Z_{n_{i},n}\in ]z_{i}-dz_{i}/2,z_{i}+dz_{i}/2[,\ \ 1\leq i\leq r)
\end{equation*}

\noindent occurs with $n_{1}-1$ observations of the sample $%
Z_{1},...,Z_{n}$ falling at left of $z_{1}$, one point in $]z_{1}-dz_{1}/2,z_{1}+dz_{1}/2[$, $n_{2}-n_{1}-1$ between $z_{1}+dz_{1}/2$
and $z_{2}-dz_{1}/2$, one point in $]z_{2}-dz_{2}/2,z_{1}+dz_{2}/2[$,
etc., and $n-k_{k}$ points at right of $z_{r}$.\\

\begin{figure}[htbp]
	\centering
		\includegraphics[width=1.00\textwidth]{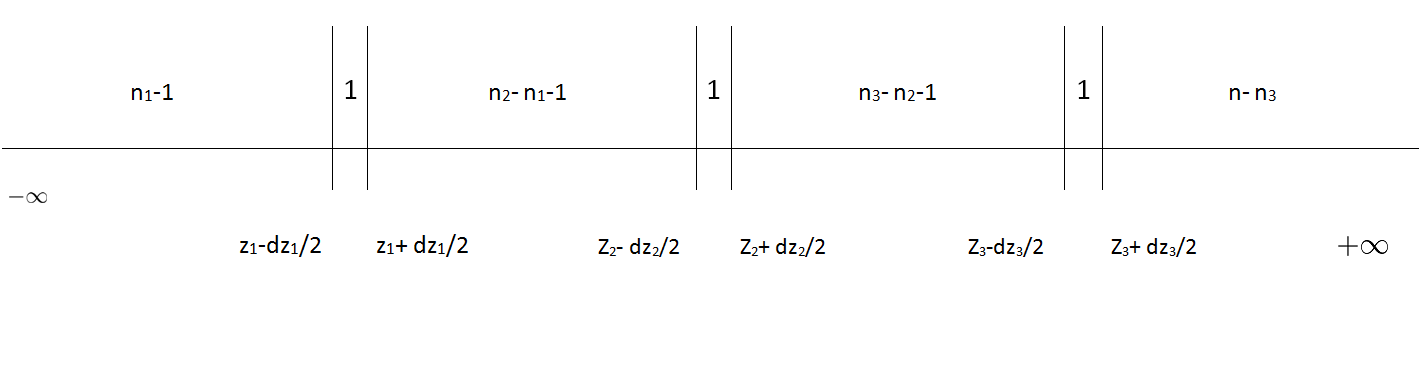}
	\caption{How are placed the observations with respect to $z_1<...<z_r$}
	\label{cv.R.FigOS}
	\end{figure}

\noindent This is illustrated in Figure \ref{cv.R.FigOS} for $r=3$.\newline

\noindent By definition, the probability density function $%
f_{(Z_{n_{1},n},...,Z_{n_{r},n})}$, whenever it exists, satisfies

\begin{eqnarray}
&&\frac{\mathbb{P}(Z_{n_{i},n}\in ]z_{1}-dz_{i}/2,z_{i}+dz_{i}/2[,\
\ 1\leq i\leq r)}{dz_{1}\times ...\times dz_{r}} \notag \\
&=&f_{(Z_{n_{1},n},...,Z_{n_{r},n})}(z_{1},...,z_{r})(1+\varepsilon(dz_{1},...,dz_{r})),  \label{cv.densityExtreme02}
\end{eqnarray}

\noindent where $\varepsilon (dz_{1},...,dz_{r})\rightarrow 0$ as each $%
dz_{i}\rightarrow 0$ $(1\leq i\leq r).$ Now, by using the independence $%
Z_{1},...,Z_{n},$ $\mathbb{P}(Z_{n_{i},n}\in ]z_{1}-dz_{i}/2,z_{i}+dz_{i}/2[,\ \ 1\leq i\leq r)$ is obtained as a multinomial
probability. Using in addition the fact  that $h$ is the common probability
density function of $Z$, we get 

\begin{eqnarray*}
&&\frac{\mathbb{P}(Z_{n_{i},n}\in ]z_{1}-dz_{i}/2,z_{i}+dz_{i}/2[,\ \ 1\leq i\leq r)}{dz_{1}\times ...\times dz_{r}} \\
&=&n!\times \frac{h(z_{1})^{n_{1}-1}}{(n_{1}-1)!}\times \frac{%
(F(z_{2})-F(z_{1}))^{n_{2}-n_{1}-1}}{(n_{2}-n_{1}-1)!}\\
&\times& ....\times \frac{(F(z_{j})-F(z_{j-1}))^{n_{j}-n_{j-1}-1}}{((n_{j}-n_{j-1}-1)!)}\\
&\times&...\times \frac{(F(z_{r})-F(z_{r-1}))^{n_{r}-n_{r-1}-1}}{(n_{r}-n_{r-1}-1)!}\\
&\times& \frac{(1-F(z_{r}))^{n-n_{r}}}{(n-n_{r})!} \\
&\times& \prod \frac{\mathbb{P}(Z_{n_{i},n}\in ]z_{i}-dz_{i}/2,z_{i}+dz_{i}/2[)}{1! \ dz_{i}}
\end{eqnarray*}

\noindent The last factor in the latter product is

$$
\prod_{i=1}^{r} h(z_{i})(1+dz_{i}).
$$
 
\noindent By setting $n_{0}=0$ and $n_{r}=n+1$ and for $-\infty
=z_{0}<z_{1}<...<z_{r}<z_{r+1}=+\infty ,$%
\begin{equation*}
f_{(Z_{n_{1},n},...,Z_{n_{r},n})}(z_{1},...,z_{r})=n!\prod_{j=1}^{r+1}\frac{%
h(z_{j})(F(z_{j})-F(z_{j-1}))^{n_{j}-n_{j-1}-1}}{(n_{j}-n_{j-1}-1)!},
\end{equation*}

\noindent we see that $f_{(Z_{n_{1},n},...,Z_{n_{r},n})}$ satisfies (\ref%
{cv.densityExtreme02}). $\blacksquare$\\

\bigskip \noindent Now, let us apply this lemma to the whole order statistics. We get this proposition.

\begin{proposition} \label{cv.R.malmquist01} Let $Z_{1},Z_{2},...,Z_{n}$ \ be $n$ independent
copies of an absolutely continuous random variable $Z$ with common probability
density function $h$, and defined on the same probability space $\left( \Omega ,%
\mathcal{A},\mathbb{P}\right)$. The associated order statistic
$$
(Z_{1,n}, Z_{2,n}, ...,  Z_{n,n})
$$ 

\noindent has the joint probability density function
 
\begin{equation*}
h_{\left( Z_{1,n},...,Z_{n,n}\right) }\left( z_{1,...,}z_{n}\right)
=n!\prod\limits_{i=1}^{n}h\left( z_{i}\right) 1_{\left( z_{1}\leq ...\leq
z_{n}\right) .}
\end{equation*} 
\end{proposition}

\bigskip \noindent \textbf{Proof}. Let us apply Lemma \ref{cv.R.malmquist01P} with $r=n$ and $n_{1}=1,n_{2}=2,...,n_{n}=n.$
Since the numbers $n_{j}-n_{j-1}-1$ vanish in (\ref{cv.densityExtreme03}),
It comes that $Z_{1,n}<Z_{2,n}<_{.}..<Z_{n,n}$ have the joint probability
density
\begin{equation*}
n!\prod_{j=1}^{n}h(z_{j})1_{(z_{1}<...<z_{r}),}
\end{equation*}

\noindent $\blacksquare$\\

\bigskip \noindent Now, we are focusing on the relation between standard uniform and
exponential order statistics.

\begin{proposition} \label{cv.R.malmquist02} Let $n\geq 1$ be a fixed integer and \ $U_{1,n}\leq
U_{2,n}\leq ...\leq U_{n,n}$ be the order statistics associated with $%
U_{1},U_{2},...,U_{n},$ which are $n$\ independent \ random variables
uniformly \ distributed on $\left( 0,1\right) .$ Let \ $%
E_{1},E_{2},...,E_{n},E_{n+1},$ $\left( n+1\right) $ independent random \
variables following the standard exponential law, that is 
\begin{equation*}
\forall x\in \mathbb{R},\text{ \ }\mathbb{P}\left( E_{i}\leq n\right)
=\left( 1-e^{-x}\right) 1_{\left( x\geq 0\right) },i=1,...,n+1.
\end{equation*}

\noindent Let $S_{j}=E_{1}+...+E_{j},$ $1\leq j\leq n+1$. Then we have the
following equality in distribution%
\begin{equation*}
\left( U_{1,n},U_{2,n},...,U_{n,n}\right) =^{d}\left( \frac{S_{1}}{S_{n+1}}%
,...,\frac{S_{n}}{S_{n+1}}\right).
\end{equation*}
\end{proposition}

\bigskip \noindent \textbf{Proof}. On one hand, by (\ref{cv.R.malmquist01}), the probability density function (\textit{pdf}) of 
$U=(U_{1,n}, U_{2,n}, U_{n,n})$ is given by 
\begin{equation*}
\forall \left( u_{1},...,u_{n}\right) \in \mathbb{R}^{n},f_{U}\left(
u_{1},...,u_{n}\right) =n!1_{\left( 0\leq u_{1}\leq ...\leq u_{n}\leq
1\right) }.
\end{equation*}

\noindent We are going to find the \ distribution \ of \ $Z_{n+1}^{\ast }=\left(
S_{1},S_{2},...,S_{n},S_{n+1}\right) $ given $S_{n+1}=t$, $t>0.$ We have \
for $y=\left( y_{1},y_{2},...,y_{n}\right) \in \mathbb{R}^{n},$%
\begin{eqnarray}
f_{Z_{n}^{\ast }}^{S_{n+1}=t}\left( y\right)  &=&\frac{f_{\left( Z_{n}^{\ast
},S_{n+1}\right) }\left( y,t\right) }{f_{S_{n+1}}\left( t\right) } \label{cv.m01} \\
&=&\frac{f_{Z_{n+1}^{\ast }}\left( y,t\right) }{f_{S_{n+1}}\left( t\right) }%
\boldsymbol{1}_{\left( 0\leq y_{1}\leq ...\leq y_{n}\leq t\right) }. \notag
\end{eqnarray}

\bigskip \noindent But $S_{n+1}$ follows a gamma law of parameters $n+1$ and $1,$ that is $%
S_{n+1}\sim \gamma \left( n+1,1\right) ,$ and its \ probability density
function is  
\begin{equation}
f_{S_{n+1}}\left( t\right) =\frac{t^{n}e^{-t}}{\Gamma \left( n+1\right) }%
1_{\left( t\geq 0\right) }=\frac{t^{n}}{n!}e^{-t}\boldsymbol{1}_{\left(
t\geq 0\right) }.  \label{cv.m02}
\end{equation}

\bigskip \noindent The distribution function of $\left( S_{1},...,S_{n+1}\right) $ comes \ from
the transformation
\begin{equation*}
\left( 
\begin{array}{c}
E_{1} \\ 
. \\ 
. \\ 
. \\ 
. \\ 
E_{n+1}%
\end{array}%
\right) =\left( 
\begin{array}{cccccc}
1 &  &  &  &  &  \\ 
-1 & 1 &  &  &  &  \\ 
0 & -1 & 1 &  &  &  \\ 
. &  & . & . &  &  \\ 
. &  &  & . & . &  \\ 
. & . & . & . & -1 & 1%
\end{array}%
\right) \left( 
\begin{array}{c}
S_{1} \\ 
. \\ 
. \\ 
. \\ 
. \\ 
S_{n+1}%
\end{array}%
\right) 
\end{equation*}

\bigskip \noindent 
Let $B$ be the matrix on the formula above. The Jacobian determinant in absolute value is $\left\vert B\right\vert =1$
and 
\begin{equation*}
B\left( 
\begin{array}{c}
y_{1} \\ 
. \\ 
. \\ 
y_{n+1}%
\end{array}%
\right) =\left( y_{1},y_{2}-y_{1},...,y_{n+1}-y_{n}\right). 
\end{equation*}

\bigskip \noindent Thus, the \ density of \ $\left( S_{1},...,S_{n+1}\right) $ \ is then given by

\begin{eqnarray*}
f_{Z_{n+1}^{n}}(y_{1},...,y_{n+1}) &=&f_{\left( E_{1},...,E_{n+1}\right)
}\left( B\left( y_{1},...,y_{n+1}\right) \right) \boldsymbol{1}_{\left(
0\leq y_{1}\leq ...\leq y_{n}\leq y_{n+1}\right) } \\
&=&\prod\limits_{i=1}^{n+1}e^{-\left( y_{i}-y_{i-1}\right) }\boldsymbol{1}%
_{\left( 0\leq y_{1}\leq ...\leq y_{n}\leq y_{n+1}\right) } \\
&=&e^{-y_{n+1}}\boldsymbol{1}_{\left( 0\leq y_{1}\leq ...\leq y_{n}\leq
y_{n+1}\right) }.
\end{eqnarray*}

\noindent where $y_{0}=0$ by convention. Going back to (\ref{cv.m01}) and (\ref{cv.m02}), we get , with $y=\left(
y_{1},y_{2},...,y_{n}\right) ,$%
\begin{equation}
f_{Z_{n}^{\ast }}^{S_{n+1}=t}\left( y\right) =\frac{n!}{t^{n}}\boldsymbol{1}%
_{\left( 0\leq y_{1}\leq ...\leq y_{n}\leq t\right) }.  \label{cv.m03}
\end{equation}

\noindent Now, for $u=\left( u_{1},u_{2},...,u_{n}\right) \in \mathbb{R}^{n},$%
\begin{equation*}
f_{\left( \frac{S_{1}}{S_{n+1}},...,\frac{S_{n}}{S_{n+1}}\right)
}^{S_{n+1}=t}\left( u\right) =f_{\left( \frac{S_{1}}{t},...,\frac{S_{n}}{t}%
\right) }^{S_{n+1}=t}\left( u_{1},u_{2},...,u_{n}\right) .
\end{equation*}

\noindent This  density \ probability function \ is obtained from (\ref{cv.m03}) by
the transform 
\begin{equation*}
\left( y_{1},y_{2},...,y_{n}\right) =t\left( u_{1},u_{2},...,u_{n}\right)
\Longleftrightarrow \left( u_{1},u_{2},...,u_{n}\right) =\frac{1}{t}\left(
y_{1},y_{2},...,y_{n}\right) 
\end{equation*}

\noindent with Jacobian determinant $t^{n}.$ Then 
\begin{eqnarray*}
f_{\left( \frac{S_{1}}{S_{n+1}},...,\frac{S_{n}}{S_{n+1}}\right)
}^{S_{n+1}=t}\left( u_{1},u_{2},...,u_{n}\right)  &=&f_{Z_{n}^{\ast
}}^{S_{n+1}=t}\left( t\left( u_{1},u_{2},...,u_{n}\right) \right) t^{n}%
\boldsymbol{1}_{\left( 0\leq tu_{1}\leq ...\leq tu_{n}\leq t\right) } \\
&=&n!\boldsymbol{1}_{\left( 0\leq u_{1}\leq ...\leq u_{n}\leq 1\right) }.
\end{eqnarray*}

\bigskip \noindent This is exactly (\ref{cv.m01}). Then the \ conditional distribution of \ $%
Z=\left( \frac{S_{1}}{S_{n+1}},...,\frac{S_{n}}{S_{n+1}}\right)$ given $%
S_{n+1}=t$ does not depend on $t$. So, its conditional distribution is also its unconditional
distribution function and then 
$$\left( \frac{S_{1}}{S_{n+1}},...,\frac{S_{n}}{S_{n+1}}\right)
$$ 

\noindent has the same law as $U=(U_{1,n}, U_{2,n}, ..., U_{n,n})$ and it is independent of $S_{n+1}$. This puts an end to the proof.\\

\bigskip \noindent We formalize the last conclusion in the following lemma.

\begin{lemma} \label{cv.R.malmquist03} Let \ $E_{1},E_{2},...,E_{n},E_{n+1}$, $n\geq 1$,  be
independent standard exponential random variables defined on the same probability
space. Let \ $S_{i}=E_{1}+E+...+E_{i},1\leq i\leq n+1,$ then 
$$
\left( 
\frac{S_{1}}{S_{n+1}},...,\frac{S_{n}}{S_{n+1}}\right)
$$

\noindent is independent of $S_{n+1}.$
\end{lemma}

\bigskip \noindent The latter proposition exposed representations of order statistics of standard
uniform random variables into that of standard exponential random variables.
The following proposition reverses the situation.

\begin{proposition} \label{cv.R.malmquist03} Assume the notations of Proposition \ref%
{cv.R.malmquist02} hold. Then for any $n\geq 1,$ 
\begin{equation*}
\left( -\log U_{1,n},...,-\log U_{n,n}\right) =_{d}\left(
E_{1,n},...,E_{n,n}\right) ,
\end{equation*}

\noindent where $E_{1,n}\leq ...\leq E_{n,n}$ are the order statistics \ of \ $%
E_{1},E_{2},...,E_{n},$ which are $n$ independent and exponentially
distributed with intensity one.
\end{proposition}

\bigskip \noindent \textbf{Proof}. By Proposition \ref{cv.R.malmquist01}, the \textit{pdf} of $E_{1,n}\leq
...\leq E_{n,n}$ is
 
\begin{equation}
f_{Z}\left( z\right) =n!e^{-\sum_{i=1}^{n}z_{i}}\boldsymbol{1}_{\left( 0\leq
z_{1}\leq ...\leq z_{n}\right) },%
\text{ }z=\left( z_{1},...,z_{n}\right) \in \mathbb{R}^{n}, \label{cv.m05}
\end{equation}

\noindent where $Z=\left( E_{1,n},...,E_{n,n}\right) .$ The distribution \ of $Z^{\ast
}=\left( -\log U_{1,n},...,-\log U_{n,n}\right) $ comes from that of $%
U=\left( U_{1,n},...,U_{n,n}\right) $ by the diffeomorphism $\left(
z_{1},...,z_{n}\right) =\left( -\log u_{1},...,-\log u_{n}\right) $ which
preserves the order of the arguments and has a jacobian determinant in
absolute value equal to

\begin{eqnarray*}
\left\vert \frac{\partial U_{i}}{\partial z_{j}}\right\vert  &=&\left\vert
\partial \frac{\partial e^{-z_{i}}}{\partial z_{j}}\right\vert =\left\vert
diag\left( -e^{-z_{1}},...,-e^{-z_{n}}\right) \right\vert  \\
&=&e^{-\sum_{i=1}^{n}z_{i}}.
\end{eqnarray*}%

\noindent Then , the \textit{pdf} of $Z^{\ast}$ is  
\begin{eqnarray*}
f_{Z^{\ast }}\left( z_{1},...,z_{n}\right)  &=&f_{U}\left(
-e^{-z_{1}},...,-e^{-z_{n}}\right) e^{-\sum_{i=1}^{n}z_{i}}\boldsymbol{1}%
_{\left( 0\leq z_{1}\leq ...\leq z_{n}\right) } \\
&=&n!e^{-\sum_{i=1}^{n}z_{i}}\boldsymbol{1}_{\left( 0\leq z_{1}\leq ...\leq
z_{n}\right) }.
\end{eqnarray*}

\noindent This \textit{pdf} is that of $\left( E_{1,n},...,E_{n,n}\right) $ by (\ref{cv.m05}). The proof ends here.\\

\bigskip \noindent \textbf{Another version}. Let us give another version of the previous result. It is \ clear that for any \
standard uniform random variable $U,$ we have \ $U=^{d}1-U.$ Then for any $%
n\geq 1$, we have
\begin{equation*}
\left( U_{1,n},...,U_{n,n}\right) =^{d}\left(
1-U_{1,n},...,1-U_{n-i+1,n},...,1-U_{n,n}\right) .
\end{equation*}

\noindent The equality in distribution in Proposition \ref{cv.R.malmquist03} becomes :
for any $n\geq 1$, 
$$
\left( -\log \left( 1-U_{n,n}\right) ,...,-\log \left(1-U_{1,n}\right) \right) =^{d}\left( E_{1,n},...,E_{n,n}\right)
$$ 

\bigskip \noindent Let us go further and denote 
\begin{equation*}
\alpha _{i,n}=-\log \left( 1-U_{i,n}\right) ,1\leq i\leq n.
\end{equation*}

\bigskip \noindent Consider the transformation for $n\geq 1,$%
\begin{equation*}
\left( 
\begin{array}{l}
n\alpha _{1,n} \\ 
\left( n-1\right) \left( \alpha _{2,n}-\alpha _{1,n}\right)  \\ 
. \\ 
. \\ 
\left( n-i+1\right) \left( \alpha _{i,n}-\alpha _{i-1,n}\right)  \\ 
. \\ 
. \\ 
1\left( \alpha _{n,n}-\alpha _{n-1,n}\right) 
\end{array}%
\right) =\left( 
\begin{array}{l}
V_{1} \\ 
V_{2} \\ 
. \\ 
. \\ 
V_{i} \\ 
. \\ 
\\ 
V_{n}%
\end{array}%
\right) .
\end{equation*}

\noindent We have 
\begin{equation*}
\left( 
\begin{array}{l}
\alpha _{1,n} \\ 
\alpha _{2,n} \\ 
. \\ 
. \\ 
. \\ 
. \\ 
\\ 
\alpha _{n,n}%
\end{array}%
\right) =\left( 
\begin{array}{l}
V_{1}/n \\ 
V_{1}/n+V_{2}/\left( n-1\right)  \\ 
. \\ 
. \\ 
. \\ 
. \\ 
\\ 
V_{1}/n+V_{2}/\left( n-1\right) +...+V_{n-1}/2+V_{1}/1%
\end{array}%
\right). 
\end{equation*}

\noindent The probability density function of $\left( V_{1},...,V_{n}\right) $ is
given \ by 
\begin{eqnarray*}
&&f_{V}\left( v_{1},...,v_{n}\right)\\
 &=&f_{(\alpha _{1,n},...,\alpha
_{n,n})}\left( v_{1}/n,v_{1}/n+v2/(n-1),...,v_{1}/n+v_{2}/\left( n-1\right)
+...+v_{n}\right)  \\
&&\times \left\vert J\left( v\right) \right\vert \times \boldsymbol{1}%
_{D_{V}}\left( v\right) .
\end{eqnarray*}

\noindent The Jacobian determinant in absolute value of this transform is
\begin{equation*}
\left\vert J\left( v\right) \right\vert =\frac{1}{n!}
\end{equation*}

\noindent and the domain of $V$ is 
\begin{equation*}
D_{V}=\mathbb{R}_{+}^{n}.
\end{equation*}

\noindent We conclude by using (\ref{cv.m05}) which gives the joint pdf of $(\alpha
_{1,n},...,\alpha _{n,n})$, and by denoting $s_{i}=v_{1}/n+v2/(n-1)+...+v_{i}/(n-i+1),$ $i=1,...,n$. We get
 
\begin{eqnarray*}
f_{V}\left( v_{1},...,v_{n}\right)  &=&\frac{1}{n!}\times
n!e^{-\sum_{i=1}^{n}s_{i}}\boldsymbol{1}_{\left( v_{1}\geq 0,..,v_{n}\geq
0\right) } \\
&=&e^{-\sum_{i=1}^{n}s_{i}}\boldsymbol{1}_{\left( v_{1}\geq 0,..,v_{n}\geq
0\right) }.
\end{eqnarray*}

\noindent We may check that $s_{1}+...+s_{n}=v_{1}+...+v_{n}.$ We arrive at 
\begin{equation*}
f_{V}\left( v_{1},...,v_{n}\right) =\prod\limits_{i=1}^{n}e^{-v_{i}}%
\boldsymbol{1}_{\left( v_{i}\geq 0\right) }.
\end{equation*}

\noindent This says that $\left( V_{1},...,V_{n}\right) $ has independent standard
exponential coordinates. We summarize our finding in :

\begin{proposition} \label{cv.R.malmquist04} Let $\alpha _{i,n}=-\log \left( 1-U_{i,n}\right) ,$ 
$i=1,...,n.$ Then the random \ variables 
$$
n\alpha _{1,n}, \ \left( n-1\right)\left( \alpha _{2,n}-\alpha _{1,n}\right),\ldots, \left( n-i+1\right) \left(\alpha _{i,n}-\alpha _{i-1,n}\right), \dots, \left( \alpha _{n,n}-\alpha_{n-1,n}\right)
$$ 

\noindent are independent standard exponential random variables.
\end{proposition}

\bigskip \noindent Let us do more and put for any $1\leq i\leq n$, 
\begin{equation*}
\left( n-i+1\right) \left( \alpha _{i,n}-\alpha _{i-1,n}\right) =\left(
n-i+1\right) \log \left( 
\frac{1-U_{i-1,n}}{1-U_{i,n}}\right) .
\end{equation*}

\noindent By our previous results we have that the random variables 
\begin{equation*}
E_{n-i+1}^{\ast }=\left( n-i+1\right) \left( \alpha _{i,n}-\alpha
_{i-1,n}\right) =\log \left( \frac{1-U_{n-i,n}}{1-U_{n-i+1,n}}\right)
^{\left( n-i+1\right) }
\end{equation*}

\noindent are independent and standard exponential random variables. We may and do change  $U_{n-i,n}$ to
$U_{i+1,n}$ to arrive at this celebrated representation.

\begin{proposition} \label{cv.R.malmquist05} (Malmquist representation). Let \ $U_{1},U_{2},...,U_{n}$ be standard uniform random variables for $n\geq 1$.
Let  $0\leq U_{1,n}<U_{2,n}<,...,<U_{n,n}\leq 1$ be their associated order
statistics. Then the \ random variables 
\begin{equation*}
\log \left( \frac{U_{i+1,n}}{U_{i,n}}\right) ^{i},i=1,...,n
\end{equation*}

\noindent are independent standard exponential random variables.
\end{proposition}

 

%% file: asymptotics_cv_04_en.tex
\chapter[Function empirical process]{The functional Empirical Process As a General
Tools in Asymptotic Statistics} \label{cv.empTool}
\section{Using the small o's and the big O's}

In this chapter, we will show how to combine all the concepts we have studied
so far to get yet simple but powerful tools that may be systematically
used to find asymptotic normal laws in a great variety of problems, even in
current research problems. We will first study the manipulations of the $o_{%
\mathbb{P}}$ and the $O_{\mathbb{P}}$ symbols concerning limits in
probability. Next, we present the functional empirical process that is used
here only in the frame of the finite distributions case. Then, we will give
some cases as illustrations.\\

It is important to notice for once that the methods given are valid for sequences of random variables and limit random variables defined on the same probability space. In consequence, we treat sequences of random variables $(X_{n})_{n\geq 1},$ $(Y_{n})_{n\geq 1},(Z_{n})_{n\geq 1}$,... defined on the
same probability space $(\Omega, \mathcal{A},\mathbb{P})$ with values in $\mathbb{R}^{k}$,$k\geq 1$ and $%
(a_{n})_{n\geq 1},$ $(b_{n})_{n\geq 1},(c_{n})_{n\geq 1}$ are positive
random numbers.\\

\section{Stochastic $o$'s and $O$'s} \label{cv.empTool.sec1}
\noindent \textbf{I - Big O's and small o's almost surely}.\\

\noindent \textbf{DEFINITIONS}.\newline

\noindent \textbf{(a)} The sequence of real random variables \noindent $%
(X_{n})_{n\geq 1}$ is said to be an $o$ (read the name of the letter o) of $%
a_{n}$ almost surely as n$\rightarrow +\infty$, denoted by 
\begin{equation*}
X_{n}=o(a_{n}),a.s.\text{ as }n\rightarrow +\infty ,
\end{equation*}

\noindent  if and onl if 
\begin{equation}
\lim_{n\rightarrow +\infty }X_{n}/a_{n}=0\text{ }a.s.  \label{eqso}
\end{equation}

\bigskip \noindent \textbf{(b)} The sequence \noindent of real random variables $\
(X_{n})_{n\geq 1}$ is said to be a big $O$ of $a_{n}$ almost surely as n$%
\rightarrow +\infty ,$ denoted by 
\begin{equation*}
X_{n}=O(a_{n}),a.s.\text{ as }n\rightarrow +\infty ,
\end{equation*}

\noindent  if and only if the sequence $\{\left\vert X_{n}\right\vert
/a_{n},n\geq 1\}$ is almost surely bounded, that is

\begin{equation}
\lim_{n\rightarrow +\infty }\sup \left\vert X_{n}\right\vert /a_{n}<+\infty ,%
\text{ }a.s.  \label{eqbO}
\end{equation}

\bigskip \noindent \textbf{BE CAREFUL}. The equality signs used in (\ref{eqso})
and (\ref{eqbO}) are to be read in one direction only in the sense : the
left member is a small $o$ of $a_{n}$\ or a big $O$ of $a_{n}.$ Do not reverse
the equality from left to right. For example, if $X_{n}$ is an $o(n),$ it is
also an $o(n^{2})$ and we may write $o(n)=o(n^{2})$ $a.s.$ but you cannot
write $o(n^{2})=o(n)$ $a.s.$ An other example : $X_{n}=n^{3/2}$ is an $o(n^{2})$
but is not an $o(n).$ This remark will extend to the notations of small $o$'s
and big $O$'s in probability to be defined below.\newline

\bigskip \noindent \textbf{Particular cases concerning the constants}. If $%
a_{n}=C>0$ for any n$\geq 1, $ denoted $a_{n}\equiv C,$ we have :\newline

\bigskip \noindent \textbf{(i)} $X_{n}=O(C)$ a.s. if and only if $X_{n}/C$
is bounded $a.s.$ if and only if $X_{n}$ is bounded $a.s.$ and we write 
\begin{equation*}
X_{n}=O(1)\text{ }a.s.
\end{equation*}

\bigskip \noindent \textbf{(ii)} $X_{n}=o(C)$ $a.s.$ if and only if $X_{n}/C$
$\rightarrow 0$ $a.s.$ if and only if $X_{n}/C\rightarrow 0$ $a.s.\ $and we
write 
\begin{equation*}
X_{n}=o(1)\text{ }a.s.
\end{equation*}

\bigskip \noindent \textbf{(iii)} For any constant $C>0$, we may write $%
C=O(1)$.\\

\bigskip \noindent \textbf{PROPERTIES}.\newline

\noindent The properties are very numerous and the user has often
to check new ones depending on his undergoing work. But a few of them must be known
and ready to be used. Let us list them in three groups.\newline

\bigskip \noindent \textbf{Group A}. Properties of small o's.\newline

\noindent \textbf{(1)} $o(a_{n})o(b_{n})=o(a_{n}b_{n})$ $a.s.$\newline

\noindent \textbf{(2)} (1) $o(o(a_{n}))=o(a_{n})$ $a.s.$\newline

\noindent \textbf{(3)} If $b_{n}\geq a_{n}$ for all $n\geq
1,o(a_{n})=o(b_{n})$ $a.s.$\newline

\noindent \textbf{(4)} $o(a_{n})+o(a_{n})=o(a_{n})$ $a.s.$\newline

\noindent \textbf{(5)} $o(a_{n})+o(b_{n})=o(a_{n}+b_{n})$ $a.s.$ and $%
o(a_{n})+o(b_{n})=o(a_{n}\vee b_{n})$ $a.s.$ where $a_{n}\vee b_{n}=\max
(a_{n},b_{n}).$\newline

\noindent \textbf{(6)} $o(a_{n})=a_{n}o(1)$ $a.s.$ and $a_{n}o(1)=o(a_{n})$ $%
a.s.$\newline

\bigskip \noindent \textbf{PROOFS}. Each of these properties is quickly
proved in :\newline

\noindent \textbf{(1)} If \bigskip $X_{n}=o(a_{n})$ and $Y_{n}=o(b_{n}),$
then 
\begin{equation*}
\lim_{n\rightarrow +\infty }\frac{\left\vert X_{n}Y_{n}\right\vert }{%
a_{n}b_{n}}=\lim_{n\rightarrow +\infty }\frac{\left\vert X_{n}\right\vert }{%
a_{n}}\times \lim_{n\rightarrow +\infty }\frac{\left\vert Y_{n}\right\vert }{%
b_{n}}=0\text{ }a.s
\end{equation*}

\noindent and then $X_{n}Y_{n}=o(a_{n}b_{n})$ $a.s.$\newline

\noindent \textbf{(2)} If $Y_{n}=o(a_{n}),$ $a.s.$ and $X_{n}=o(Y_{n}),$ $%
a.s.,$%
\begin{equation*}
\lim_{n\rightarrow +\infty }\frac{\left\vert X_{n}\right\vert }{a_{n}}%
=\lim_{n\rightarrow +\infty }\left\vert \frac{X_{n}}{Y_{n}}\right\vert
\times \frac{\left\vert Y_{n}\right\vert }{a_{n}}=\lim_{n\rightarrow +\infty
}\left\vert \frac{X_{n}}{Y_{n}}\right\vert \times \lim_{n\rightarrow +\infty
}\frac{\left\vert Y_{n}\right\vert }{a_{n}}=0,
\end{equation*}

\noindent that is $X_{n}=o(a_{n})$ $a.s.$\newline

\noindent \textbf{(3)} If \bigskip $X_{n}=o(a_{n})$ and $b_{n}\geq a_{n}$
for all n$\geq 1,$ then 
\begin{equation*}
0\leq \lim_{n\rightarrow +\infty }\sup \frac{\left\vert X_{n}\right\vert }{%
b_{n}}=\lim_{n\rightarrow +\infty }\sup \frac{\left\vert X_{n}\right\vert }{%
a_{n}}\frac{a_{n}}{b_{n}}\leq \lim_{n\rightarrow +\infty }\sup \frac{%
\left\vert X_{n}\right\vert }{a_{n}}=0\text{ }a.s
\end{equation*}

\noindent and $\left\vert X_{n}/b_{n}\right\vert \rightarrow 0$ $a.s.,$ that
is $X_{n}=o(b_{n}),$ $a.s.$\newline

\noindent \textbf{(4)} If \bigskip $X_{n}=o(a_{n})$ and $Y_{n}=o(a_{n}),$
then%
\begin{equation*}
\lim_{n\rightarrow +\infty }\sup \frac{\left\vert X_{n}+Y_{n}\right\vert }{%
a_{n}}\leq \lim_{n\rightarrow +\infty }\frac{\left\vert X_{n}\right\vert }{%
a_{n}}+\lim_{n\rightarrow +\infty }\frac{\left\vert Y_{n}\right\vert }{b_{n}}%
=0\text{ }a.s
\end{equation*}

\noindent and then $X_{n}+Y_{n}=o(a_{n})$ $a.s.$\newline

\noindent \textbf{(5)} To prove that $o(a_{n})+o(b_{n})=o(a_{n}+b_{n})$ $a.s.
$, use Point (3) to see that $o(a_{n})=o(a_{n}+b_{n})$ $a.s.$ since $%
a_{n}+b_{n}\geq a_{n}$ for all $n\geq 1$ and as well $\
o(b_{n})=o(a_{n}+b_{n})$ $a.s.$ and then use Point (3) to conclude. We
prove that $o(a_{n})+o(b_{n})=o(a_{n}\vee b_{n})$ $a.s.$ in the very same
manner.\newline

\noindent \textbf{(6)} This is a simple rephrasing of the definition.\newline

\bigskip \noindent \textbf{Group B}. Properties of big $O$'s.\newline

\noindent \textbf{(1)} $O(a_{n})O(b_{n})=O(a_{n}b_{n})$ $a.s.$\newline

\noindent \textbf{(2)} $O(O(a_{n}))=O(a_{n})$ $a.s.$\newline

\noindent \textbf{(3)} If $b_{n}\geq a_{n}$ for all $n\geq
1,O(a_{n})=O(b_{n})$ $a.s.$\newline

\noindent \textbf{(4)} $O(a_{n})+O(a_{n})=O(a_{n})$ $a.s.$\newline

\noindent \textbf{(5)} $O(a_{n})+O(b_{n})=O(a_{n}+b_{n})$ $a.s.$ and $%
O(a_{n})+O(b_{n})=O(a_{n}\vee b_{n})$ $a.s.$, where $a_{n}\vee b_{n}=\max
(a_{n},b_{n}).$\newline

\noindent \textbf{(6)} $O(a_{n})=a_{n}O(1)$ $a.s.,$ and $a_{n}O(1)=O(a_{n})$ 
$a.s.$\newline

\bigskip \noindent \textbf{PROOFS}. These properties are proved exactly as
those of \textbf{Group A}, where superior limits are used at the place of limits.\newline

\noindent \textbf{Group C}. Properties of combinations of small o's and big
O's.\newline

\noindent \textbf{(1)} $o(a_{n})O(b_{n})=o(a_{n}b_{n})$ $a.s.$\newline

\noindent \textbf{(2)} $o(O(a_{n}))=o(a_{n})$ $a.s.$ and $%
O(o(a_{n}))=o(a_{n}),$ $a.s.$\newline

\noindent \textbf{(3a)} If $a_{n}=O(b_{n}),$ $a.s.,$ then $%
o(a_{n})+O(b_{n})=O(b_{n})$ $a.s.$\newline

\noindent \textbf{(3b)} If $b_{n}=O(a_{n}),$ $a.s.,$ then $%
o(a_{n})+O(b_{n})=O(a_{n})$ $a.s.$\newline

\noindent \textbf{(3c)} If $b_{n}=o(a_{n}),$ $a.s.,$ then $%
o(a_{n})+O(b_{n})=o(b_{n})$ $a.s.$\newline

\noindent \textbf{(4)} $(1+o(a_{n}))^{-1}-1=O(a_{n})$, $a.s.$\newline

\bigskip \noindent \textbf{PROOFS}.\newline

\noindent \textbf{(1)} If $X_{n}=o(a_{n})$ and $Y_{n}=O(b_{n}),$ then%
\begin{equation*}
\lim_{n\rightarrow +\infty }\sup \frac{\left\vert Y_{n}\right\vert }{b_{n}}%
=C<+\infty \text{ }a.s.
\end{equation*}

\noindent and 
\begin{eqnarray*}
\lim_{n\rightarrow +\infty }\sup \frac{\left\vert X_{n}Y_{n}\right\vert }{%
a_{n}b_{n}} &=&\lim_{n\rightarrow +\infty }\left( \left\vert \frac{X_{n}}{%
a_{n}}\right\vert \times \frac{\left\vert Y_{n}\right\vert }{b_{n}}\right)
=\lim_{n\rightarrow +\infty }\sup \left\vert \frac{X_{n}}{a_{n}}\right\vert
\times \lim_{n\rightarrow +\infty }\sup \frac{\left\vert Y_{n}\right\vert }{%
b_{n}} \\
&\leq &C\lim_{n\rightarrow +\infty }\sup \left\vert \frac{X_{n}}{a_{n}}%
\right\vert =0\text{ }a.s.
\end{eqnarray*}

\noindent \textbf{(2)} Use Points (6) of Groups A and B to say 
\begin{equation*}
o(O(a_{n}))=o(1)\times O(a_{n})=a_{n}\times o(1)\times O(1)=a_{n}\times
o(1)=o(a_{n})
\end{equation*}

\noindent and 
\begin{equation*}
O(o(a_{n}))=o(a_{n})O(1)=a_{n}\times o(1)\times O(1)=a_{n}\times
o(1)=o(a_{n}).
\end{equation*}

\bigskip
\noindent \textbf{(3a-b-c)} These three points are proved in similar ways.
Let us give the details of (3b) for example. Let $X_{n}=o(a_{n})$ and $%
Y_{n}=O(b_{n})$ and $b_{n}=O(a_{n}).$ Then 
\begin{eqnarray*}
o(a_{n})+O(b_{n}) &=&o(a_{n})+O(O(a_{n}))=o(a_{n})+O(a_{n}) \\
&=&a_{n}(o(1)+O(1))=a_{n}\times O(1)=O(a_{n}).
\end{eqnarray*}

\noindent \textbf{(4)} We have 
\begin{eqnarray*}
(1+o(a_{n}))^{-1}-1 &=&\frac{o(a_{n})}{1+o(a_{n})}=o(a_{n})O(1) \\
&=&a_{n}o(1)O(1)=a_{n}o(1)=o(a_{n}).
\end{eqnarray*}

\bigskip \noindent \textbf{II - Big $O$'s and small $o$'s in probability}.\\

\noindent \textbf{DEFINITIONS}.\newline

\noindent \textbf{(a)} The sequence of real random variables \noindent $%
(X_{n})_{n\geq 1}$ is said to be an o (read the name of the letter o) of $%
a_{n}$ in probability as n$\rightarrow +\infty ,$ denoted by%
\begin{equation*}
X_{n}=o_{\mathbb{P}}(a_{n}),\text{ as }n\rightarrow +\infty ,
\end{equation*}

\bigskip \noindent if and only if
\begin{equation}
\lim_{n\rightarrow +\infty }X_{n}/a_{n}=0\text{ in probability},
\label{eqso}
\end{equation}

\noindent that is for any $\lambda>0$ 
\begin{equation*}
\lim_{n\rightarrow +\infty }\mathbb{P}(\left\vert X_{n}\right\vert >\lambda
a_{n})=0.
\end{equation*}

\noindent \textbf{(b)} The sequence \noindent of real random variables $\
(X_{n})_{n\geq 1}$ is said to be a big O of $a_{n}$ in probability as n$%
\rightarrow +\infty ,$ denoted by

\begin{equation*}
X_{n}=O_{\mathbb{P}}(a_{n}),\text{ as }n\rightarrow +\infty ,
\end{equation*}

\bigskip \noindent if and only if the sequence $\{\left\vert X_{n}\right\vert
/a_{n},n\geq 1\}$ is bounded in probability, that is \ : For any $%
\varepsilon >0,$ there exists a constant $\lambda >0$, such that%
\begin{equation}
\inf_{n\geq 1}\mathbb{P}(\left\vert X_{n}\right\vert \leq \lambda a_{n})\geq
1-\varepsilon,  \label{big-O1}
\end{equation}

\noindent which is equivalent to

\begin{equation}
\liminf_{\lambda \uparrow +\infty }\limsup_{n\rightarrow +\infty }%
\mathbb{P}(\left\vert X_{n}\right\vert >\lambda a_{n})=0.  \label{big-O2}
\end{equation}

\noindent Before we go further, let us prove the following Lemma.

\begin{lemma}
\label{oO.01} Each of (\ref{big-O1}) and (\ref{big-O2}) is equivalent to :
For any $\varepsilon >0,$ there exists an integer $N\geq 1$ a constant $%
\lambda >0$, such that
\begin{equation}
\inf_{n\geq N}\mathbb{P}(\left\vert X_{n}\right\vert \leq \lambda a_{n})\geq
1-\varepsilon .  \label{big-O3}
\end{equation}
\end{lemma}

\bigskip \noindent \textbf{PROOF}. To prove that (\ref{big-O1}) and (\ref%
{big-O3}) are equivalent, it will be enough to show that (\ref{big-O3}) $\Longrightarrow $
(\ref{big-O1}) since the reverse implication is obvious. Suppose that (\ref%
{big-O3}), that is, for $\varepsilon >0,$ there exist $N\geq 1$ and a real
number $\lambda _{0}>0$ such that

\begin{equation}
\forall (n\geq N),\text{ }\mathbb{P}(\left\vert X_{n}/a_{n}\right\vert \leq
\lambda _{0})\geq 1-\varepsilon .
\end{equation}

\bigskip
\noindent If $N=1,$ then (\ref{big-O1})\ holds. If not, we have for $j\in
\{1,...,N-1\}$ fixed$,$ $(\left\vert X_{j}/a_{j}\right\vert \leq \lambda
)\uparrow \Omega $ as $\lambda \uparrow +\infty .$ So by the Monotone
Convergence Theorem, there exists for each $j\in \{1,...,N-1\}$ a real
number $\lambda _{j}>0$ such that $\mathbb{P}(\left\vert
X_{j}/a_{j}\right\vert \leq \lambda _{j})>1-\varepsilon$. We take $\lambda
=\max (\lambda _{0},\lambda _{1},...,\lambda _{N-1})$ and get
\begin{equation*}
\forall (n\geq 1),\text{ }\mathbb{P}(\left\vert X_{n}/a_{n}\right\vert \leq
\lambda )\geq 1-\varepsilon ,
\end{equation*}

\bigskip
\noindent which is (\ref{big-O1}). Now, let us prove that (\ref{big-O2})$%
\Longleftrightarrow $(\ref{big-O3}). First (\ref{big-O2}) means%
\begin{equation*}
\lim_{\lambda \uparrow +\infty }\limsup_{n\rightarrow +\infty }\mathbb{P}%
(\left\vert X_{n}\right\vert >\lambda a_{n})=0,
\end{equation*}

\noindent since $\limsup_{n\rightarrow +\infty }\mathbb{P}(\left\vert
X_{n}\right\vert >\lambda a_{n})$ in non-increasing as $\lambda \uparrow
+\infty $ on $[0,1]$.\\

\bigskip \noindent We get for any $\varepsilon >0,$ there exists a real
number $\lambda >0$ such that 
\begin{equation*}
\limsup_{n\rightarrow +\infty }\mathbb{P}(\left\vert X_{n}\right\vert
>\lambda a_{n})=\lim_{N\uparrow +\infty }\sup_{n\geq N}\mathbb{P}(\left\vert
X_{n}\right\vert >\lambda a_{n})\leq \varepsilon /2.
\end{equation*}

\noindent Then for some $N>0$, 
\begin{equation*}
\sup_{n\geq N}\mathbb{P}(\left\vert X_{n}\right\vert >\lambda a_{n})\leq
\varepsilon ,
\end{equation*}

\noindent that is 
\begin{equation*}
\inf_{n\geq N}\mathbb{P}(\left\vert X_{n}\right\vert \leq \lambda a_{n})\geq
1-\varepsilon ,
\end{equation*}

\noindent which is (\ref{big-O3}). Now, a rephrasing of this gives :
for any $\varepsilon >0,$ there exists $N_{0}>0$ and a real number $\lambda
_{0}>0$ such that

\begin{equation}
\inf_{n\geq N}\mathbb{P}(\left\vert X_{n}\right\vert \leq \lambda
_{0}a_{n})\geq 1-\varepsilon ,
\end{equation}

\noindent that is

\begin{equation*}
\sup_{n\geq N}\mathbb{P}(\left\vert X_{n}\right\vert >\lambda
_{0}a_{n})<\varepsilon ,
\end{equation*}

\noindent which leads to 
\begin{equation*}
\inf_{N\geq 1}\sup_{n\geq N}\mathbb{P}(\left\vert X_{n}\right\vert >\lambda
_{0}a_{n})<\varepsilon ,
\end{equation*}

\noindent and next 
\begin{equation*}
\inf_{\lambda >0}\inf_{N\geq 1}\sup_{n\geq N}\mathbb{P}(\left\vert
X_{n}\right\vert >\lambda a_{n})<\varepsilon ,
\end{equation*}

\noindent which is (\ref{big-O2}).\newline

\bigskip 

\bigskip \noindent \textbf{COMMENTS, NOTATIONS AND SOMME LEMMAS}.\newline

\noindent \textbf{(a)} From Chapter \ref{cv.tensRk}, an $O_{\mathbb{P}}(1)$
is simply a tight sequence of random variables. From Theorem \ref{tensTheo2}
of Chapter \ref{cv.tensRk}, we have that any sequence $X_{n}=O_{\mathbb{P}%
}(a_{n})$ contains a sub-sequence $(X_{n_{k}})_{k\geq 1}$ such that $%
(X_{n_{k}}/a_{n_{k}})_{k\geq 1}$ weakly converges in $\mathbb{R}$.\newline

\noindent \textbf{(b)} It may be convenient to rephrase (\ref{big-O1}) into
the following sentence.\\

\noindent For any $\varepsilon >0$ there exists a real number 
$\lambda >0$ such that $\left\vert X_{n}\right\vert \leq \lambda a_{n}$ With
Probability At Least Equal to $1-\varepsilon$ for all $n\geq 1$.\\

\noindent By using the complementary events, we will say : for any $\varepsilon >0$ there exists a real number 
$\lambda >0$ such that $\left\vert X_{n}\right\vert > \lambda a_{n}$ With Probability At Most Equal to $\varepsilon$ for all $n\geq 1$.\\

\noindent With Probability At Least Equal to $1-\varepsilon$ will be abbreviated by $WPALE(1-\varepsilon)$.\\

\noindent As well, $WPAME(\varepsilon)$  is an abbreviation of With Probability At Most Equal to $\varepsilon$.\\

\noindent For lengthy demonstrations, using these types of sentences described above may be handy.\newline

\bigskip \noindent We will need two other lemmas.

\begin{lemma} \label{oO-02} We have the following properties :

\noindent (a) If  $X_{n}$ is a sequence of $k-$random vectors weakly converging (say, to a $k$-random vector $X$), then $\left\Vert X_{n}\right\Vert =O_{\mathbb{P}}(1)$.\\

\noindent (b) Let $X_{n}$ be a sequence of random vectors with values in the metric
space $(S,d)$ converging in probability to a constant $C\in S$ and let $g$ be a measurable mapping from $(S,d)$ to another metric space $(E,r)$. If $g$ is continuous at $C$, then $g(X_{n})$ converges in probability to $C$.\\

\bigskip \noindent (c) Consider a sequence of $k-$random vectors $%
(X_{n})_{n\geq 1}$ converging to zero in probability. Let $R(x)$ be a real
function of $x\in \mathbb{R}^{k}$ continuous at zero and such that  $R(0)=0$. Let $p>0$ be a fixed integer. If $R(x)=o(\Vert
x\Vert ^{p})$ as $x\rightarrow 0$, then $R(X_{n})=o_{\mathbb{P%
}}(\Vert X_{n}\Vert ^{p}).$ If $R(x)=O(\Vert x\Vert ^{p})$ as $x\rightarrow
0\ $, then $R(X_{n})=O_{\mathbb{P}}(\Vert X_{n}\Vert ^{p})$.
\end{lemma}

\bigskip \noindent \textbf{Proof}.\\

\noindent \textbf{Proof of Point (a)}. If $X_{n}\rightarrow _{\mathbb{P}}X$, then by Proposition \ref%
{cv.CvCp.01}, $X_{n}\rightarrow _{w}X$ and by the continuous mapping Theorem %
\ref{cv.mappingTh} of Chapter \ref{cv},  $\left\Vert X_{n}\right\Vert
\rightarrow _{w}\left\Vert X\right\Vert .$ Then by Theorem \ref%
{cv.theo.portmanteau.rk}, we have for any continuity point of $F_{\left\Vert
X\right\Vert }(\lambda )=P(\left\Vert X\right\Vert \leq \lambda ),$

\begin{equation*}
\lim_{n\rightarrow +\infty} \mathbb{P}\left( \left\vert X_{n}\right\vert >\lambda \right) =\mathbb{P%
}\left( \left\Vert X\right\Vert >\lambda \right) =F_{\left\Vert X\right\Vert
}(\lambda ).
\end{equation*}

\bigskip
\noindent Since the set of discontinuity points of $F_{\left\Vert X\right\Vert }$ is
at most countable (see Point 6 of Chapter \ref{cv.R}, Section \ref{cv.sec2}), apply the formula above for $\lambda
\rightarrow +\infty $ while $\lambda $ are continuity points. Since $%
1-F_{\left\Vert X\right\Vert }(\lambda )\rightarrow 0$ as $\lambda
\rightarrow +\infty ,$ then for any $\varepsilon >0,$ we are able to pick
one value of $\lambda (\varepsilon )$ which is a continuity point of $%
F_{\left\Vert X\right\Vert }$ satisfying $1-F_{\left\Vert X\right\Vert
}(\lambda )<\varepsilon .$ For any $\varepsilon >0,$ we have found $\lambda
(\varepsilon )>0$ such that

\begin{equation*}
\limsup_{n\rightarrow +\infty }\mathbb{P}\left( \left\vert
X_{n}\right\vert >\lambda (\varepsilon )\right) \leq \varepsilon ,
\end{equation*}

\noindent which implies
\begin{equation*}
\lim_{\lambda \rightarrow +\infty }\limsup_{n\rightarrow +\infty }\mathbb{P%
}\left( \left\vert X_{n}\right\vert >\lambda \right) =0.
\end{equation*}

\noindent Point (a) is proved.\\

\noindent \textbf{Proof of Point (b)}. Assume the notations of this point and suppose that $g$ is continous at $%
C$. Let $\varepsilon >0.$ By the continuity of $g$ at $C$, there exists $%
\eta >0$ such that

\begin{equation*}
d(x,C)<\eta \Longrightarrow r(g(x),g(C))<\varepsilon /2.
\end{equation*}%

\noindent Now

\begin{eqnarray*}
\mathbb{P}(r(g(X_{n}),g(C) &>&\varepsilon )=\mathbb{P}(\left\{
r(g(X_{n}),g(C)>\varepsilon \right\} \cap \left\{ d(X_{n},C)\geq \eta
\right\} )\\
&+&\mathbb{P}(\left\{ r(g(X_{n}),g(C)>\varepsilon \right\} \cap
\left\{ d(X_{n},C)<\eta \right\} ) \\
&\leq &\mathbb{P}(d(X_{n},C)\geq \eta ),
\end{eqnarray*}

\noindent since $(\left\{ r(g(X_{n}),g(C)>\lambda \right\} \cap \left\{
d(X_{n},C)<\eta \right\} )\subset $ $(\left\{ r(g(X_{n}),g(C)>\varepsilon
\right\} \cap \left\{ d(X_{n},C)<\eta \right\} \cap \left\{ \left\{
r(g(X_{n}),g(C)<\varepsilon /2\right\} \right\} )=\emptyset$.\\

\noindent Then, since $X_{n}\rightarrow _{\mathbb{P}}C$, we have 

\begin{equation*}
\limsup_{n\rightarrow +\infty }\mathbb{P}(r(g(X_{n}),g(C)>\varepsilon
)\leq \limsup_{n\rightarrow +\infty }\mathbb{P}(d(X_{n},C)\geq \eta )=0.
\end{equation*}

\noindent So Point (b) is true.\\

\noindent \textbf{Proof of Point (c-1)}. Let $R(x)=o(\Vert x\Vert ^{p})$ as $x\rightarrow 0$. Then 
$$
g(x)=\left\vert R(x)/\Vert x\Vert ^{p}\right\vert \rightarrow 0
$$

\noindent as $x\rightarrow 0$. This proves that $g(0)=0$, and that $g$ is continuous at zero. By continuity of $g$ at zero and by Point (a), 
we have, as $n\rightarrow +\infty$, that $g(X_{n})=\left\vert R(X_{n})\right\vert /\Vert X_{n}\Vert ^{p}\rightarrow _{%
\mathbb{P}}0$  whenever $X_{n}\rightarrow _{\mathbb{P}}0$.\\

\noindent Hence $R(X_{n})=o_{\mathbb{P}}(\Vert X_{n}\Vert ^{p})$.\\

\noindent \textbf{Proof of Point (c-2)}. Let $R(x)=O(\Vert x\Vert ^{p})$ as $x\rightarrow 0.$ Then for any $%
\varepsilon >0,$\ there exist $\eta >0$ and $C>0$ such  $\left\vert
R(x)\right\vert /\Vert x\Vert ^{p}\leq C$ for all $\left\Vert x\right\Vert
<\eta .$ Then for $\lambda >C,$

\begin{eqnarray*}
\mathbb{P}(\left\vert R(X_{n})\right\vert /\Vert X_{n}\Vert ^{p} &>&\lambda
)=\mathbb{P}(\left\{ \left\vert R(X_{n})\right\vert /\Vert X_{n}\Vert
^{p}>\lambda \right\} \cap \left\{ \left\Vert X_{n}\right\Vert \geq \eta
\right\} )\\
&+&\mathbb{P}(\left\{ \left\vert R(X_{n})\right\vert /\Vert
X_{n}\Vert ^{p}>\lambda \right\} \cap \left\{ \left\Vert X_{n}\right\Vert
<\eta \right\} ) \\
&\leq &\mathbb{P}(\left\Vert X_{n}\right\Vert \geq \eta ),
\end{eqnarray*}

\noindent since $(\left\{ \left\vert R(X_{n})\right\vert /\Vert X_{n}\Vert
^{p}>\lambda \right\} \cap \left\{ \left\Vert X_{n}\right\Vert <\eta
\right\} )\subset (\left\{ \left\vert R(X_{n})\right\vert /\Vert X_{n}\Vert
^{p}>\lambda \right\} \cap \left\{ \left\Vert X_{n}\right\Vert <\eta
\right\} \cap \left\{ \left\{ \left\vert R(X_{n})\right\vert /\Vert
X_{n}\Vert ^{p}<C\right\} \right\} )=\emptyset .$ Then for all $\lambda >C,$

\begin{equation*}
\limsup_{n\rightarrow +\infty }\mathbb{P}(\left\vert R(X_{n})\right\vert
/\Vert X_{n}\Vert ^{p}>\lambda )=0
\end{equation*}

\noindent and then

\begin{equation*}
\lim_{\lambda \rightarrow +\infty }\limsup_{n\rightarrow +\infty }\mathbb{P%
}(\left\vert R(X_{n})\right\vert /\Vert X_{n}\Vert ^{p}>\lambda )=0.
\end{equation*}

\noindent We conclude that $\left\vert R(X_{n})\right\vert /\Vert X_{n}\Vert ^{p}=O_{%
\mathbb{P}}(1)$.\\

\bigskip In some situations, we would be able to work with convergence in
probability while we are not sure of measurability of some sequences. For
example, using the Mean Value Theorem with real random sequences $X_{n}$ and 
$Z_{n}$ and real function $g$ of class $C^{1}$ may lead to this kind of
formula \begin{equation*}
g(X_{n})-g(Y_{n})=(X_{n}-Y_{n})g^{\prime }(Z_{n}),
\end{equation*}

\noindent with $\min (X_{n},Y_{n})\leq Z_{n}\leq \max (X_{n},Y_{n}).$ Here, we know
that $g^{\prime }(Z_{n})$ is measurable but we do not know if $Z_{n}$ is. In
such a situation, we may need the notion of outer probability.

\begin{definition} \label{cv.outerprob} Let $(u_{n})_{n\geq 1}$ be a sequence of real-valued applications defined on $%
\Omega .$ It has a measurable covering sequence if and only if there exists
a sequence of non-negative real random variables $(v_{n})_{n\geq 1}$\ 
defined on $(\Omega ,\mathcal{A})$ such that
\begin{equation*}
\forall (n\geq 1),u_{n}\leq v_{n}.
\end{equation*}

\noindent Next, $(u_{n})_{n\geq 1}$ converges in outer probability to a real-valued application $u$ defined
on $\Omega $, as $n\rightarrow +\infty ,$ if and only the sequence $(u_{n}-u)_{n\geq 1}$ has a measurable covering sequence $(v_{n})_{n\geq 1}$ which converges to
zero in probability, and we denote

\begin{equation*}
u_{n}\rightarrow _{\mathbb{P}^{\ast }}u\text{ as }n\rightarrow +\infty. 
\end{equation*}
\end{definition}

\bigskip \noindent  We are going to see that the result of Point (b) of the lemma above
still holds for convergence in outer probability in the special case of $\mathbb{R}$.

\begin{lemma} \label{oO-03} Let  $X_{n}$ be a sequence of real-valued applications defined on $%
\Omega $ converging in outer probability to $c\in \mathbb{R}$. Let $g$ be a
real-valued function defined on $\mathbb{R},$ continuous at $c$ and such
that for each $n\geq 1,$  $g(X_{n})$ is measurable. Then $g(X_{n})$
converges in probability to $g(c)$.
\end{lemma}

\bigskip \textbf{Proof}.  Assume the notations of the lemma. Let $Y_{n}$ be a sequence of
random variables such that $\left\vert X_{n}-c\right\vert \leq Y_{n}$ for
all $n\geq 1$ and $Y_{n}\rightarrow 0$ in probability.\\

\noindent Now, by the continuity of $g$ at $c$, there exists $\eta >0$ such that%
\begin{equation*}
\left\vert x-c\right\vert <\eta \Longrightarrow \left\vert
g(x)-g(c)\right\vert <\varepsilon /2.
\end{equation*}

\noindent Next 
\begin{eqnarray*}
\mathbb{P}(\left\vert g(X_{n})-g(c)\right\vert  >\varepsilon )&=&\mathbb{P}%
(\left\{ r(g(X_{n}),g(C)>\varepsilon \right\} \cap \left\{ Y_{n}\geq \eta
\right\} )\\
&+&\mathbb{P}(\left\{ \left\vert g(X_{n})-g(c)\right\vert>\varepsilon \right\} \cap \left\{ Y_{n}<\eta \right\} ) \\
&\leq &\mathbb{P}(\left\{ Y_{n}\geq \eta \right\} )
\end{eqnarray*}

\noindent since $\left\{ \left\vert g(X_{n})-g(c)\right\vert >\varepsilon \right\}
\cap \left\{ Y_{n}<\eta \right\} =\emptyset .$ The reason on this is that on 
$\left\{ \left\vert g(X_{n})-g(c)\right\vert >\varepsilon \right\} \cap
\left\{ Y_{n}<\eta \right\} ,$ we have $\left\vert X_{n}-c\right\vert \leq
Y_{n}<\eta $ and then $\left\vert g(X_{n})-g(c)\right\vert <\varepsilon /2.$
This is impossible. Next,  since $Y_{n}\rightarrow _{\mathbb{P}}0,$ we have 
\begin{equation*}
\lim \sup_{n\rightarrow +\infty }\mathbb{P}(\left\vert
g(X_{n})-g(c)\right\vert >\varepsilon )\leq \lim \sup_{n\rightarrow +\infty }%
\mathbb{P}(Y_{n}\geq \eta )=0.
\end{equation*}

\noindent The proof is complete.\\

\bigskip
\bigskip 
\noindent Now, we may give some important properties of the small $o^{\prime
}s$ and the big $O^{\prime }s$ in probability.\newline

\noindent \textbf{MAIN PROPERTIES}.\\ 

\noindent  \textbf{(1)} If $X_{n}=o(1)$ $a.s.,$ then $X_{n}=o_{\mathbb{P}}(1)$.\newline

\noindent \textbf{(2)} $o_{\mathbb{P}}(a_{n})=a_{n}o_{\mathbb{P}}(1)$ and $%
a_{n}o_{\mathbb{P}}(1)=o_{\mathbb{P}}(a_{n})$.\newline

\noindent \textbf{(3)} $o_{\mathbb{P}}(a_{n})o_{\mathbb{P}}(b_{n})=o_{%
\mathbb{P}}(a_{n}b_{n})$.\newline

\noindent \textbf{(4)} $o_{\mathbb{P}}(o_{\mathbb{P}}(a_{n}))=o_{\mathbb{P}%
}(a_{n})$.\newline

\noindent \textbf{(5)} If $b_{n}\geq a_{n}$ for all $n\geq 1,o_{\mathbb{P}%
}(a_{n})=o_{\mathbb{P}}(b_{n})$.\newline

\noindent \textbf{(6)} $o_{\mathbb{P}}(a_{n})+o_{\mathbb{P}}(a_{n})=o_{%
\mathbb{P}}(a_{n})$.\newline

\noindent \textbf{(7)} $o_{\mathbb{P}}(a_{n})+o_{\mathbb{P}}(b_{n})=o_{%
\mathbb{P} }(a_{n}+b_{n})$ and $o_{\mathbb{P}}(a_{n})+o_{\mathbb{P}%
}(b_{n})=o_{\mathbb{P}}(a_{n}\vee b_{n})$, where $a_{n}\vee b_{n}=\max
(a_{n},b_{n})$.\newline

\noindent \textbf{(8)} If $X_{n}=O(1)$ $a.s.,$ then $X_{n}=O_{\mathbb{P}}(1)$%
.\newline

\noindent \textbf{(9)} $O_{\mathbb{P}}(a_{n})=a_{n}O_{\mathbb{P}}(1)$.\\

\noindent \textbf{(10)} $O_{\mathbb{P}}(a_{n})O_{\mathbb{P}}(b_{n})=O_{%
\mathbb{P}}(a_{n}b_{n})$.\newline

\noindent \textbf{(11)} $O_{\mathbb{P}}(O_{\mathbb{P}}(a_{n}))=O_{\mathbb{P}%
}(a_{n})$.\newline

\noindent \textbf{(12)} If $b_{n}\geq a_{n}$ for all $n\geq 1,O_{\mathbb{P}%
}(a_{n})=O_{\mathbb{P}}(b_{n})$.\newline

\noindent \textbf{(13)} $O_{\mathbb{P}}(a_{n})+O_{\mathbb{P}}(a_{n})=O_{%
\mathbb{P}}(a_{n})$.\newline

\noindent \textbf{(14)} $O_{\mathbb{P}}(a_{n})+O_{\mathbb{P}}(b_{n})=O_{%
\mathbb{P}}(a_{n}+b_{n})$ and $O_{\mathbb{P}}(a_{n})+O_{\mathbb{P}%
}(b_{n})=O_{\mathbb{P}}(a_{n}\vee b_{n})$ where $a_{n}\vee b_{n}=\max
(a_{n},b_{n})$.\newline

\noindent \textbf{(15)} $o_{\mathbb{P}}(a_{n})O_{\mathbb{P}}(b_{n})=o_{%
\mathbb{P}}(a_{n}b_{n})$.\newline

\noindent \textbf{(16)} $o_{\mathbb{P}}(O_{\mathbb{P}}(a_{n}))=o_{\mathbb{P}%
}(a_{n})$ and $O_{\mathbb{P}}(o_{\mathbb{P}}(a_{n}))=o_{\mathbb{P}}(a_{n})$.%
\newline

\noindent \textbf{(17a)} If $a_{n}=O_{\mathbb{P}}(b_{n})$, then $o(a_{n})+O_{%
\mathbb{P}}(b_{n})=O_{\mathbb{P}}(b_{n})$.\newline

\noindent \textbf{(17b)} If $b_{n}=O_{\mathbb{P}}(a_{n})$, then $o_{\mathbb{P%
}}(a_{n})+O_{\mathbb{P}}(b_{n})=O_{\mathbb{P}}(a_{n})$.\newline

\noindent \textbf{(17c)} If $b_{n}=o_{\mathbb{P}}(a_{n}),$ $a.s.$, then $%
o(a_{n})+O_{\mathbb{P}}(b_{n})=o(b_{n})$.\newline

\noindent \textbf{(18)} $(1+o_{\mathbb{P}}(a_{n}))^{-1}-1=o_{\mathbb{P}%
}(a_{n})$.\\

\noindent \textbf{(19)} An $o_{\mathbb{P}}(1)$ is an $O_{\mathbb{P}}(1)$.

\bigskip \noindent \textbf{PROOFS}.\\

\bigskip \noindent (1) This derived from the implication : $X_{n}\rightarrow 0$ a.s. $%
\Longrightarrow X_{n}\rightarrow _{P}0$ (See Proposition \ref{cv.CvCp.01}).\\

\bigskip \noindent (2) We have :\\

\noindent $X_{n}=o_{\mathbb{P}}(a_{n})\Longleftrightarrow \left\vert
X_{n}/a_{n}\right\vert \longrightarrow _{P}0\Longleftrightarrow
X_{n}/a_{n}=o_{\mathbb{P}}(1)\Longleftrightarrow X_{n}=a_{n}o_{\mathbb{P}%
}(1).$\\

\bigskip \noindent (3) By Point (2) above, $o_{\mathbb{P}}(a_{n})o_{\mathbb{P}%
}(b_{n})=a_{n}b_{n}\times o_{\mathbb{P}}(1)o_{\mathbb{P}}(1)=a_{n}b_{n}%
\times o_{\mathbb{P}}(1)=o_{\mathbb{P}}(a_{n}b_{n})$ (Check $o_{\mathbb{P}%
}(1)o_{\mathbb{P}}(1)=o_{\mathbb{P}}(1)$ in Property (A1) in the Appendix
subsection below).\\

\bigskip \noindent (4) $o_{\mathbb{P}}(o_{\mathbb{P}}(a_{n}))=o_{\mathbb{P}}(a_{n})o_{\mathbb{P}%
}(1)=a_{n}\times o_{\mathbb{P}}(1)o_{\mathbb{P}}(1)=a_{n}\times o_{\mathbb{P}%
}(1)=o_{\mathbb{P}}(a_{n})$ (Again, use Property (A1) of the Appendix subsection
below).\\

\bigskip \noindent (5) Let $b_{n}\geq a_{n}$ for all $n\geq 1,X_{n}=o_{\mathbb{P}}(a_{n}).$ For
any $\eta >0,0\leq \lim_{n\rightarrow +\infty }\sup P(\left\vert
X_{n}/b_{n}\right\vert >\eta )\leq \lim_{n\rightarrow +\infty }P(\left\vert
X_{n}/a_{n}\right\vert >\eta )=0.$\\

\bigskip \noindent (6) Let  $X_{n}=o_{\mathbb{P}}(a_{n})$ and $Y_{n}=o_{\mathbb{P}}(a_{n}).$
Use the classical stuff, for $\eta >0,$%
\begin{equation*}
\left( \frac{\left\vert X_{n}\right\vert }{a_{n}}>\eta /2\right) \cap \left( 
\frac{\left\vert Y_{n}\right\vert }{a_{n}}>\eta /2\right) \subset \left( 
\frac{\left\vert X_{n}+Y_{n}\right\vert }{a_{n}}>\eta \right) .
\end{equation*}

\noindent Then for $\eta >0,$%
\begin{eqnarray}
\limsup_{n\rightarrow +\infty }\mathbb{P}\left( \frac{\left\vert
X_{n}+Y_{n}\right\vert }{a_{n}}>\eta \right) &\leq& \limsup_{n\rightarrow
+\infty }\mathbb{P}\left( \frac{\left\vert X_{n}\right\vert }{a_{n}}>\eta
/2\right) \notag \\
&+&\limsup_{n\rightarrow +\infty }\mathbb{P}\left( \frac{\left\vert Y_{n}\right\vert }{a_{n}}>\eta /2\right) =0.  \label{decompSum}
\end{eqnarray}

\bigskip \noindent (7) To prove this point, combine Points (5) and (6) above.\\

\bigskip \noindent (8) $X_{n}=O(1)$ $a.s.$ as $n\rightarrow +\infty $ means there exists $%
\Omega _{0}$ measurable such that $\mathbb{P}(\Omega _{0})=1$ and for any $%
\omega \in \Omega _{0},$ 
\begin{equation*}
\limsup_{n\rightarrow +\infty }\left\vert X_{n}(\omega )\right\vert
=\inf_{n\geq 1}\sup_{p\geq n}\left\vert X_{p}\right\vert =M(\omega )<+\infty.
\end{equation*}

\noindent We have for all $n\geq 1,$%
\begin{equation*}
\mathbb{P}\left( \left\vert X_{n}\right\vert >\lambda \right) \leq \mathbb{P}%
\left( \sup_{p\geq n}\left\vert X_{p}\right\vert >\lambda \right). 
\end{equation*}

\noindent We have $Y_{n}=\sup_{p\geq n}\left\vert X_{p}\right\vert \downarrow M$ $a.s.$
Then  $Y_{n}1_{\Omega _{0}}\rightarrow _{\mathbb{P}}M1_{\Omega _{0}}$. We
are dealing with real random variables and we may apply the weak convergence
results. We get $Y_{n}1_{\Omega _{0}}\rightarrow _{w}M1_{\Omega _{0}}$ by
Proposition \ref{cv.CvCp.01}$.$ By Theorem \ref{cv.theo.portmanteau.rk}, we
have for any continuity point of $F_{M}(\lambda )=P(M1_{\Omega _{0}}\leq
\lambda )$. Use the Monotone Convergence Theorem to get

\begin{eqnarray*}
\limsup_{n\rightarrow +\infty }\mathbb{P}\left( \left\vert
X_{n}\right\vert >\lambda \right)  &=&\limsup_{n\rightarrow +\infty }%
\mathbb{P}\left( \sup_{p\geq n}\left\vert X_{p}\right\vert 1_{\Omega
_{0}}>\lambda \right)  \\
&=&\limsup_{n\rightarrow +\infty }\mathbb{P}\left( \sup_{p\geq
n}\left\vert X_{p}\right\vert 1_{\Omega _{0}}>\lambda \right)\\
& =&\mathbb{P}\left( M1_{\Omega _{0}}>\lambda \right)\\
&=&1-F_{M}(\lambda).
\end{eqnarray*}

\bigskip \noindent Since the set of discontinuity points of $F_{M}$ is at most countable (see Point 6 of Chapter \ref{cv.R}), we  apply the formula above for 
$\lambda \rightarrow +\infty$ while $\lambda $ are continuity points. Since $1-F_{M}(\lambda )\rightarrow 0$ as $\lambda \rightarrow +\infty ,$ then for any $\varepsilon >0,$ we are able to pick one value of $\lambda (\varepsilon )$ which is a continuity point of $F_{M}$ satisfying $1-F_{M}(\lambda )<\varepsilon$. For any $\varepsilon >0$, we have found $\lambda (\varepsilon )>0$ such that

\begin{equation*}
\limsup_{n\rightarrow +\infty }\mathbb{P}\left( \left\vert
X_{n}\right\vert >\lambda (\varepsilon )\right) \leq \varepsilon ,
\end{equation*}

\noindent which implies
\begin{equation*}
\lim_{\lambda \rightarrow +\infty }\limsup_{n\rightarrow +\infty }\mathbb{P%
}\left( \left\vert X_{n}\right\vert >\lambda \right) =0.
\end{equation*}

\noindent Then $X_{n}=O_{\mathbb{P}}(1).$\\

\bigskip \noindent (9) Let $X_{n}=O_{\mathbb{P}}(a_{n}).$ Then

\begin{equation*}
\lim_{\lambda \rightarrow +\infty }\limsup_{n\rightarrow +\infty }\mathbb{P%
}\left( \left\vert X_{n}/a_{n}\right\vert >\lambda \right) =0.
\end{equation*}

\noindent This is the definition that $X_{n}/a_{n}=O_{\mathbb{P}}(1)$ and then $%
X_{n}=a_{n}O_{\mathbb{P}}(1).$\\

\bigskip \noindent (10) $O_{\mathbb{P}}(a_{n})O_{\mathbb{P}}(b_{n})=a_{n}b_{n}O_{\mathbb{P}%
}(1)O_{\mathbb{P}}(1)=a_{n}b_{n}O_{\mathbb{P}}(1)$ (Check that $O_{\mathbb{P}}(1)O_{\mathbb{P}}(1)=O_{\mathbb{P}}(1)$ in Property (A2) in the Appendix subsection).\\

\bigskip \noindent \textbf{(11)} $O_{\mathbb{P}}(O_{\mathbb{P}}(a_{n}))=O_{\mathbb{P}%
}(a_{n})O_{\mathbb{P}}(1)=a_{n}O_{\mathbb{P}}(1)O_{\mathbb{P}}(1)=O_{\mathbb{%
P}}(a_{n})$.\\

\bigskip \noindent \textbf{(12)} Let $b_{n}\geq a_{n}$ for all $n\geq 1$ and $%
X_{n}=O_{\mathbb{P}}(a_{n})$. Then%
\begin{equation*}
\lim_{\lambda \rightarrow +\infty }\limsup_{n\rightarrow +\infty }\mathbb{P%
}\left( \left\vert X_{n}/b_{n}\right\vert >\lambda \right) \leq
\lim_{\lambda \rightarrow +\infty }\limsup_{n\rightarrow +\infty }\mathbb{P%
}\left( \left\vert X_{n}/a_{n}\right\vert >\lambda \right) =0.
\end{equation*}

\noindent Then $X_{n}=O_{\mathbb{P}}(b_{n})$.\\

\bigskip \noindent \textbf{(13)} Let $X_{n}=O_{\mathbb{P}}(a_{n})$ and $X_{n}=O_{%
\mathbb{P}}(a_{n}).$ Use the same technique as in Formula \ref{decompSum}\
below to get  \newline

\begin{eqnarray}
\lim_{\lambda \rightarrow +\infty }\limsup_{n\rightarrow +\infty }\mathbb{P%
}\left( \frac{\left\vert X_{n}+Y_{n}\right\vert }{a_{n}}>\lambda \right)
&\leq& \lim_{\lambda \rightarrow +\infty }\limsup_{n\rightarrow +\infty }%
\mathbb{P}\left( \frac{\left\vert X_{n}\right\vert }{a_{n}}>\lambda
/2\right) \notag\\
&+&\lim_{\lambda \rightarrow +\infty }\limsup_{n\rightarrow
+\infty }\mathbb{P}\left( \frac{\left\vert Y_{n}\right\vert }{a_{n}}>\lambda
/2\right) =0.
\end{eqnarray}

\bigskip \noindent  \textbf{(14)} Combine Points (12) and (13) to get this one.\\

\bigskip \noindent  \textbf{(15)} $o_{\mathbb{P}}(a_{n})O_{\mathbb{P}%
}(b_{n})=a_{n}b_{n}o_{\mathbb{P}}(1)O_{\mathbb{P}}(1)=a_{n}b_{n}o_{\mathbb{P}%
}(1)=o_{\mathbb{P}}(a_{n}b_{n})$. (Check that $o_{\mathbb{P}}(1)O_{\mathbb{P}%
}(1)$ in Property (A3) in the Appendix subsection below)$.$\newline

\bigskip \noindent  \textbf{(16)} $o_{\mathbb{P}}(O_{\mathbb{P}}(a_{n}))=O_{\mathbb{P}%
}(a_{n})o_{\mathbb{P}}(1)=a_{n}O_{\mathbb{P}}(1)o_{\mathbb{P}}(1)=a_{n}o_{%
\mathbb{P}}(1)=o_{\mathbb{P}}(a_{n})$ and $O_{\mathbb{P}}(o_{\mathbb{P}%
}(a_{n}))=o_{\mathbb{P}}(a_{n})O_{\mathbb{P}}(1)=o_{\mathbb{P}}(a_{n})$.\\

\bigskip \noindent  \textbf{(17a)} Let $a_{n}=O(b_{n})$ and $X_{n}=o_{\mathbb{P}%
}(a_{n})$ and $Y_{n}=O_{\mathbb{P}}(b_{n}).$ There exists $C>0$ such that $%
a_{n}\leq Cb_{n}$ for any $n\geq 1.$ Then $X_{n}=o_{\mathbb{P}}(a_{n})=o_{%
\mathbb{P}}(Cb_{n})$ by Point (5). But obviously $X_{n}=o_{\mathbb{P}%
}(Cb_{n})=o_{\mathbb{P}}(b_{n})$ and then $X_{n}=O_{\mathbb{P}}(b_{n})$ by
Point (19) below. Finally $X_{n}+Y_{n}=O_{\mathbb{P}}(b_{n})+O_{%
\mathbb{P}}(b_{n})=O_{\mathbb{P}}(b_{n})$.\\

\bigskip \noindent \textbf{(17b)} Let $b_{n}=O(a_{n})$, $X_{n}=o_{\mathbb{P}}(a_{n})$ and $Y_{n}=O_{\mathbb{P}}(b_{n}).$ We exchange the roles of $a_{n}$
and $b_{n}$ to get $b_{n}\leq Ca_{n}$ and $Y_{n}=O_{\mathbb{P}}(Ca_{n})=O_{%
\mathbb{P}}(a_{n})$ by Point (1) and finally,
$$
X_{n}+Y_{n}=o_{\mathbb{P}%
}(a_{n})+O_{\mathbb{P}}(a_{n})=O_{\mathbb{P}}(a_{n})+O_{\mathbb{P}}(a_{n})=O_{\mathbb{P}}(a_{n}).
$$

\bigskip \noindent  \textbf{(17c)} Let  $b_{n}=o_{\mathbb{P}}(a_{n})$, $X_{n}=o_{%
\mathbb{P}}(a_{n})$ and $Y_{n}=O_{\mathbb{P}}(b_{n}).$ We have $%
X_{n}+Y_{n}=o_{\mathbb{P}}(a_{n})+O_{\mathbb{P}}(o_{\mathbb{P}}(a_{n}))=o_{%
\mathbb{P}}(a_{n})+o_{\mathbb{P}}(a_{n})$ by Point (16). Finally $%
X_{n}+Y_{n}=o_{\mathbb{P}}(a_{n}).$\\

\bigskip \noindent \textbf{(18)} We have
\begin{equation*}
(1+o_{\mathbb{P}}(a_{n}))^{-1}-1=\frac{o_{\mathbb{P}}(a_{n})}{1+o_{\mathbb{P}%
}(a_{n})}.
\end{equation*}

\noindent By Point (b) of Lemma \ref{oO-02}, $(1+o_{\mathbb{P}}(a_{n}))^{-1}\rightarrow _{P}1$ and by
Point (a) of the same lemma, $(1+o_{\mathbb{P}}(a_{n}))^{-1}=O_{\mathbb{P}%
}(1).$ Then 
\begin{equation*}
(1+o_{\mathbb{P}}(a_{n}))^{-1}-1=O_{\mathbb{P}}(1)o_{\mathbb{P}}(a_{n})=o_{%
\mathbb{P}}(a_{n}),
\end{equation*}

\noindent by Point (15).\\

\bigskip \noindent $(19)$ By Lemma \ref{oO-02}, an $o_{\mathbb{P}}(1)$ converges to $0$ in
probability and then is an $O_{\mathbb{P}}(1)$.\\

\subsection{Extensions} $ $\\

\noindent The concepts of \textit{small o's} and \textit{big O's} are extended to $\mathbb{R}^{k}$ in the following
way :\\

\noindent \textbf{(a)} The sequence of random vectors $(X_{n})_{n\geq 1}$ of $\mathbb{R}^{k},
$ is an $o(a_{n})$ $a.s.$ if and only if $\left\Vert X_{n}\right\Vert
/a_{n}=o(1)$ $a.s.$, and is an $o_{\mathbb{P}}(a_{n})$ if and only if $%
\left\Vert X_{n}\right\Vert /a_{n}=o_{\mathbb{P}}(1)$.\\

\noindent \textbf{(b)} The sequence of random vectors $(X_{n})_{n\geq 1}$  of $\mathbb{R}^{k},
$ is an $O(a_{n})$ $a.s.$ if and only if $\left\Vert X_{n}\right\Vert
/a_{n}=O(1)$ $a.s.$, and is an $O_{\mathbb{P}}(a_{n})$ if and only if $%
\left\Vert X_{n}\right\Vert /a_{n}=O_{\mathbb{P}}(1)$.\\

\noindent From there, handling these concepts is easy by combining their properties in 
$\mathbb{R}$ and those of the norms in $\mathbb{R}^{k}.$

\subsection{Balanced sequences}  $ $\\

\noindent It may help in some cases to have sequences $X_{n}$ such that both $\left\Vert
X_{n}\right\Vert $ and $1/\left\Vert X_{n}\right\Vert $ are bounded in
probability. Let us give some notations for real sequences.\\

\noindent \textbf{(1)} For $0\leq a<b<+\infty ,$ we denote by $X_{n}=O_{\mathbb{P}%
}(a,b,a_{n},b_{n})$ the property that for any $\varepsilon >0,$ there exists 
$\lambda >0$ such that we have $(a+\lambda \leq \left\vert X_{n}\right\vert
/a_{n},\left\vert X_{n}\right\vert /a_{n}\leq b-\lambda )$ \textit{WPALE}$%
(1-\varepsilon ),$ for large values of $n.$ If $a_{n}=b_{n}$ for all $n\geq
1,$ we simply write $X_{n}=O_{\mathbb{P}}(a,b,a_{n}).$\\

\noindent \textbf{(1)} For $0\leq a,$ we denote by $X_{n}=O_{\mathbb{P}}(a,+\infty ,a_{n},b_{n})
$ the property that for ant $\varepsilon >0,$ there exists $\lambda >0$ such
that we have 
$$
(a+\lambda ^{-1}\leq \left\vert X_{n}\right\vert /a_{n},\left\vert X_{n}\right\vert /a_{n}\leq \lambda )
$$ 

\noindent \textit{WPALE}$(1-\varepsilon )$, for large values of $n$. If $a_{n}=b_{n}$ for all $n\geq
1,$ we simply write $X_{n}=O_{\mathbb{P}}(a,+\infty ,a_{n})$.\\

\noindent An example of a sequence of random variables satisfying $X_{n}=O_{\mathbb{P}}(0,+\infty ,1)$ is a sequence $X_{n}$
weakly converging to $X>0$ $a.s.$ In this case $1/X_{n}\rightsquigarrow 1/X$
finite $a.s.$ and then $X_{n}=O_{\mathbb{P}}(1)$ and $1/X_{n}=O_{\mathbb{P}%
}(1).$ Combining these two points leads to $X_{n}=O_{\mathbb{P}}(0,+\infty
,1).$

\subsection{Appendix}  $ $\\

\noindent \textbf{(A1)} If $X_{n}\rightarrow _{\mathbb{P}}a\in \mathbb{R}$ and $%
X_{n}\rightarrow _{\mathbb{P}}b\in \mathbb{R},$ then $X_{n}Y_{n}\rightarrow
_{\mathbb{P}}ab.$\\

\noindent \textbf{Proof}. We have $(\eta +\left\vert b\right\vert )\eta +\left\vert
a\right\vert \eta \rightarrow 0$ $\ as$ $\eta \rightarrow 0.$ For any $%
\varepsilon >0,$ for any $\delta >0,$ choose a value of $\eta >0$ such that $%
(\eta +\left\vert b\right\vert )\eta +\left\vert a\right\vert \eta <\delta .$
We apply the definition of the convergences $X_{n}\rightarrow _{\mathbb{P}}a$
and $X_{n}\rightarrow _{\mathbb{P}}b$ to get a value $N_{0}\geq 1$ such that
for any $n\geq N_{0},$

\begin{equation*}
\mathbb{P}(\left\vert X_{n}-a\right\vert \geq \eta )\leq \varepsilon /2\text{
and }\mathbb{P}(\left\vert Y_{n}-b\right\vert \geq \eta )\leq \varepsilon /2.
\end{equation*}

\noindent But
\begin{eqnarray*}
\left\vert X_{n}Y_{n}-ab\right\vert  &=&\left\vert
X_{n}Y_{n}-aY_{n}+aY_{n}-ab\right\vert  \\
&\leq &\left\vert Y_{n}\right\vert \left\vert X_{n}-a\right\vert +\left\vert
a\right\vert \left\vert Y_{n}-b\right\vert  \\
&\leq &(\left\vert Y_{n}-b\right\vert +\left\vert b\right\vert )\text{ }%
\left\vert X_{n}-a\right\vert +\left\vert a\right\vert \left\vert
Y_{n}-b\right\vert 
\end{eqnarray*}

\noindent On $(\left\vert X_{n}-a\right\vert \geq \eta )^{c}\cap (\left\vert
Y_{n}-b\right\vert \geq \eta )^{c},$

\begin{equation*}
\left\vert X_{n}Y_{n}-ab\right\vert \leq (\eta +\left\vert b\right\vert
)\eta +\left\vert a\right\vert \eta \leq \delta .
\end{equation*}%

\noindent Then for $n\geq N_{0},$%
\begin{equation*}
(\left\vert X_{n}-a\right\vert \geq \eta )^{c}\cap (\left\vert
Y_{n}-b\right\vert \geq \eta )^{c}\subset (\left\vert
X_{n}Y_{n}-ab\right\vert \leq \delta ),
\end{equation*}

\noindent that is
\begin{equation*}
\mathbb{P}(\left\vert X_{n}-a\right\vert \geq \eta )^{c}\cap (\left\vert
Y_{n}-b\right\vert \geq \eta )^{c}\leq \mathbb{P}(\left\vert
X_{n}Y_{n}-ab\right\vert \leq \delta ),
\end{equation*}

\bigskip \noindent and by taking complements,
\begin{equation*}
\mathbb{P}(\left\vert X_{n}Y_{n}-ab\right\vert >\delta )\leq \mathbb{P}%
((\left\vert X_{n}-a\right\vert \geq \eta )\cup (\left\vert
Y_{n}-b\right\vert \geq \eta ))\leq \varepsilon /2+\varepsilon
/2=\varepsilon .
\end{equation*}

\noindent Thus,
\begin{equation*}
X_{n}Y_{n}\rightarrow _{\mathbb{P}}ab.
\end{equation*}

\bigskip \noindent \textbf{Property (A2)}. If $X_{n}=O_{\mathbb{P}}(1)$ and $X_{n}=O_{\mathbb{P}}(1)$ then $%
X_{n}Y_{n}=O_{\mathbb{P}}(1).$\\

\noindent \textbf{Proof}. By applying the definition of an $O_{\mathbb{P}}(1)$, we may find for any $%
\varepsilon >0,$ two integer numbers $N_{1}$ and $N_{2}$ and two positive
numbers $\lambda _{1}>0$ and $\lambda _{2}>0$ such that

\begin{equation*}
\forall (n\geq N_{1}),\mathbb{P(}\left\vert X_{n}\right\vert \leq \lambda
_{1})\geq 1-\varepsilon /2\text{ and }\forall (n\geq N_{2}),\mathbb{P(}%
\left\vert Y_{n}\right\vert \leq \lambda _{2})\geq 1-\varepsilon /2.
\end{equation*}

\bigskip
\noindent For $n\geq \max (N_{1,}N_{2}),$ 
\begin{equation*}
\mathbb{(}\left\vert X_{n}\right\vert \leq \lambda _{1})\cap \mathbb{(}%
\left\vert Y_{n}\right\vert \leq \lambda _{2})\subset \mathbb{(}\left\vert
X_{n}Y_{n}\right\vert \leq \lambda _{1}\lambda _{2}),
\end{equation*}

\noindent which is equivalent to
\begin{equation*}
\mathbb{(}\left\vert X_{n}Y_{n}\right\vert >\lambda _{1}\lambda _{2})\subset 
\mathbb{(}\left\vert X_{n}\right\vert >\lambda _{1})\cup \mathbb{(}%
\left\vert Y_{n}\right\vert >\lambda _{2}),
\end{equation*}

\noindent which implies for $n\geq \max (N_{1,}N_{2}),$%
\begin{equation*}
\mathbb{P(}\left\vert X_{n}Y_{n}\right\vert >\lambda _{1}\lambda _{2})\leq 
\mathbb{P(}\left\vert X_{n}\right\vert >\lambda _{1})+\mathbb{P(}\left\vert
Y_{n}\right\vert >\lambda _{2})\leq \varepsilon /2+\varepsilon
/2=\varepsilon.
\end{equation*}

\bigskip
\noindent Thus, for any $\varepsilon >0,$ there exists a non-negative $N$ $(=\max
(N_{1,}N_{2})),$ there exists $\lambda >0$ $(=\lambda _{1}\lambda _{2})$ and
for any $n\geq N,$
\begin{equation*}
\mathbb{P(}\left\vert X_{n}Y_{n}\right\vert \leq \lambda )\geq 1-\varepsilon.
\end{equation*}

\noindent Hence $X_{n}Y_{n}=O_{\mathbb{P}}(1)$.\\

\begin{equation*}
X_{n}Y_{n}\rightarrow _{\mathbb{P}}ab.
\end{equation*}

\bigskip \noindent \textbf{Property (A3)}. If $X_{n}=o_{\mathbb{P}}(1)$ and $Y_{n}=O_{\mathbb{P}}(1)$ then $%
X_{n}Y_{n}=o_{\mathbb{P}}(1).$\\

\noindent \textbf{Proof}. Fix $\varepsilon >0.$ By applying the definition of an $O_{\mathbb{P}}(1)$
there exist an integer number $N_{1}$ and a positive number $\lambda >0$
such that  
\begin{equation*}
\mathbb{P(}\left\vert Y_{n}\right\vert \leq \lambda )\geq 1-\varepsilon /2.
\end{equation*}

\noindent Now let $\eta >0.$ Let us apply the definition of $X_{n}=o_{\mathbb{P}}(1)$
to get that there exists a positive integer $N_{2}$ such that

\begin{equation*}
\forall (n\geq N_{2}),\text{ }\mathbb{P(}\left\vert X_{n}\right\vert >\eta
/\lambda )\leq \varepsilon /2.
\end{equation*}

\noindent Thus for  $n\geq \max (N_{1,}N_{2}),$ 
\begin{equation*}
\mathbb{(}\left\vert X_{n}\right\vert \leq \eta /\lambda )\cap \mathbb{(}%
\left\vert Y_{n}\right\vert \leq \lambda )\subset \mathbb{(}\left\vert
X_{n}Y_{n}\right\vert \leq \eta )
\end{equation*}

\noindent which is equivalent to
\begin{equation*}
\mathbb{(}\left\vert X_{n}Y_{n}\right\vert >\eta )\subset \mathbb{(}%
\left\vert X_{n}\right\vert >\eta /\lambda )\cup \mathbb{(}\left\vert
Y_{n}\right\vert >\lambda ),
\end{equation*}

\noindent which implies for $n\geq \max (N_{1,}N_{2}),$%
\begin{equation*}
\mathbb{P(}\left\vert X_{n}Y_{n}\right\vert >\eta )\leq \mathbb{P}(\left\vert X_{n}\right\vert >\eta /\lambda )+\mathbb{P(}\left\vert
Y_{n}\right\vert >\lambda )\leq \varepsilon /2+\varepsilon /2=\varepsilon .
\end{equation*}

\bigskip \noindent Thus for any $\varepsilon >0,$ for any $\eta >0,$ there exists a non negative 
$N$ $(=\max (N_{1,}N_{2})),$ for any $n\geq N,$%
\begin{equation*}
\mathbb{P(}\left\vert X_{n}Y_{n}\right\vert >\eta )\leq \varepsilon .
\end{equation*}

\noindent Hence $X_{n}Y_{n}=o_{\mathbb{P}}(1)$.

\newpage

\section{Delta Methods} \label{cv.empTool.sec2}

\bigskip The Delta method is a quick way to derive new asymptotic laws for
sequences of random variables defined on the same probability measure $%
(\Omega ,\mathcal{A},\mathbb{P})$. Here, we present the univariate and multivariate
case. Here, we will see the usefulness of the results in Section \ref{cv.CvCp}
of Chapter \ref{cv} combined with the manipulations of the $o^{\prime }s$ and the $%
O^{\prime }s$ in probability we just exposed in the first section of this
chapter.\\

\noindent We begin by Delta Methods in $\mathbb{R}$.\\ 

\subsection{Univariate Version}  $ $\\

\begin{proposition} \label{delta01} Let $(X_{n})_{n\geq 1}$ be a sequence real random variables defined on the
same probability space $(\Omega ,A,\mathbb{P})$ and let $\theta $ be a real
number and $(a_{n}>0)_{n\geq 1}$ be a sequence of real numbers such that $%
a_{n}\rightarrow +\infty$ as $n\rightarrow +\infty $.\\

\noindent Let $g:D\rightarrow \mathbb{R}$ be a function of class $C^{1}$, such that $D$ is a domain of $\mathbb{R},$ $%
\theta $ is in the interior $\overset{o}{D}$ of $D,$ $\{X_{n},n\geq
1\}\subset \overset{o}{D}$.\\

\noindent  If $a_{n}(X_{n}-\theta )$ wealky converges to a random variable $Z$ as $n\rightarrow +\infty$, then 
$$
a_{n}(g(X_{n})-g(\theta ))\rightsquigarrow g^{\prime }(\theta )Z \ \ as \ \  n\rightarrow +\infty,
$$

\noindent where $\nabla g(a)=g^{\prime }(\theta )$ is the derivative of $g$ at $\theta$.
\end{proposition}

\bigskip \noindent \textbf{Proof of Proposition \ref{delta01}}. Assume that all the hypotheses of the proposition are true. By Point (a) of
Lemma \ref{oO-02}, we have $a_{n}(X_{n}-\theta )=O_{P}(1)$ and then%
\begin{equation*}
X_{n}=\theta +O_{P}(1)a_{n}^{-1}\rightarrow _{\mathbb{P}}\theta 
\end{equation*}

\noindent which by Proposition \ref{cv.CvECp} in Section \ref{cv.CvCp} of Chapter \ref%
{cv}, is equivalent to the weak convergence
\begin{equation*}
X_{n}\rightsquigarrow \theta .
\end{equation*}

\noindent Now the Mean Value Theorem implies
\begin{equation}
g(X_{n})-g(\theta )=g^{\prime }(Y_{n})(X_{n}-\theta ),  \label{appliMVT}
\end{equation}

\noindent where
\begin{equation*}
\min (X_{n},\theta )\leq Y_{n}\leq \max (X_{n},\theta ),
\end{equation*}

\noindent that is

\begin{equation*}
\left\vert Y_{n}-\theta \right\vert \leq \left\vert X_{n}-\theta \right\vert.
\end{equation*}

\bigskip
\noindent It follows that $Y_{n}\rightarrow _{\mathbb{P}}\theta $ and since $g^{\prime }$
is continuous, we have $g^{\prime }(Y_{n})\rightarrow _{\mathbb{P}}g(\theta )
$ by Point (b) of  Lemma \ref{oO-02}. Then by using  Proposition \ref{cv.CvECp} in Section \ref{cv.CvCp} of Chapter \ref{cv}, we see that is
equivalent to 
\begin{equation*}
g^{\prime }(Y_{n})\rightsquigarrow g^{\prime }(\theta ).
\end{equation*}

\noindent By the property of Slutsky given in \ \ref{cv.slutsky} in Section \ref%
{cv.CvCp} of Chapter \ref{cv}, we have the weak convergence%
\begin{equation*}
(g^{\prime }(Y_{n}),a_{n}(X_{n}-\theta ))\rightsquigarrow (g^{\prime
}(\theta ),Z)
\end{equation*}

\noindent and by the continuous mapping Theorem \ref{cv.mappingTh} in Chapter \ref{cv}
combined with (\ref{appliMVT}), we get the final conclusion 
\begin{equation*}
a_{n}(g(X_{n})-g(\theta ))=(g^{\prime }(Y_{n})\times a_{n}(X_{n}-\theta
))\rightsquigarrow g^{\prime }(\theta )Z.
\end{equation*}

\bigskip \noindent \textbf{Remark}. If we use the derivative map (total derivative)

\begin{equation*}
h\rightarrow g_{\theta }^{\prime }(h)=g^{\prime }(\theta )h,
\end{equation*}

\noindent in Proposition  \ref{delta01}, we may write the conclusion in the form
\begin{equation*}
a_{n}(g(X_{n})-g(\theta ))=g_{\theta }^{\prime }(Z).
\end{equation*}

\noindent This writing suggests we may have this kind of results in more general
spaces. Let us move to the multivariate case.

\subsection{Multivariate version}  $ $\\

\noindent  The first statement concerns the transformation of the converging sequence of $k$ components by a real function of 
$k$ arguments.

\begin{proposition} \label{delta02}
Let $(X_{n})_{n\geq 1}$ be a sequence $k$-random vectors, $k\geq 1,$ defined
on the same probability space $(\Omega ,A,\mathbb{P})$ and let $\theta \in 
\mathbb{R}^{k}$ and $(a_{n}>0)_{n\geq 1}$ be sequence of real numbers such
that $a_{n}\rightarrow +\infty$ as $n\rightarrow +\infty$.\\

\noindent Let $g:D\rightarrow \mathbb{R}$ be a function of class $C^{1}$, such that $D$ is a domain of $%
\mathbb{R}^{k},$  $\theta $ is in $\overset{o}{D}$, the interior of $D,$ $%
\{X_{n},n\geq 1\}\subset \overset{o}{D}$.\\

\noindent If $a_{n}(X_{n}-\theta )$ wealky converges to a $k-$random vector $Z$ as $n\rightarrow +\infty$, then

\begin{equation*}
a_{n}(g(X_{n})-g(\theta ))\rightsquigarrow \text{ }^{t}\nabla g(\theta
)Z=<\nabla g(\theta ),Z>\text{ as }n\rightarrow +\infty ,
\end{equation*}

\noindent where 
\begin{equation*}
^{t}\nabla g(\theta )=\left(\frac{\partial g(\theta )}{\partial \theta _{1}},...,%
\frac{\partial g(\theta )}{\partial \theta _{k}}\right)
\end{equation*}

\noindent is the gradient vector of $g$ at $\theta$.
\end{proposition}

\bigskip \noindent The second statement is the most general in the finite dimension frame, in which the converging sequence of $k$ components is transformed by a multicomponent function of $k$ arguments.

\begin{proposition} \label{delta03}
Let $(X_{n})_{n\geq 1}$ be a sequence $k$-random vectors, $k\geq 1,$ defined
on the same probability space $(\Omega ,A,\mathbb{P})$ and let $\theta $ be
a real number and $(a_{n}>0)_{n\geq 1}$ be a sequence of real numbers such
that $a_{n}\rightarrow +\infty$ as $n\rightarrow +\infty$.\\

\noindent Let $g:D\rightarrow \mathbb{R}^{m}$ be a function of class $C^{1}$, such that $\theta $ is $%
\overset{o}{D},$ the interior $\overset{o}{D}$ of $D,$ $\{X_{n},n\geq
1\}\subset \overset{o}{D}$. Denote by $g_j$, $1\leq j \leq m$, the components of the function $g$.\\

\noindent If $a_{n}(X_{n}-\theta )$ weakly converges to a $%
k-$random vector $Z$ as $n\rightarrow +\infty$, then
\begin{equation*}
a_{n}(g(X_{n})-g(\theta ))\rightsquigarrow \text{ }g_{\theta }^{\prime }Z=%
\text{ as }n\rightarrow +\infty ,
\end{equation*}

\noindent where $g_{\theta }^{\prime }$ is the matrix of partial derivatives of first
order

\begin{equation*}
g_{\theta }^{\prime }=\left( 
\begin{tabular}{lllll}
$\frac{\partial g_{1}}{\partial \theta _{1}}$ & ... & $\frac{\partial g_{1}}{%
\partial \theta _{j}}$ & .. & $\frac{\partial g_{1}}{\partial \theta _{k}}$
\\ 
...  & ... & ... & ... & ... \\ 
$\frac{\partial g_{i}}{\partial \theta _{1}}$ & ... & $\frac{\partial g_{i}}{%
\partial \theta _{j}}$ & ... & $\frac{\partial g_{j}}{\partial \theta _{k}}$
\\ 
... & ... & ... & .... & ... \\ 
$\frac{\partial g_{m}}{\partial \theta _{1}}$ & ... & $\frac{\partial g_{m}}{%
\partial \theta _{j}}$ & ... & $\frac{\partial g_{m}}{\partial \theta _{k}}$%
\end{tabular}%
\right) ,
\end{equation*}
\end{proposition}

\bigskip \noindent \textbf{Proof of Proposition \ref{delta02}}. Assume that the hypotheses of the proposition hold.\\

 \noindent Let us use the expansion of $g$ of first order at $\theta =$ $^{t}(\theta
_{1},...,\theta _{k})$ for $x=$ $(x_{1},...,x_{k})^{T}$

\begin{equation}
g(x)-g(\theta )=(x_{1}-\theta _{1})\frac{\partial g}{\partial \theta _{1}}%
(\theta )+...+(x_{1}-\theta _{k})\frac{\partial g}{\partial \theta _{k}}%
(\theta )+o(\left\Vert x-\theta \right\Vert ).  \label{develop01}
\end{equation}

 \noindent Since $a_{n}(T_{n}-\theta )=$ $a_{n}((T_{1,n},...,T_{k,n})^{T}-(\theta
_{1},...,\theta _{k})^{T})\rightsquigarrow Z=(Z_{1},...,Z_{k})^{T}$, it
follows from the continuous mapping theorem \ref{cv.mappingTh} of Chapter \ref{cv}, that for each $1\leq i\leq k,$ $%
a_{n}(T_{j,n}-\theta _{j})$ converges to $Z_{j}$ and then Point (a) of lemme \ref{oO-02}, we get
\begin{equation*}
1\leq j\leq k,(T_{j,n}-\theta _{j})=O_{\mathbb{P}}(a_{n}^{-1}).
\end{equation*}

 \noindent Thus by Points (10) and (13) of the main properties in Part II of the above
section,

\begin{equation}
\left\Vert T_{n}-\theta \right\Vert =\left\{ \sum_{j=1}^{k}(T_{j,n}-\theta
_{j})^{2}\right\} ^{1/2}=O_{\mathbb{P}}(a_{n}^{-1})=o_{\mathbb{P}}(1).
\label{develop02}
\end{equation}

\noindent The term $o(\left\Vert x-\theta \right\Vert )$ in \ref{develop01} is
continuous as a difference of two continuous functions and takes the value $0
$ for $\left\Vert x-\theta \right\Vert =0.$ By using Part (c) of Lemma \ref%
{oO-02}, a combination of (\ref{develop01}) and (\ref{develop02}) leads to
 
\begin{eqnarray*}
a_{n}(g(x)-g(\theta )) &=&a_{n}(T_{1,n}-\theta _{1})\frac{\partial g}{%
\partial \theta _{1}}(\theta )+...+a_{n}(T_{k,n}-\theta _{k})\frac{\partial g%
}{\partial \theta _{k}}(\theta )\\
&+&a_{n}o_{\mathbb{P}}(O_{\mathbb{P}}(a_{n}^{-1})) \\
&=&\text{ }^{t}\nabla g(\theta) (a_{n}(T_{n}-\theta ))+o_{\mathbb{P}}(1).
\end{eqnarray*}

 \noindent This says that $a_{n}(g(x)-g(\theta ))$ and $^{t}\nabla g(\theta
)(a_{n}(T_{n}-\theta )$ are equivalent in probability. Since $^{t}\nabla
g(\theta )(a_{n}(T_{n}-\theta )\rightsquigarrow $ $^{t}\nabla g(\theta )Z$
by the continuous mapping Theorem, we get by Proposition \ref{cvcp.prop4}\
in Section \ref{cv.CvCp} of Chapter \ref{cv}, 

\begin{equation*}
a_{n}(g(x)-g(\theta ))\rightsquigarrow ^{t}\nabla g(\theta )Z.
\end{equation*}

\bigskip \noindent \textbf{Proof of Proposition \ref{delta03}}. Assume that the hypotheses of the proposition hold.\\

 \noindent The function $g$ has $m$ components $g_{j} \in \mathbb{R}^{m}$ so that we
write $g=$$^{t}(g_{1},...,g_{m})$. Each component is of class $C^{1}$. Let us use
the conclusion of Proposition \ref{delta02} for each of these components at $\theta =$ $%
^{t}(\theta _{1},...,\theta _{k})$ for $x=(x_{1},...,x_{k})^{T}$ to get
\begin{equation}
g_{j}(x)-g_{j}(\theta )=(x_{1}-\theta _{1})\frac{\partial g_{j}}{\partial
\theta _{1}}(\theta )+...+(x_{1}-\theta _{k})\frac{\partial g_{j}}{\partial
\theta _{k}}(\theta )+o(\left\Vert x-\theta \right\Vert ).
\end{equation}

 \noindent This can be written using matrices as 
\begin{equation*}
g(x)-g(\theta )=g_{\theta }^{\prime }(x-\theta )+o^{(m)}(\left\Vert x-\theta
\right\Vert ),
\end{equation*}

 \noindent where  $o^{(m)}(\left\Vert x-\theta \right\Vert )$ is a vector of $m$
coordinates such that each of them is a continuous function which is also an $%
o(\left\Vert x-\theta \right\Vert ).$ A similar notation is also used for $o_{\mathbb{P}}(\circ).$ \ By applying the method used in Proposition \ref{delta02}, we
get 
\begin{equation*}
g(T_{n})-g(\theta )=g_{\theta }^{\prime }(T_{n}-\theta )+o_{\mathbb{P}%
}^{(m)}(a_{n}^{-1})
\end{equation*}

 \noindent and 
\begin{equation*}
a_{n}(g(T_{n})-g(\theta ))=g_{\theta }^{\prime }a_{n}(T_{n}-\theta )+o_{%
\mathbb{P}}^{(m)}(1).
\end{equation*}

 \noindent We have
\begin{equation*}
\left\Vert a_{n}(g(T_{n})-g(\theta ))-g_{\theta }^{\prime
}a_{n}(T_{n}-\theta )\right\Vert _{\mathbb{R}^{m}}=\left\Vert o_{\mathbb{P}%
}^{(m)}(1)\right\Vert _{\mathbb{R}^{m}}=o_{\mathbb{P}}(1).
\end{equation*}

 \noindent Then $a_{n}(g(T_{n})-g(\theta ))$ has the same weak limit as $g_{\theta
}^{\prime }a_{n}(T_{n}-\theta )$ which is $g_{\theta }^{\prime }Z$ by the
continuous mapping.

\newpage

\section{Using the Functional Empirical Process in Asymptotic Statistics} \label{cv.empTool.sec3}

\subsection{The Functional empirical process}  $ $\\

\noindent The functional empirical process (FEP) is a powerful tool for deriving asymptotic
limit distributions. It is similar to the multivariate delta method. But
the PEF has an advantage we describe below.\\

\noindent Given a sequence $Z_{1}$, $Z_{2}$, ..., if independent and identically distributed random variables, of common probability law 
$\mathbb{P}_{0}$, we will be able\\

\noindent \textbf{(1)} to find a Gaussian stochastic process $\mathbb{G}_{\mathbb{P}_{0}}$\\

\noindent and\\

\noindent \textbf{(2)} to express the asymptotic distributions of statistics which are functions of  $Z_{1},Z_{2},...,Z_{n}$ with respect to 
$\mathbb{G}_{\mathbb{P}_{0}}$.\\

\bigskip \noindent This allows to separately study all statistics based on $Z_{1},Z_{2},...,Z_{n}$ and, each time we want it, to get the joint 
asymptotic distributions of any finite number of them. \textbf{We say that we place the asymptotic distributions of statistics based on $Z_{1},Z_{2},...,Z_{n}$ in the Gaussian field of $G_{\mathbb{P}_{0}}$}.\\

\noindent Another interesting point is that the joint distributions we obtain by using the \textit{FEP} tool, have their covariance functions expressed in functional forms. Whatever be complicated these covariances, we do not have to worry about their form since the powerful computers of modern times are able to compute them in very short times.\\

\noindent The Delta method does not have this unified frame. Instead, each work is done for once. When we need to add or drop any statistic, we have to do the job again.\\

\noindent  Before we present the functional empirical process, we want to reassure the reader that we will only use finite distributions of the functional empirical process, that is, we remain in $\mathbb{R}^k$ and we will not use the heavy tools of functional topologies or 
Vapnick-Cervonenkis classes.\\

\noindent Let $Z_{1}$, $Z_{2}$, ... be a sequence of independent copies of a random variable $Z$ defined on the same probability space $(\Omega ,\mathcal{A},\mathbb{P})$ with values on some metric space $(S,d)$. Define for each $n\geq 1,$ the functional empirical process by 
\begin{equation*}
\mathbb{G}_{n}(f)=\frac{1}{\sqrt{n}}\sum_{i=1}^{n}(f(Z_{i})-\mathbb{E}
f(Z_{i})),
\end{equation*}

\bigskip \noindent where $f$ is a real and measurable function defined on $\mathbb{R}$ such that

\begin{equation}
\mathbb{V}_{Z}(f)=\int \left( f(x)-\mathbb{P}_{Z}(f)\right)
^{2}dP_{Z}(x)<\infty ,  \label{fep.var}
\end{equation}

\noindent which entails

\begin{equation}
\mathbb{P}_{Z}(\left\vert f\right\vert )=\int \left\vert f(x)\right\vert
dP_{Z}(x)<\infty \text{.}  \label{fep.esp}
\end{equation}

\bigskip \noindent Denote by $\mathcal{F}(S)$ - $\mathcal{F}$ for short -
the class of real-valued measurable functions that are defined on $S$ such
that (\ref{fep.var}) holds. The space $\mathcal{F}$, when endowed with the
addition and the external multiplication by real scalars, is a linear space.
Next, it is remarkable that $\mathbb{G}_{n}$ is linear on $\mathcal{F}$, that
is for $f$ and $g$ in $\mathcal{F}$ and for $(a,b)\in \mathbb{R}{^{2}}$, we
have

\begin{equation*}
a\mathbb{G}_{n}(f)+b\mathbb{G}_{n}(g)=\mathbb{G}_{n}(af+bg).
\end{equation*}

\bigskip \noindent We have this result

\begin{lemma} \label{fep.lemma.tool.1} Given the notation above, then for any finite
number of elements $f_{1},...,f_{k}$ of $\mathcal{S},k\geq 1,$ we have

\begin{equation*}
(\mathbb{G}_{n}(f_{1}),...,\mathbb{G}_{n}(f_{k}))^{T}\rightsquigarrow 
\mathcal{N}_{k}(0
\end{equation*}

\bigskip \noindent where 
\begin{equation*}
\Gamma (f_{i},f_{j})=\int \left( f_{i}-\mathbb{P}_{Z}(f_{i})\right) \left(
f_{j}-\mathbb{P}_{Z}(f_{j})\right) d\mathbb{P}_{Z}(x),1\leq ,j\leq k.
\end{equation*}
\end{lemma}

\bigskip \noindent This lemma says that the weak limit of the sequence $^{t}$ $(\mathbb{G}_{n}(f_1), \mathbb{G}_{n}(f_2), ...,\mathbb{G}_{n}(f_k))$ has the same law than the vector $^{t}$ $(\mathbb{G}(f_1), \mathbb{G}(f_2), ...,\mathbb{G}(f_k))$, where $\{\mathbb{G}(f),f\in \mathcal{F}\}$ is a Gaussian process of variance-covariance function

\begin{equation}
\Gamma (f,g)=\int \left( f-\mathbb{P}_{Z}(f)\right) \left(
g-\mathbb{P}_{Z}(g)\right) d\mathbb{P}_{Z}(x), \ \ (f,g)\in \mathcal{F}^2. \label{cv.pefTools.cov01}
\end{equation}

\noindent By applying the Skorohod-Wichura Theorem \label{cv.skorohodWichura} (See Chapter \ref{cv}), we may suppose that we are on the saùe probability space on which we have the following approximation : 

\begin{equation}
\mathbb{G}_{n}(f_{1})=\mathbb{G}_{n}(f_{1})+o_{\mathbb{P}}(1). 1\leq i \leq p. \label{cv.pefTools.rep01}
\end{equation}

\noindent We will come back later on the application of the formula.\\

\bigskip \noindent \textbf{PROOF}. It is enough to use the Cram\'{e}r-Wold
Criterion (see Proposition \ref{cv.wold} in Chapter \ref{ChapRevCvRk}), that
is to show that for any $a=(a_{1},...,a_{k})^{T}\in \mathbb{R}^{k}$, we have

$$
<a,T_{n}>\rightsquigarrow <a,T>
$$

\noindent where we have used the notation $T_{n}=(\mathbb{G}_{n}(f_{1}),...,\mathbb{G}_{n}(f_{k}))^{T}$, and where $T$ follows the $\mathcal{N}%
_{k}(0,\Gamma (f_{i},f_{j})_{1\leq i,j\leq k})$ law, and $<\circ ,\circ >$
stands for the usual scalar product in $\mathbb{R}^{k}$.\\

\noindent  But, by the standard central limit theorem in $\mathbb{R}$, we have
\begin{equation*}
<a,T_{n}>=\mathbb{G}_{n}\left( \sum\limits_{i=1}^{k}a_{i}f_{i}\right)
\rightsquigarrow \mathcal{N}(0,\sigma _{\infty }^{2}),
\end{equation*}

\bigskip \noindent where, for $g=\sum_{1\leq i\leq k}a_{i}f_{i},$%
\begin{equation*}
\sigma _{\infty }^{2}=\int \left( g(x)-\mathbb{P}_{Z}(g)\right) ^{2}dP_{Z}(x)
\end{equation*}

\bigskip \noindent and this easily gives
\begin{equation*}
\sigma _{\infty }^{2}=\sum\limits_{1\leq i,j\leq k}a_{i}a_{j}\Gamma
(f_{i},f_{j}),
\end{equation*}

\noindent so that $N(0,\sigma _{\infty }^{2})$ is the law of $<a,T>.$ The proof is finished.

\subsection{How to use the FEP tool?}  $ $\\

\noindent The usual statistics we are working with in Asymptotic Statistics are based on univariate or multivariate samples, meaning we usually work on $\mathbb{R}^{k}$. Once we have our sample $Z_{1},Z_{2},...$ as random
variables defined in the same probability space with values in $\mathbb{R}%
^{k}$, the studied statistic, say $T_{n}$, is usually a combination of
expressions of the form
\begin{equation*}
H_{n}=\frac{1}{n}\sum\limits_{i=1}^{k}H(Z_{i})
\end{equation*}

\bigskip \noindent for $H\in \mathcal{F}$. We use the results of Lemma \ref{fep.lemma.tool.1} and  Point (a) of 
Lemma \ref{oO-02}, to have this very sample expansion
$\mu (H)=\mathbb{E}H(Z),$ 
\begin{equation}
H_{n}=\mu (H)+n^{-1/2}\mathbb{G}_{n}(H).  \label{fep.expan}
\end{equation}

\bigskip \noindent We have that $\mathbb{G}_{n}(H)$ is asymptotically
bounded in probability since $\mathbb{G}_{n}(H)$ weakly converges to, say $%
M(H)$ and then by the continuous mapping theorem $\left\Vert \mathbb{G}%
_{n}(H)\right\Vert \rightsquigarrow \left\Vert M(H)\right\Vert .$ Since all
the $\mathbb{G}_{n}(H)$ are defined on the same probability space, we get
for all $\lambda >0,$ by the assertion of the Portmanteau Theorem for
concerning open sets,%
\begin{equation*}
\lim \sup_{n\rightarrow \infty }P(\left\Vert \mathbb{G}_{n}(H)\right\Vert
>\lambda )\leq P(\left\Vert M(H)\right\Vert >\lambda )
\end{equation*}

\bigskip \noindent and then 
\begin{equation*}
\liminf_{\lambda \rightarrow \infty}\limsup_{n\rightarrow \infty
}\mathbb{P}(\left\Vert \mathbb{G}_{n}(H)\right\Vert >\lambda )\leq \liminf_{\lambda \rightarrow \infty } \mathbb{P}(\left\Vert M(H)\right\Vert >\lambda )=0.
\end{equation*}

\bigskip \noindent From this, we use the big $O_{\mathbb{P}}$ notation, that
is $\mathbb{G}_{n}(H)=O_{\mathbb{P}}(1).$ Formula (\ref{fep.expan}) becomes 
\begin{equation*}
H_{n}=\mu (H)+n^{-1/2}\mathbb{G}_{n}(H)=\mu (H)+O_{\mathbb{P}}(n^{-1/2})
\end{equation*}

\bigskip \noindent and we will be able to use the delta method. Indeed, let $%
g:\mathbb{R} \rightarrow \mathbb{R}$ be continuously differentiable in a
neighborhood of $\mu (H).$ The mean value theorem leads to%
\begin{equation}
g(H_{n})=g(\mu (H))+g^{\prime }(\mu _{n}(H))\text{ }n^{-1/2}\mathbb{G}_{n}(H),
\label{fep.expan01}
\end{equation}

\bigskip \noindent where 
\begin{equation}
\mu _{n}(H)\in \lbrack (\mu (H)+n^{-1/2}\mathbb{G}_{n}(H))\wedge \mu
(H),(\mu (H)+n^{-1/2}\mathbb{G}_{n}(H))\vee \mu (H)],  \label{cv.pef.TAF01}
\end{equation}

\bigskip \noindent so that 
\begin{equation*}
\left\vert \mu _{n}(H)-\mu (H)\right\vert \leq n^{-1/2}\mathbb{G}_{n}(H)=O_{\mathbb{P}}(n^{-1/2}).
\end{equation*}

\bigskip \noindent \textbf{Warning}. Some people, if not many, would wrongly use $\mu_{n}(H)$ as a measurable random variable. But although $g^{\prime}(\mu _{n}(H))$ is measurable as a fraction of measurable applications, we cannot say that $\mu_{n}(H)$ without further information.\\

\noindent We may proceed as $\mu_{n}(H)$ is measurable, just for demonstration purposes. Next, we will give the most correct way. In that case, we have :\\

\bigskip \noindent The sequence $\mu _{n}(H)$ converges to $\mu _{n}(H)$ in outer probability, (denoted $\mu _{n}(H)\rightarrow _{\mathbb{P}}\mu (H)).$ But the convergence in probability to a constant is equivalent to the weak convergence. Then $\mu _{n}(H)\rightsquigarrow \mu (H).$ Using again the continuous mapping theorem, $g^{\prime }(\mu _{n}(H))\rightsquigarrow g^{\prime }(\mu (H))$, which in turn yields $g^{\prime}(\mu_{n}(H))\rightarrow _{\mathbb{P}}g^{\prime }(\mu (H))$ by the characterization of the weak convergence to a constant.\\

\noindent However, the correct way uses outer probability and Lemma \ref{cv.outerprob} above. We say : 

\bigskip \noindent Based on definition \ref{oO-03} and on Formula (\ref{cv.pef.TAF01}), we may see that rhe sequence $\mu _{n}(H)$ converges to $\mu _{n}(H)$ in outer probability, denoted $\mu _{n}(H)$ $\rightarrow _{\mathbb{P^{\ast}}}\mu (H))$. Then, in vertue of Lemma \ref{cv.outerprob} above, we get $g^{\prime }(\mu_{n}(H))\rightarrow _{\mathbb{P}}g^{\prime }(\mu (H))$ by the
characterization of the weak convergence to a constant.\\

\noindent Now, (\ref{fep.expan01}) becomes : 

\begin{eqnarray*}
g(H_{n}) &=&g(\mu (H))+(g^{\prime }(\mu (H)+o_{P}(1))\text{ }n^{-1/2}\mathbb{%
G}_{n}(H) \\
&=&g(\mu (H))+g^{\prime }(\mu (H)\times \text{ }n^{-1/2}\mathbb{G}%
_{n}(H)+o_{P}(1))\text{ }n^{-1/2}\mathbb{G}_{n}(H) \\
&=&g(\mu (H))+\text{ }n^{-1/2}\mathbb{G}_{n}(g^{\prime }(\mu
(H)H)+o_{P}(n^{-1/2}).
\end{eqnarray*}

\bigskip \noindent We obtain at the final expansion
\begin{equation}
g(H_{n})=g(\mu (H))+\text{ }n^{-1/2}\mathbb{G}_{n}(g^{\prime }(\mu
(H)H)+o_{P}(n^{-1/2}).  \label{fep.expanFinal}
\end{equation}

\noindent By using the Skorohod-Wichura representation, we get by Formula \label{cv.pefTools.rep01}, that 

\begin{equation}
g(H_{n})=g(\mu (H))+\text{ }n^{-1/2}\mathbb{G}(g^{\prime }(\mu
(H)H)+o_{P}(n^{-1/2}).  \label{fep.expanFinal01}
\end{equation}
\bigskip \noindent The method consists of using the expansion (\ref{fep.expanFinal}) as many times as needed and next to do some algebra on these
expansions.\\

\noindent The algebraic computations we refereed above are based on the application of the following lemma.

\begin{lemma}
\label{fep.lemma.tool.2} Let ($A_{n})$ and ($B_{n})$ be two sequences of real
valued random variables defined on the same probability space holding the
sequence $Z_{1}$, $Z_{2}$, ...

\noindent Let $A$ and $B$ be two real numbers and Let $L(z)$
and $H(z)$ be two real-valued functions$\ of$ $z\in S$, with $(L,H)\in \mathcal{F}^2$.\\

\noindent Suppose that 
$$
A_{n}=A+n^{-1/2}\mathbb{G}_{n}(L)+o_{P}(n^{-1/2})
$$ 

\noindent and 

$$
A_{n}=B+n^{-1/2}\mathbb{G}_{n}(H)+o_{P}(n^{-1/2}).
$$ 

\noindent Then, we have

\begin{equation*}
A_{n}+B_{n}=A+B+n^{-1/2}\mathbb{G}_{n}(L+H)+o_{P}(n^{-1/2}),
\end{equation*}

\noindent and

\begin{equation*}
A_{n}B_{n}=AB+n^{-1/2}\mathbb{G}_{n}(BL+AH)
\end{equation*}

\noindent and if $B\neq 0$, we also have
\begin{equation*}
\frac{A_{n}}{B_{n}}=\frac{A}{B}+n^{-1/2}\mathbb{G}_{n}(\frac{1}{B}L-\frac{A}{%
B^{2}}H)+o_{P}(n^{-1/2}).
\end{equation*}
\end{lemma}

\bigskip \noindent By putting together all the previous described steps in a smart
way, the methodology will lead us to a final result of the form%
\begin{equation*}
T_{n}=t+n^{-1/2}\mathbb{G}_{n}(h)+o_{P}(n^{-1/2}),
\end{equation*}

\noindent which entails the weak convergence

\begin{eqnarray*}
\sqrt{n}(T_{n}-t)&=&\mathbb{G}_{n}(h)+o_{P}(1)\rightsquigarrow \mathcal{N}(0,\Gamma(h,h))\\
&=&\mathbb{G}(h)+o_{P}(1).
\end{eqnarray*}

\bigskip \noindent Now, we are going show how to apply the methodology on the empirical linear correlation coefficient.

\subsection{An Example} \label{fep.subsec4}  $ $\\

\noindent We are going to illustrate our tool on the plug-in estimator of
the linear correlation coefficient of \ two random variables $(X,Y)$, with
neither $X$ nor $Y$ is degenerated, defined as follows

\begin{equation*}
\rho =\frac{\sigma _{xy}}{\sigma _{x}^{2}\sigma _{y}^{2}},
\end{equation*}

\noindent where

\begin{equation*}
\mu _{x}=\int x\text{ }dP_{X}(x),\text{ }\mu _{y}=\int x\text{ }dP_{X}(x),%
\text{ }\sigma _{xy}=\int (x-\mu _{x})(y-\mu _{y})dP_{(X,Y)}(x,y)
\end{equation*}

\noindent and

\begin{equation*}
\sigma _{x}^{2}=\int (x-\mu _{x})^{2}dP_{X}(x),\text{ }\sigma _{y}^{2}=\int
(x-\mu _{x})(y-\mu _{y})dP_{X}(y).
\end{equation*}

\noindent We also dismiss the case where $\left\vert \rho \right\vert =1$,
for which one of $X$ and $Y$ is an affine function of the other, meaning for example that we have $X=aY+b$ for some $(a,b)\mathbb{R}^2$.\\

\noindent It is clear that centering the variables $X$ and $Y$ at their expectations and normalizing them by their standard deviations $\sigma _{x}$ and $\sigma _{y}$ do not change the correlation coefficients $\rho$. So we may
and do center $X$ and $Y$ at their expectations and normalize them so that
we can and do assume that

\begin{equation*}
\mu _{x}=\text{ }\mu _{y}=0,\text{ }\sigma _{x}=\sigma _{y}=1.
\end{equation*}

\bigskip \noindent However, we will let these coefficients appear with their
names and we only use their particular values at the conclusion stage.%
\newline

\bigskip \noindent Let us construct the plug-in estimator of $\rho $. To this
end, let $(X_{1},Y_{1}),$ $(X_{2},Y_{2}),...$ be a sequence of independent
observations of $(X,Y).$ For each $n\geq 1,$ the plug-in estimator is the
following

\begin{equation*}
\rho _{n}=\left\{ \frac{1}{n}\sum_{i=1}^{n}(X_{i}-\overline{X})(Y_{i}-%
\overline{Y})\right\} \left\{ \frac{1}{n^{2}}\sum_{i=1}^{n}(X_{i}-\overline{X%
})^{2}\times \sum_{i=1}^{n}(X_{i}-\overline{X})^{2}\right\} ^{-1/2}.
\end{equation*}

\bigskip \noindent We are going to give the asymptotic theory of $\rho _{n}$
as an estimator of $\rho$. Introduce the notation

\begin{equation*}
\mu _{(p,x),(q,y)}=E((X-\mu _{x})^{p}(Y-\mu _{y})^{q}),\mu _{4,x}=E(X-\mu
_{x})^{4}\text{, }\mu _{4,y}=E(Y-\mu _{x})^{4}.
\end{equation*}

\bigskip \noindent Here is the outcome of the application of the method.

\begin{theorem}
\label{fep.theo1}

\bigskip \noindent Suppose that neither of $X$ and $Y$ is degenerated and
both have finite fourth moments and that $X^{3}Y$ and $XY^{3}$ have finite
expectations. Then, as $n\rightarrow \infty ,$%
\begin{equation*}
\sqrt{n}(\rho _{n}-\rho )\rightsquigarrow N(0,\sigma ^{2}),
\end{equation*}

\bigskip \noindent where 
\begin{eqnarray*}
\sigma ^{2} &=&\sigma _{x}^{-2}\sigma _{y}^{-2}(1+\rho ^{2}/2)\mu
_{(2,x),(2,y)}+\rho ^{2}(\sigma _{x}^{-4}\mu _{4,x}+\sigma _{y}^{-4}\mu
_{4,y})/4 \\
&&-\rho (\sigma _{x}^{-3}\sigma _{y}^{-1}\mu _{(3,x),(1,y)}+\sigma
_{x}^{-1}\sigma _{y}^{-3}\mu _{(1,x),(3,y)}).
\end{eqnarray*}
\end{theorem}

\bigskip \noindent This result enables to test independence between $X$ and $%
Y$, or to test non linear correlation in the following sense.

\bigskip

\begin{theorem}
\label{fep.theo2} Suppose that the assumptions of Theorem \ref{fep.theo1} hold. Then%
\newline
\bigskip \noindent \textbf{(1)} If $X$ and $Y$ are not linearly correlated,
that is $\rho =0,$we have 
\begin{equation*}
\sqrt{n}\rho _{n}\rightsquigarrow N(0,\sigma _{1}^{2}),
\end{equation*}

\bigskip \noindent where 
\begin{equation*}
\sigma _{1}^{2}=\sigma _{x}^{-2}\sigma _{y}^{-2}\mu _{(2,x),(2,y)}.
\end{equation*}

\bigskip \noindent \textbf{(2)} If $X$ and $Y$ are independent, then $\rho
=0,$ and\bigskip 
\begin{equation*}
\sqrt{n}\rho _{n}\rightsquigarrow N(0,1)
\end{equation*}
\end{theorem}

\bigskip \noindent \textbf{Proofs}. We are going to use the functional empirical process based on the
observations $(X_{i},Y_{i}),i=1,2,...$ that are independent copies of $(X,Y)$%
. Write

\begin{equation*}
\rho _{n}^{2}=\frac{\frac{1}{n}\sum_{i=1}^{n}X_{i}Y_{i}-\overline{X}\text{ }%
\overline{Y}}{\left\{ \frac{1}{n}\sum_{i=1}^{n}X_{i}^{2}-\overline{X}%
^{2}\right\} ^{1/2}\left\{ \frac{1}{n}\sum_{i=1}^{n}Y_{i}^{2}-\overline{Y}%
^{2}\right\} ^{1/2}}=\frac{A_{n}}{B_{n}}.
\end{equation*}

\bigskip \noindent Let us say for once that all the functions of $Z=(X,Y)$
that will appear below are measurable and have finite second moments. Let us
handle the numerator and denominator separately. To treat $A_{n},$ using the
empirical process implies that

\begin{equation}
\left\{ 
\begin{tabular}{l}
$\frac{1}{n}\sum_{i=1}^{n}X_{i}Y_{i}=\mu _{xy}+n^{-1/2}G_{n}(p),$ \\ 
$\overline{X}=\mu _{x}+n^{-1/2}G_{n}(\pi _{1}),$ \\ 
$\overline{Y}=\mu _{y}+n^{-1/2}G_{n}(\pi _{2}),$%
\end{tabular}%
\right.  \label{fep.casAN}
\end{equation}

\bigskip \noindent where $p(x,y)=xy$, $\pi _{1}(x,y)=x$ and $\pi _{2}(x,y)=y$.
From there we use the fact that $G_{n}(g)=O_{P}(1)$ for $\mathbb{E}%
(g(X,Y)^{2})<+\infty $ and get

\begin{equation}
A_{n}=\mu _{xy}+n^{-1/2}G_{n}(p)-(\mu _{x}+n^{-1/2}G_{n}(\pi _{1}))(\mu
_{y}+n^{-1/2}G_{n}(\pi _{2})).  \label{fep.haut}
\end{equation}

\bigskip \noindent This leads to 
\begin{equation*}
A_{n}=\sigma _{xy}+n^{-1/2}G_{n}(H_{1})+o_{P}(n^{-1/2})
\end{equation*}

\bigskip \noindent with 
\begin{equation*}
H_{1}(x,y)=p(x,y)-\mu _{x}\pi _{2}-\mu _{y}\pi _{1.}
\end{equation*}

\bigskip \noindent Next, we have to handle $B_{n}.$ Since the roles of $%
\left\{ \frac{1}{n}\sum_{i=1}^{n}X_{i}^{2}-\overline{X}^{2}\right\} ^{1/2}$
and of $\left\{ \frac{1}{n}\sum_{i=1}^{n}Y_{i}^{2}-\overline{Y}^{2}\right\}
^{1/2}$ are symmetrical, we treat one of them and extend the results to the
other. Let us handle $\left\{ \frac{1}{n}\sum_{i=1}^{n}X_{i}^{2}-\overline{X}%
^{2}\right\} ^{1/2}$. The combination of (\ref{fep.casAN}) and the Delta method enables to get%
\begin{equation*}
\overline{X}^{2}=\left( \mu _{x}+n^{-1/2}G_{n}(\pi _{1})\right) ^{2}=\mu
_{x}^{2}+2\mu _{x}n^{-1/2}G_{n}(\pi _{1})+o_{P}(n^{-1/2}),
\end{equation*}

\bigskip \noindent that is,
\bigskip
\begin{equation*}
\overline{X}^{2}=\left( \mu _{x}+n^{-1/2}G_{n}(\pi _{1})\right) ^{2}=\mu
_{x}^{2}+n^{-1/2}G_{n}(2\mu _{x}\pi _{1})+o_{P}(n^{-1/2}).
\end{equation*}

\bigskip \noindent From there, we get
\begin{eqnarray*}
\frac{1}{n}\sum_{i=1}^{n}X_{i}^{2}-\overline{X}^{2}
&=&m_{2,x}+n^{-1/2}G_{n}(\pi _{1}^{2})-\overline{X}^{2} \\
&=&m_{2,x}-\mu _{x}^{2}+n^{-1/2}G_{n}(\pi _{1}^{2}-2\mu _{x}\pi
_{1})+o_{P}(n^{-1/2}) \\
&=&\sigma _{x}^{2}+n^{-1/2}G_{n}(\pi _{1}^{2}-2\mu _{x}\pi
_{1})+o_{P}(n^{-1/2}).
\end{eqnarray*}

\bigskip \noindent Using the Delta-method once again leads to 
\begin{equation*}
\left\{ \frac{1}{n}\sum_{i=1}^{n}X_{i}^{2}-\overline{X}^{2}\right\}
^{1/2}=\sigma _{x}+n^{-1/2}G_{n}(\frac{1}{2\sigma _{x}}\left\{ \pi
_{1}^{2}-2\mu _{x}\pi _{1}\right\} )+o_{P}(n^{-1/2}).
\end{equation*}

\bigskip \noindent In a similar way, we get
\begin{equation*}
\left\{ \frac{1}{n}\sum_{i=1}^{n}Y_{i}^{2}-\overline{Y}^{2}\right\}
^{1/2}=\sigma _{y}+n^{-1/2}G_{n}(\frac{1}{2\sigma _{y}}\left\{ \pi
_{2}^{2}-2\mu _{y}\pi _{2}\right\} )+o_{P}(n^{-1/2}).
\end{equation*}

\bigskip \noindent We get
\begin{eqnarray*}
B_{n} &=&\left\{ \frac{1}{n}\sum_{i=1}^{n}X_{i}^{2}-\overline{X}^{2}\right\}
^{1/2}\left\{ \frac{1}{n}\sum_{i=1}^{n}Y_{i}^{2}-\overline{Y}^{2}\right\}
^{1/2} \\
&=&\sigma _{x}\sigma _{y}+n^{-1/2}G_{n}(\frac{\sigma _{y}}{2\sigma _{x}}%
\left\{ \pi _{1}^{2}-2\mu _{x}\pi _{1}\right\} +\frac{\sigma _{x}}{2\sigma
_{y}}\left\{ \pi _{2}^{2}-2\mu _{y}\pi _{2}\right\} )+o_{P}(n^{-1/2}).
\end{eqnarray*}

\bigskip \noindent By setting
\begin{equation*}
H_{2}(x,y)=\frac{\sigma _{y}}{2\sigma _{x}}\left\{ \pi _{1}^{2}-2\mu _{x}\pi
_{1}\right\} +\frac{\sigma _{x}}{2\sigma _{y}}\left\{ \pi _{2}^{2}-2\mu
_{y}\pi _{2}\right\},
\end{equation*}

\bigskip \noindent we have

\begin{equation}
B_{n}=\sigma _{x}\sigma _{y}+n^{-1/2}G_{n}(H_{2})+n^{-1/2}.  \label{fep.bas}
\end{equation}

\bigskip \noindent Now, combining (\ref{fep.haut}) and (\ref{fep.bas}) and using
Lemma \ref{fep.lemma.tool.2} yields

\begin{equation*}
\sqrt{n}(\rho _{n}^{2}-\rho ^{2})=n^{-1/2}G_{n}(\frac{1}{\sigma _{x}\sigma
_{y}}H_{1}-\frac{\sigma _{xy}}{\sigma _{x}^{2}\sigma _{y}^{2}}%
H_{2})+o_{P}(1).
\end{equation*}

\bigskip \noindent Put
\begin{equation*}
H=\frac{1}{\sigma _{x}\sigma _{y}}(p(x,y)-\mu _{x}\pi _{2}-\mu _{y}\pi _{1})-%
\frac{\rho }{\sigma _{x}\sigma _{y}}\left\{ \frac{1}{2\sigma _{x}^{2}}%
\left\{ \pi _{1}^{2}-2\mu _{x}\pi _{1}\right\} +\frac{1}{2\sigma _{y}^{2}}%
\left\{ \pi _{2}^{2}-2\mu _{y}\pi _{2}\right\} \right\}.
\end{equation*}

\bigskip \noindent Now we continue with the centered and normalized case to get

\begin{equation*}
H(x,y)=p(x,y)-\frac{\rho }{2}(\pi _{1}^{2}+\pi _{2}^{2})
\end{equation*}

\bigskip \noindent and 
\begin{equation*}
H(X,Y)=XY-\frac{\rho }{2}(X^{2}+Y^{2}).
\end{equation*}

\bigskip \noindent Denote

\begin{equation*}
\mu _{(p,x),(q,y)}=E((X-\mu _{x})^{p}(Y-\mu _{y})^{q}).
\end{equation*}

\bigskip \noindent We have 
\begin{equation*}
\mathbb{E}H(X,Y)=\sigma _{xy}-\rho =0
\end{equation*}

\bigskip \noindent and $var(H(X,Y))$ is equal to 
\begin{equation*}
\mu _{(2,x),(2,y)}+\rho ^{2}(\mu _{4,x}+\mu _{4,y})/4-\rho (\mu
_{(3,x),(1,y)}+\mu _{(1,x),(3,y)})+\rho ^{2}\mu _{(2,x),(2,y)}/2
\end{equation*}

\noindent and finally 

$$
var(H(X,Y))=\sigma _{0}^{2}
$$ 

\noindent with
\begin{equation*}
\sigma _{0}^{2}=(1+\rho ^{2}/2)\mu _{(2,x),(2,y)}+\rho ^{2}(\mu _{4,x}+\mu
_{4,y})/4-\rho (\mu _{(3,x),(1,y)}+\mu _{(1,x),(3,y)}).
\end{equation*}

\bigskip \noindent This gives the conclusion that for centered and
normalized $X$ and $Y$,

\begin{equation*}
\sqrt{n}(\rho _{n}-\rho )\rightsquigarrow N(0,\sigma _{0}^{2}).
\end{equation*}

\bigskip \noindent Next, if we use the normalizing coefficients in $\sigma
_{0}$, we get

\begin{eqnarray*}
\sigma ^{2} &=&\sigma _{x}^{2}\sigma _{y}^{2}(1+\rho ^{2}/2)\mu
_{(2,x),(2,y)}+\rho ^{2}(\sigma _{x}^{4}\mu _{4,x}+\sigma _{y}^{4}\mu
_{4,y})/4 \\
&&-\rho (\sigma _{x}^{3}\sigma _{y}\mu _{(3,x),(1,y)}+\sigma _{x}\sigma
_{y}^{3}\mu _{(1,x),(3,y)})
\end{eqnarray*}

\bigskip \noindent and we conclude in the general case that

\begin{equation*}
\sqrt{n}(\rho _{n}-\rho )\rightsquigarrow N(0,\sigma ^{2}).
\end{equation*}

\bigskip \noindent The proof of Theorem \ref{fep.theo2} follows by easy
computations under the particular conditions of $\rho $ and under
independence.

 

%% file: asymptotics_analysis_01_en.tex
\chapter{Elements of Theory of Functions and Real Analysis} \label{funct}

\section{Review on limits in $\overline{\mathbb{R}}$. What should not be ignored on limits.} \label{funct.sec.1}

\noindent \textbf{Definition} $\ell \in \overline{\mathbb{R}}$ is an accumulation point of a sequence 
 $(x_{n})_{n\geq 0}$ of real numbers finite or infinite, in $\overline{\mathbb{R}}$, if and only if there exists a sub-sequence $(x_{n(k)})_{k\geq 0}$ of
 $(x_{n})_{n\geq 0}$ such that $%
x_{n(k)}$ converges to $\ell $, as $k\rightarrow +\infty $.\newline

\noindent \textbf{Exercise 1}.\\

\noindent  Set $y_{n}=\inf_{p\geq n}x_{p}$ and $%
z_{n}=\sup_{p\geq n}x_{p} $ for all $n\geq 0$. Show that :\newline

\noindent \textbf{(1)} $\forall n\geq 0,y_{n}\leq x_{n}\leq z_{n}$.\newline

\noindent \textbf{(2)} Justify the existence of the limit of $y_{n}$ called limit inferior of the sequence $(x_{n})_{n\geq 0}$, denoted by $%
\liminf x_{n}$ or $\underline{\lim }$ $x_{n},$ and that it is equal to the following%
\begin{equation*}
\underline{\lim }\text{ }x_{n}=\lim \inf x_{n}=\sup_{n\geq 0}\inf_{p\geq
n}x_{p}.
\end{equation*}

\noindent \textbf{(3)} Justify the existence of the limit of $z_{n}$ called limit superior of the sequence $(x_{n})_{n\geq 0}$ denoted by $%
\lim \sup x_{n}$ or $\overline{\lim }$ $x_{n},$ and that it is equal%
\begin{equation*}
\overline{\lim }\text{ }x_{n}=\lim \sup x_{n}=\inf_{n\geq 0}\sup_{p\geq
n}x_{p}x_{p}.
\end{equation*}

\bigskip

\noindent \textbf{(4)} Establish that 
\begin{equation*}
-\liminf x_{n}=\limsup (-x_{n})\noindent \text{ \ \ and \ }-\limsup
x_{n}=\liminf (-x_{n}).
\end{equation*}

\newpage \noindent \textbf{(5)} Show that the limit superior is sub-additive and the limit inferior is super-additive, i.e. :  for two sequences
$(s_{n})_{n\geq 0}$ and $(t_{n})_{n\geq 0}$ 
\begin{equation*}
\limsup (s_{n}+t_{n})\leq \limsup s_{n}+\limsup t_{n}
\end{equation*}

\noindent and
\begin{equation*}
\lim \inf (s_{n}+t_{n})\geq \lim \inf s_{n}+\lim \inf t_{n}.
\end{equation*}

\noindent \textbf{(6)} Deduce from (1) that if%
\begin{equation*}
\lim \inf x_{n}=\lim \sup x_{n},
\end{equation*}%
then $(x_{n})_{n\geq 0}$ has a limit and 
\begin{equation*}
\lim x_{n}=\lim \inf x_{n}=\lim \sup x_{n}
\end{equation*}

\bigskip

\noindent \textbf{Exercise 2.} Accumulation points of $%
(x_{n})_{n\geq 0}$.\newline

\noindent \textbf{(a)} Show that if $\ell _{1}$=$\lim \inf x_{n}$ and $\ell
_{2}=\lim \sup x_{n}$ are accumulation points of $(x_{n})_{n\geq 0}.
$ Show one case and deduce the second one and by using Point (3) of Exercise 1.\newline

\noindent \textbf{(b)} Show that $\ell _{1}$ is the smallest accumulation point of $(x_{n})_{n\geq 0}$ and $\ell _{2}$ is the biggest.
(Similarly, show one case and deduce the second one and by using Point (3) of Exercise 1).\newline

\noindent \textbf{(c)} Deduce from (a) that if $(x_{n})_{n\geq 0}$ has
a limit $\ell ,$ then it is equal to the unique accumulation point and so,%
\begin{equation*}
\ell =\overline{\lim }\text{ }x_{n}=\lim \sup x_{n}=\inf_{n\geq
0}\sup_{p\geq n}x_{p}.
\end{equation*}

\noindent \textbf{(d)} Combine this result with Point \textbf{(6)} of Exercise 1 to show that a sequence $(x_{n})_{n\geq 0}$ of $\overline{\mathbb{R}}
$ has a limit $\ell $ in $\overline{\mathbb{R}}$ if and only if\ $\lim \inf
x_{n}=\lim \sup x_{n}$ and then%
\begin{equation*}
\ell =\lim x_{n}=\lim \inf x_{n}=\lim \sup x_{n}.
\end{equation*}

\newpage

\noindent \textbf{Exercise 3. } Let $(x_{n})_{n\geq 0}$ be a non-decreasing sequence
of $\overline{\mathbb{R}}$. Study its limit superior and its limit inferior and deduce that%
\begin{equation*}
\lim x_{n}=\sup_{n\geq 0}x_{n}.
\end{equation*}

\noindent Deduce that for a non-increasing sequence $(x_{n})_{n\geq 0}$
of $\overline{\mathbb{R}},$%
\begin{equation*}
\lim x_{n}=\inf_{n\geq 0}x_{n}.
\end{equation*}

\bigskip

\noindent \textbf{Exercise 4.} (Convergence criteria)\newline

\noindent \textbf{Prohorov Criterion} Let $(x_{n})_{n\geq 0}$ be a sequence of $\overline{%
\mathbb{R}}$ and a real number $\ell \in \overline{\mathbb{R}}$ such that: Every subsequence of $(x_{n})_{n\geq 0}$ 
also has a subsequence ( that is a subssubsequence of $(x_{n})_{n\geq 0}$ ) that converges to $\ell .$
Then, the limit of $(x_{n})_{n\geq 0}$ exists and is equal $\ell .$\newline

\noindent \textbf{Upcrossing or Downcrossing Criterion}. \newline

\noindent Let $(x_{n})_{n\geq 0}$ be a sequence in $\overline{\mathbb{R}}$ and two real numbers $a$ and $b$ such that $a<b.$
We define%
\begin{equation*}
\nu _{1}=\left\{ 
\begin{array}{cc}
\inf  & \{n\geq 0,x_{n}<a\} \\ 
+\infty  & \text{if (}\forall n\geq 0,x_{n}\geq a\text{)}%
\end{array}%
\right. .
\end{equation*}%
If $\nu _{1}$ is finite, let%
\begin{equation*}
\nu _{2}=\left\{ 
\begin{array}{cc}
\inf  & \{n>\nu _{1},x_{n}>b\} \\ 
+\infty  & \text{if (}n>\nu _{1},x_{n}\leq b\text{)}%
\end{array}%
\right. .
\end{equation*}%
.

\noindent As long as the $\nu _{j}'s$ are finite, we can define for $\nu
_{2k-2}(k\geq 2)$

\begin{equation*}
\nu _{2k-1}=\left\{ 
\begin{array}{cc}
\inf  & \{n>\nu _{2k-2},x_{n}<a\} \\ 
+\infty  & \text{if (}\forall n>\nu _{2k-2},x_{n}\geq a\text{)}%
\end{array}%
\right. .
\end{equation*}%
and for $\nu _{2k-1}$ finite, 
\begin{equation*}
\nu _{2k}=\left\{ 
\begin{array}{cc}
\inf  & \{n>\nu _{2k-1},x_{n}>b\} \\ 
+\infty  & \text{if (}n>\nu _{2k-1},x_{n}\leq b\text{)}%
\end{array}%
\right. .
\end{equation*}

\noindent We stop once one $\nu _{j}$ is $+\infty$. If $\nu
_{2j}$ is finite, then 
\begin{equation*}
x_{\nu _{2j}}-x_{\nu _{2j-1}}>b-a. 
\end{equation*}

\noindent We then say : by that moving from $x_{\nu _{2j-1}}$ to $x_{\nu
_{2j}},$ we have accomplished a crossing (toward the up) of the segment $[a,b]$
called \textit{up-crossings}. Similarly, if one $\nu _{2j+1}$
is finite, then the segment $[x_{\nu _{2j}},x_{\nu _{2j+1}}]$ is a crossing downward (down-crossing) of the segment $[a,b].$ Let%
\begin{equation*}
D(a,b)=\text{ number of up-crossings of the sequence of the segment }[a,b]\text{.}
\end{equation*}

\bigskip

\noindent \textbf{(a)} What is the value of $D(a,b)$ if \ $\nu _{2k}$ is finite and $\nu
_{2k+1}$ infinite.\newline

\noindent \textbf{(b)} What is the value of $D(a,b)$ if \ $\nu _{2k+1}$ is finite and $\nu
_{2k+2}$ infinite.\newline

\noindent \textbf{(c)} What is the value of $D(a,b)$ if \ all the $\nu _{j}'s$ are finite.%
\newline

\noindent \textbf{(d)} Show that $(x_{n})_{n\geq 0}$ has a limit iff
for all $a<b,$ $D(a,b)<\infty.$\newline

\noindent \textbf{(e)} Show that $(x_{n})_{n\geq 0}$ has a limit iff
for all $a<b,$ $(a,b)\in \mathbb{Q}^{2},D(a,b)<\infty .$\newline

\bigskip

\noindent \textbf{Exercise 5. } (Cauchy Criterion). Let $%
(x_{n})_{n\geq 0}$ $\mathbb{R}$ be a sequence of (\textbf{real numbers}).\newline

\noindent \textbf{(a)} Show that if $(x_{n})_{n\geq 0}$ is Cauchy,
then it has a unique accumulation point $\ell \in 
\mathbb{R}$ which is its limit.\newline

\noindent \textbf{(b)} Show that if a sequence $(x_{n})_{n\geq 0}\subset 
\mathbb{R}$ \ converges to $\ell \in \mathbb{R},$ then, it is Cauchy.%
\newline

\noindent \textbf{(c)} Deduce the Cauchy criterion for sequences of real numbers.

\newpage

\begin{center}
\textbf{SOLUTIONS}
\end{center}

\noindent \textbf{Exercise 1}.\newline

\noindent \textbf{Question (1)}. It is obvious that :%
\begin{equation*}
\underset{p\geq n}{\inf }x_{p}\leq x_{n}\leq \underset{p\geq n}{\sup }x_{p},
\end{equation*}

\noindent since $x_{n}$ is an element of $\left\{
x_{n},x_{n+1},...\right\} $ on which we take the supremum or the infimum.%
\newline

\noindent \textbf{Question (2)}. Let $y_{n}=\underset{p\geq 0}{\inf }%
x_{p}=\underset{p\geq n}{\inf }A_{n},$ where $A_{n}=\left\{
x_{n},x_{n+1},...\right\} $ is a non-increasing sequence of sets : $\forall n\geq 0$,
\begin{equation*}
A_{n+1}\subset A_{n}.
\end{equation*}

\noindent So the infimum on $A_{n}$ increases. If $y_{n}$ increases in $%
\overline{\mathbb{R}},$ its limit is its upper bound, finite or infinite. So%
\begin{equation*}
y_{n}\nearrow \underline{\lim }\text{ }x_{n},
\end{equation*}

\noindent is a finite or infinite number.\newline

\noindent \textbf{Question (3)}. We also show that $z_{n}=\sup A_{n}$ decreases and $z_{n}\downarrow \overline{\lim }$ $x_{n}$.\newline

\noindent \textbf{Question (4) \label{qst4}}. We recall that 
\begin{equation*}
-\sup \left\{ x,x\in A\right\} =\inf \left\{ -x,x\in A\right\}, 
\end{equation*}

\noindent which we write 
\begin{equation*}
-\sup A=\inf (-A).
\end{equation*}

\noindent Thus,

\begin{equation*}
-z_{n}=-\sup A_{n}=\inf (-A_{n}) = \inf \left\{-x_{p},p\geq n\right\}.
\end{equation*}

\noindent The right hand term tends to $-\overline{\lim}\ x_{n}$ and the left hand to $\underline{\lim} (-x_{n})$ and so 

\begin{equation*}
-\overline{\lim}\ x_{n}=\underline{\lim }\ (-x_{n}).
\end{equation*}

\bigskip \noindent Similarly, we show:
\begin{equation*}
-\underline{\lim } \ (x_{n})=\overline{\lim} \ (-x_{n}).
\end{equation*}

\noindent 

\noindent \textbf{Question (5)}. These properties come from the formulas, where $A\subseteq \mathbb{R},B\subseteq \mathbb{R}$ :%
\begin{equation*}
\sup \left\{ x+y,A\subseteq \mathbb{R},B\subseteq \mathbb{R}\right\} \leq
\sup A+\sup B.
\end{equation*}

\noindent In fact : 
\begin{equation*}
\forall x\in \mathbb{R},x\leq \sup A
\end{equation*}

\noindent and
\begin{equation*}
\forall y\in \mathbb{R},y\leq \sup B.
\end{equation*}

\noindent Thus 
\begin{equation*}
x+y\leq \sup A+\sup B,
\end{equation*}

\noindent where 
\begin{equation*}
\underset{x\in A,y\in B}{\sup }x+y\leq \sup A+\sup B.
\end{equation*}%
Similarly,%
\begin{equation*}
\inf (A+B\geq \inf A+\inf B.
\end{equation*}

\noindent In fact :

\begin{equation*}
\forall (x,y)\in A\times B,x\geq \inf A\text{ and }y\geq \inf B.
\end{equation*}

\noindent Thus 
\begin{equation*}
x+y\geq \inf A+\inf B,
\end{equation*}

\noindent and so
\begin{equation*}
\underset{x\in A,y\in B}{\inf }(x+y)\geq \inf A+\inf B
\end{equation*}

\noindent \textbf{Application}.\newline

\begin{equation*}
\underset{p\geq n}{\sup } \ (x_{p}+y_{p})\leq \underset{p\geq n}{\sup } \ x_{p}+\underset{p\geq n}{\sup } \ y_{p}.
\end{equation*}

\noindent All these sequences are non-increasing. By taking the infimum, we obtain the limits superior :

\begin{equation*}
\overline{\lim }\text{ }(x_{n}+y_{n})\leq \overline{\lim }\text{ }x_{n}+%
\overline{\lim }\text{ }x_{n}.
\end{equation*}

\bigskip

\noindent \textbf{Question (6)}. Set

\begin{equation*}
\underline{\lim } \ x_{n}=\overline{\lim } \ x_{n}.
\end{equation*}

\noindent Since : 
\begin{equation*}
\forall x\geq 1,\text{ }y_{n}\leq x_{n}\leq z_{n},
\end{equation*}%

\begin{equation*}
y_{n}\rightarrow \underline{\lim} \ x_{n}
\end{equation*}%

\noindent and 

\begin{equation*}
z_{n}\rightarrow \overline{\lim } \ x_{n},
\end{equation*}

\noindent we apply the Sandwich Theorem to conclude that the limit of $x_{n}$ exists and :

\begin{equation*}
\lim \text{ }x_{n}=\underline{\lim }\text{ }x_{n}=\overline{\lim }\text{ }%
x_{n}.
\end{equation*}

\bigskip 
\noindent \textbf{Exercice 2}.\newline

\noindent \textbf{Question (a).}\\

\noindent Thanks to Question (4) of Exercise 1, it suffices to show this property for one of the limits. Consider the limit superior and the three cases:\\

\noindent \textbf{The case of a finite limit superior} :

\begin{equation*}
\underline{\lim} x_{n}=\ell \text{ finite.}
\end{equation*}

\noindent By definition, 
\begin{equation*}
z_{n}=\underset{p\geq n}{\sup }x_{p}\downarrow \ell .
\end{equation*}

\noindent So: 
\begin{equation*}
\forall \varepsilon >0,\exists (N(\varepsilon )\geq 1),\forall p\geq
N(\varepsilon ),\ell -\varepsilon <x_{p}\leq \ell +\varepsilon .
\end{equation*}

\noindent Take less than that:

\begin{equation*}
\forall \varepsilon >0,\exists n_{\varepsilon }\geq 1:\ell -\varepsilon
<x_{n_{\varepsilon }}\leq \ell +\varepsilon.
\end{equation*}

\noindent We shall construct a sub-sequence converging to $\ell$.\\

\noindent Let $\varepsilon =1:$%
\begin{equation*}
\exists N_{1}:\ell -1<x_{N_{1}}=\underset{p\geq n}{\sup }x_{p}\leq \ell +1.
\end{equation*}

\noindent But if 
\begin{equation}
z_{N_{1}}=\underset{p\geq n}{\sup }x_{p}>\ell -1, \label{cc}
\end{equation}

\noindent there surely exists an $n_{1}\geq N_{1}$ such that%
\begin{equation*}
x_{n_{1}}>\ell -1.
\end{equation*}

\noindent If not, we would have 
\begin{equation*}
( \forall p\geq N_{1},x_{p}\leq \ell -1\ ) \Longrightarrow \sup \left\{
x_{p},p\geq N_{1}\right\} =z_{N_{1}}\geq \ell -1,
\end{equation*}

\noindent which is contradictory with (\ref{cc}). So, there exists $n_{1}\geq N_{1}$ such that
\begin{equation*}
\ell -1<x_{n_{1}}\leq \underset{p\geq N_{1}}{\sup }x_{p}\leq \ell -1.
\end{equation*}

\noindent i.e.

\begin{equation*}
\ell -1<x_{n_{1}}\leq \ell +1.
\end{equation*}

\noindent We move to step $\varepsilon =\frac{1}{2}$ and we consider the sequence%
 $(z_{n})_{n\geq n_{1}}$ whose limit remains $\ell$. So, there exists $N_{2}>n_{1}:$%
\begin{equation*}
\ell -\frac{1}{2}<z_{N_{2}}\leq \ell -\frac{1}{2}.
\end{equation*}

\noindent We deduce like previously that $n_{2}\geq N_{2}$ such that%
\begin{equation*}
\ell -\frac{1}{2}<x_{n_{2}}\leq \ell +\frac{1}{2}
\end{equation*}

\noindent with $n_{2}\geq N_{1}>n_{1}$.\\

\noindent Next, we set $\varepsilon =1/3,$ there will exist $N_{3}>n_{2}$ such that%
\begin{equation*}
\ell -\frac{1}{3}<z_{N_{3}}\leq \ell -\frac{1}{3}
\end{equation*}

\noindent and we could find an $n_{3}\geq N_{3}$ such that%

\begin{equation*}
\ell -\frac{1}{3}<x_{n_{3}}\leq \ell -\frac{1}{3}.
\end{equation*}

\noindent Step by step, we deduce the existence of $%
x_{n_{1}},x_{n_{2}},x_{n_{3}},...,x_{n_{k}},...$ with $n_{1}<n_{2}<n_{3}%
\,<...<n_{k}<n_{k+1}<...$ such that

$$
\forall k\geq 1, \ell -\frac{1}{k}<x_{n_{k}}\leq \ell -\frac{1}{k},
$$

\noindent i.e.

\begin{equation*}
\left\vert \ell -x_{n_{k}}\right\vert \leq \frac{1}{k},
\end{equation*}

\noindent which will imply: 
\begin{equation*}
x_{n_{k}}\rightarrow \ell 
\end{equation*}

\noindent Conclusion : $(x_{n_{k}})_{k\geq 1}$ is very well a subsequence since $n_{k}<n_{k+1}$ for all $k \geq 1$ 
and it converges to $\ell$, which is then an accumulation point.\\

\noindent \textbf{Case of the limit superior equal $+\infty$} : 
$$
\overline{\lim} \text{ } x_{n}=+\infty.
$$
\noindent Since $z_{n}\uparrow +\infty ,$ we have : $\forall k\geq 1,\exists
N_{k}\geq 1,$ 
\begin{equation*}
z_{N_{k}}\geq k+1.
\end{equation*}

\noindent For $k=1$, let $z_{N_{1}}=\underset{p\geq N_{1}}{\inf }%
x_{p}\geq 1+1=2.$ So there exists 
\begin{equation*}
n_{1}\geq N_{1}
\end{equation*}

\noindent such that :
\begin{equation*}
x_{n_{1}}\geq 1.
\end{equation*}

\noindent For $k=2$, consider the sequence $(z_{n})_{n\geq n_{1}+1}.$
We find in the same manner 
\begin{equation*}
n_2 \geq n_{1}+1
\end{equation*}%
\noindent and 
\begin{equation*}
x_{n_{2}}\geq 2.
\end{equation*}

\noindent Step by step, we find for all $k\geq 3$, an $n_{k}\geq n_{k-1}+1$ such that
\begin{equation*}
x_{n_{k}}\geq k,
\end{equation*}

\noindent which leads to $x_{n_{k}}\rightarrow +\infty $ as $k\rightarrow +\infty $.\\

\noindent \textbf{Case of the limit superior equal $-\infty$} : 

$$
\overline{\lim }x_{n}=-\infty.
$$

\noindent This implies : $\forall k\geq 1,\exists N_{k}\geq 1,$ such that%
\begin{equation*}
z_{n_{k}}\leq -k.
\end{equation*}

\noindent For $k=1$, there exists $n_{1}$ such that%
\begin{equation*}
z_{n_{1}}\leq -1.
\end{equation*}
But 
\begin{equation*}
x_{n_{1}}\leq z_{n_{1}}\leq -1.
\end{equation*}

\noindent Let $k=2$. Consider $\left( z_{n}\right) _{n\geq
n_{1}+1}\downarrow -\infty .$ There will exist $n_{2}\geq n_{1}+1:$%
\begin{equation*}
x_{n_{2}}\leq z_{n_{2}}\leq -2
\end{equation*}

\noindent Step by step, we find $n_{k1}<n_{k+1}$ in such a way that $x_{n_{k}}<-k$ for all $k$ bigger than $1$. So
\begin{equation*}
x_{n_{k}}\rightarrow +\infty 
\end{equation*}

\bigskip

\noindent \textbf{Question (b).}\\

\noindent Let $\ell$ be an accumulation point of $(x_n)_{n \geq 1}$, the limit of one of its sub-sequences $(x_{n_{k}})_{k \geq 1}$. We have

$$
y_{n_{k}}=\inf_{p\geq n_k} \ x_p \leq x_{n_{k}} \leq  \sup_{p\geq n_k} \ x_p=z_{n_{k}}.
$$

\noindent The left hand side term is a sub-sequence of $(y_n)$ tending to the limit inferior and the right hand side is a 
sub-sequence of $(z_n)$ tending to the limit superior. So we will have:

$$
\underline{\lim} \ x_{n} \leq \ell \leq \overline{\lim } \ x_{n},
$$

\noindent which shows that $\underline{\lim} \ x_{n}$ is the smallest accumulation point and $\overline{\lim } \ x_{n}$ is the largest.\\

\noindent \textbf{Question (c).} If the sequence $(x_n)_{n \geq 1}$ has a limit $\ell$, it is the limit of all its sub-sequences,
so subsequences tending to the limits superior and inferior. Which answers question (b).\\

\noindent \textbf{Question (d).} We answer this question by combining point (d) of this exercise and Point 6) of the Exercise 1.\\

\noindent \textbf{Exercise 3}. Let $(x_{n})_{n\geq 0}$ be a non-decreasing sequence, we have:%
\begin{equation*}
z_{n}=\underset{p\geq n}{\sup} \ x_{p}=\underset{p\geq 0}{\sup} \ x_{p},\forall
n\geq 0.
\end{equation*}

\noindent Why? Because by increasingness,%
\begin{equation*}
\left\{ x_{p},p\geq 0\right\} =\left\{ x_{p},0\leq p\leq n-1\right\} \cup
\left\{ x_{p},p\geq n\right\}.
\end{equation*}

\bigskip

\noindent Since all the elements of $\left\{ x_{p},0\leq p\leq
n-1\right\} $ are smaller than than those of $\left\{ x_{p},p\geq n\right\} ,$
the supremum is achieved on $\left\{ x_{p},p\geq n\right\} $ and so 
\begin{equation*}
\ell =\underset{p\geq 0}{\sup } \ x_{p}=\underset{p\geq n}{\sup }x_{p}=z_{n}.
\end{equation*}

\noindent Thus
\begin{equation*}
z_{n}=\ell \rightarrow \ell .
\end{equation*}

\noindent We also have $y_n=\inf \left\{ x_{p},0\leq p\leq n\right\}=x_n$, which is a non-decreasing sequence and so converges to
$\ell =\underset{p\geq 0}{\sup } \ x_{p}$. \\

\bigskip

\noindent \textbf{Exercise 4}.\\

\noindent Let $\ell \in \overline{\mathbb{R}}$ having the indicated property. Let $\ell ^{\prime }$ be a given accumulation point.%
\begin{equation*}
 \left( x_{n_{k}}\right)_{k\geq 1} \subseteq \left( x_{n}\right) _{n\geq 0}%
\text{ such that }x_{n_{K}}\rightarrow \ell ^{\prime}.
\end{equation*}

\noindent By hypothesis this sub-sequence $\left( x_{n_{K}}\right) $
has in turn a sub-sub-sequence $\left( x_{n_{\left( k(p)\right) }}\right)_{p\geq 1} $ such that $x_{n_{\left( k(p)\right) }}\rightarrow
\ell $ as $p\rightarrow +\infty $.\newline

\noindent But as a sub-sequence of $\left( x_{n_{\left( k\right)
}}\right) ,$ 
\begin{equation*}
x_{n_{\left( k(\ell )\right) }}\rightarrow \ell ^{\prime }.
\end{equation*}%
Thus
\begin{equation*}
\ell =\ell ^{\prime}.
\end{equation*}

\noindent Applying that to the limit superior and limit inferior, we have:%
\begin{equation*}
\overline{\lim} \ x_{n}=\underline{\lim}\ x_{n}=\ell.
\end{equation*}

\noindent And so $\lim x_{n}$ exists and equals $\ell$.\\

\noindent \textbf{Exercise 5}.\\

\noindent \textbf{Question (a)}. If $\nu _{2k}$ is finite and if $\nu _{2k+1}$ is infinite, then there ate exactly $k$ up-crossings : 
$[x_{\nu_{2j-1}},x_{\nu _{2j}}]$, $j=1,...,k$, that is, we have $D(a,b)=k$.\\

\noindent \textbf{Question (b)}. If $\nu _{2k+1}$ is finite and $\nu _{2k+2}$ is infinite, then there are exactly $k$ up-crossings:
$[x_{\nu_{2j-1}},x_{\nu_{2j}}]$, $j=1,...,k$, that is we have $D(a,b)=k$.\\

\noindent \textbf{Question (c)}. If all the $\nu_{j}'s$ are finite, then there are an infinite number of up-crossings : 
$[x_{\nu_{2j-1}},x_{\nu_{2j}}]$, $j\geq 1k$ : $D(a,b)=+\infty$.\\

\noindent \textbf{Question (d)}. Suppose that there exist $a < b$ rationals such that $D(a,b)=+\infty$. 
Then all the $\nu _{j}'s$ are finite. The subsequence $x_{\nu_{2j-1}}$ is strictly below $a$. 
So its limit inferior is below $a$. This limit inferior is an accumulation point of the sequence $(x_n)_{n\geq 1}$, 
so is more than $\underline{\lim}\ x_{n}$, which is below $a$.\\

\noindent Similarly, the subsequence $x_{\nu_{2j}}$ is strictly below $b$. So the limit superior is above $a$. 
This limit superior is an accumulation point of the sequence $(x_n)_{n\geq 1}$, so it is below $\overline{\lim}\ x_{n}$, 
which is directly above $b$. This leads to :
$$
\underline{\lim}\ x_{n} \leq a < b \leq \overline{\lim}\ x_{n}. 
$$

\noindent That implies that the limit of $(x_n)$ does not exist. In contrary, we just proved that the limit of $(x_n)$ exists, 
meanwhile for all the real numbers $a$ and $b$ such that $a<b$, $D(a,b)$ is finite.\\

\noindent Now, suppose that the limit of $(x_n)$ does not exist. Then,

$$
\underline{\lim}\ x_{n} < \overline{\lim}\ x_{n}. 
$$

\noindent We can then find two rationals $a$ and $b$ such that $a<b$ and a number $\epsilon$ such that $0<\epsilon$, such that 

$$
\underline{\lim}\ x_{n} < a-\epsilon < a < b < b+\epsilon <  \overline{\lim}\ x_{n}. 
$$

\noindent If $\underline{\lim}\ x_{n} < a-\epsilon$, we can return to Question \textbf{(a)} of Exercise 2 and construct a sub-sequence of $(x_n)$
which tends to $\underline{\lim}\ x_{n}$ while remaining below $a-\epsilon$. Similarly, if $b+\epsilon < \overline{\lim}\ x_{n}$, 
we can create a sub-sequence of $(x_n)$ which tends to $\overline{\lim}\ x_{n}$ while staying above $b+\epsilon$. 
It is evident with these two sequences that we could define with these two sequences all $\nu_{j}$ finite and so $D(a,b)=+\infty$.\\

\noindent We have just shown by contradiction that if all the $D(a,b)$ are finite for all rationals $a$ and $b$ such that $a<b$, 
then, the limit of $(x_n)_{n\geq 0}$ exists.\\

\noindent \textbf{Exercise 5}. Cauchy criterion in $\mathbb{R}$.\\

\noindent Suppose that the sequence is Cauchy, $i.e.$,
$$
\lim_{(p,q)\rightarrow (+\infty,+\infty)} \ (x_p-x_q)=0.
$$

\noindent Then let $x_{n_{k,1}}$ and $x_{n_{k,2}}$ be two sub-sequences converging respectively to $\ell_1=\underline{\lim}\ x_{n}$ and $\ell_2=\overline{\lim}\ x_{n}$. So

$$
\lim_{(p,q)\rightarrow (+\infty,+\infty)} \ (x_{n_{p,1}}-x_{n_{q,2}})=0.
$$

\noindent, By first letting $p\rightarrow +\infty$, we have

$$
\lim_{q\rightarrow +\infty} \ \ell_1-x_{n_{q,2}}=0,
$$

\noindent which shows that $\ell_1$ is finite, else $\ell_1-x_{n_{q,2}}$ would remain infinite and would not tend to $0$. 
By interchanging the roles of $p$ and $q$, we also have that $\ell_2$ is finite.\\

\noindent Finally, by letting $q\rightarrow +\infty$, in the last equation, we obtain
$$
\ell_1=\underline{\lim}\ x_{n}=\overline{\lim}\ x_{n}=\ell_2.
$$

\noindent which proves the existence of the finite limit of the sequence $(x_n)$.\\

\noindent Now suppose that the finite limit $\ell$ of $(x_n)$ exists. Then

$$
\lim_{(p,q)\rightarrow (+\infty,+\infty)} \ (x_p-x_q)=\ell-\ell=0,
$$0
\noindent which shows that the sequence is Cauchy.\\

%% file: asymptotics_analysis_02_en.tex
\section{Miscellaneous facts} \label{funct.facts}

\bigskip \noindent \textbf{FACT 1}. For any $a\in \mathbb{R},$ 
\begin{equation*}
\left\vert e^{ia}-1\right\vert =\sqrt{2(1-\cos a)}\leq 2\left\vert \sin
(a/2)\right\vert \leq 2\left\vert a/2\right\vert ^{\delta }.
\end{equation*}

\bigskip \noindent This is easy for $\left\vert a/2\right\vert >1.$ Indeed
for $\delta >0,\left\vert a/2\right\vert ^{\delta }>0$ and%
\begin{equation*}
2\left\vert \sin (a/2)\right\vert \leq 2\leq 2\left\vert a/2\right\vert
^{\delta }
\end{equation*}

\bigskip \noindent Now for $\left\vert a/2\right\vert >1,$ we have the
expansion

\begin{eqnarray*}
2(1-\cos a) &=&a^{2}-\sum\limits_{k=2}^{\infty }(-1)^{2}\frac{a^{2k}}{(2k)!}%
=x^{2}-2\sum\limits_{k\geq 2,k\text{ }even}^{\infty }\frac{a^{2k}}{(2k)!}-%
\frac{a^{2(k+1)}}{(2(k+1))!} \\
&=&a^{2}-2x^{2(k+1)}\sum\limits_{k\geq 2,k\text{ }even}^{\infty }\frac{1}{%
(2k)!}\left\{ \frac{1}{a^{2}}-\frac{1}{(2k+1)((2k+2)...(2k+k)}\right\} .
\end{eqnarray*}

\bigskip \noindent For each $k\geq 2,$ for $\left\vert a/2\right\vert <1,$%
\begin{equation*}
\left\{ \frac{1}{a^{2}}-\frac{1}{(2k+1)((2k+2)...(2k+k)}\right\} \geq
\left\{ \frac{1}{4}-\frac{1}{(2k+1)((2k+2)...(2k+k)}\right\} \geq 0.
\end{equation*}

\bigskip \noindent Hence 
\begin{equation*}
2(1-\cos a)\leq a^{2}.
\end{equation*}

\bigskip \noindent But for $\left\vert a/2\right\vert ,$ the function $%
\delta \hookrightarrow \left\vert a/2\right\vert ^{\delta }$ is
non-increasing in $\delta ,0\leq \delta \leq 1$. Then%
\begin{equation*}
\sqrt{2(1-\cos a)}\leq \left\vert a\right\vert =2\left\vert a/2\right\vert
^{1}\leq 2\left\vert a/2\right\vert ^{\delta }.
\end{equation*}

 